\documentclass[a4paper,10pt]{amsart}
\usepackage{amssymb}

\def\bfit{\bfseries\itshape}

\input xypic
\xyoption{all}
\xyoption{arc}

\def\theoreme#1{{\refstepcounter{theo}\label{#1}\noindent\bf Th\'eor\`eme 
\arabic{section}.\arabic{theo}.}}
\def\proposition#1{{\refstepcounter{theo}\label{#1}\noindent\bf Proposition 
\arabic{section}.\arabic{theo}.}}
\def\lemme#1{{\refstepcounter{theo}\label{#1}\noindent\bf Lemme  
\arabic{section}.\arabic{theo}.}}
\def\corollaire#1{{\refstepcounter{theo}\label{#1}\noindent\bf Corollaire 
\arabic{section}.\arabic{theo}.}}

\def\definition#1{{\refstepcounter{theo}\label{#1}\noindent\bf D\'efinition 
\arabic{section}.\arabic{theo}.}}
\def\Theoreme#1#2{{\refstepcounter{theo}\label{#2}\noindent\bf Th\'eor\`eme 
\arabic{section}.\arabic{theo} (#1).}}
\def\Proposition#1#2{{\refstepcounter{theo}\label{#2}\noindent\bf Proposition 
\arabic{section}.\arabic{theo} (#1).}}
\def\Lemme#1#2{{\refstepcounter{theo}\label{#2}\noindent\bf Lemme  
\arabic{section}.\arabic{theo} (#1).}}
\def\Corollaire#1#2{{\refstepcounter{theo}\label{#2}\noindent\bf Corollaire 
\arabic{section}.\arabic{theo} (#1).}}

\def\remarque#1{{\refstepcounter{theo}\label{#1}\noindent\sc Remarque   
\arabic{section}.\arabic{theo} - }}
\def\remarques#1{{\refstepcounter{theo}\label{#1}\noindent\sc Remarques   
\arabic{section}.\arabic{theo} - }}
\def\exemple#1{{\refstepcounter{theo}\label{#1}\noindent\sc Exemple 
\arabic{section}.\arabic{theo} - }}
\def\exemples#1{{\refstepcounter{theo}\label{#1}\noindent\sc Exemples 
\arabic{section}.\arabic{theo} - }}

\def\equat{\refstepcounter{theo}$$~}
\def\endequat{\leqno{\boldsymbol{(\arabic{section}.\arabic{theo})}}~$$}

\newcounter{soussection}[section]

\def\soussection#1{\refstepcounter{soussection}
\noindent{\bf \arabic{section}.\Alph{soussection}. #1.}}

\def\FM{{\mathbb{F}}}

\def\NM{{\mathbb{N}}}

\def\QM{{\mathbb{Q}}}
\def\RM{{\mathbb{R}}}

\def\ZM{{\mathbb{Z}}}

          \def\eG{{\mathfrak e}}
          
\def\GG{{\mathfrak G}}

\def\RG{{\mathfrak R}}          
\def\SG{{\mathfrak S}}          \def\sG{{\mathfrak s}}

\def\GGB{{\boldsymbol{\mathfrak G}}}

\def\a{\alpha}
\def\b{\beta}

\def\g{\gamma}
\def\G{\Gamma}
\def\d{\delta}
\def\D{\Delta}
\def\e{\varepsilon}
\def\ph{\varphi}
\def\ch{\chi}
\def\l{\lambda}
\def\L{\Lambda}
\def\m{\mu}
\def\n{\nu}
\def\o{\omega}

\def\r{\rho}
\def\s{\sigma}

\def\th{\theta}
\def\Th{\Theta}
\def\t{\tau}
\def\x{\xi}

\def\z{\zeta}

\def\lamb{{\boldsymbol{\lambda}}}

\def\mub{{\boldsymbol{\mu}}}

\def\Sigb{{\boldsymbol{\Sig}}}

\def\AC{{\mathcal{A}}}

\def\EC{{\mathcal{E}}}
\def\FC{{\mathcal{F}}}
\def\GC{{\mathcal{G}}}
\def\HC{{\mathcal{H}}}
\def\IC{{\mathcal{I}}}

\def\LC{{\mathcal{L}}}
\def\MC{{\mathcal{M}}}
\def\NC{{\mathcal{N}}}

\def\PC{{\mathcal{P}}}

\def\RC{{\mathcal{R}}}
\def\SC{{\mathcal{S}}}

\def\UC{{\mathcal{U}}}

\def\XC{{\mathcal{X}}}

\def\ZC{{\mathcal{Z}}}

\def\ECt{{\tilde{\mathcal{E}}}}
\def\FCt{{\tilde{\mathcal{F}}}}

\def\HCt{{\tilde{\mathcal{H}}}}

\def\LCt{{\tilde{\mathcal{L}}}}

\def\BCB{{\boldsymbol{\mathcal{B}}}}

\def\UCB{{\boldsymbol{\mathcal{U}}}}

          \def\ab{{\mathbf a}}
\def\Bb{{\mathbf B}}          
\def\Cb{{\mathbf C}}          
\def\Db{{\mathbf D}}          
          \def\eb{{\mathbf e}}
          
\def\Gb{{\mathbf G}}          
\def\Hb{{\mathbf H}}

\def\Lb{{\mathbf L}}          
\def\Mb{{\mathbf M}}          
          \def\nb{{\mathbf n}}
\def\Ob{{\mathbf O}}          
\def\Pb{{\mathbf P}}          \def\pb{{\mathbf p}}
\def\Qb{{\mathbf Q}}          
\def\Rb{{\mathbf R}}          
\def\Sb{{\mathbf S}}          
\def\Tb{{\mathbf T}}          \def\tb{{\mathbf t}}
\def\Ub{{\mathbf U}}          \def\ub{{\mathbf u}}
\def\Vb{{\mathbf V}}

\def\Yb{{\mathbf Y}}          
\def\Zb{{\mathbf Z}}

          \def\ati{{\tilde{a}}}

          \def\fti{{\tilde{f}}}
\def\Gti{{\tilde{G}}}          \def\gti{{\tilde{g}}}

          \def\lti{{\tilde{l}}}
\def\Mti{{\tilde{M}}}          
\def\Nti{{\tilde{N}}}

\def\Rti{{\tilde{R}}}          
          \def\sti{{\tilde{s}}}
\def\Tti{{\tilde{T}}}          \def\tti{{\tilde{t}}}
          
\def\Vti{{\tilde{V}}}          
          \def\wti{{\tilde{w}}}
          \def\xti{{\tilde{x}}}
          \def\yti{{\tilde{y}}}
          \def\zti{{\tilde{z}}}

          \def\fha{{\hat{f}}}
          
          \def\hha{{\hat{h}}}

\def\Rha{{\hat{R}}}          
          \def\sha{{\hat{s}}}

          \def\wha{{\hat{w}}}

          \def\zha{{\hat{z}}}

          \def\wdo{{\dot{w}}}

          \def\zdo{{\dot{z}}}

\def\Wba{{\bar{W}}}          \def\wba{{\bar{w}}}
          \def\xba{{\bar{x}}}
          
          \def\zba{{\bar{z}}}

\def\Bbt{{\tilde{\Bb}}}

\def\Gbt{{\tilde{\Gb}}}

\def\Lbt{{\tilde{\Lb}}}

\def\Pbt{{\tilde{\Pb}}}

\def\Tbt{{\tilde{\Tb}}}

\def\Ybt{{\tilde{\Yb}}}          
\def\Zbt{{\tilde{\Zb}}}

          \def\Gbh{{\hat{\Gb}}}

          \def\Ybh{{\hat{\Yb}}}

\def\lamb{{\boldsymbol{\lambda}}}

\def\mub{{\boldsymbol{\mu}}}

\def\Sigb{{\boldsymbol{\Sigma}}}

\def\alpt{{\tilde{\alpha}}}
\def\bett{{\tilde{\beta}}}
\def\gamt{{\tilde{\gamma}}}

\def\Delt{{\tilde{\Delta}}}

\def\phit{{\tilde{\phi}}}
\def\Phit{{\tilde{\Phi}}}

\def\chit{{\tilde{\chi}}}
\def\lamt{{\tilde{\lambda}}}

\def\mut{{\tilde{\mu}}}
\def\nut{{\tilde{\nu}}}

\def\sigt{{\tilde{\s}}}

\def\taut{{\tilde{\t}}}
\def\thet{{\tilde{\theta}}}

\def\etat{{\tilde{\eta}}}

\def\phh{{\hat{\varphi}}}

\def\chih{{\hat{\chi}}}

\def\omeh{{\hat{\omega}}}

\def\rhoh{{\hat{\rho}}}
\def\sigh{{\hat{\s}}}

\def\Gambh{{\hat{\boldsymbol{\Gamma}}}}

\def\ad{\mathop{\mathrm{ad}}\nolimits}

\def\Aut{\mathop{\mathrm{Aut}}\nolimits}

\def\Cent{\mathop{\mathrm{Cent}}\nolimits}
\def\Class{\mathop{\mathrm{Class}}\nolimits}

\def\cus{{\mathrm{cus}}}
\def\Cus{{\mathrm{Cus}}}

\def\Ext{\mathop{\mathrm{Ext}}\nolimits}

\def\Hom{\mathop{\mathrm{Hom}}\nolimits}

\def\Id{\mathop{\mathrm{Id}}\nolimits}

\def\Im{\mathop{\mathrm{Im}}\nolimits}

\def\Ind{\mathop{\mathrm{Ind}}\nolimits}
\def\Inf{\mathop{\mathrm{Inf}}\nolimits}
\def\Int{\mathop{\mathrm{Int}}\nolimits}
\def\INT{\mathop{\mathrm{int}}\nolimits}
\def\Irr{\mathop{\mathrm{Irr}}\nolimits}

\def\Ker{\mathop{\mathrm{Ker}}\nolimits}

\def\mod{\mathop{\mathrm{mod}}\nolimits}

\def\Out{\mathop{\mathrm{Out}}\nolimits}

\def\reg{{\mathrm{r\acute{e}g}}}

\def\res{\mathop{\mathrm{res}}\nolimits}
\def\Res{\mathop{\mathrm{Res}}\nolimits}
\def\rang{\mathop{\mathrm{rg}}\nolimits}

\def\Tr{\mathop{\mathrm{Tr}}\nolimits}

\def\uni{{\mathrm{uni}}}
\def\faisceau{{\mathrm{FCar}}}

\def\tete#1{\par\leavevmode\makebox[0.7cm]{$(\mathrm{#1})$}}

\def\imp{\Rightarrow}

\def\to{\rightarrow}
\def\longto{\longrightarrow}
\def\injto{\hookrightarrow}

\def\mapright#1{\hspace{0.2em}\smash{
     \mathop{\rightarrow}\limits^{\SS#1}}\hspace{0.2em}}
\def\doublefleche#1{\hspace{0.2em}\smash{
     \mathop{\longleftrightarrow}\limits^{\SS#1}}\hspace{0.2em}}

\def\longmapright#1{\hspace{0.3em}\smash{
     \mathop{\longrightarrow}\limits^{#1}}\hspace{0.3em}}

\def\fonction#1#2#3#4#5{\begin{array}{rccc}
{#1} : & {#2} & \longto & {#3} \\
& {#4} & \longmapsto & {#5} 
\end{array}}

\def\fonctio#1#2#3#4{\begin{array}{ccc}
{#1} & \longto & {#2} \\
{#3} & \longmapsto & {#4} 
\end{array}}

\def\ci{\circ}
\def\pr{\prime}
\def\ve{\vee}
\def\we{\wedge}

\def\incl{\hspace{0.05cm}{\subset}\hspace{0.05cm}}
\def\notincl{\hspace{0.05cm}{\not\subset}\hspace{0.05cm}}

\def\vide{\varnothing}

\def\fq{\FM_q}
\def\fp{\FM_p}
\def\ql{{\QM_\el}}
\def\qlb{{\overline{\QM}_\el}}

\def\DS{\displaystyle}
\def\SS{\scriptstyle}
\def\SSS{\scriptscriptstyle}

\def\fin{~$\SS \blacksquare$}
\def\finl{~$\SS \square$}
\def\el{\ell}

\def\lexp#1#2{\kern\scriptspace\vphantom{#2}^{#1}\kern-\scriptspace#2}
\def\le{\hspace{0.1em}\mathop{\leqslant}\nolimits\hspace{0.1em}}
\def\ge{\hspace{0.1em}\mathop{\geqslant}\nolimits\hspace{0.1em}}

\mathchardef\lllllll="3278
\def\SEC{$\lllllll$}
\mathchardef\inferieur="321E
\mathchardef\superieur="321F

\def\proof{\noindent{\sc{D\'emonstration}~-} }

\def\eqna{\begin{eqnarray*}}
\def\endeqna{\end{eqnarray*}}

\def\borel{sous-groupe de Borel }

\def\cad{c'est-\`a-dire }
\def\car{caract\`ere }
\def\cars{caract\`eres }
\def\ele{\'el\'ement }
\def\eles{\'el\'ements }

\def\irr{irr\'eductible }
\def\irrs{irr\'eductibles }
\def\levi{sous-groupe de Levi }
\def\levis{sous-groupes de Levi }
\def\levic{compl\'ement de Levi }
\def\levics{compl\'ements de Levi }
\def\mor{morphisme }
\def\mors{morphismes }

\def\iso{isomorphisme }
\def\isos{isomorphismes }
\def\auto{automorphisme }

\def\para{sous-groupe parabolique }
\def\paras{sous-groupes paraboliques }

\def\resp{respectivement }
\def\ssi{si et seulement si }

\def\tor{tore maximal }
\def\tors{tores maximaux }

\def\tors{{\mathrm{tors}}}
\def\mini{{\mathrm{min}}}

\def\phan{{\phantom{\bullet}}}

\def\sem{{\mathrm{sem}}}
\def\pro{{\prime *}}

\def\OF{{*F^*}}
\def\RES{{\RG\eG\sG}}
\def\ORES{\lexp{*}{\RES}}
\def\minus{{\mathrm{minus}}}
\def\compose{\mathop{\circ}}

\renewcommand{\today}{%
  \number\day\space\ifcase\month\or
  Janvier\or F\'evrier\or Mars\or Avril\or Mai\or Juin\or
  Juillet\or Ao\^ut\or Septembre\or Octobre\or Novembre\or D\'ecembre\fi
  \space\number\year}

\def\rhodot{{\dot{\rho}}}

\begin{document}

\begin{centerline}
{\Large \bf Sur les caract\`eres des groupes r\'eductifs finis \`a centre non connexe~:}
\end{centerline}

\bigskip

\begin{centerline}
{\Large \bf applications aux groupes sp\'eciaux lin\'eaires et unitaires}
\end{centerline}

\bigskip

\begin{centerline}
{\sc C\'edric Bonnaf\'e\footnote{CNRS, UMR 6623, \\ D\'epartement de Math\'ematiques, \\
Universit\'e de Franche-Comt\'e\\ 16 Route de Gray \\ 25030 BESAN\c{C}ON Cedex \\
FRANCE \\ {\tt bonnafe@math.univ-fcomte.fr}}}
\end{centerline}

\bigskip

\begin{centerline}{\today}\end{centerline}

\vskip1cm

\begin{quotation}
\noindent{\bf R\'esum\'e :} 
Un premier but de cet article est de pr\'esenter une synth\`ese des 
r\'esultats de plusieurs auteurs concernant les caract\`eres des groupes 
r\'eductifs finis \`a centre non connexe. Nous nous int\'eressons 
particuli\`erement aux probl\`emes directement li\'es \`a la 
non connexit\'e du centre. Nous insistons notamment sur les 
caract\`eres de Gelfand-Graev et les caract\`eres semisimples.

Un deuxi\`eme but est d'\'etudier l'influence de la non connexit\'e 
du centre sur la th\'eorie des faiceaux-caract\`eres. Nous nous concentrons 
notamment sur la famille des faisceaux-caract\`eres dont le support 
rencontre la classe unipotente r\'eguli\`ere~: ce sont les analogues 
naturels des caract\`eres semisimples. 

Le dernier but est l'application de ces r\'esultats aux groupes 
r\'eductifs finis de type $A$, d\'eploy\'es ou non (comme par 
exemple les groupes sp\'eciaux lin\'eaires ou unitaires). 
Lorsque le cardinal du corps fini de r\'ef\'erence 
est assez grand, nous obtenons un param\'etrage des caract\`eres 
irr\'eductibles, calculons explicitement le foncteur d'induction de 
Lusztig dans la base des caract\`eres irr\'eductibles, param\'etrons 
les faisceaux-caract\`eres et montrons que les fonctions caract\'eristiques 
de ces faisceaux-caract\`eres sont des transform\'ees de Fourier des 
caract\`eres irr\'eductibles (conjecture de Lusztig). Ces r\'esultats 
permettent de construire un algorithme th\'eorique pour 
calculer la table de caract\`eres de ces groupes.
\end{quotation}

\vskip1cm

\begin{quotation}
\noindent{\bf Abstract :} 
A first aim of this paper is to present an overview of results 
obtained by several authors on the characters of finite reductive groups 
with non-connected centre. We are particularly interested in problems 
directly linked to the non-connectedness of the centre. We insist 
on Gelfand-Graev and semisimple characters. 

A second aim is to study the influence of the non-connectedness 
of the centre on the theory of character sheaves. We study 
more precisely the family of character sheaves whose support meets 
the regular unipotent class: these are analogues of the semisimple 
characters. 

The last aim is the application of these results to finite reductive 
groups of type $A$, split or not (as for instance the special linear 
or special unitary groups). Whenever the cardinality of the 
finite field is large enough, we obtain a parametrization of 
the irreducible characters, a parametrization of the character sheaves, 
and we show that the characteristic functions of character sheaves 
are Fourier transforms of the irreducible characters (Lusztig's conjecture). 
This gives a theoretical algorithm for computing the 
character table of these groups. 
\end{quotation}

\newpage

\def\proof{\noindent{\sc D\'emonstration - }} 
\def\ZCC{{\SSS{\ZC}}}
\def\fromgg{\begin{quotation}{\tiny \noindent{\sc 
Comparaison de $\Gb'$ et $\Gb$.\hskip0.2cm~}}}
\def\endfromgg{\end{quotation}}
\def\aff{{\mathrm{aff}}}
\def\diagr{{\mathrm{diag}}}
\def\acompleter{$${\boldsymbol{\bullet \# ! ****** ! \#\bullet}}$$}
\def\itemth#1{\item[${\mathrm{(#1)}}$]}

{\Large \part*{Introduction\label{intro intro}}}

\vskip1cm

En 1907, Schur \cite{schur} et Jordan \cite{jordan} d\'eterminaient 
les tables de caract\`eres du groupe g\'en\'eral 
lin\'eaire $GL(2,q)$ et du groupe sp\'ecial lin\'eaire $SL(2,q)$. 
En 1951, Steinberg \cite{steinberg gl3} d\'eterminait la table 
de caract\`eres de $GL(3,q)$ et $GL(4,q)$, en utilisant entre autres 
des constructions g\'en\'erales de repr\'esentations (que nous appelons 
de nos jours {\it unipotentes}) de $GL(n,q)$ (voir \cite{steinberg uni}). 
Finalement, au prix d'un tour de force combinatoire remarquable, Green \cite{green} 
d\'eterminait en 1955 la table de caract\`eres de $GL(n,q)$ (au moins 
algorithmiquement). Par contre, les progr\`es concernant le groupe 
sp\'ecial lin\'eaire furent beaucoup plus long. 
En 1971, Lehrer \cite{lehrer these} 
d\'eterminait une partie de la table de caract\`eres de 
$SL(4,q)$, celle correspondant aux s\'eries discr\`etes. Citons 
\'egalement les travaux de Lehrer \cite{lehrer 1973} 
puis Digne, Lehrer et Michel \cite{DLM1} 
qui donnent des informations partielles pour le calcul de la table 
de caract\`eres de $SL(n,q)$. Ces informations sont suffisantes 
pour compl\'eter la table de caract\`eres lorsque $n$ est premier. 
Un des buts de cet article est de fournir un algorithme 
th\'eorique pour calculer la table de caract\`eres de $SL(n,q)$~: 
cependant, nous ne sommes capable de montrer la validit\'e de 
cet algorithme que lorsque $q$ est assez grand. 

Plus g\'en\'eralement, si $\Gb$ est un groupe r\'eductif connexe d\'efini 
sur une cl\^oture alg\'ebrique $\FM$ du corps fini 
\`a $p$ \'el\'ements $\FM_p$ ($p$ premier) et si $F : \Gb \to \Gb$ est 
une isog\'enie dont une puissance est un endomorphisme de Frobenius 
relatif \`a une $\FM_q$-structure sur $\Gb$ ($q=p^?$), le calcul de 
la table de caract\`eres du groupe fini $\Gb^F$ (appel\'e 
{\it groupe r\'eductif fini}) est loin d'\^etre 
r\'esolu en toute g\'en\'eralit\'e. Rappelons quand m\^eme que le 
cas du groupe $Sp(4,q)$ a \'et\'e r\'esolu, pour $q$ impair, 
par Srinivasan \cite{srinivasan} en 1968. D'autres r\'esultats 
ont \'et\'e obtenus sur les petits groupes (groupes de Suzuki, groupes 
de Ree...). 

Pourtant, en 1976, l'article fondateur de Deligne et Lusztig \cite{delu} 
permettait \`a la th\'eorie des caract\`eres des groupes r\'eductifs finis 
de faire des progr\`es consid\'erables. Leur id\'ee, inspir\'ee par des calculs 
de Drinfeld montrant que la s\'erie discr\`ete de $SL(2,q)$ apparaissait 
dans la cohomologie $\ell$-adique de la vari\'et\'e d\'efinie par l'\'equation 
$xy^q-yx^q=1$, \'etait d'utiliser la structure de vari\'et\'e de $\Gb$ 
pour produire des sous-vari\'et\'es de $\Gb$ sur lequel le groupe fini 
$\Gb^F$ agit et de r\'ecup\'erer ainsi des repr\'esentations de $\Gb^F$ 
dans la cohomologie $\ell$-adique de ces vari\'et\'es. Poursuivant dans 
cette voix, Lusztig \cite{lubook}, apr\`es une s\'erie impressionnante d'articles, 
obtenait en 1984 le param\'etrage des caract\`eres irr\'eductibles de 
$\Gb^F$ (dans l'esprit du programme de Langlands) lorsque le centre de 
$\Gb$ est connexe. En plus de ce param\'etrage, 
il obtenait une formule explicite pour le degr\'e de ces caract\`eres 
ainsi qu'un algorithme (th\'eorique) permettant de calculer les valeurs 
de ces caract\`eres en les \'el\'ements semi-simples. 

Au cours de sa d\'emarche, Lusztig introduisait une nouvelle base orthonormale de 
l'espace des fonctions centrales, la base des {\it caract\`eres fant\^omes}, 
obtenue \`a partir de la base des caract\`eres irr\'eductibles 
par une matrice diagonale par blocs, les blocs \'etant des matrices de 
transform\'ees de Fourier associ\'ees \`a des petits groupes finis 
(dont la taille ne d\'epend pas de $q$). En 1984-1986, 
Lusztig (voir \cite{luicc} et \cite{lucs}) 
d\'eveloppait une nouvelle th\'eorie, la th\'eorie des 
{\it faisceaux-caract\`eres}, dans le but de comprendre l'intrusion de ces 
petits groupes finis et de ces caract\`eres fant\^omes. Un faisceau-caract\`ere 
est un faisceau pervers $\Gb$-\'equivariant irr\'eductible sur $\Gb$ satisfaisant 
\`a certaines conditions. Si $A$ est un faisceau-caract\`ere $F$-stable, 
on peut lui associer une fonction centrale sur $\Gb^F$, appel\'ee 
{\it fonction caract\'eristique de $A$}~; cette fonction n'est d\'efinie 
qu'\`a une constante multiplicative pr\`es mais Lusztig a d\'efini des 
normalisations qui en font des fonctions de norme $1$. De plus, Lusztig a montr\'e
que ces fonctions caract\'eristiques de faisceaux-caract\`eres $F$-stables 
forment une base orthonormale de l'espace des fonctions centrales. 
Il a fait la conjecture suivante~:

\bigskip

\begin{quotation}
\noindent{\bf Conjecture de Lusztig :} {\it Si le centre de $\Gb$ est 
connexe, la matrice de passage entre la base des caract\`eres 
fant\^omes et la base des fonctions caract\'eristiques de faisceaux-caract\`eres 
$F$-stables sur $\Gb$ est diagonale.}
\end{quotation}

\bigskip

\noindent D'autre part, Lusztig a aussi d\'ecrit un algorithme th\'eorique permettant 
de calculer les fonctions caract\'eristiques de faisceaux-caract\`eres 
$F$-stables, m\^eme lorsque le centre de $\Gb$ n'est pas connexe. 
En 1995, Shoji \cite{shoji almost} d\'emontrait cette conjecture. 

\medskip

Il appara\^\i t ainsi au cours de l'\'evolution de la th\'eorie que 
la non-connexit\'e du centre de $\Gb$ entra\^\i ne de nombreuses 
complications. Certaines sont techniques (comme par exemple 
le param\'etrage des classes de conjugaison), d'autres sont th\'eoriques 
(comme par exemple la non-connexit\'e du centralisateur des \'el\'ements 
semi-simples du dual de $\Gb$ ou l'augmentation significative du nombre de 
faisceaux-caract\`eres cuspidaux). 
Cette non-connexit\'e du centre explique les difficult\'es 
qui ont \'emaill\'e la recherche d'une table de caract\`eres pour 
le groupe sp\'ecial lin\'eaire. 

Plusieurs auteurs ont \'etudi\'e les groupes r\'eductifs finis 
\`a centre non connexe (Asai \cite{asai}, Lusztig \cite{luznc}, 
Digne et Michel \cite{dm}, Digne, Lehrer et Michel \cite{DLM1}, 
\cite{DLM2}, Shoji \cite{shoji} ou l'auteur \cite{cedthese}, 
\cite{bonnafe mackey}, \cite{bonnafe torsion}). Le param\'etrage des 
caract\`eres irr\'eductibles a pu ainsi \^etre achev\'ee \cite{luznc} 
et des nouvelles informations sur la table de caract\`eres du groupe $\Gb^F$ 
(par exemple la valeur en les \'el\'ements unipotents r\'eguliers \cite{DLM1}) 
ont \'et\'e obtenues. 

Concernant l'analogue de la conjecture de Lusztig, un des premiers 
probl\`emes vient de ce qu'il n'y a pas de d\'efinition 
indiscutable de la notion de caract\`ere fant\^ome. On 
peut alors consid\'erer comme une r\'eponse positive 
\`a la conjecture de Lusztig pour les groupes \`a centre 
non connexe un th\'eor\`eme qui montrerait 
que la base des fonctions caract\'eristiques de faisceaux-caract\`eres 
$F$-stables est obtenue \`a partir de la base des caract\`eres 
irr\'eductibles par une matrice diagonale par blocs, les blocs 
\'etant des matrices de transform\'ees de Fourier associ\'ees 
\`a des petits groupes finis. Dans cette acceptation, 
la conjecture de Lusztig a \'et\'e d\'emontr\'ee pour les groupes 
sp\'eciaux orthogonaux et symplectiques par Waldspurger \cite{waldspurger} 
lorsque $q$ est assez grand.  

Shoji \cite{shoji banff} a d\'emontr\'e la conjecture 
de Lusztig pour le groupe $SL(n,q)$ pour $p > 3n$ 
et $q$ une puissance quelconque de $p$. Il a aussi propos\'e une d\'efinition 
int\'eressante de caract\`ere fant\^ome~: un caract\`ere fant\^ome 
devrait \^etre la descente de Shintani de $\Gb^{F^n}$ \`a $\Gb^F$ d'un 
caract\`ere irr\'eductible $F$-stable de $\Gb^{F^n}$ (pour $n$ suffisamment 
divisible). Cette d\'efinition a le m\'erite d'\^etre correcte lorsque 
le centre de $\Gb$ est connexe et de rendre vraie la conjecture de 
Lusztig dans le groupe sp\'ecial lin\'eaire.

Un des buts du pr\'esent article est de d\'emontrer la conjecture 
de Lusztig (sans prendre la d\'efinition de Shoji de caract\`ere 
fant\^ome) pour le groupe sp\'ecial lin\'eaire et le groupe 
sp\'ecial unitaire lorsque $p$ est quelconque et $q$ est assez grand. 
Plus pr\'ecis\'ement, 
nous obtenons un param\'etrage des caract\`eres de ces groupes 
r\'eductifs finis et d\'efinissons a priori, sans r\'ef\'erence 
\`a la th\'eorie des faisceaux-caract\`eres, des transform\'ees de Fourier 
naturelles de ces caract\`eres irr\'eductibles (nous nous inspirons de 
\cite[\SEC 5 et 6]{dm}). Nous montrons alors 
que la matrice de passage entre la base des transform\'ees de Fourier 
et la base des fonctions caract\'eristiques de faisceaux-caract\`eres 
$F$-stables est diagonale et calculons explicitement les coefficients 
diagonaux. Dans l'optique 
d'obtenir un algorithme pour calculer la table de caract\`eres 
de ces groupes, notre r\'esultat est satisfaisant. 
Il faut cependant \^etre r\'ealiste~: la mise en \oe uvre de 
cet algorithme n\'ecessite encore un travail consid\'erable. 
Une derni\`ere remarque~: notre r\'esultat est valide pour tous 
les groupes de type $A$, quel que soit l'endomorphisme de Frobenius 
consid\'er\'e. M\^eme dans le cas {\it d\'eploy\'e}, il s'applique 
aux groupes interm\'ediaires de la forme $\Sb\Lb_n/\mub_d$, o\`u 
$d$ divise $n$~: ces groupes sont des extensions non triviales 
du groupe fini $SL(n,q)/\mub_d(\FM_q)$ qui ne sont pas contenus 
dans le travail de Shoji \cite{shoji banff}.

Dans le cas du groupe sp\'ecial lin\'eaire, il serait int\'eressant 
de relier plus finement le param\'etrage de Shoji et le notre pour 
d\'eterminer la matrice de passage entre nos transform\'ees de 
Fourier et les caract\`eres fant\^omes de Shoji~: lorsque $q$ est 
assez grand, cette matrice de passage est diagonale mais 
nous n'en connaissons pas les coefficients. Il serait 
aussi int\'eressant, dans le cas du groupe sp\'ecial unitaire, 
de savoir si nos transform\'ees de Fourier sont des caract\`eres 
fant\^omes (\`a une constante pr\`es) au sens de Shoji. 

\bigskip

Cet article a aussi un autre but~: pr\'esenter une synth\`ese 
des r\'esultats sur les groupes r\'eductifs finis 
directement li\'es \`a la non connexit\'e du centre. 
Notons $\ZC(\Gb)=\Zb(\Gb)/\Zb(\Gb)^\circ$ le groupe 
des composantes connexes du centre de $\Gb$. Nous montrons 
comment relier $\ZC(\Gb)$ au syst\`eme de racines de $\Gb$, 
\'etudions le morphisme $\ZC(\Gb)\to\ZC(\Lb)$ 
(o\`u $\Lb$ est un \levi de $\Gb$) en lien avec 
les automorphismes du diagramme de Dynkin affine, 
relions la structure de $\ZC(\Gb)$ avec la non-connexit\'e 
du centralisateur des \'el\'ements semi-simples du dual 
$\Gb^*$ de $\Gb$, \'etudions la distinction qu'elle 
entra\^\i ne entre s\'eries de Lusztig g\'eom\'etriques 
et rationnelles, \'etudions l'action de $H^1(F,\ZC(\Gb))$ 
sur les caract\`eres de $\Gb^F$ \`a travers la th\'eorie de 
Harish-Chandra, calculons les composantes irr\'eductibles 
des caract\`eres de Gelfand-Graev, \'etudions l'action de $\ZC(\Gb)$ (par 
conjugaison ou par translation) sur les faisceaux-caract\`eres, 
avant d'appliquer tout ceci aux caract\`eres des groupes de type $A$. 
Parmi ces r\'esultats, beaucoup sont bien connus 
et d\^us \`a d'autres auteurs, mais nous avons 
souhait\'e les pr\'esenter ensemble, notamment pour les relier 
entre eux et quelquefois pour en am\'eliorer l\'eg\`erement 
le degr\'e de g\'en\'eralit\'e. 

Pour \'etudier les groupes \`a centre non connexe, 
nous reprenons une technique courante \cite{delu}~: elle consiste \`a voir 
$\Gb$ comme un sous-groupe ferm\'e distingu\'e d'un groupe 
$\Gbt$ \`a centre connexe tel que $\Gbt/\Gb$ soit ab\'elien 
(c'est toujours possible~; par exemple, plonger $\Sb\Lb_n(\FM)$ dans 
$\Gb\Lb_n(\FM)$). La th\'eorie de Clifford permet alors, 
par restriction de $\Gbt^F$ \`a $\Gb^F$, d'utiliser ce que l'on sait 
de $\Gbt^F$, par exemple par les avantages li\'es \`a la connexit\'e 
du centre de $\Gbt$. 

Cet article est organis\'e comme suit. Dans le chapitre 
\ref{chapitre preliminaire}, nous introduisons les notations g\'en\'erales 
en vigueur dans tout l'article, pr\'esentons le contexte et 
\'etablissons quelques r\'esultats pr\'eliminaires \`a la suite. 
Dans le chapitre \ref{chapitre ZG}, nous montrons comment calculer 
$\ZC(\Gb)$ et $\ZC(\Lb)$ (pour un \levi $\Lb$ de $\Gb$) de 
plusieurs mani\`eres. Nous rappelons les diff\'erentes constructions 
d'un morphisme entre le groupe $A_{\Gb^*}(s)$ des composantes connexes 
du centralisateur d'un \'el\'ement semi-simple $s$ de $\Gb^*$ et le groupe 
$\ZC(\Gb)^\wedge$ des caract\`eres lin\'eaires de $\ZC(\Gb)$. Nous 
y construisons une action de $H^1(F,\ZC(\Gb))$ sur les fonctions 
centrales sur $\Gb^F$. Nous rappelons aussi les notions de {\it cuspidalit\'e} 
introduites dans \cite{bonnafe cras}, \cite{bonnafe mackey}, 
\cite{bonnafe torsion} et 
\cite{bonnafe regulier}.
Un des buts du chapitre \ref{chapitre lusztig} est de d\'emontrer 
la disjonction des s\'eries de Lusztig rationnelles. L'essentiel de 
cette preuve est contenu dans \cite{luirr} ou \cite{dmbook}. 
Dans le chapitre \ref{chapitre harish-chandra}, nous \'etudions 
l'action de $H^1(F,\ZC(\Gb))$ \`a travers la th\'eorie de Harish-Chandra. 
Nous ne pensons pas que ceci soit trait\'e ailleurs dans ce degr\'e 
de g\'en\'eralit\'e. Le chapitre \ref{chapitre gelfand}, largement inspir\'e 
par \cite{asai}, \cite{DLM1}, \cite{DLM2} et \cite{cedthese}, 
traite des \'el\'ements unipotents r\'eguliers, 
des caract\`eres de Gelfand-Graev, de leurs composantes 
irr\'eductibles (les caract\`eres dits {\it r\'eguliers}) et de 
leur dual de Curtis (les caract\`eres dits {\it semi-simples}). 
Nous obtenons notamment une d\'ecomposition des 
caract\`eres semi-simples comme combinaison lin\'eaire d'induits de 
fonctions absolument cuspidales. Dans le chapitre \ref{chapitre faisceaux}, 
nous \'etudions les diff\'erentes actions de $\ZC(\Gb)$ sur les 
faisceaux-caract\`eres. Tout d'abord, si $A$ est un faisceau-caract\`ere, 
la $\Gb$-\'equivariance de $A$ induit une action de $\ZC(\Gb)$ sur $A$ 
via un caract\`ere lin\'eaire. De plus, via l'action par translation, 
$\ZC(\Gb)$ permute les 
faisceaux-caract\`eres~: nous d\'eterminons l'action de cette permutation 
\`a travers le proc\'ed\'e d'induction \`a partir des faisceaux-caract\`eres 
cuspidaux. Pour finir, nous \'etudions les faisceaux-caract\`eres 
apparaissant dans l'induit de faisceaux-caract\`eres cuspidaux 
dont le support rencontre la classe unipotente r\'eguli\`ere 
et d\'ecrivons leur fonction caract\'eristique en termes 
d'induction de Lusztig. 

Dans le chapitre \ref{chapitre a}, nous supposons que toutes les composantes 
quasi-simples de $\Gb$ sont  
de type $A$ et que $\Gb$ est muni d'un endomorphisme de Frobenius quelconque, 
d\'eploy\'e ou non. Nous montrons comment les fonctions centrales 
introduites dans le chapitre \ref{chapitre gelfand} permettent, lorsque 
$q$ est grand, de construire les caract\`eres irr\'eductibles 
comme combinaisons lin\'eaires d'induits de fonctions caract\'eristiques 
de fonctions absolument cuspidales. En utilisant le fait 
que le support de tout faisceau-caract\`ere cuspidal sur $\Gb$ 
rencontre (\`a translation pr\`es par $\ZC(\Gb)$) la classe unipotente 
r\'eguli\`ere et les formules de la derni\`ere section du chapitre 
\ref{chapitre faisceaux}, nous obtenons la conjecture de Lusztig. 
Donnons-en un \'enonc\'e sommaire~: soit $s$ un \'el\'ement semi-simple 
de $\Gb^{*F^*}$, soit $W^\circ(s)$ le groupe de Weyl de $C_{\Gb^*}^\circ(s)$, 
soit $A_{\Gb^*}(s)=C_{\Gb^*}(s)/C_{\Gb^*}^\circ(s)$ (c'est un groupe 
ab\'elien), soit $\EC(\Gb^F,(s))$ 
la s\'erie de Lusztig g\'eom\'etrique de $\Gb^F$ associ\'ee \`a la classe 
de conjugaison de $s$ dans $\Gb^*$ et soit $\faisceau(\Gb,(s))^F$ 
la s\'erie g\'eom\'etrique des faisceaux-caract\`eres $F$-stables associ\'ee 
\`a $s$. Notons $\IC(\Gb,s)$ l'ensemble des triplets $(\chi,\xi,\a)$ 
tels que $\chi$ parcourt un ensemble de repr\'esentants des 
$A_{\Gb^*}(s)$-orbites $F^*$-stables de caract\`eres 
irr\'eductibles de $W^\circ(s)$, $\xi\in (A_{\Gb^*}(s,\chi)^{F^*})^\we$ 
et $\z \in H^1(F^*,A_{\Gb^*}(s,\chi))$ 
(o\`u $A_{\Gb^*}(s,\chi)$ est le stabilisateur de $\chi$ dans 
$A_{\Gb^*}(s)$). Notons $\IC^\vee(\Gb,s)$ l'ensemble des triplets 
$(\chi,a,\t)$ o\`u $\chi$ parcourt un ensemble de repr\'esentants des 
$A_{\Gb^*}(s)$-orbites $F^*$-stables 
de caract\`eres irr\'eductibles de $W^\circ(s)$, $a \in A_{\Gb^*}(s,\chi)^{F^*}$ et 
$\t \in H^1(F^*,A_{\Gb^*}(s,\chi))^\we$. Si $A$ est un faisceau-caract\`ere 
$F$-stable sur $\Gb$, nous noterons $\XC_A$ sa fonction caract\'eristique 
(explicitement normalis\'ee comme dans l'article). 

\bigskip

\noindent{\bf Th\'eor\`eme.} 
{\it Supposons $q$ assez grand. Alors il existe deux bijections 
$$\fonctio{\IC(\Gb,s)}{\EC(\Gb^F,(s))}{(\chi,\xi,\a)}{R_\chi(s)_{\xi,\a}}$$
$$\fonctio{\IC^\vee(\Gb,s)}{\faisceau(\Gb,(s))^F}{(\chi,a,\t)}{A_\chi(s)_{a,\t}}
\leqno{\text{et}}$$
telles que, si $(\chi,a,\t) \in \IC^\vee(\Gb,s)$, alors
$$\XC_{A_\chi(s)_{a,\t}}=\frac{\z_{s,\chi,a,\t}}{|A_{\Gb^*}(s,\chi)|} 
\sum_{\SS{\xi \in (A_{\Gb^*}(s,\chi)^{F^*})^\we} 
\atop \SS{\a \in H^1(F^*,A_{\Gb^*}(s,\chi))}}
\overline{\xi(a) \t(\a)} R_\chi(s)_{\xi,\a},$$
o\`u $\z_{s,\chi,a,\t}$ est une racine de l'unit\'e explicitement 
d\'etermin\'ee.}

\bigskip

\noindent Il est \`a noter que, comme cons\'equence des travaux effectu\'es, 
on obtient une description explicite du foncteur d'induction de 
Lusztig en termes du groupe de Weyl (lorsque $q$ est assez grand). 

Dans la section \ref{chapitre sln}, nous \'etudions plus pr\'ecis\'ement 
le cas o\`u $\Gb$ est un sous-groupe de Levi d'un groupe d\'eploy\'e 
de type $A$. Nous obtenons par exemple, en utilisant uniquement 
la th\'eorie de Harish-Chandra, un param\'etrage des caract\`eres 
irr\'eductibles de $\EC(\Gb^F,(s))$ par $\IC(\Gb,s)$ dont nous montrons 
qu'il co\"\i ncide avec le param\'etrage du th\'eor\`eme pr\'ec\'edent 
lorsque $q$ est assez grand. Cela nous permet de retrouver, comme cas 
particulier des th\'eor\`emes du chapitre \ref{chapitre a}, 
le r\'esultat de notre th\`ese \cite[th\'eor\`eme 16.2.1]{cedthese} 
sur le calcul de l'induction de Lusztig dans le groupe sp\'ecial lin\'eaire. 
Nous rappelons aussi comment fonctionne la d\'ecomposition de Jordan. 

Dans l'appendice A, nous rassemblons les r\'esultats techniques sur les 
caract\`eres de produits en couronne que nous utilisons dans 
les deux derniers chapitres. Dans l'appendice B, nous rappelons des 
r\'esultats classiques sur les sommes de Gauss et montrons 
comment ils permettent d'obtenir les valeurs des racines de 
l'unit\'e $\z_{s,\chi,a,\t}$ intervenant dans l'\'enonc\'e 
du th\'eor\`eme pr\'ec\'edent.

\vskip1cm

\noindent{\sc Remarque - } Cet article est largement inspir\'e de notre 
th\`ese \cite{cedthese}, notamment des parties qui n'ont fait l'objet 
d'aucune publication. Nous en avons cependant am\'elior\'e et enrichi 
le traitement. L'appendice est essentiellement contenu dans 
\cite[partie 1, chapitre I]{cedthese}~: il est \`a noter que 
le corollaire \ref{frobenius}, qui correspond \`a 
\cite[proposition 1.9.1]{cedthese}, 
est ici affubl\'e d'une preuve correcte, contrairement \`a 
ce qui est \'ecrit dans \cite{cedthese}~! Le chapitre \ref{chapitre lusztig} 
est une version tr\`es enrichie de \cite[\SEC 6]{cedthese}. 
Le chapitre \ref{chapitre harish-chandra} correspond \`a 
\cite[\SEC 7]{cedthese}~: remarquons que le groupe not\'e ici 
$W_{\Gb^F}'(\Lbt,\lamt)$ et not\'e $\Wba_{\Gb^F}(\Lbt,\lamt)$ dans 
\cite[\SEC 7.4]{cedthese} est d\'efini ici de mani\`ere intrins\`eque 
et non par un produit semi-direct peu canonique. Le chapitre 
\ref{chapitre gelfand} correspond \`a \cite[\SEC 12 et 13]{cedthese}~: 
ici, l'am\'elioration consiste \`a utiliser la version pr\'ecis\'ee 
du th\'eor\`eme de Digne, Lehrer et Michel sur la restriction de Lusztig des 
caract\`eres de Gelfand-Graev que l'auteur a obtenue dans 
\cite{bonnafe action}. La section \ref{chapitre sln} correspond 
\`a \cite[\SEC 15 et 16]{cedthese}~: compte tenu de la remarque 
pr\'ec\'edente, le r\'esultat sur l'induction de Lusztig est 
ici plus pr\'ecis. Il faut aussi noter que la convention 
dans le param\'etrage des caract\`eres irr\'eductibles de $\Gb^F$ ayant 
\'et\'e l\'eg\`erement modifi\'e, les formules obtenues ici se 
retrouvent all\'eg\'ees de certains signes.

\vskip1cm

\noindent{\sc Remerciements - } L'auteur tient \`a remercier tr\`es chaleureusement 
Jean Michel pour l'avoir lanc\'e dans ce sujet lors de sa th\`ese, pour 
l'avoir initi\'e \`a la th\'eorie des faisceaux-caract\`eres et pour 
les innombrables et fructueuses discussions que nous avons eues sur ce sujet 
depuis. 

\newpage

\tableofcontents

\newpage

{\Large \part{Pr\'eliminaires, notations, 
d\'efinitions\label{chapitre preliminaire}}}

\bigskip

Dans la premi\`ere section de ce chapitre, nous introduisons les notations 
et conventions g\'en\'erales valables dans tout cet article. Dans la deuxi\`eme 
section, nous introduisons les objets que nous allons \'etudier 
(groupes r\'eductifs finis) tout en \'etablissant quelques r\'esultats 
pr\'eliminaires. La troisi\`eme section est une collection de r\'esultats, 
notamment sur les centralisateurs de sous-tores de groupes r\'eductifs, 
que nous utiliserons dans la suite de l'article.

\bigskip

\section{Notations g\'en\'erales}~

\medskip

\soussection{Notations usuelles}
Nous notons $\NM=\{0,1,2,3,\dots\}$ l'ensemble des entiers naturels 
et $\NM^*=\NM\setminus\{0\}=\{1,2,3,\dots\}$ 
l'ensemble des entiers naturels non nuls. 
Comme il est d'usage, $\ZM$, $\QM$ et $\RM$ d\'esignent respectivement 
l'anneau des entiers relatifs, le corps des nombres rationnels 
et le corps des nombres r\'eels. 
Si $r$ est un nombre premier, $\ZM_r$ d\'esigne l'anneau des entiers 
$r$-adiques et nous notons $\QM_r$ son corps des fractions. 
Le corps r\'esiduel de $\ZM_r$ est not\'e $\FM_r$. Si $x \in \QM$, 
nous notons $\n_r(x) \in \ZM \cup \{+ \infty\}$ sa valuation 
$r$-adique. Si $x \not= 0$, nous d\'efinissons $x_r=r^{\n_r(x)}$ 
et $x_{r'}=x x_r^{-1}$. 
Nous notons $\ZM_{(r)}=\{x \in \QM~|~\n_r(x) \ge 0\}$. 

\bigskip

\soussection{Groupes, anneaux, corps}\label{sous groupes}
Nous fixons dans cet article un nombre premier $p$. Soit $\FM$ une cl\^oture 
alg\'ebrique du corps fini \`a $p$ \'el\'ements $\fp$. 
Si $q$ est une puissance de $p$, nous notons $\fq$ le sous-corps 
de $\FM$ de cardinal $q$. Par {\it vari\'et\'e} (ou {\it groupe 
alg\'ebrique}), nous entendons une vari\'et\'e (respectivement un groupe 
alg\'ebrique) sur $\FM$. Si $n \in \NM^*$, nous posons
$$\mub_n(\FM)=\{\x \in \FM^\times~|~\x^n=1\}.$$
Alors $|\mub_n(\FM)|=n_{p'}$.

Nous nous fixons aussi un nombre premier $\el$ diff\'erent de $p$ 
et nous notons $\qlb$ une cl\^oture alg\'ebrique du corps des nombres 
$\el$-adiques 
$\ql$. Nous fixons une fois pour toutes un automorphisme involutif 
$\qlb \to \qlb$, $x \mapsto \xba$ tel que $\bar{\o}=\o^{-1}$ pour toute 
racine de l'unit\'e $\o$ dans $\qlb^\times$. 

\medskip

Si $E$ est un ensemble et si $\sim$ est une relation d'\'equivalence 
sur $E$, on notera $E/\sim$ l'ensemble des classes d'\'equivalence de $\sim$ 
dans $E$ et $[E/\sim]$ un ensemble de repr\'esentants de ces classes 
d'\'equivalence. Le lecteur pourra v\'erifier que, chaque fois 
que cette notation sera employ\'ee (par exemple dans une somme 
$\sum_{x \in [E/\sim]} f(x)$), le r\'esultat sera ind\'ependant 
du choix des repr\'esentants. 
Si $X$ est une partie de $E$, nous noterons $1_X$ (ou $1_X^E$ s'il est n\'ecessaire 
de pr\'eciser l'ensemble de r\'ef\'erence) la fonction caract\'eristique de $X$ 
\`a valeurs dans $\qlb$. 

\medskip

Si $G$ est un groupe, nous notons $|G| \in \NM^* \cup \{+\infty\}$ son ordre. 
Si $X$ est un sous-ensemble de $G$, $<X>$ d\'esigne le sous-groupe de $G$ 
engendr\'e par $X$, $N_G(X)$ le normalisateur de $X$ dans $G$ et $C_G(X)$ le 
centralisateur de $X$ dans $G$. Si $g \in G$, nous notons $o(g)=|<g>|$ son 
ordre. Nous notons $G_\tors$ l'ensemble des \'el\'ements de $G$ d'ordre fini, 
$G_p$ l'ensemble des \'el\'ements de $G$ d'ordre fini \'egal \`a une puissance 
de $p$ et $G_{p'}$ l'ensemble des \'el\'ements de $G$ d'ordre fini premier \`a 
$p$. Si $g \in G_\tors$, nous notons $g_p$ et $g_{p'}$ les uniques 
\'el\'ements de $G_p$ et $G_{p'}$ respectivement tels que 
$g=g_p g_{p'} =g_{p'} g_p$~: $g_p$ est appel\'e la $p$-partie de $g$ tandis 
que $g_{p'}$ est appel\'e la $p'$-partie de $g$. Remarquons que $o(g_p)=o(g)_p$ 
et $o(g_{p'})=o(g)_{p'}$. 

Si $\AC$ est un groupe agissant sur $G$ et si $\ph \in \AC$, nous noterons 
$H^1(\ph,A)$ l'ensemble des classes de {\it $\ph$-conjugaison} de $G$ 
(deux \'el\'ements $g$ et $g'$ de $G$ sont dits $\ph$-conjugu\'es 
s'il existe $x \in G$ tel que $g' = x^{-1} g \ph(x)$). En d'autres 
termes, $H^1(\ph,G) = H^1(\ZM,G)$, o\`u le g\'en\'erateur $1$ de $\ZM$ 
agit sur $G$ via $\ph$. Si $\ph$ est d'ordre fini, on a en g\'en\'eral 
$H^1(\ph,G) \not= H^1(<\ph>,G)$. 

\medskip

\exemple{H1 abelien}
Supposons $G$ ab\'elien. Alors $H^1(\ph,G) = G/\Im(\ph-1)$ o\`u on 
note $\ph-1 : G \to G$, $g \mapsto g^{-1} \ph(g)$~: en particulier, 
$H^1(\ph,G)$ h\'erite naturellement d'une structure de groupe. 
Si de plus $G$ est fini, alors $|A^\ph|=|H^1(\ph,G)|$.\finl

\medskip

Nous fixons une fois pour toutes un isomorphisme de groupes 
$$\imath : (\QM/\ZM)_{p'} \longto \FM^\times$$
et un morphisme injectif de groupes 
$$\jmath : \QM/\ZM \longto \qlb^\times.$$
On obtient alors un morphisme injectif 
$$\kappa : \FM^\times \longto \qlb^\times$$ 
d\'efini par $\kappa=\jmath \compose \imath^{-1}$. 
Pour finir ce paragraphe, nous d\'efinissons le morphisme surjectif 
$\tilde{\imath} : \QM \to \FM^\times$ 
comme \'etant la composition de $\imath$ avec le morphisme 
$\QM \to \QM/\ZM \to (\QM/\ZM)_{p'}$. Notons que 
$\Ker \tilde{\imath} =\ZM[1/p]=\{a p^r~|~a \in \ZM$ et $r \in \ZM\}$. 
De m\^eme, nous notons $\tilde{\jmath} : \QM \to \qlb^\times$ le compos\'e de 
$\jmath$ et du morphisme canonique $\QM \to \QM/\ZM$~; on a 
$\Ker \tilde{\jmath}=\ZM$. 

\bigskip

\soussection{Caract\`eres des groupes finis\label{car sous}}
Si $G$ est un groupe fini, nous notons $\Irr G$ l'ensemble de ses caract\`eres 
irr\'eductibles sur $\qlb$ et $G^\wedge$ le groupe de ses caract\`eres 
lin\'eaires \`a valeurs dans ${\overline{\QM}}_\ell^\times$. 
On a $G^\we \incl \Irr G$~; 
on a $G^\we=\Irr G$ si et seulement si $G$ est ab\'elien. Si $f : G \to H$ est 
un morphisme de groupes finis, nous noterons $\fha : H^\we \to G^\we$, 
$\th \mapsto \th \circ f$ le morphisme dual de $f$.

Si $G$ est un sous-groupe 
distingu\'e d'un groupe $\Gti$ et si $\phi \in \Gti$, nous noterons 
$\Cent(G\phi)$ le $\qlb$-espace vectoriel des fonctions 
$G\phi \to \qlb$ invariantes par $G$-conjugaison. 
Nous d\'efinissons sur $\Cent(G\phi)$ le produit scalaire
$$\fonction{\langle , \rangle_G}{\Cent(G\phi) \times 
\Cent(G\phi)}{\qlb}{(\g,\g')}{\DS{\frac{1}{|G|} \sum_{g \in G} \g(g\phi) 
\overline{\g'(g\phi).}}}$$
Si $H$ est un sous-groupe de $G$ et si $g \in G$ est tel que $\lexp{g\phi}{H}=H$, 
nous noterons $\Ind_{Hg\phi}^{G\phi} : \Cent(Hg\phi) \to \Cent(G\phi)$ 
l'adjoint, pour les produits scalaires $\langle,\rangle_{Hg\phi}$ et 
$\langle,\rangle_{G\phi}$, de l'application de restriction 
naturelle $\Res_{Hg\phi}^{G\phi} : \Cent(G\phi) \to \Cent(Hg\phi)$. 
En fait, si $f \in \Cent(Hg\phi)$, on a 
\equat\label{induit tordu formule}
(\Inf_{Hg\phi}^{G\phi} f)(x\phi)=
\sum_{\SS{y \in [G/H]} \atop \SS{y^{-1}x \phi y \in Hg\phi}} f(y^{-1}x \phi y).
\endequat
On en d\'eduit que  
$$\Ind_{Hg\phi}^{G\phi}=\Res_{G\phi}^{G<\phi>} \circ \Ind_{H <g\phi>}^{G<\phi>} 
\fti,$$
o\`u $\fti$ est l'extension par z\'ero de $f$ \`a $H<g\phi>$. 
Avec ces notations, $\Cent(G)$ est le $\qlb$-espace vectoriel des fonctions 
centrales $G \to \qlb$ et $\Irr G$ en est une base orthonormale. 
Nous identifions $\ZM \Irr G$ avec le groupe de Grothendieck de la cat\'egorie des 
$\qlb G$-modules de type fini.

Si $g \in G$, nous noterons $\g_g^G$ la fonction centrale sur $G$ 
d\'efinie par
$$\g_g^G(g')=\begin{cases}
             0 & {\text{si $g$ et $g'$ ne sont pas conjugu\'es dans $G$,}} \\
	     |C_G(g)| & {\text{sinon.}}
	     \end{cases}$$
En fait, 
$$\g_g^G=\sum_{\g \in \Irr G} \overline{\g(g)} \g.$$
De plus, si $g \in H$, alors
$$\Ind_H^G \g_g^H = \g_g^G.$$

Pour finir cette sous-section, nous allons donner une formule permettant 
de calculer l'induction $\Ind_{Hg\phi}^{G\phi}$ dans un cas particulier. 
Supposons maintenant que $\Gti=G=H<g>$ et $\phi=1$. En particulier, $H$ est 
distingu\'e 
dans $G$. Pour tout \car \irr $\chi$ de $H$ invariant par $G$ (\cad par $g$), 
on fixe une extension $\chit$ de $\chi$ \`a $G$ (l'existence de $\chit$ est 
assur\'ee par la cyclicit\'e de $G/H$). On note $\chit_g$ la restriction de 
$\chit$ \`a $Hg$. Alors $(\chit_g)_{\chi \in (\Irr H)^g}$ 
est une base orthonormale de $\Cent(Hg)$ et
\equat\label{une formule tordue}
\Ind_{Hg}^G \chit_g = \sum_{\xi \in  (G/H)^\we} \xi(g)^{-1} (\chit \otimes \xi).
\endequat
Ici, un \car lin\'eaire de $G/H$ est aussi vu comme 
un \car lin\'eaire de $G$. 

\bigskip

\noindent{\sc D\'emonstration de \ref{une formule tordue} - } 
Posons $\g=\Ind_{Hg}^G \chit_g$ et 
$\g'=\sum_{\xi \in  (G/H)^\we} \xi(g)^{-1} (\chit \otimes \xi)$. 
D'apr\`es \ref{induit tordu formule}, on a, pour $x \in G$, 
\eqna
\g(x)&=&\sum_{\SS{y \in [G/H]} \atop \SS{y^{-1}xy \in Hg}} 
\chit(y^{-1}xy) \\
&=&\sum_{\SS{y \in [G/H]} \atop \SS{y^{-1}xy \in Hg}} \chit(x) \\
&=& \begin{cases}
    |G/H| \chit(x) & \text{si } x \in Hg \\
    0 & \text{sinon.}
    \end{cases}
\endeqna
D'autre part,
\eqna
\g'(x) &=& \chit(x) \DS{(\sum_{\xi \in (G/H)^\we} \xi(g)^{-1}\xi(x))} \\
&=& \begin{cases}
    |G/H| \chit(x) & \text{si } x \in Hg \\
    0 & \text{sinon.}~\SS{\blacksquare}
    \end{cases}
\endeqna

\bigskip

\soussection{Groupes alg\'ebriques\label{sous groupes algebriques}}
Si $\Hb$ est un groupe alg\'ebrique lin\'eaire, nous notons $\Hb^\ci$ 
sa composante connexe contenant l'\ele neutre, $\Hb_\uni$ sa sous-vari\'et\'e 
ferm\'ee form\'ee des \eles unipotents, $\Hb_\sem$ l'ensemble de 
ses \eles semi-simples, $\Zb(\Hb)$ son centre, $\Db(\Hb)$ 
son groupe d\'eriv\'e et $\Rb_u(\Hb)$ son radical unipotent. 
Nous posons $\ZC(\Hb)=\Zb(\Hb)/\Zb(\Hb)^\circ$. 
Nous notons $\rang(\Hb)$ le rang de $\Hb$ \cad la dimension d'un de ses tores 
maximaux. Nous posons $\rang_\sem(\Hb)=\rang(\Db(\Hb))$. 
Si $h \in \Hb$, nous posons $C_\Hb^\circ(h)=C_\Hb(h)^\ci$ et 
$A_\Hb(h)=C_\Hb(h)/C_\Hb^\ci(h)$. Remarquons que $\Hb_\sem=\Hb_{p'}$ et que 
$\Hb_\uni=\Hb_p$. 

Nous appellerons {\it compl\'ement de Levi} de $\Hb$ tout 
sous-groupe $\Lb$ de $\Hb$ tel que $\Hb=\Lb \ltimes \Rb_u(\Hb)$ 
(il est \`a noter qu'il n'existe pas toujours de compl\'ement de 
Levi). Nous appellerons {\it \levi} de $\Hb$ tout compl\'ement de 
Levi d'un \para de $\Hb$.

Si $F : \Hb \to \Hb$ est une isog\'enie dont une puissance $F^\d$ est 
un endomorphisme de Frobenius pour une structure rationnelle sur $\Hb$, on a, 
d'apr\`es le th\'eor\`eme de Lang, $H^1(F,\Hb)=H^1(F,\Hb/\Hb^\ci)$. 
D'autre part, si $h \in \Hb^F$, nous noterons, pour all\'eger les notations 
lorsqu'il n'y aura pas d'ambigu\"\i t\'e sur $F$, $\g_h^\Hb$ 
la fonction centrale $\g_h^{\Hb^F}$. 

\bigskip

\soussection{Caract\`eres rationnels et sous-groupes 
\`a un param\`etre}\label{sous un} 
Soit $\Hb$ un groupe alg\'ebrique lin\'eaire. Nous notons $X(\Hb)$ le groupe 
(ab\'elien, not\'e additivement) des caract\`eres rationnels $\Hb \to \FM^\times$ et 
$Y(\Hb)$ l'ensemble des sous-groupes \`a un param\`etre $\FM^\times \to \Hb$. 
On a bien s\^ur $Y(\Hb)=Y(\Hb^\circ)$. Si $\pi : \Hb \to \Hb'$ 
est un morphisme de groupes alg\'ebriques, nous posons 
$\pi_X : X(\Hb') \to X(\Hb)$, $x \mapsto x \compose \pi$ 
et $\pi_Y : Y(\Hb) \to Y(\Hb')$, $y \mapsto \pi \compose y$. 
Remarquons que $\pi_X$ est un morphisme de groupes. Si $F : \Hb \to \Hb$ est 
l'isog\'enie pr\'ec\'edente, les applications $F_X$ et $F_Y$ seront not\'ees 
par la m\^eme lettre $F$.
 
Si $\Hb$ est commutatif, alors $Y(\Hb)$ est un groupe ab\'elien (que nous 
noterons additivement)~: c'est un $\ZM$-module libre de rang 
$\rang(\Hb)$ et nous d\'efinissons alors une forme bilin\'eaire 
$$< , >_\Hb : X(\Hb) \times Y(\Hb) \longto \ZM$$ 
par la condition 
$$x(y(\x))=\x^{<x,y>_\Hb}$$
pour tous $x \in X(\Hb)$, $y \in Y(\Hb)$ et 
$\x \in \FM^\times$. Cette forme bilin\'eaire 
est \'etendue par lin\'earit\'e en une forme bilin\'eaire 
$$(X(\Hb) \otimes_\ZM \QM) \times (Y(\Hb) \otimes_\ZM \QM) \longto \QM$$
que l'on notera encore $< , >_\Hb$ par abus de notation. Cette derni\`ere 
forme bilin\'eaire est non d\'eg\'en\'er\'ee.

Si $\pi : \Hb \to \Hb'$ est un morphisme 
de groupes alg\'ebriques commutatifs, alors $\pi_Y$ est un morphisme 
de groupes et $\pi_X$ et $\pi_Y$ 
sont adjoints par rapport \`a $< , >_\Hb$ et $< , >_{\Hb'}$. En d'autres 
termes, 
$$<\pi_X(x) , y >_{\Hb} = < x, \pi_Y(y) >_{\Hb'}$$ 
pour tous $x \in X(\Hb')$ et $y \in Y(\Hb)$. 
Nous d\'efinissons par ailleurs le morphisme 
$$\fonction{\tilde{\imath}_\Hb}{Y(\Hb) 
\otimes_\ZM \QM}{\Hb}{y \otimes_\ZM r}{y(\tilde{\imath}(r)).}$$
Soient $x \in X(\Hb)$ et $y \in Y(\Hb) \otimes_\ZM \QM$. Alors

\equat\label{scalaire xy}
x(\tilde{\imath}_\Hb(y))=\tilde{\imath}(<x , y >_\Hb).
\endequat

\bigskip

\noindent{\it D\'emonstration de \ref{scalaire xy}.} 
\'Ecrivons $y=y' \otimes_\ZM r$ avec $y' \in Y(\Hb)$ et $r \in \QM$. Alors 
\eqna
x(\tilde{\imath}_\Hb(y))&=& x(y'(\tilde{\imath}(r))) \\
&=&\tilde{\imath}(r)^{<x,y'>_\Hb} \\
&=&\tilde{\imath}(r <x,y'>_\Hb) \\
&=&\tilde{\imath}(<x , y >_\Hb), \\
\endeqna
ce qui est le r\'esultat annonc\'e.\fin

\bigskip

\noindent{\sc Remarque -} Si $\Hb$ est un tore, 
alors $<,>_\Hb$ est une dualit\'e parfaite et $\tilde{\imath}_\Hb$ 
est surjective (en effet, l'application $Y(\Hb) \otimes_\ZM \FM^\times$, 
$y \otimes \x \mapsto y(\x)$ est un isomorphisme de groupes).\finl

\bigskip

Si $A$ est un groupe agissant sur le groupe commutatif $\Hb$, nous 
d\'efinissons une action de $A$ sur les $\ZM$-modules $X(\Hb)$ et 
$Y(\Hb)$ par les formules suivantes~:
$$\fonctio{A \times X(\Hb)}{X(\Hb)}{(\s,x)}{\s_X^{-1}(x)=x \compose \s^{-1}}$$
$$\fonctio{A \times Y(\Hb)}{Y(\Hb)}{(\s,y)}{\s_Y(y)=\s \compose y.}
\leqno{\mathrm{et}}$$
Il est alors facile de v\'erifier que 
$$< \s(x),\s(y) >_{\Hb}=<x,y>_\Hb$$
pour tous $x \in X(\Hb)$, $y \in Y(\Hb)$ et $\s \in A$. 

\bigskip

\soussection{Groupes diagonalisables} 
Nous terminons cette section en rappelant quelques faits \'el\'ementaires 
sur les groupes diagonalisables. 
Premi\`erement, si $\Db'$ est un sous-groupe ferm\'e d'un groupe 
diagonalisable $\Db$ (en particulier, $\Db'$ est aussi un groupe 
diagonalisable), alors la suite de $\ZM$-modules 
\equat\label{X exact}
0 \longto X(\Db/\Db') \longto X(\Db) \longto X(\Db') \longto 0
\endequat
induite par l'inclusion $\Db' \injto \Db$ est exacte. Si $X'$ est un 
sous-groupe de $X(\Db)$ tel que 
$$\Db'=\{d \in \Db~|~\forall~\chi \in X',~\chi(d)=1\},\leqno{(*)}$$
alors l'application naturelle $X(\Db) \to X(\Db')$ induit un 
isomorphisme de groupes 
\equat\label{X/X'}
X(\Db') \simeq (X(\Db)/X')/(X(\Db)/X')_p.
\endequat
Puisque $X(\Db'/\Db^{\prime \circ})\simeq X(\Db')_\tors$, on d\'eduit 
de \ref{X/X'} un isomorphisme de groupes ab\'eliens finis 
\equat\label{xd'}
X(\Db'/\Db^{\prime \circ}) \simeq (X(\Db)/X')_{p'}.
\endequat

Nous supposons maintenant, et ce jusqu'\`a la fin de cette sous-section, 
que $\Db$ est {\it connexe} (\cad que $\Db$ est un tore) et que $\Db'$ 
est {\it fini}. Tout d'abord, l'application 
\equat\label{X irr}
\fonctio{X(\Db')}{\Db^{\prime \we}}{x}{\kappa \ci x}
\endequat
est un isomorphisme de groupes ab\'eliens finis. 
L'isomorphisme \ref{xd'} montre que $X'$ est 
d'indice fini dans $X(\Db)$. Posons maintenant 
$$Y'=\{y \in Y(\Db) \otimes_\ZM \QM~|~\forall~x \in X',~<x,y>_\Db \in \ZM\}.$$
Alors $Y'$ contient $Y(\Db)$ et 
la restriction de $\tilde{\imath}_\Db$ \`a $Y'$ a pour image 
$\Db'$ et induit un isomorphisme de groupes ab\'eliens finis

\medskip 
 
\equat\label{Y/Y'}
(Y'/Y(\Db))_{p'} \simeq \Db'.
\endequat

\bigskip

\noindent{\sc D\'emonstration de \ref{Y/Y'} - } D'apr\`es $(*)$ et \ref{scalaire xy}, 
on a $\tilde{\imath}_\Db(Y') \incl \Db'$. D'autre part, d'apr\`es 
\cite[\SEC 4, ${\text{n}}^\circ$ 8]{bourbaki algebre}, 
l'application 
$$\fonctio{X(\Db)/X' \times
Y'/Y(\Db)}{\qlb^\times}{(x+X',
y+Y(\Db))}{\tilde{\jmath}(<x,y>_\Db)}$$
est une dualit\'e parfaite. Cela montre que $\tilde{\imath}_\Db(Y')=\Db'$ 
et que l'on a bien un isomorphisme $(Y'/Y(\Db))_{p'} \simeq \Db'$ 
(toujours gr\^ace \`a \ref{scalaire xy}).\fin

\bigskip

\lemme{divisibilite}
{\it Soit $y \in Y(\Db)$ et soit $n \in \ZM$, premier \`a $p$, tels que 
$y(\tilde{\imath}(1/n))=1$. Alors il existe $y_0 \in Y(\Db)$ tel que 
$y=ny_0$.}

\bigskip

\proof On a $|\mub_n(\FM)|=n$ car $n$ est premier \`a $p$. 
De plus, $\mub_n(\FM) \incl \Ker y$ par hypoth\`ese. 
Donc $y$ induit un \mor de groupes alg\'ebriques 
$\bar{y} : \FM^\times/\mub_n(\FM) \to \Db$. 

Mais l'application $\FM^\times \to \FM^\times$, $\x \mapsto \x^n$ est s\'eparable car 
$p$ ne divise pas $n$. Donc elle induit un isomorphisme de groupes 
alg\'ebriques $\a : \FM^\times/\mub_n(\FM) \to \FM^\times$. Soit 
$y_0=\bar{y} \ci \a^{-1}$. 
Alors $y_0 \in Y(\Db)$ et $y=ny_0$ par construction.\fin

\bigskip

\proposition{curieux}
{\it Soient $\Tb$ et $\Tb'$ deux tores et soit $\pi : \Tb \to \Tb'$ un \mor 
de groupes alg\'ebriques de noyau fini. Alors le \mor 
$\pi_Y : Y(\Tb) \to Y(\Tb')$ est injectif et on a un isomorphisme naturel 
$$(Y(\Tb')/\Im \pi_Y)_{p'} \simeq \Ker \pi.$$}

\bigskip

\proof Soit $y \in \Ker\pi_Y$. Alors l'image de $y$ est contenue dans $\Ker \pi$ 
mais est aussi connexe et contient $1$. 
Donc $y=0$ car $\Ker \pi$ est fini~: l'injectivit\'e de $\pi_Y$ est prouv\'ee.

Notons $Y'$ le sous-groupe de $Y(\Tb')$ fom\'e des \'el\'ements $y' \in Y(\Tb')$ 
tels que $ny' \in \Im \pi_Y$ pour un $n \in \ZM$ premier \`a $p$. 
Alors $\Im \pi_Y \incl Y'$ et, par construction, 
$Y'/\Im \pi_Y=(Y(\Tb')/\Im \pi_Y)_{p'}$. 

Soit $y' \in Y'$ et soit $n \in \ZM$, non 
divisible par $p$, tels que $ny' \in \Im \pi_Y$. Soit $y$ l'unique 
\'el\'ement de $Y(\Tb)$ tel que $\pi_Y(y)=ny'$. On pose 
$$\s(y')=y(\tilde{\imath}(\frac{1}{n})) \in \Tb.$$
Alors $\pi(\s(y'))=(ny')(\tilde{\imath}(1/n))=1$ donc $\s(y') \in \Ker \pi$. 
Il est facile de v\'erifier que l'application 
$$\s : Y' \longto \Ker \pi$$
est bien d\'efinie et est un morphisme de groupes.

Il est aussi facile de voir que $\Im \pi_Y \incl \Ker \s$. 
R\'eciproquement, soit $y' \in \Ker \s$. Soient $n \in \ZM$, premier \`a $p$, 
et $y \in Y(\Tb)$ tels que $ny'=\pi_Y(y)$. Alors 
$$y(\imath^{-1}(\frac{1}{n}))=1.$$
D'apr\`es le lemme \ref{divisibilite}, il existe $y_0 \in Y(\Tb)$ tel que $y=ny_0$. 
En particulier, $y'=\pi_Y(y_0)$ donc $\Ker \s \incl \Im \pi_Y$. cela montre 
que $\Ker \s=\Im \pi_Y$.

Il reste \`a prouver que $\s$ est surjectif. Soit $t \in \Ker \pi$ et soit $n$ l'ordre 
de $t$. Alors $n$ est premier \`a $p$ et, puisque $\Tb$ est connexe, 
il existe $y \in Y(\Tb)$ tel que $t=y(\imath^{-1}(1/n))$ (voir par exemple 
\cite[Proposition 0.20]{dmbook}). 
Alors $\pi_Y(y)(\imath^{-1}(1/n)) = 1$ donc, d'apr\`es le lemme \ref{divisibilite}, 
il existe $y' \in Y(\Tb')$ tel que $ny'=\pi_Y(y)$. Alors 
$y' \in Y'$ et $\s(y')=t$ par construction.\fin

\bigskip

\soussection{Dualit\'e entre tores} 
Soient $\Tb$ et $\Tb^*$ deux tores. On suppose qu'ils sont munis respectivement 
d'isog\'enies $F$ et $F^*$ dont une puissance est un endomorphisme 
de Frobenius. Nous dirons que $(\Tb,F)$ et $(\Tb^*,F^*)$ sont 
{\it duaux} (ou simplement que $\Tb$ et $\Tb^*$ sont {\it duaux}) 
s'il existe un isomorphisme $\nu : X(\Tb) \longmapright{\sim} Y(\Tb^*)$ tel que 
$F^* \circ \nu = \nu \circ F$. Bien s\^ur, la relation de dualit\'e 
est sym\'etrique car, en utilisant les dualit\'es parfaites donn\'ees 
par $<,>_\Tb$ et $<,>_{\Tb^*}$, le morphisme adjoint 
$\nu^* : X(\Tb^*) \to Y(\Tb)$ de $\nu$ est un isomorphisme. 

\medskip

Supposons donc que $(\Tb,F)$ et $(\Tb^*,F^*)$ sont duaux et soit 
$\nu : X(\Tb) \longmapright{\sim} Y(\Tb^*)$ un isomorphisme tel que 
$F^* \circ \nu = \nu \circ F$. Nous identifierons $X(\Tb)$ et $Y(\Tb^*)$ 
via $\nu$ et nous identifierons $X(\Tb^*)$ et $Y(\Tb)$ via $\nu^*$. On a 
une suite exacte
$$1 \longto \Tb^F \longto \Tb \longmapright{F-1} \Tb \longto 1.$$
De plus, le morphisme $F-1 : \Tb \to \Tb$ est \'etale, donc induit 
un isomorphisme entre $\Tb$ et $\Tb/\Tb^F$. 
D'apr\`es \ref{X exact}, on obtient donc une suite exacte 
\equat\label{X F}
0 \longto X(\Tb) \longmapright{F-1} X(\Tb) \longto X(\Tb^F) \longto 0.
\endequat
Par dualit\'e, on obtient une suite exacte
\equat\label{Y X F}
0 \longto Y(\Tb^*) \longmapright{F^*-1} Y(\Tb^*) \longto X(\Tb^F) \longto 0.
\endequat
Soit maintenant $n$ un entier naturel non nul tel que $F^n$ soit un endomorphisme 
de Frobenius d\'eploy\'e pour une structure sur un corps fini 
\`a $q$ \'el\'ements. L'application $\Tb^{*F^{*n}} \to \Tb^{*F^*}$, 
$t \mapsto N_{F^{*n}/F^*}(t)$ est surjective. De plus l'application 
$Y(\Tb^*) \to \Tb^{*F^{*n}}$, $y \mapsto y(\tilde{\imath}(1/(q-1)))$ est 
surjective. Il est alors facile de v\'erifier que l'on a une suite exacte
\equat\label{Y F}
0 \longto Y(\Tb^*) \longmapright{F^*-1} Y(\Tb^*) \longmapright{f}
\Tb^{*F^*} \longto 0,
\endequat
o\`u $f : Y(\Tb^*) \to \Tb^{*F^*}$, 
$y \mapsto N_{F^{*n}/F^*}(y(\tilde{\imath}(1/(q-1))))$. Il est \`a noter 
que l'application $f$ ne d\'epend pas du choix de $n$. La comparaison 
des suites exactes \ref{Y X F} et \ref{Y F} et l'isomorphisme \ref{X irr} 
fournit un isomorphisme $\Tb^{*F^*} \simeq (\Tb^F)^\we$. Cet 
isomorphisme ne d\'epend que du choix de $\imath$ et $\jmath$. 

Nous allons l'expliciter. Soit $s \in \Tb^{*F^*}$. 
Notons $\sha$ le caract\`ere lin\'eaire de $\Tb^F$ d\'efini par $s$ via 
l'isomorphisme pr\'ec\'edent. Pour calculer $\sha$, il faut tout 
d'abord trouver un \'el\'ement $y \in Y(\Tb^*)\simeq X(\Tb)$ 
tel que $s=N_{F^{*n}/F^*}(y(\tilde{\imath}(1/(q-1))))$. Alors 
\equat\label{egalite sha}
\sha=\kappa \circ \Res_{\Tb^F}^\Tb y.
\endequat

\bigskip

\section{Le contexte}~

\medskip

\soussection{Le probl\`eme}
Nous nous int\'eressons dans cet article \`a la th\'eorie des caract\`eres 
d'un groupe fini de la forme 
$\Gb^F$ (param\'etrage des caract\`eres, table de caract\`eres, 
th\'eorie de Deligne-Lusztig, th\'eorie de Harish-Chandra, conjecture 
de Lusztig sur les faisceaux-caract\`eres...), o\`u $\Gb$ est un groupe 
r\'eductif connexe et $F : \Gb \to \Gb$ est une isog\'enie telle que $F^\d$ est 
l'endomorphisme de Frobenius de $\Gb$ relatif \`a une structure 
rationnelle sur $\Gb$.

Nous nous concentrons plus particuli\`erement sur les probl\`emes reli\'es \`a la 
non connexit\'e du centre de $\Gb$. En d'autres termes, nous essayons 
de r\'esoudre les questions concernant les groupes \`a centre non connexe 
en supposant que la m\^eme question est r\'esolue pour les groupes \`a 
centre connexe. 
La strat\'egie habituelle est la suivante. Il est possible \cite{delu} 
de construire un groupe r\'eductif connexe $\Gbt$ (muni lui aussi 
d'une isog\'enie encore not\'ee $F$) dont $\Gb$ est un sous-groupe 
ferm\'e $F$-stable contenant $\Db(\Gbt)$. L'\'etude pr\'ecise du foncteur 
de restriction $\Res_{\Gb^F}^{\Gbt^F}$ fournit alors des \'el\'ements 
de r\'eponse (th\'eorie de Clifford).

\bigskip

\noindent{\sc Remarque - } Il n'est pas d\'eraisonnable de supposer 
que beaucoup de choses sont connues pour les groupes \`a centre connexe. 
Par exemple, la conjecture de Lusztig sur les faisceaux-caract\`eres 
a \'et\'e r\'esolue par T. Shoji \cite{shoji almost}. 

C'est d'autant moins d\'eraisonnable que notre but est d'\'etudier 
le groupe sp\'ecial lin\'eaire. En effet, ce groupe est inclus 
dans le groupe g\'en\'eral lin\'eaire et pratiquement tout 
ce qui concerne la table de caract\`eres de ce dernier 
est connu (aussi bien du point de vue \'el\'ementaire de J.A. Green \cite{green}, 
que du point de vue de la th\'eorie de Deligne-Lusztig \cite{lusr}, 
voire m\^eme du point de vue de la th\'eorie des faisceaux-caract\`eres~: 
le lien entre ces trois th\'eories est lui aussi bien compris).\finl 

\bigskip

\soussection{Plongements}
Comme expliqu\'e ci-dessus, l'un des buts de cet article est 
d'\'etudier les foncteurs de restriction entre groupes 
de m\^eme type. Pour cela, nous nous fixons 
un groupe r\'eductif connexe $\Gbt$ muni d'une isog\'enie 
$F : \Gbt \to \Gbt$ telle que $F^\d$ est un endomorphisme 
de Frobenius de $\Gb$ relatif \`a une structure sur le corps 
fini $\fq$ (ici, $\d$ est un entier naturel non nul et $q$ 
est une puissance de $p$ fix\'es une fois pour toutes~: bien 
qu'ils ne soient pas uniquement d\'etermin\'es par la donn\'ee 
de $(\Gb,F)$, le nombre r\'eel positif $q^{1/\d}$ l'est). 

Nous nous fixons aussi un sous-groupe ferm\'e connexe $F$-stable 
$\Gb$ de $\Gbt$. Tout au long de cet article, nous supposerons 
que les hypoth\`eses suivantes sont satisfaites~:

\begin{quotation}
{\it \noindent 
\begin{itemize}
\itemth{1} Le centre de $\Gbt$ est connexe~;

\itemth{2} Le groupe $\Gb$ contient le groupe d\'eriv\'e de $\Gbt$.
\end{itemize}}
\end{quotation}

\bigskip

\remarque{tous les groupes} Puisque tout groupe r\'eductif peut-\^etre 
plong\'e dans un groupe \`a centre connexe de m\^eme type \cite{delu}, les 
r\'esultats que nous allons d\'emontrer concernant le groupe $\Gb$ 
seront vrais pour tous les groupes r\'eductifs connexes.\finl

\bigskip

Il r\'esulte de ces hypoth\`eses que $\Zb(\Gb)=\Zb(\Gbt) \cap \Gb$ et 
$\Gbt=\Gb.\Zb(\Gbt)$. De plus, $\Db(\Gbt)=\Db(\Gb)$. 
Nous notons
$$i : \Gb \injto \Gbt$$
l'inclusion canonique. 

Nous fixons aussi dans cet article un \borel $F$-stable $\Bbt_0$ de $\Gbt$ ainsi 
qu'un \tor $F$-stable $\Tbt_0$ de $\Bbt_0$. Nous notons $\Ub_0$ le radical unipotent 
de $\Bbt_0$. On pose
$$\Bb_0=\Bbt_0 \cap \Gb \qquad\qquad{\mathrm{et}}\qquad\qquad
\Tb_0=\Tbt_0 \cap \Gb.$$
Alors $\Bb_0$ est un \borel $F$-stable de $\Gb$, $\Tb_0$ est un \tor 
$F$-stable de $\Bb_0$ et $\Ub_0$ est le radical unipotent 
de $\Bb_0$. 

\bigskip

\soussection{Syst\`eme de racines\label{soussection racines}} 
Nous notons $W_0$ le groupe de Weyl de $\Gb$ relativement \`a $\Tb_0$~;
remarquons que $W_0$ est canoniquement isomorphe au groupe de Weyl de 
$\Gbt$ relativement \`a $\Tbt_0$. 
Nous notons $\Phi_0$ (\resp $\Phit_0$) le syst\`eme de racines 
de $\Gb$ (\resp $\Gbt$) relativement \`a $\Tb_0$ 
(\resp $\Tbt_0$). Le morphisme  
$i_X : X(\Tbt_0) \injto X(\Tb_0)$ associ\'e \`a $i$ 
induit une bijection entre $\Phi_0$ et $\Phit_0$.
Nous notons $\D_0$ (\resp $\Delt_0$) 
la base de $\Phi_0$ (\resp $\Phit_0$) associ\'ee \`a $\Bb_0$ 
(\resp $\Bbt_0$). Si $\a \in \Phi_0$, nous notons 
$\Ub_\a$ le sous-groupe unipotent de dimension $1$ de $\Gb$ normalis\'e 
par $\Tb_0$ et associ\'e \`a $\a$. 

Soit $\phi_0 : X(\Tb_0) \otimes_\ZM \RM \to X(\Tb_0) \otimes_\ZM \RM$ 
l'automorphisme d'ordre fini \'egal \`a $q^{-1/\d} F$. Puisqu'il est d'ordre 
fini, on a $\det \phi_0 \in \{1,-1\}$. Nous poserons 
$$\e_\Gb=\det \phi_0 \qquad{\text{et}}\qquad \eta_\Gb=\e_{\Db(\Gb)}.$$
D'autre part, $\phi_0$ normalise $W_0$ et induit sur $W_0$ le 
m\^eme automorphisme que celui induit par $F$. Nous noterons 
aussi $\phi_0 : Y(\Tb_0) \otimes_\ZM \RM \to Y(\Tb_0) \otimes_\ZM \RM$ 
l'automorphisme d'ordre fini \'egal \`a $q^{-1/\d} F$. Pour finir, 
nous noterons $\phit_0 : \Phi_0 \to \Phi_0$ la bijection telle 
que, pour toute racine $\a \in \Phi_0$, il existe un entier naturel $\d_\a$ 
tel que $F(\a)= p^{\d_\a} \phit_0(\a)$ (si $\o$ est une orbite 
sous l'action de $\phit_0$, alors $\sum_{\a \in \o} \d_\a > 0$). 
Bien s\^ur, $\phit_0$ stabilise $\D_0$ et $\Phi_0^+$. 

\medskip

Si $I$ est une partie de $\D_0$, nous noterons $<\Phi_I>$ le 
sous-syst\`eme de $\Phi$ de base $I$, $W_I$ le groupe de Weyl 
de $\Phi_I$, $\Pb_I$ le \para $\Bb_0 W_I \Bb_0$ de $\Gb$, $\Ub_I$ son 
radical unipotent et $\Lb_I$ 
le \levic de $\Pb_I$ contenant $\Tb_0$. Alors $\Pb_I$ (ou $\Lb_I$) 
est $F$-stable \ssi $\phit_0(I)=I$. 

\bigskip

\soussection{Dualit\'e} Nous fixons un triplet $(\Gbt^*,\Tbt_0^*,F^*)$ dual de  
$(\Gbt,\Tbt_0,F)$ au sens de \cite[d\'efinition 13.10]{dmbook}. 
Nous fixons aussi un triplet $(\Gb^*,\Tb_0^*,F^*)$ 
dual de $(\Gb,\Tb_0,F)$. Le morphisme $i$ induit un morphisme 
$i^* : \Gbt^* \to \Gb^*$ commutant avec $F^*$ 
et tel que $i^*(\Tbt_0^*)=\Tb_0^*$. Il faut cependant 
faire attention~: $i^*$ n'est pas uniquement d\'etermin\'e par $i$. 
On peut le composer avec n'importe quel automorphisme int\'erieur 
induit par un \'el\'ement $F^*$-stable du tore maximal. 
Notons que $i^*$ est surjectif.

\bigskip

\soussection{Un r\'esultat \`a la Borel-Tits} 
Borel et Tits ont montr\'e, pour un groupe r\'eductif d\'efini 
sur un corps quelconque $K$, que tout $K$-\levic 
d'un $K$-\para est le centralisateur d'un $K$-tore d\'eploy\'e 
\cite[th\'eor\`eme 4.15]{borel tits}.
 
Ici, $F$ ne d\'efinit pas forc\'ement une structure sur un corps 
fini, mais il est tout de m\^eme possible de donner une caract\'erisation 
similaire des \levics $F$-stables de \paras $F$-stables. 
Un \levi $F$-stable de $\Gb$ est dit {\it $\Gb$-d\'eploy\'e} 
(ou {\it $(\Gb,F)$-d\'eploy\'e} s'il peut y avoir ambigu\"\i t\'e 
sur l'isog\'enie) s'il existe un \para $F$-stable de $\Gb$ dont c'est 
un compl\'ement de Levi.

\smallskip

Soit $\Tb$ un \tor $F$-stable de $\Gb$ et soit $\Lb$ un \levi $F$-stable 
de $\Gb$ contenant $\Tb$. On note $\Phi$ et $\Phi_\Lb$ les 
syst\`emes de racines respectifs de $\Gb$ et $\Lb$ relativement \`a $\Tb$. 

\bigskip

\proposition{prop borel tits}
{\it Avec les notations ci-dessus, les assertions suivantes sont 
\'equivalentes~:
\begin{itemize}
\itemth{1} $\Lb$ est $\Gb$-d\'eploy\'e.

\itemth{2} On a 
$$\Phi_\Lb=\{\a \in \Phi~|~
\forall v \in \Ker(F-q^{1/\d},Y(\Zb(\Lb)^\circ)\otimes_\ZM \QM(q^{1/\d})),
~<\a,v>_\Tb=0\}.$$

\itemth{3} Il existe un sous-$\QM(q^{1/\d})$-espace vectoriel $E$ de 
$\Ker(F-q^{1/\d},Y(\Tb) \otimes_\ZM \QM(q^{1/\d}))$ tel que 
$$\Phi_\Lb=\{\a \in \Phi~|~\forall v \in E,~ <\a,v>_\Tb=0\}.$$

\itemth{4} Il existe $v \in \Ker(F-q^{1/\d},Y(\Tb) \otimes_\ZM \QM(q^{1/\d}))$ tel que 
$$\Phi_\Lb=\{\a \in \Phi~|~<\a,v>_\Tb=0\}.$$
\end{itemize}}

\bigskip

\noindent{\sc Remarque - } Lorsque $F$ est un endomorphisme de Frobenius 
(par exemple lorsque $\d=1$), alors la proposition pr\'ec\'edente 
est une cons\'equence imm\'ediate du th\'eor\`eme de Borel-Tits.\finl

\bigskip

\proof Il est clair que (4) $\imp$ (3). Le fait que (3) $\imp$ (4) 
r\'esulte de la finitude de $\Phi$ et du fait que $\QM(q^{1/\d})$ est 
un corps infini. 

\medskip

Montrons maintenant que (4) $\imp$ (1). Notons $\phit : \Phi \to \Phi$ 
la bijection telle que $F(\a)$ soit un multiple positif de $\phit(\a)$ 
pour tout $\a \in \Phi$. Posons 
$$\Psi=\{\a \in \Phi~|~<\a,v>_\Tb \ge 0\}.$$
Alors $\phit(\Psi)=\Psi$. De plus, $\Psi$ est close, 
$\Psi \cap -\Psi=\Phi_\Lb$ et $\Psi \cup -\Psi=\Phi$. 
Donc il existe un \para $\Pb$ de $\Gb$ dont $\Lb$ est un \levic 
et dont $\Psi$ est le ``syst\`eme de racines'' relativement \`a $\Tb$. 
Puisque $\phit(\Psi)=\Psi$, $\Pb$ est $F$-stable. 

\medskip

Montrons maintenant que (1) $\imp$ (4). Soit $\Pb$ 
un \para $F$-stable de $\Gb$ dont $\Lb$ est un compl\'ement de Levi. 
Fixons un entier naturel non nul $n$ tel que $F^{n\d}(t)=t^{q^n}$ pour tout 
$t \in \Tb$. 
Notons $\m$ l'endomorphisme de $Y(\Zb(\Lb)^\circ)\otimes_\ZM \QM(q^{1/\d})$ 
\'egal \`a $\sum_{k=0}^{n\d-1} q^{k/\d} F^{n\d-1-k}$. Alors 
$\m \circ (F-q^{1/\d})=0$ et 
$$Y(\Zb(\Lb)^\circ)\otimes_\ZM \QM(q^{1/\d}) = \Ker(F-q^{1/\d}) \oplus \Ker \m.
\leqno{(*)}$$
Notons $\Psi$ le ``syst\`eme de racines'' de $\Pb$ relativement 
\`a $\Tb$. Alors il existe $\l \in Y(\Zb(\Lb)^\circ)$ tel que 
$$\Psi=\{\a \in \Phi~|~<\a,\l>_\Tb \ge 0\}.$$
\'Ecrivons $\l=v_1+v_2$, o\`u $v_1$ et $v_2$ appartiennent \`a 
$Y(\Zb(\Lb)^\circ)\otimes_\ZM \QM(q^{1/\d})$ et v\'erifient $F(v_1)=q^{1/\d} v_1$ 
et $\m(v_2)=0$ (voir $(*)$). 
Posons 
$$\Psi'=\{\a \in \Phi~|~<\a,v_1>_\Tb \ge 0\}.$$
Nous allons montrer que $\Psi=\Psi'$. Soit $\a \in \Psi$. Puisque 
$\Psi$ est $\phit$-stable, on a, pour tout $k \in \NM$, 
$<F^k(\a),\l>_\Tb =<\a,F^k(\l)>_\Tb \ge 0$. Par cons\'equent $<\a,\m(\l)>_\Tb \ge 0$ ou, 
en d'autres termes, $<\a,n\d v_1 >_\Tb \ge 0$. Donc $\a \in \Psi'$. 
R\'eciproquement, soit $\a \in \Psi'$. Supposons que $<\a,\l>_\Tb < 0$. 
Alors, puisque $\Phi\setminus\Psi$ est $\phit$-stable, on obtient 
comme pr\'ec\'edemment que $<\a,\m(\l)>_\Tb < 0$, \cad $<\a,n\d v_1 >_\Tb < 0$.

\medskip

Puisque (3) $\imp$ (1), on en d\'eduit que (2) $\imp$ (1). 
Pour finir, montrons que (1) $\imp$ (2). Notons 
$$\Phi'=\{\a \in \Phi~|~\forall v \in 
\Ker(F-q^{1/\d},Y(\Zb(\Lb)^\circ)\otimes_\ZM \QM(q^{1/\d})),~<\a,v>_\Tb=0\}$$
et soit $v_0 \in \Ker(F-q^{1/\d},Y(\Tb) \otimes \QM(q^{1/\d}))$ tel que 
$$\Phi_\Lb=\{\a \in \Phi~|~<\a,v_0>_\Tb=0\}.$$
Alors $\Phi' \incl \Phi_\Lb \incl \Phi'$.\fin

\bigskip

\corollaire{coro borel tits}
{\it Notons $\EC$ l'ensemble des sous-espaces vectoriels 
de $\Ker(F-q^{1/\d},Y(\Tb) \otimes_\ZM \QM(q^{1/\d}))$ et $\LC$ 
l'ensemble des \levis $F$-stables $\Gb$-d\'eploy\'es de $\Gb$ contenant $\Tb$. 
Si $E \in \EC$, notons $\Lb_E$ le \levi $F$-stable 
de $\Gb$ dont le syst\`eme de racines relativement \`a $\Tb$ est 
$$\{\a \in \Phi~|~\forall v \in E,~ <\a,v> = 0\}.$$
Alors l'application
$$\fonctio{\EC}{\LC}{E}{\Lb_E}$$
est une surjection d\'ecroissante.}

\bigskip

\remarque{plus petit levi}
Soient $\Pb_1$ et $\Pb_2$ deux \paras $F$-stables de $\Gb$ 
contenant $\Tb$ et soient $\Lb_1$ et $\Lb_2$ les \levics respectifs de 
$\Pb_1$ et $\Pb_2$ contenant $\Tb$ (ils sont donc $F$-stables). Notons 
$\Ub_1$ et $\Ub_2$ les radicaux unipotents respectifs de $\Pb_1$ et $\Pb_2$. 
Alors, d'apr\`es par exemple \cite[proposition 2.1]{dmbook}, $\Lb_1 \cap \Lb_2$ 
est un \levic $F$-stable du \para $F$-stable $(\Pb_1 \cap \Pb_2).\Ub_1$ de $\Gb$. 
Cela montre qu'il existe un \levi $F$-stable $\Gb$-d\'eploy\'e 
minimal contenant $\Tb$. 

Cela aurait pu se voir gr\^ace au corollaire \ref{coro borel tits} 
dont nous reprenons les notations~: ce \levi $F$-stable $\Gb$-d\'eploy\'e 
minimal minimal est $\Lb_{\Ker(F-q^{1/\d})}$.\finl

\bigskip

\soussection{Quelques propri\'et\'es du morphisme ${\boldsymbol{i^*}}$}
Soit $\Tbt$ un \tor $F$-stable de $\Gbt$ et soit $\Tbt^*$ un \tor 
$F^*$-stable de $\Gbt^*$ dual de $\Tbt$. On pose
$$\Tb=\Tbt \cap \Gb 
\quad\quad\quad {\mathrm{et}}\quad\quad\quad \Tb^*=j^*(\Tbt^*).$$
Alors, d'apr\`es \ref{X exact}, on a une suite exacte 
$$0 \longto X(\Tbt/\Tb) \longto X(\Tbt) \longto X(\Tb) \longto 0.$$
Tous les groupes impliqu\'es dans cette suite exacte sont sans torsion donc, 
par dualit\'e, on obtient que la suite 
$$0 \longto X(\Tb^*) \longto X(\Tbt^*) \longto \Hom(X(\Tbt/\Tb),\ZM) 
\longto 0$$
est exacte. Ici, l' application $X(\Tb^*) \to X(\Tbt^*)$ 
est induite par le morphisme $i^* : \Tb^\pro \to \Tb^*$. En utilisant 
\`a nouveau \ref{X exact}, on obtient un isomorphisme de groupes 
$$\Hom(X(\Tbt/\Tb),\ZM) \simeq X(\Ker i^*).$$
Cela prouve la proposition suivante~:

\bigskip

\proposition{dualite}
{\it Le groupe $\Ker i^*$ est un tore central 
$F^*$-stable de $\Gbt^*$ qui est dual de 
$\Tbt/\Tb$. De plus, cette dualit\'e est compatible avec les isog\'enies 
$F$ et $F^*$.}

\bigskip

\corollaire{coro dual}
{\it Les tores $\Gbt/\Gb$ et $\Ker i^*$ sont duaux et cette dualit\'e est 
compatible avec les isog\'enies $F$ et $F^*$.}

\bigskip

\proof En effet, l'injection $\Tbt \injto \Gbt$ induit un isomorphisme 
$\Tbt/\Tb \simeq \Gbt/\Gb$.\fin

\bigskip

Si $z \in (\Ker i^*)^{F^*}$, nous notons $\zha^{\Gbt}$ le caract\`ere 
lin\'eaire de $\Gbt^F/\Gb^F$ d\'efini par la dualit\'e du corollaire 
\ref{coro dual}. Nous identifions $\zha^{\Gbt}$ avec le caract\`ere 
lin\'eaire de $\Gbt^F$ qu'il induit.

\bigskip

\corollaire{surjecte}
{\it Le morphisme $i^* : \Gbt^{*F^*} \to \Gb^{*F^*}$ est surjectif.}

\bigskip

\proof C'est une cons\'equence imm\'ediate de la connexit\'e de $\Ker i^*$ et 
du th\'eor\`eme de Lang.\fin

\bigskip

\soussection{Action de ${\boldsymbol{\Zb(\Gb)^F}}$ sur 
${\boldsymbol{\Cb\eb\nb\tb(\Gb^F)}}$} 
Si $z \in \Zb(\Gb)^F$ et si $\g \in \Cent \Gb^F$, on pose
$$\fonction{t_z^\Gb \g}{\Gb^F}{\qlb}{g}{\g(zg).}$$
Alors $t_z^\Gb \g \in \Cent \Gb^F$ et l'application 
$$t_z^\Gb : \Cent \Gb^F \to \Cent \Gb^F$$ 
est une isom\'etrie. D'autre part, l'application
$$t^\Gb : \Zb(\Gb)^F \longto \Gb\Lb_\qlb(\Cent \Gb^F)$$ 
est un morphisme de groupes, \cad que l'on a d\'efini ainsi 
une action de $\Zb(\Gb)^F$ par isom\'etries sur $\Cent \Gb^F$. 

\bigskip

\section{Fourre-tout}~

\medskip

\soussection{Morphismes isotypiques} 
Un morphisme $\pi : \Gbh \to \Gb$ entre groupes r\'eductifs 
est dit {\it isotypique} si $\Ker \pi$ est central et $\pi(\Gbh)$ contient 
le groupe d\'eriv\'e de $\Gb$.

\bigskip

\exemple{exemples isotypique} 
Les morphismes $j$ et $j^*$ sont isotypiques.\finl

\bigskip

\soussection{Sous-groupes de Levi auto-oppos\'es} 
Soit $\Lb$ un \levi de $\Gb$. On dit que $\Lb$ est 
($\Gb$-){\it auto-oppos\'e} si, 
pour tout \levi $\Mb$ de $\Gb$ contenant $\Lb$ strictement, 
on a $|N_\Mb(\Lb)/\Lb| \ge 2$. 
D'apr\`es \cite{howlett}, $\Lb$ est $\Gb$-auto-oppos\'e 
\ssi tout \para de $\Gb$ dont $\Lb$ est un \levi est conjugu\'e 
\`a $\Pb$. 

Le groupe $\Gb$ est dit {\it universellement auto-oppos\'e} si, 
pour tout morphisme isotypique $\pi : \Gbh \to \Gb$ et pour 
tout groupe r\'eductif $\Gambh$ dont $\Gbh$ est un 
sous-groupe de Levi, $\Gbh$ est $\Gambh$-auto-oppos\'e.

\bigskip

\exemples{cuspidal auto}
(a) S'il existe une classe unipotente de $\Gb$ supportant un syst\`eme 
local cuspidal (au sens de \cite[introduction]{luicc}), alors $\Gb$ est 
universellement auto-oppos\'e \cite[th\'eor\`eme 9.2]{luicc}.

\medskip

(b) Comme nous le verrons dans la proposition \ref{cuspidal prop} (b), 
un groupe cuspidal (voir \SEC\ref{section cuspidal} pour 
la d\'efinition) est universellement auto-oppos\'e.\finl

\bigskip

Soit maintenant $I$ une partie de $\D$. Alors $I$ est dite 
($W$-){\it auto-oppos\'ee} si $\Lb_I$ est $\Gb$-auto-oppos\'e. 
Posons 
$$W(I)=\{w \in W~|~w(I)=I\},$$
$$W^I=\{w \in W~|~w(I) \incl \D\}$$
$$I^{(1)}=\bigcap_{w \in W^I} w(I).\leqno{\mathrm{et}}$$
Il est clair que $W(I) \incl W^I$ et il est bien connu que 
$N_W(W_I)=W(I) \ltimes W_I$. La deuxi\`eme d\'efinition de 
sous-groupe de Levi auto-oppos\'e montre que $I$ est 
$W$-auto-oppos\'ee \ssi $W(I)=W^I$, \cad \ssi $I=I^{(1)}$.

D\'efinissons une suite d\'ecroissante $(I^{(n)})_{n \in \NM}$ 
de parties de $\D$ par r\'ecurrence de la fa\c{c}on suivante~:
$$\left\{\begin{array}{l}
I^{(0)}=I \\
I^{(n+1)}=(I^{(n)})^{(1)}{\mathrm{~pour~tout~}}n \in \NM.
\end{array}\right.$$
Posons $I^{(\infty)}=\cap_{n \in \NM} I^{(n)}$. Alors 
$I^{(\infty)}$ est la plus grande partie $W$-auto-oppos\'e 
de $\D$ contenue dans $I$.

\bigskip

\soussection{Centralisateurs de sous-tores de $\Gb$} 
Le r\'esultat suivant est une g\'en\'eralisation de 
\cite[corollaire 4.2.3]{bonnafe couro}~:

\bigskip

\lemme{centralisateur tore}
{\it Soit $\Bb$ un \borel de 
$\Gb$ et soit $\Tb$ un \tor de $\Bb$. Notons 
$\Lb$ un \levic d'un \para $\Pb$ de $\Gb$ tels que $\Tb \incl \Lb$ 
et $\Bb \incl \Pb$. Soit $\Phi^+$ (\resp $\Phi_\Lb^+$) le syst\`eme de racines 
positives de $\Gb$ (\resp $\Lb$) relativement \`a $\Tb$ associ\'e \`a $\Bb$. 
Soit $A$ un sous-groupe de $\Aut(\Tb)$ stabilisant 
$\Phi^+$ et $\Phi_\Lb$ et notons $W_\Lb$ le groupe de Weyl de 
$\Lb$ relativement \`a $\Tb$. Alors 
$$C_\Gb\bigl((\Tb^{W_\Lb \rtimes A})^\ci\bigr)=\Lb.$$}

\bigskip

\noindent{\sc Remarque - } Gardons les notations du lemme \ref{centralisateur tore}. 
Alors $W_\Lb$ est un sous-groupe de $\Aut(\Tb)$ normalis\'e par $A$ et 
$W_\Lb \cap A=\{1\}$ car $A$ stabilise $\Phi^+$~: le produit semi-direct 
$W_\Lb \rtimes A$ est donc bien d\'efini. De plus, 
$(\Tb^{W_\Lb})^\ci=\Zb(\Lb)^\ci$ et ce dernier groupe diagonalisable est 
stable sous l'action de $A$. Le lemme \ref{centralisateur tore} dit donc que 
\equat\label{ce t}
C_\Gb\bigl((\Zb(\Lb)^A)^\ci\bigr)=\Lb.~{\SS{\square}}
\endequat

\bigskip

\proof Le groupe 
$(\Tb^{W_\Lb \rtimes A})^\ci$ 
est un sous-tore de $\Gb$ donc $\Mb=C_\Gb((\Tb^{W_\Lb \rtimes A})^\ci)$ est un \levi 
de $\Gb$. De plus, $\Lb \incl \Mb$ par la remarque pr\'ec\'edente. 

Soit $\a \in \Phi^+$ une racine de $\Mb$ relativement \`a $\Tb$. Alors, d'apr\`es  
\cite[Proposition 4.2.1]{bonnafe couro}, on a 
$$\sum_{g \in A \ltimes W_\Lb} g(\a)=0.$$
Soit 
$$\b=\sum_{w \in W_\Lb} w(\a).$$
Supposons que $w(\a)$ soit positive pour tout $w \in W_\Lb$. Alors 
$\b$ est une somme de racines positives et 
$$\sum_{a \in A} a(\b)=0$$
ce qui contredit le fait que $A$ stabilise $\Phi^+$. Donc il existe 
$w \in W_\Lb$ tel que $w(\a)$ n'est pas positive. 
Par cons\'equent, $\a \in \Phi_\Lb$ car $\Tb \incl \Lb$ 
et $\Bb \incl \Pb$.\fin

\bigskip

\soussection{Centralisateur d'\'el\'ements semi-simples} 
Nous rappelons ici le th\'eor\`eme de Steinberg 
\cite[th\'eor\`eme 8.1]{steinberg endo}~:

\bigskip

\Theoreme{Steinberg}{cgs connexe}
{\it Si $\sti \in \Gbt^*$ est semi-simple, alors $C_{\Gbt^*}(\sti)$ est connexe.}

\bigskip

\noindent{\sc Remarque - } Le th\'eor\`eme de Steinberg n\'ecessite l'hypoth\`ese 
de connexit\'e du centre de $\Gbt$. Il n'y a pas de r\'esultat analogue 
pour le groupe $\Gb^*$.\finl

\newpage

{\Large \part{Le groupe ${\boldsymbol{\ZC(\Gb)}}$\label{chapitre ZG}}}

\bigskip

Nous \'etudions ici en d\'etails le groupe des composantes du centre 
$\ZC(\Gb)$. Les r\'esultats de ce chapitre sont pour la plupart 
classiques, sauf les propositions \ref{stabilisateur I} et 
\ref{cuspidal las}. Dans la section \ref{calcul zcg}, nous 
\'etudions le morphisme surjectif $\ZC(\Gb) \to \ZC(\Lb)$ lorsque 
$\Lb$ est un \levi de $\Gb$. Dans la section \ref{sssc}, 
nous montrons comment calculer le groupe $\ZC(\Gb)$ et 
le noyau du morphisme pr\'ec\'edent en termes du diagramme 
de Dynkin affine, du moins lorsque $\Gb$ est simplement connexe. 
La section \ref{h1 section} est consacr\'ee aux multiples 
r\'ealisations du groupe $H^1(F,\ZC(\Gb))$ ainsi qu'\`a ses liens 
avec les \'el\'ements semi-simples de $\Gb^*$. Dans la section 
\ref{section cuspidal}, nous rappelons les diff\'erentes notions 
de cuspidalit\'es introduites par l'auteur (\cite{bonnafe cras}, 
\cite{bonnafe regulier} et \cite{bonnafe mackey}). 

\bigskip

\begin{centerline}{\sc Notations}\end{centerline}

\medskip

Nous nous fixons dans ce chapitre un \tor $F$-stable $\Tbt$ de 
$\Gbt$ et nous posons $\Tb=\Tbt \cap \Gb$. 
Nous notons $W$ le groupe de Weyl de $\Gb$ relativement \`a $\Tb$~;
remarquons que $W$ est canoniquement isomorphe au groupe de Weyl de 
$\Gbt$ relativement \`a $\Tbt$. 
Nous notons $\Phi$ (\resp $\Phit$) le syst\`eme de racines 
de $\Gb$ (\resp $\Gbt$) relativement \`a $\Tb$ 
(\resp $\Tbt$). Le morphisme  
$i_X : X(\Tbt) \injto X(\Tb)$ associ\'e \`a $i$ 
induit une bijection entre $\Phi$ et $\Phit$.

Nous fixons aussi un \borel $\Bbt$ de $\Gbt$ contenant $\Tbt$ et nous posons 
$\Bb=\Bbt \cap \Gb$. Il est \`a noter que $\Bbt$ n'est 
pas n\'ecessairement $F$-stable. Nous notons $\D$ (\resp $\Delt$) 
la base de $\Phi$ (\resp $\Phit$) associ\'ee \`a $\Bb$ 
(\resp $\Bbt$). Alors $\D$ est une base du $\QM$-espace vectoriel 
$X(\Tb/\Zb(\Gb)^\circ) \otimes_\ZM \QM$~: 
nous notons $(\varpi_\a^\vee)_{\a \in \D}$ la base de 
$Y(\Tb/\Zb(\Gb)^\circ) \otimes_\ZM \QM$ duale de $\D$ (pour la dualit\'e 
induite par $< , >_{\Tb/\Zb(\Gb)^\circ}$). 

\bigskip

Si $I$ est une partie de $\D$, nous notons $\Phi_I$ le sous-syst\`eme de 
racines parabolique ayant $\Phi$ comme base et nous notons $W_I$ le groupe 
de Weyl de $\Phi_I$. Nous posons $\Pb_I=\Bb W_I \Bb$ et nous notons $\Lb_I$ 
l'unique \levi de $\Pb_I$ contenant $\Tb$. Alors $\Phi_I$ est le syst\`eme 
de racines de $\Lb_I$ relativement \`a $\Tb$ et $W_I$ est le groupe de Weyl 
de $\Lb_I$ relativement \`a $\Tb$.

\bigskip

\section{Calcul de $\ZC(\Gb)$\label{calcul zcg}}

\medskip

\soussection{Le groupe ${\boldsymbol{\ZC(\Gb)}}$ et le syst\`eme 
de racines de $\Gb$} 
Calculer $\ZC(\Gb)$ et calculer $X(\ZC(\Gb)) \simeq \ZC(\Gb)^\we$ sont des 
probl\`emes \'equivalents. Puisque 
$$\Zb(\Gb)=\{t \in \Tb~|~\forall \a \in \Phi,~\a(t)=1\},$$
on a, d'apr\`es \ref{xd'} et \ref{Y/Y'}~:

\bigskip

\proposition{XZ}
{\it Le morphisme canonique $X(\Tb) \to X(\Zb(\Gb))$ induit un isomorphisme 
de groupes ab\'eliens  
$$X(\ZC(\Gb)) \simeq (X(\Tb)/<\Phi>)_{p'}.$$
De plus, $\tilde{\imath}_{\Tb/\Zb(\Gb)^\circ}$ induit un isomorphisme 
$$(\mathop{\oplus}_{\a \in \D} \ZM \varpi_\a^\vee / 
Y(\Tb/\Zb(\Gb)^\circ))_{p'} \simeq \ZC(\Gb).$$}

\bigskip

\soussection{Sous-groupes de Levi} 
Soit $\Lb$ un \levi de $\Gb$. Alors~:

\bigskip

\proposition{ZG ZL}
{\it Le morphisme $\ZC(\Gb) \to \ZC(\Lb)$ induit par l'inclusion  
$\Zb(\Gb) \injto \Zb(\Lb)$ est surjectif.}

\bigskip

\proof On peut supposer (et nous le ferons) que $\Lb=\Lb_I$ pour une 
partie $I$ de $\D$. Alors $<\Phi_I>$ est un facteur direct de $<\Phi>$, 
donc $(X(\Tb)/<\Phi_I>)_{p'} \to (X(\Tb)/<\Phi>)_{p'}$ est injectif. 
Donc, d'apr\`es la proposition \ref{XZ}, le morphisme naturel 
$X(\ZC(\Lb)) \to X(\ZC(\Gb))$ est injectif.\fin

\bigskip

\corollaire{NGL}
{\it Le groupe $N_\Gb(\Lb)$ agit trivialement sur $\ZC(\Lb)$.}

\bigskip

\corollaire{ZL tilde}
{\it Soit $\Lbt$ un \levi de $\Gbt$. Alors $\Zb(\Lbt)$ est connexe.}

\bigskip

Le morphisme surjectif $\ZC(\Gb) \to \ZC(\Lb)$ sera not\'e $h_\Lb^\Gb$ 
(ou bien $h_\Lb$ lorsqu'il n'y a pas d'ambigu\"\i t\'e). Son 
morphisme dual sera not\'e $\hha_\Lb : \ZC(\Lb)^\we \injto \ZC(\Gb)^\we$. 
Une fois que $\ZC(\Gb)$ 
est d\'etermin\'e, le calcul de $\ZC(\Lb)$ est \'equivalent au calcul de 
$\Ker h_\Lb$. La proposition suivante compl\`ete la proposition \ref{XZ}. 

\bigskip

\Proposition{Digne-Lehrer-Michel}{varpi noyau}
{\it Soit $I$ une partie de $\D$. Alors $\Ker h_{\Lb_I}$ est engendr\'e par  
$(\tilde{\imath}_{\Tb/\Zb(\Gb)^\circ}(\varpi_\a^\vee))_{\a \in \D-I}$.}

\bigskip

\proof En rempla\c{c}ant $\Gb$ par $\Gb/\Zb(\Gb)^\circ$ si c'est n\'ecessaire, 
on peut supposer que $\Zb(\Gb)^\circ=1$. 
Soit $X_I=X(\Tb) \cap ( <\Phi_I> \otimes_\ZM \QM)$. Alors  
$X_I=X(\Tb/\Zb(\Lb_I)^\circ)$. Mais, si $\a \in \D-I$ et si  
$x \in X_I$, alors $<x,\varpi_\a^\vee >_\Tb = 0$. 
Donc $\tilde{\imath}_\Tb(\varpi_\a^\vee) \in \Ker h_{\Lb_I}$. 

R\'eciproquement, soit $z \in \Ker h_{\Lb_I}=\Zb(\Lb_I)^\circ \cap \Zb(\Gb)$. 
Il existe donc $y \in Y(\Zb(\Lb_I)^\circ) \otimes_\ZM \QM$ tel que 
$z=\tilde{\imath}_\Tb(y)$. Mais, $(\varpi_\a^\vee)_{\a \in \D-I}$ est une base 
du $\QM$-espace vectoriel $Y(\Zb(\Lb_I)^\circ) \otimes_\ZM \QM$. Donc 
$$y=\sum_{\a \in \D-I} r_\a \varpi_\a^\vee,$$
avec $r_\a \in \QM$ (pour tout $\a \in \D-I$). Si on pose
$$y'=\sum_{\a \in \D-I} (r_\a)_{p'} \varpi_\a^\vee,$$
alors $z=\tilde{\imath}(y')$. Mais, puisque $z \in \Zb(\Gb)$, on a, 
pour tout $\a \in \D-I$, $<\a,y'>=(r_\a)_{p'} \in \ZM[1/p]$. 
Donc $(r_\a)_{p'} \in \ZM$, ce qui montre que
$$z=\prod_{\a \in \D-I} \tilde{\imath}(\varpi_\a^\vee)^{(r_\a)_{p'}}.$$
Cela termine la preuve de la proposition \ref{varpi noyau}.\fin

\bigskip

La proposition \ref{varpi noyau} entra\^{\i}ne le r\'esultat 
suivant (dont une preuve diff\'erente peut \^etre trouv\'ee par exemple 
dans \cite[proposition 2.4]{bonnafe regulier}).

\bigskip

\corollaire{inter ker}
{\it Soient $I$ et $J$ deux parties de $\D$. Alors  
$\Ker h_{\Lb_{I\cap J}} = (\Ker h_{\Lb_I}). (\Ker h_{\Lb_J})$.}

\bigskip

\corollaire{I infty}
{\it Soit $I$ une partie de $\D$. Alors 
$\Ker h_{\Lb_I}=\Ker h_{\Lb_{I^{(\infty)}}}$.}

\bigskip

\soussection{Les groupes ${\boldsymbol{\ZC(\Gb)}}$ et 
${\boldsymbol{\Ker' i^*}}$} 
Dans cette sous-section, nous construisons un isomorphisme entre les groupes 
$\Ker' i^*$ et $\ZC(\Gb)^\wedge$. Pour cela, remarquons tout d'abord que 
$$X(\Tbt/\Zb(\Gbt)) \simeq (<\Phit> \otimes_\ZM \ZM[1/p]) \cap X(\Tbt/\Zb(\Gbt)).$$
Identifions $i_X : X(\Tbt/\Zb(\Gbt)) \to X(\Tb)$ et 
$i_Y^* : Y(\Tbt^* \cap \Db(\Gbt^*)) \to Y(\Tb^*)$. Alors l'image de 
$i_X$ contient $<\Phi>$ et est contenue dans 
$(<\Phi> \otimes_\ZM \ZM[1/p]) \cap X(\Tb)$. D'apr\`es la proposition 
\ref{curieux}, on obtient alors un isomorphisme de groupes ab\'eliens 
finis
$$(X(\Tb)/<\Phi>)_{p'} \simeq \Ker' i^*,$$
ce qui, en composant avec l'isomorphisme de la proposition \ref{XZ}, 
fournit un isomorphisme
$$X(\ZC(\Gb)) \simeq \Ker' i^*.$$
Gr\^ace \`a l'isomorphisme \ref{X irr}, on obtient finalement 
un isomorphisme
\equat
\o : \Ker' i^* \longto \ZC(\Gb)^\wedge.
\endequat
Cet isomorphisme peut \^etre d\'ecrit de la fa\c{c}on suivante. Soit 
$a \in \Ker' i^*$ et soit $n \in \ZM$, premier \`a $p$, tel que $a^n=1$. 
Alors il existe $\xti \in X(\Tbt/\Zb(\Gbt)) \simeq Y(\Tbt^* \cap \Db(\Gbt^*))$ 
tel que $\xti(\tilde{\imath}(1/n))=a$ et il existe $x \in X(\Tb)$ tel que 
$i_X(\xti)=nx$. En fait, $x \in X(\Tb/\Zb(\Gb)^\circ)$ et 
\equat\label{definition omega}
\o(a)=\kappa \circ \Res_{\ZC(\Gb)}^{\Tb/\Zb(\Gb)^\circ} x.
\endequat
Le r\'esultat suivant est imm\'ediat~:

\bigskip

\lemme{omega commute}
{\it L'isomorphisme $\o$ ne d\'epend pas du choix de $\Tbt$ et $\Tbt^*$ 
$($il d\'epend uniquement du choix de $\imath$ et $\jmath$). De plus, 
$\o \circ F^* = F \circ \o$.}

\bigskip

Les groupes $(\ZC(\Gb)^\we)^F$ et $H^1(F,\ZC(\Gb))^\we$ sont canoniquement 
isomorphes. De m\^eme, les groupes $H^1(F,\ZC(\Gb)^\we)$ et $(\ZC(\Gb)^F)^\we$ 
sont canoniquement isomorphes. Donc, gr\^ace \`a $\o$, nous pouvons 
construire des isomorphismes de groupes
\equat\label{les omegas}
\begin{array}{c}
\o^0 : (\Ker' i^*)^{F^*} \longto H^1(F,\ZC(\Gb))^\we \\
\o^1 : H^1(F^*,\Ker' i^*) \longto (\ZC(\Gb)^F)^\we.
\end{array}
\endequat
Par dualit\'e, nous obtenons des isomorphismes 
\equat\label{les autres omegas}
\begin{array}{c}
\omeh : \ZC(\Gb) \longto (\Ker' i^*)^\we \\
\omeh^0 : H^1(F,\ZC(\Gb)) \longto ((\Ker' i^*)^{F^*})^\we \\
\omeh^1 : \ZC(\Gb)^F \longto H^1(F^*,\Ker' i^*)^\we.
\end{array}
\endequat
Tous ces isomorphismes ne d\'ependent pas du choix de $\Tbt$ et $\Tbt^*$.

\bigskip

\remarque{omega 0} L'isomorphisme $\o^0$ recevra une autre interpr\'etation 
dans la section \ref{h1 section} (voir diagramme \ref{h1 commutatif}).\finl 

\bigskip

\remarque{levi omega}
Soit $\Lb$ un \levi de $\Gb$ et soit $\Lb^*$ un \levi de $\Gb^*$ dual de $\Lb$. 
Soient $\Lbt=\Lb.\Zb(\Gbt)$ et $\Lbt^*=i^{*-1}(\Lb^*)$. 
Notons $\Ker_\Lb' i^* = \Db(\Lbt^*) \cap \Ker i^*$ et 
$\omega_\Lb : \Ker_\Lb' i^* \to \ZC(\Lb)^\we$ l'isomorphisme analogue de $\o$ 
obtenu en rempla\c{c}ant $\Gb$ par $\Lb$. Alors 
le diagramme
$$\diagram
\Ker_\Lb' i^* \rrto^{\DS{\o_\Lb}} \ddto && \ZC(\Lb)^\we \ddto^{\DS{\hha_\Lb}} \\
&&\\
\Ker' i^* \rrto^{\DS{\o}} && \ZC(\Gb)^\we
\enddiagram$$
est commutatif. L'application verticale de gauche
 est bien s\^ur l'inclusion 
canonique et rappelons que $\hha_\Lb$ est l'application duale de $h_\Lb$.\finl

\bigskip

\section{Groupes simplement connexes\label{sssc}}~

\medskip

{\it Nous supposons dans cette section, et uniquement dans cette 
section, que $\Gb$ est semi-simple, quasi-simple et simplement connexe et 
que $p$ ne divise pas le cardinal de $X(\Tb)/<\Phi>$.} 
Nous allons aborder le calcul explicite des groupes $\ZC(\Gb)$ 
et $\ZC(\Lb)$ en utilisant uniquement le syst\`eme de racines de $\Gb$. 
Le calcul du groupe $\ZC(\Gb)$ en termes des poids minuscules 
est fait dans \cite[chapitre VI, \SEC 2, corollaire de la proposition 5]{bourbaki}. 
Le calcul du groupe 
$\ZC(\Lb)$ peut alors \^etre fait gr\^ace \`a la proposition 
\ref{varpi noyau}. En revanche, la proposition \ref{stabilisateur I} 
nous semble nouvelle~: elle contient une description de 
$\Ker h_{\Lb_I}$ lorsque $I$ est auto-oppos\'ee en termes 
du groupe d'automorphismes du diagramme de Dynkin affine 
de $\Gb$.   

\bigskip

\soussection{Poids minuscules} 
Notons $\alpt$ la plus grande racine de $\D$ (elle est bien d\'efinie 
car, puisque $\Gb$ est quasi-simple, $\Phi$ est irr\'eductible). 
Posons 
$$\D^\aff=\D \cup \{-\alpt\}$$
$$W^\aff=W \rtimes Y(\Tb).\leqno{\mathrm{et}}$$
Puisque $\Gb$ est simplement connexe, $W^\aff$ est le groupe de Weyl 
affine de $\Phi$. On d\'efinit 
$$\Aut_W(\D^\aff)=\{w \in W~|~w(\D^\aff)=\D^\aff\}.$$
Alors $\Aut_W(\D^\aff)$ est le groupe des automorphismes du diagramme 
de Dynkin affine de $\Gb$ induits par un \'el\'ement de $W$. 
Pour finir, nous aurons besoin des notations suivantes~:
$$\D_\minus=\{\a \in \D~|~<\alpt,\varpi_\a^\vee>_\Tb = 1\}$$
$$\D_\minus^\aff=\D_\minus \cup \{-\alpt\}.\leqno{\mathrm{et}}$$
Posons conventionnellement $\varpi_{-\alpt}^\vee=0$. 

\bigskip

\proposition{zg quasi}
{\it Supposons $\Gb$ semi-simple, simplement connexe et quasi-simple. Alors 
l'application $\D_\minus^\aff \to \ZC(\Gb)$, 
$\a \mapsto \tilde{\imath}(\varpi_\a^\vee)$ est bijective.}

\bigskip

\proof Voir \cite[chapitre VI, \SEC 2, corollaire de la proposition 5]{bourbaki}.\fin

\bigskip

\soussection{Automorphismes du diagramme de Dynkin affine}
Notons $C_0$ l'alc\^ove
$$C_0=\{y \in Y(\Tb) \otimes_\ZM \RM~|~ \Bigl(\forall \a \in \D,~
<\a, y>_\Tb \ge 0 \Bigr)~ {\mathrm{et}}~<-\alpt,y>_\Tb \le 1~ \}.$$
Soit $z \in \Zb(\Gb)$. Soit $y \in Y(\Tb) \otimes_\ZM \ZM_{(p)}$ tel que  
$\tilde{\imath}_\Tb(y)=z$. Puisque $< \a , y >_\Tb \in \ZM$ pour tout  
$\a \in \Phi$, il existe un unique $w \in W^\aff$ tel que  
$w(C_0)=y + C_0$. Nous notons $w_z$ la projection de $w$ sur $W$. Alors  
$w_z$ ne d\'epend que de $z$ et non du choix de $y$. 

Nous allons maintenant rappeler la formule explicite de 
$w_{\tilde{\imath}_\Tb(\varpi_\a^\vee)}$ pour $\a \in \D_\minus^\aff$. 
Si $\a \in \D^\aff$, nous notons $w_\a = w_\D w_{\D\setminus\{\a\}}$ (remarquons 
que $w_{-\alpt}=1$). Alors on a, pour tout 
$\a \in \Delt$, 
\equat\label{wz=oa}
w_{\tilde{\imath}_\Tb(\varpi_\a^\vee)}=w_\a
\endequat
(voir \cite[chapitre VI, \SEC 2, proposition 6]{bourbaki}). 

La proposition suivante est d\'emontr\'ee dans \cite[\SEC 2.3]{bourbaki}

\bigskip

\proposition{zg affine}
{\it Supposons $\Gb$ semi-simple, quasi-simple et simplement 
connexe. Alors l'application $\ZC(\Gb) \to \Aut_W(\D^\aff)$, 
$z \mapsto w_z$ est bien d\'efinie~; c'est un isomorphisme de groupes.}

\bigskip

Compte tenu du corollaire \ref{I infty}, on peut, pour calculer 
$\Ker h_{\Lb_I}$, se ramener au cas o\`u $I$ est auto-oppos\'ee. 
Dans ce cas, la proposition suivante en fournit une description en termes 
du groupe d'automorphismes du diagramme de Dynkin affine de $\Gb$. 

\bigskip

\proposition{stabilisateur I}
{\it Supposons $\Gb$ semi-simple, quasi-simple et simplement 
connexe. Si $I$ est une partie auto-oppos\'ee de $\D$, alors  
$$\Ker h_{\Lb_I} = \{z \in \ZC(\Gb)~|~w_z(I)=I\}.$$}

\bigskip

\proof Soit $\a \in \D$. Compte tenu de \ref{wz=oa} et de la proposition 
\ref{varpi noyau}, il suffit de montrer que 
$w_\a(I)=I$ si et seulement si $\a \in \D-I$.
Tout d'abord, si $\a \in \D-I$, alors $w_{\D-\{\a\}}(I)=-I$ et $w_\D(I)=-I$ 
car $I$ est auto-oppos\'ee. Donc $w_\a(I)=I$. 
R\'eciproquement, si $\a \in I$, alors $w_{\D-\{\a\}}(\a) \in \Phi^+$ 
et donc $w_\a(\a) \in \Phi^-$. Par suite, $w_\a(\a) \not\in I$, ce 
qui montre que $w_\a(I) \not= I$.\fin

\bigskip

\section{Le groupe $H^1(F,\ZC(\Gb))$\label{h1 section}}~

\medskip

\soussection{Morphisme vers ${\boldsymbol{\Ob\ub\tb(\Gb^F)}}$} 
Commen\c{c}ons par un rappel \'el\'ementaire~:

\bigskip

\lemme{cgg}
{\it On a $C_\Gb(\Gb^F)=\Zb(\Gb)$.}

\bigskip

\proof Nous verrons dans \SEC\ref{unireg soussection} qu'il existe un 
\'el\'ement unipotent $u \in \Bb_0^F$ tel que $C_\Gb(u)=\Zb(\Gb).C_{\Ub_0}(u)$. 
Donc $C_\Gb(\Gb^F) \incl C_\Gb(u)=\Zb(\Gb).C_{\Ub_0}(u)$ 
Si on note $w_0$ un \'el\'ement de $\Gb^F$ repr\'esentant l'\'el\'ement 
de plus grande longueur de $W_0$, alors 
$C_\Gb(\Gb^F) \incl C_\Gb(w_0 u w_0^{-1})=\Zb(\Gb). \lexp{w_0}{C_{\Ub_0}(u)}$. 
Donc $C_\Gb(\Gb^F) \incl C_\Gb(u) \cap C_\Gb(w_0 u w_0^{-1}) = \Zb(\Gb)$. 
D'autre part, il est clair que $\Zb(\Gb) \incl C_\Gb(\Gb^F)$. 
D'o\`u le r\'esultat.\fin

\bigskip

\remarque{zgf} Le lemme \ref{cgg} montre en particulier que 
le centre de $\Gb^F$ est $\Zb(\Gb)^F$.\finl

\bigskip

Nous notons $\Aut(\Gb,F)$ le groupe des automorphismes de $\Gb$ 
commutant avec $F$. Alors le groupe $\Int(\Gb^F)$ des automorphismes 
de $\Gb$ induits par la conjugaison par un \'el\'ement de $\Gb^F$ est 
un sous-groupe distingu\'e de $\Aut(\Gb,F)$. D'apr\`es 
la remarque \ref{zgf}, le groupe 
$\Int(\Gb^F)$ est isomorphe au groupe des automorphismes 
int\'erieurs de $\Gb^F$, ce qui justifie la notation utilis\'ee. 
Nous notons $\Out(\Gb,F)$ le groupe quotient $\Aut(\Gb,F)/\Int(\Gb^F)$. 
On a un morphisme canonique $\Out(\Gb,F) \to \Out(\Gb^F)$.

Si $z \in H^1(F,\ZC(\Gb))$, nous notons $g_z$ un \'el\'ement de $\Gb$ tel 
que $g_z^{-1} F(g_z)$ appartient \`a $\Zb(\Gb)$ et repr\'esente $z$. 
Alors l'automorphisme int\'erieur $\INT g_z$ appartient 
\`a $\Aut(\Gb,F)$. On note $\t_z^\Gb$ son image dans $\Out(\Gb^F)$. 

\bigskip

\proposition{h1 sur g}
{\it L'application $\t^\Gb : H^1(F,\ZC(\Gb)) \to \Out(\Gb^F)$, $z \mapsto \t_z^\Gb$ 
est bien d\'efinie~: c'est un morphisme injectif de groupes.}

\bigskip

\proof Soient $z \in H^1(F,\ZC(\Gb))$ et soient $g$ et $h$ deux 
\'el\'ements de $\Gb$ tels que $h^{-1}F(h)$ et $g^{-1}F(g)$ appartiennent 
\`a $\Zb(\Gb)$ et repr\'esentent $z$. Alors il existe $x \in \Zb(\Gb)$ 
tel que $x^{-1} F(x) h^{-1} F(h) = g^{-1} F(g)$. Puisque 
$x$ est central dans $\Gb$, on a $F(gh^{-1}x^{-1})=gh^{-1}x^{-1}$. 
Posons $y = gh^{-1}x^{-1}$. Alors $y \in \Gb^F$ et 
$\INT g=\INT(y) \compose \INT(hx)=\INT(y) \compose \INT(h)$. 
Donc l'image de $\INT g$ dans $\Out(\Gb,F)$ co\"{\i}ncide 
avec l'image de $\INT h$. Cela montre que $\t^\Gb$ est bien d\'efinie. 

\medskip

Montrons maintenant 
que c'est un morphisme de groupes. Soient $z$ et $z'$ deux \'el\'ements 
de $H^1(F,\ZC(\Gb))$ et soient $g$ et $g'$ deux \'el\'ements de 
$\Gb$ tels que $g^{-1} F(g)$ et $g^{\prime -1} F(g')$ appartiennent 
\`a $\Zb(\Gb)$ et repr\'esentent respectivement $z$ et $z'$. 
Alors, si on pose $a=(gg')^{-1}F(gg')$, on a 
$a=g^{\prime -1} g^{-1} F(g) F(g')= g^{\prime -1} F(g') g^{-1} F(g)$ car 
$g^{-1}F(g)$ est central. Donc $a \in \Zb(\Gb)$ et $a$ repr\'esente $zz'$. 
Mais $\INT(gg')=\INT(g) \compose \INT(g')$, donc 
$\t_{zz'}^\Gb=\t_z^\Gb \compose \t_{z'}^\Gb$. Par cons\'equent, $\t^\Gb$ est 
bien un morphisme de groupes.

\medskip

Pour finir, montrons que $\t^\Gb$ est injectif. Soit $z \in H^1(F,\ZC(\Gb))$ 
et soit $g \in \Gb$ tel que $g^{-1}F(g)$ appartient \`a $\Zb(\Gb)$ et 
repr\'esente $z$. Supposons que $\INT(g)$ induit un 
automorphisme int\'erieur de $\Gb^F$. Alors il existe $h \in \Gb^F$ 
tel que $h^{-1}g \in C_\Gb(\Gb^F)=\Zb(\Gb)$ (voir lemme \ref{cgg}). 
Donc, si on pose $a=h^{-1}g$, alors $g^{-1}F(g)=(ha)^{-1}F(ha)=a^{-1}F(a)$ 
car $F(h)=h$. Donc $z=1$.\fin

\bigskip

\noindent{\sc Remarque - } Le morphisme de groupes 
$\Out(\Gb,F) \to \Out(\Gb^F)$ est en g\'en\'eral non injectif. 
Par exemple, si $\Gb$ est un tore, si $\d=1$, et 
si $\Gb$ est d\'eploy\'e sur $\fq$, alors tout automorphisme 
de $\Gb$ commute avec $F$, donc $\Out(\Gb,F) = \Aut(\Gb)$ car $\Gb$ 
est ab\'elien. Mais, si $n=\dim \Gb$, on a $\Aut(\Gb) \simeq \Gb\Lb_n(\ZM)$ 
et $\Aut(\Gb^F) \simeq \Gb\Lb_n(\ZM/(q-1)\ZM)$, donc le morphisme 
$\Aut(\Gb) \to \Aut(\Gb^F)$ a un noyau infini lorsque $n \ge 2$.\finl

\bigskip

Le groupe $\Out(\Gb^F)$ agit sur $\Cent \Gb^F$ et $\Irr \Gb^F$ de la fa\c{c}on 
suivante~: si $\t \in \Out(\Gb^F)$ et si $\ch$ appartient \`a $\Cent \Gb^F$ ou 
$\Irr \Gb^F$, on pose $\t(\ch)=\ch \compose \taut^{-1}$, o\`u $\taut$ est un 
automorphisme de $\Gb^F$ repr\'esentant $\t$. Cela nous d\'efinit donc, 
\`a travers le morphisme $\t^\Gb$, une action de $H^1(F,\ZC(\Gb))$ 
sur $\Cent \Gb^F$ et $\Irr \Gb^F$. Si $V$ est un sous-espace de 
$\Cent(\Gb^F)$ stable sous l'action de $H^1(F,\ZC(\Gb))$ et si 
$\z \in H^1(F,\ZC(\Gb))^\we$, nous noterons $V_\z$ la composante 
$\z$-isotypique de $V$. On a donc $V_\z=V \cap \Cent(\Gb^F)_\z$. 

\bigskip

\soussection{Le groupe ${\boldsymbol{\Gbt^F/\Gb^F.\Zb(\Gbt)^F}}$} 
Soit $\Lbt$ un \levi de $\Gbt$. On 
suppose ici que $\Lbt$ est $F$-stable. On pose $\Lb = \Gb \cap \Lbt$. 
Alors $\Lb$ est un \levi $F$-stable de $\Gb$. 
Le morphisme de groupes $h_\Lb^1 : H^1(F,\ZC(\Gb)) \to H^1(F,\ZC(\Lb))$ 
induit par $h_\Lb$ est surjectif. S'il y a ambigu\"\i t\'e, nous le noterons 
$h_\Lb^{\Gb,1}$. 

\bigskip

Soit $\lti \in \Lbt^F$. Alors il existe $l \in \Lb$ et $\zti \in \Zb(\Gbt)$ 
tels que $\lti=l\zti$. Puisque $F(\lti)=\lti$, on en d\'eduit que 
$l^{-1}F(l) = F(\zti)^{-1} \zti$. Donc $l^{-1}F(l) \in \Zb(\Gb)$. 
On note $\s_\Lb^\Gb(\lti)$ sa classe dans $H^1(F,\ZC(\Gb))$. 
Il est facile de v\'erifier que $\s_\Lb^\Gb(\lti)$ ne d\'epend que 
de $\lti$ et non du choix de $l$ et $\zti$. Il est tout aussi imm\'ediat 
que $\s_\Lb^\Gb$ induit un isomorphisme de groupes (toujours not\'e $\s_\Lb^\Gb$) 
$$\Lbt^F/\Lb^F.\Zb(\Gbt)^F \longmapright{\sim} H^1(F,\ZC(\Gb)).$$
D'autre part, on a un morphisme canonique 
$$\Lbt^F/\Lb^F.\Zb(\Gbt)^F \longto \Out(\Lb^F)$$ 
et un simple calcul montre que le diagramme
$$\diagram
\Lbt^F/\Lb^F.\Zb(\Gbt)^F \rrto^{\DS{\s_\Lb^\Gb}} \ddto && H^1(F,\ZC(\Gb)) 
\ddto^{\DS{h_\Lb^1}} \\
&& \\
\Out(\Lb^F) && H^1(F,\ZC(\Lb)) \llto_{\DS{\t^\Lb}} 
\enddiagram$$
est commutatif. Lorsque $\Lb=\Gb$, l'isomorphisme $\s_\Lb^\Gb$ sera not\'e 
$\s_\Gb$. De plus, l'isomorphisme dual de $\s_\Lb^\Gb$ sera not\'e 
$\sigh_\Lb^\Gb : H^1(F,\ZC(\Gb))^\we \longmapright{\sim} 
(\Lbt^F/\Lb^F.\Zb(\Gbt)^F)^\we$ 
(ou $\sigh_\Gb$ si $\Lb=\Gb$).

\medskip

En utilisant l'isomorphisme 
$\o^0 : (\Ker' i^*)^{F^*} \longmapright{\sim} H^1(F,\ZC(\Gb))^\we$ 
(voir \ref{les omegas}), on obtient un diagramme commutatif~:

\medskip

\equat\label{h1 commutatif}
\diagram
&& (\Ker' i^*)^{F^*} \dllto_{\DS{\o^0}}\rrto && (\Ker i^*)^{F^*} \ddto^{\sim}\\
H^1(F,\ZC(\Gb))^\wedge \drrto_{\DS{\sigh_\Lb^\Gb}} &&&& \\
&& (\Lbt^F/\Lb^F.\Zb(\Gbt)^F)^\we \rrto && (\Lbt^F/\Lb^F)^\wedge.
\enddiagram
\endequat

\bigskip

\noindent{\sc D\'emonstration de \ref{h1 commutatif} - } Soit 
$\Tbt$ un \tor $F$-stable de $\Lbt$. Soit $\Tbt^*$ un \tor $F^*$-stable 
de $\Gbt$ dual de $\Tbt$. Posons $\Tb=\Tbt \cap \Gb$ et $\Tb^*=i^*(\Tbt^*)$. 
Puisque $\Lbt^F/\Lb^F \simeq \Tbt^F/\Tb^F$, il suffit de montrer 
la commutativit\'e du diagramme \ref{h1 commutatif} lorsque $\Lbt=\Tbt$, 
ce que nous supposerons dor\'enavant. 
 
Soit $n \in \NM^*$. Notons $\o_n^0 : (\Ker' i^*)^{F^{*n}} \to H^1(F^n,\ZC(\Gb))$ 
induit par $\o$. Notons 
$$\b_n : (\Ker' i^*)^{F^{*n}} \to (\Tbt^{F^n})^\we$$ 
le morphisme compos\'e 
$(\Ker' i^*)^{F^{*n}} \to (\Ker i^*)^{F^{*n}} 
\mapright{\sim} (\Tbt^{F^n})^\we$ et 
$$\g_n : H^1(F^n,\ZC(\Gb))^\we \to (\Tbt^{F^n})^\we$$ 
le morphisme compos\'e 
$H^1(F^n,\ZC(\Gb))^\we \to (\Tbt^{F^n}/\Tb^{F^n}.\Zb(\Gbt)^{F^n})^\we 
\to (\Tbt^{F^n})^\we$. Il s'agit de montrer que $\g_1 \circ \o_1^0 = \b_1$. 
Mais, le diagramme
$$\diagram
&& (\Ker' i^*)^{F^*} \dllto_{\DS{\o_1^0}} \ddto \drrto^{\DS{\b_1}} && \\
H^1(F,\ZC(\Gb))^\wedge \xto[0,4]^{\DS{\g_1\qquad\qquad}} \ddto 
&&&& (\Tbt^F)^\we 
\ddto^{\DS{N_{F^n/F}^\we}} \\
&& (\Ker' i^*)^{F^{*n}} \dllto_{\DS{\o_n^0}} \drrto^{\DS{\b_n}} && \\
H^1(F^n,\ZC(\Gb))^\wedge \xto[0,4]^{\DS{\g_n\qquad\qquad}} 
&&&& (\Tbt^{F^n})^\we \\
\enddiagram$$
est commutatif. Ici, toutes les applications verticales sont injectives. 
En particulier, cela montre qu'il suffit de montrer que $\g_n \circ \o_n^0 = \b_n$ 
et donc que l'on peut remplacer $F$ par n'importe laquelle de ses puissances. 
Par exemple, et c'est ce que nous ferons par la suite, nous pouvons 
supposer que $F$ est un endomorphisme de Frobenius d\'eploy\'e de $\Tbt$ 
(sur un corps fini \`a $q$ \'el\'ements) et que $F$ 
agit trivialement sur $\ZC(\Gb)$. En particulier, $F^*$ est un endomorphisme 
de Frobenius d\'eploy\'e de $\Tbt^*$ et $F^*$ agit trivialement sur $\Ker' i^*$. 

Soit $a \in \Ker' i^* = (\Ker' i^*)^{F^*}$ et soit $\tti \in \Tbt^F$. 
Puisque $F^*$ est d\'eploy\'e, il existe 
$\yti \in Y(\Tbt^* \cap \Db(\Gbt^*)) \simeq X(\Tbt/\Zb(\Gbt)) \incl X(\Tbt)$ tel que 
$a=y(\tilde{\imath}(1/(q-1)))$. Soit alors $y \in Y(\Tb^*)\simeq X(\Tb)$ 
tel que $i_X(\yti)=(q-1) y$. 
Soient maintenant $t \in \Tb$ et $\zti \in \Zb(\Gbt)$ tels que $\tti=t\zti$. 
Posons $z=t^{-1}F(t) \in \Zb(\Gb)$. 
D'apr\`es \ref{egalite sha} et \ref{definition omega}, il suffit de 
montrer que $\yti(\tti)=y(z)$. Mais, 
$$y(z)=y(t^{-1}F(t))=y(t^{q-1})=((q-1)y)(t)=i_X(\yti)(t)=\yti(t).$$ 
Cela montre 
le r\'esultat car $t \Zb(\Gbt)=\tti\Zb(\Gbt)$ et $\yti \in X(\Tbt/\Zb(\Gbt))$.\fin

\bigskip

\section{Cuspidalit\'e\label{section cuspidal}}~

\medskip

\soussection{D\'efinition et premi\`eres propri\'et\'es}  
Le groupe r\'eductif $\Gb$ est dit {\it cuspidal} si, pour tout 
\levi $\Lb$ propre de $\Gb$, on a $\Ker h_\Lb^\Gb \not= \{1\}$ 
(voir \cite[\SEC 1]{bonnafe cras} ou \cite[\SEC 2.C]{bonnafe regulier}). 
Nous rappelons ici quelques propri\'et\'es 
des groupes cuspidaux dont le lecteur pourra trouver une preuve dans 
\cite[propositions 2.12 et 2.18 et remarque 2.14]{bonnafe regulier}.

\bigskip

\proposition{cuspidal prop}
{\it Supposons que $\Gb$ est cuspidal. Alors~:
\begin{itemize}
\itemth{a} Si $\pi : \Gbh \to \Gb$ est un morphisme isotypique, alors 
$\Gbh$ est cuspidal.

\itemth{b} Le groupe $\Gb$ est universellement auto-oppos\'e.

\itemth{c} Toutes les composantes quasi-simples de $\Gb$ sont de type $A$.
\end{itemize}}

\bigskip

\soussection{Sous-groupes de Levi cuspidaux\label{sous cuspidal levi}}
Si $K$ est un sous-groupe de $\ZC(\Gb)$, nous notons $\LC(K)$ (ou $\LC^\Gb(K)$ 
s'il y a ambigu\"\i t\'e) l'ensemble des \levis de $\Gb$ 
tels que $\Ker h_\Lb \incl K$. Nous notons $\LC_\mini(K)$ 
(ou $\LC_\mini^\Gb(K)$) l'ensemble des \'el\'ements minimaux pour l'inclusion 
de $\LC(K)$.
La preuve de la proposition suivante peut \^etre trouv\'ee dans 
\cite[lemme 2.16]{bonnafe regulier}~:

\bigskip

\proposition{cuspidal L}
{\it Si $K$ est un sous-groupe de $\ZC(\Gb)$, alors $\LC_\mini(K)$ 
est une seule classe de conjugaison de sous-groupes de $\Gb$. 
De plus, ses \'el\'ements sont cuspidaux.}

\bigskip

Il est facile, en utilisant entre autres la proposition \ref{stabilisateur I}, 
de classifier les groupes $\LC_\mini(K)$ lorsque $\Gb$ 
est semi-simple, quasi-simple et simplement connexe. Cette classification 
est faite dans \cite[table 2.17]{bonnafe regulier}. Nous la rappelons 
dans la table \ref{tabletable} (cette table ne contient pas les groupes 
$\LC_\mini(\ZC(\Gb))$ car ce sont les tores maximaux de $\Gb$). Pour pouvoir 
lire cette table, il convient de signaler que les copoids fondamentaux 
en type $D_{2r}$ sont num\'erot\'es comme dans \cite[planches]{bourbaki}). 
La classification 
dans le cas g\'en\'eral d\'ecoule de \cite[\SEC 2.B]{bonnafe regulier}. 

\bigskip

\begin{table}
\begin{tabular}{|c||c|c|c|c|c|}
\hline
$\xy (0,3) *+={\mathrm{Type}}; (0,-3) *+={{\mathrm{de~}}\Gb};\endxy$ & $\ZC(\Gb)$ & $K$ & 
${\DS{{{\mathrm{Type~de}}}\vphantom{A_N^C}} \atop 
\DS{\vphantom{A_N^C}\Lb \in \LC_\mini(K)}}$ 
& $\ZC(\Lb)$ & Diagramme de $(\Gb,\Lb)$\\
\hline
\hline
$A_r$ & $\mub_{r+1}$ & $\xy (0,6) *+={\mub_{r+1/d}}; 
(0,0) *+={\DS{{\mathrm{o\grave{u}~}}d~|~r+1}}; 
(0,-6) *+={\DS{{\mathrm{et}}~p~\not|~ d}};\endxy$ & 
$\underbrace{A_{d-1} \times \dots \times A_{d-1}}_{{{\frac{r+1}{d}} ~{\mathrm{fois}}}}$ 
&$\mub_d$ & 
{\small $\xy
(0,0) *+={~\SS{A_{d-1}}~} *\frm{-} ;(8,0) *+={\phan} *\frm{o} **@{-} ;
(16,0) *+={~\SS{A_{d-1}}~} *\frm{-}**@{-}; (24,0) *+={\phan} *\frm{o} **@{-} ;
(27,0) *+={~\SS{\dots}~} **@{-};(30,0) *+={\phan} *\frm{o} **@{-};
(38,0) *+={~\SS{A_{d-1}}~} *\frm{-}**@{-};
\endxy $}\\
\hline
\hline
$\xy (0,0) *+={B_{2r+1}}; (0,-6) *+={p\not=2};\endxy$ & $\mub_2$ & $1$ & 
$\underbrace{A_1 \times \dots \times A_1}_{r+1~{\mathrm{fois}}}$ &$\mub_2$  & 
$\vphantom{\xy (0,4) *+={a} ; (0,-4) *+={a};\endxy }$
{\small $\xy (0,0) *+={\bullet} *\frm{o} ; (7,0) *+={\phan} *\frm{o} **@{-} ;
(14,0) *+={\bullet} *\frm{o} **@{-} ; (21,0) *+={\phan} *\frm{o} **@{-} ;
(26,0)  *+={~\dots~} **@{-};(31,0) *+={\phan} *\frm{o} **@{-};
(38,0) *+={\bullet} *\frm{o} **@{=}; (34.5,0) *+={>};
\endxy$}\\
\hline
$\xy (0,0) *+={B_{2r}}; (0,-6) *+={p\not=2};\endxy$
& $\mub_2$ & $1$ & 
$\underbrace{A_1 \times \dots \times A_1}_{r~{\mathrm{fois}}}$ &$\mub_2$  & 
$\vphantom{\xy (0,4) *+={a} ; (0,-4) *+={a};\endxy }$
{\small $\xy (0,0) *+={\bullet} *\frm{o} ; (7,0) *+={\phan} *\frm{o} **@{-} ;
(14,0) *+={\bullet} *\frm{o} **@{-} ; (21,0) *+={\phan} *\frm{o} **@{-} ;
(26,0)  *+={~\dots~} **@{-};(31,0) *+={\bullet} *\frm{o} **@{-};
(38,0) *+={\phan} *\frm{o} **@{=}; (34.5,0) *+={>};
\endxy$}\\
\hline
\hline
$C_r$ 
& $\mub_2$ & $1$ & $A_1$ & $\mub_2$ & 
$\vphantom{\xy (0,4) *+={a} ; (0,-2) *+={a};\endxy }$
{\small $\xy (0,0) *+={\phan} *\frm{o};(8,0) *+={\phan} *\frm{o}**@{-}; 
(14,0) *+={~\dots~} **@{-};
(20,0) *+={\phan} *\frm{o} **@{-}; (28,0) *+={\phan} *\frm{o} **@{-};
(36,0)*+={\bullet} *\frm{o} **@{=}; (32,0) *+={<};
\endxy$}\\
$p \not= 2$ &&& (grande racine)$\vphantom{\xy (0,0) *+={a};(0,-4) *+={a};\endxy }$ 
&&\\
\hline
\hline
$\xy (0,0) *+={D_{2r+1}}; (0,-6) *+={p\not=2};\endxy$& $\mub_4$ &
$\xy (0,7) *+={1}; (0,-7) *+={\mub_2};\endxy$ & 
$\xy 
(0,4) *+={\underbrace{A_1 \times \dots \times A_1}_{r-1~{\mathrm{fois}}} 
\times A_3};
(0,-7) *+={A_1 \times A_1};\endxy $ 
& $\xy (0,7) *+={\mub_4}; (0,-7) *+={\mub_2};\endxy$ 
& 
{\small $\xy (0,7) *+={\bullet} *\frm{o} ; (7,7) *+={\phan} *\frm{o} **@{-} ;
(14,7) *+={\bullet} *\frm{o} **@{-} ; (19,7) *+={~\dots~}  **@{-} ;
(24,7)  *+={\bullet} *\frm{o}**@{-};(31,7) *+={\phan} *\frm{o} **@{-};
(38,7) *+={\bullet} *\frm{o} **@{-}; (43,11) *+={\bullet} *\frm{o} **@{-};
(38,7) *+={\bullet} *\frm{o} **@{-}; (43,3) *+={\bullet} *\frm{o} **@{-};
(0,-7) *+={\phan} *\frm{o} ; (10,-7) *+={\phan} *\frm{o} **@{-} ;
(19,-7) *+={~\dots~}  **@{-} ;
(28,-7)  *+={\phan} *\frm{o}**@{-};
(38,-7) *+={\phan} *\frm{o} **@{-}; (43,-3) *+={\bullet} *\frm{o} **@{-};
(38,-7) *+={\phan} *\frm{o} **@{-}; (43,-11) *+={\bullet} *\frm{o} **@{-};
\endxy$}\\
\hline
&& 1  & 
$\underbrace{A_1 \times \dots \times A_1}_{r+1~{\mathrm{fois}}}$ &
$\mub_2 \times \mub_2$&
$\vphantom{\xy (0,6) *+={a} ; (0,-6) *+={a};\endxy }$
{\small $\xy (0,0) *+={\bullet} *\frm{o} ; (7,0) *+={\phan} *\frm{o} **@{-} ;
(14,0) *+={\bullet} *\frm{o} **@{-} ; (21,0) *+={\phan} *\frm{o} **@{-} ;
(26,0)  *+={~\dots~} **@{-};(31,0) *+={\bullet} *\frm{o} **@{-};
(38,0) *+={\phan} *\frm{o} **@{-}; (43,4) *+={\bullet} *\frm{o} **@{-};
(38,0) *+={\phan} *\frm{o} **@{-}; (43,-4) *+={\bullet} *\frm{o} **@{-};\endxy$}\\
$D_{2r}$& $\xy (0,-5) *+={\DS{\mub_2 \times \mub_2}}; \endxy$& 
$<\tilde{\iota}_\Tb(\varpi_{2r-1}^\ve)>$  & 
$\underbrace{A_1 \times \dots \times A_1}_{r~{\mathrm{fois}}}$ & $\mub_2$ & 
$\vphantom{\xy (0,6) *+={a} ; (0,-6) *+={a};\endxy }$ 
{\small $\xy (0,0) *+={\bullet} *\frm{o} ; (7,0) *+={\phan} *\frm{o} **@{-} ;
(14,0) *+={\bullet} *\frm{o} **@{-} ; (21,0) *+={\phan} *\frm{o} **@{-} ;
(26,0)  *+={~\dots~} **@{-};(31,0) *+={\bullet} *\frm{o} **@{-};
(38,0) *+={\phan} *\frm{o} **@{-}; (43,4) *+={\phan} *\frm{o} **@{-};
(38,0) *+={\phan} *\frm{o} **@{-}; (43,-4) *+={\bullet} *\frm{o} **@{-};
\endxy$}\\
$p \not= 2$&& $<\tilde{\iota}_\Tb(\varpi_{2r}^\ve)>$ & 
$\underbrace{A_1 \times \dots \times A_1}_{r~{\mathrm{fois}}}$ & $\mub_2$ & 
$\vphantom{\xy (0,2) *+={a} ; (0,-10) *+={a};\endxy }$ 
{\small $\xy (0,0) *+={\bullet} *\frm{o} ; (7,0) *+={\phan} *\frm{o} **@{-} ;
(14,0) *+={\bullet} *\frm{o} **@{-} ; (21,0) *+={\phan} *\frm{o} **@{-} ;
(26,0)  *+={~\dots~} **@{-};(31,0) *+={\bullet} *\frm{o} **@{-};
(38,0) *+={\phan} *\frm{o} **@{-}; (43,4) *+={\bullet} *\frm{o} **@{-};
(38,0) *+={\phan} *\frm{o} **@{-}; (43,-4) *+={\phan} *\frm{o} **@{-};
\endxy$}\\
&&  $<\tilde{\iota}_\Tb(\varpi_1^\ve)>$  & $A_1 \times A_1$ & $\mub_2$ &  
{\small $\xy 
(0,0) *+={\phan} *\frm{o} ; (10,0) *+={\phan} *\frm{o} **@{-} ;
(19,0) *+={~\dots~}  **@{-} ;
(28,0)  *+={\phan} *\frm{o}**@{-};
(38,0) *+={\phan} *\frm{o} **@{-}; (43,4) *+={\bullet} *\frm{o} **@{-};
(38,0) *+={\phan} *\frm{o} **@{-}; (43,-4) *+={\bullet} *\frm{o} **@{-}; 
(43,-7) *+={~}
\endxy$}\\
\hline
\hline
$\xy (0,0) *+={E_6}; (0,-6) *+={p\not=3};\endxy$& $\mub_3$ & 
$1$ & $A_2 \times A_2$ & $\mub_3$ & 
$\vphantom{\xy (0,5) *+={a} ; (0,-12) *+={a};\endxy }$ 
{\small $\xy (0,0) *+={\bullet} *\frm{o} ; (10,0) *+={\bullet} *\frm{o} **@{-} ;
(20,0) *+={\phan} *\frm{o} **@{-}; (30,0) *+={\bullet} *\frm{o} **@{-} ; 
(40,0)*+={\bullet} *\frm{o} **@{-} ; 
(20,0) *+={\phan} *\frm{o}  ; (20,-6) *+={\phan} *\frm{o} **@{-}; 
\endxy$}\\
\hline
\hline
$\xy (0,0) *+={E_7}; (0,-6) *+={p\not=2};\endxy$
& $\mub_2$ & $1$ & $A_1 \times A_1 \times A_1$ & $\mub_2$ & 
$\vphantom{\xy (0,5) *+={a} ; (0,-12) *+={a};\endxy }$ 
{\small $\xy (0,0) *+={\phan} *\frm{o} ; (9,0) *+={\phan} *\frm{o} **@{-} ;
(18,0) *+={\phan} *\frm{o} **@{-}; (27,0) *+={\bullet} *\frm{o} **@{-} ; 
(36,0)*+={\phan} *\frm{o} **@{-} ; (45,0) *+={\bullet} *\frm{o} **@{-} ;
(18,0) *+={\phan} *\frm{o}  ; (18,-6) *+={\bullet} *\frm{o} **@{-}; 
\endxy$}\\
\hline
\end{tabular}

\bigskip

\begin{centerline}{\refstepcounter{theo}\label{tabletable}\noindent\bf Table   
\arabic{section}.\arabic{theo}}\end{centerline}
\end{table}

\bigskip

\soussection{Caract\`eres lin\'eaires cuspidaux} 
Un caract\`ere lin\'eaire $\z : \ZC(\Gb) \to {\overline{\QM}}_\ell^\times$ est dit 
{\it cuspidal} si, pour tout sous-groupe de Levi propre $\Lb$ de $\Gb$, 
on a $\Ker h_\Lb \notincl \Ker \z$. Nous noterons $\ZC_\cus^\we(\Gb)$ 
l'ensemble des caract\`eres lin\'eaires cuspidaux de $\ZC(\Gb)$. 
Si $\ZC_\cus^\we(\Gb) \not=\vide$, alors $\Gb$ est cuspidal. En particulier, 
toutes ses composantes quasi-simples sont de type $A$ (voir proposition 
\ref{cuspidal prop} (c)). 

\bigskip

\section{\'El\'ements semi-simples et non connexit\'e 
du centre\label{section semisimple}}

\medskip

\begin{quotation}
\noindent{\bf Hypoth\`ese : } {\it Nous fixons dans cette section un 
\'el\'ement semi-simple $s \in \Gb^{*F^*}$. Nous fixons aussi 
un \'el\'ement semi-simple $\sti \in \Gb^{*F^*}$ tel que $i^*(\sti)=s$.}
\end{quotation}

\medskip

L'existence de $\sti$ est assur\'ee par le th\'eor\`eme de Lang 
et la connexit\'e de $\Ker i^*$.

\bigskip

\soussection{Centralisateur de ${\boldsymbol{s}}$\label{soussection ags}}
Soit $\Bb^*_1$ un \borel $F^*$-stable de $C_{\Gb^*}^\ci(s)$ et soit  
$\Tb_1^*$ un \tor $F^*$-stable de $\Bb_1^*$. On note $W$ (\resp $W(s)$, 
\resp $W^\ci(s)$) le groupe de Weyl de $\Gb^*$ (\resp $C_{\Gb^*}(s)$, 
\resp $C_{\Gb^*}^\ci(s)$) relativement \`a $\Tb_1^*$. Alors 
$$W(s)=\{w \in W~|~ w(s)=s\}$$
et $W^\ci(s)$ est un sous-groupe distingu\'e de $W(s)$. De plus, 
$W(s)/W^\ci(s)$ est canoniquement isomorphe \`a $A_{\Gb^*}(s)$. Soit 
$A(s)=\{w \in W(s)~|~\lexp{w}{\Bb_1^*}=\Bb_1^*\}$. Alors 
$A(s)$ est un sous-groupe $F^*$-stable de $W(s)$ et 
$W(s)=W^\ci(s) \rtimes A(s)$. Donc $A(s)$ est canoniquement isomorphe 
\`a $A_{\Gb^*}(s)$. Nous identifierons par la suite 
$A_{\Gb^*}(s)$ avec $A(s)$, de sorte que 
\equat
W(s)=W^\ci(s) \rtimes A_{\Gb^*}(s).
\endequat
Soit $\Phi_1$ (\resp $\Phi_s$) le syst\`eme de racines de $\Gb^*$ 
(\resp $C_{\Gb^*}^\circ(s)$) relativement \`a $\Tb_1^*$. Alors 
$$\Phi_s=\{\a \in \Phi_1~|~\a(s)=1\}.$$
Pour tout $w \in W(s)$, l'automorphisme 
$q^{-1/\d} w F^*$ de $X(\Tb_1^*) \otimes_\ZM \RM$ 
est d'ordre fini $o_w$. Soit $N$ le plus petit commun multiple de 
$(o_w)_{w \in W(s)}$. On notera $\phi_1$ un g\'en\'erateur 
d'un groupe cyclique $<\phi_1>$ d'ordre $N$ et nous ferons 
agir $\phi_1$ sur $X(\Tb_1)^* \otimes_\ZM \RM$ comme $q^{-1/\d} F^*$. 
Ainsi, $\phi_1$ normalise $W(s)$ et, si $w \in W(s)$, on a 
$\phi_1 w \phi_1^{-1}=F^*(w)$. On peut donc d\'efinir le produit 
semi-direct $W(s) \rtimes <\phi_1>$. 

Le but de cette sous-section 
est de relier le groupe $A_{\Gb^*}(s)$ aux groupes $\Ker' i^*$, 
$\ZC(\Gb)^\we$ et $(\Gbt^F/\Gb^F)^\we$. Tout d'abord, consid\'erons l'application 
$$\fonction{\ph_s}{C_{\Gb^*}(s)}{\Ker' i^*}{g}{[\gti,\sti]=
\gti\sti\gti^{-1}\sti^{-1}}$$
o\`u, pour tout $g \in C_{\Gb^*}(s)$, $\gti$ d\'esigne un \'el\'ement de $\Gbt^*$ 
tel que $i^*(\gti)=g$. Alors $\ph_s(g)$ ne d\'epend ni du choix de $\gti$
ni du choix de $\sti$. Puisque $\Ker' i^*$ est central, $\ph_s$ 
est un \mor de groupes 
et le noyau de $\ph_s$ est $i^*(C_{\Gbt^*}(\sti))$. Mais, d'apr\`es 
le th\'eor\`eme \ref{cgs connexe} 
et par exemple \cite[Proposition 2.3]{dmbook}, 
on a $i^*(C_{\Gbt^*}(\sti))=C_{\Gb^*}^\ci(s)$ donc $\ph_s$ induit un morphisme 
injectif de groupes encore not\'e 
\equat
\ph_s : A_{\Gb^*}(s) \injto \Ker' i^*.
\endequat
Ce morphisme commute \`a l'action de $F^*$. Comme cons\'equence, on obtient le

\bigskip

\lemme{ags}
{\it Le groupe $A_{\Gb^*}(s)$ est ab\'elien et, via $\ph_s$, 
$$A_{\Gb^*}(s)\simeq\{z \in \Ker i^*~|~\sti\text{ et }\sti z
\text{ sont conjugu\'es dans $\Gbt^*$}\}.$$}

\bigskip

Par dualit\'e, $\ph_s$ induit un morphisme surjectif 
$\phh_s : (\Ker' i^*)^\we \to A_{\Gb^*}(s)^\we$ et, 
par composition avec les \isos 
$\o$, $\o^0$, $\o^1$, $\omeh$, $\omeh^0$ et $\omeh^1$, on obtient des morphismes 
\equat
\begin{array}{c}
\o_s : A_{\Gb^*}(s) \injto \ZC(\Gb)^\we, \\
\\
\o_s^0 : A_{\Gb^*}(s)^{F^*} \injto H^1(F,\ZC(\Gb))^\we, \\
\\
\o_s^1 : H^1(F^*,A_{\Gb^*}(s)) \longto (\ZC(\Gb)^F)^\we,\\
\\
\omeh_s : \ZC(\Gb) \twoheadrightarrow A_{\Gb^*}(s)^\we, \\
\\
\omeh^0_s : H^1(F,\ZC(\Gb)) \twoheadrightarrow (A_{\Gb^*}(s)^{F^*})^\we, \\
\\
\omeh^1_s : \ZC(\Gb)^F \longto H^1(F^*,A_{\Gb^*}(s))^\we.
\end{array}
\endequat
Les \mors $\o_s$ et $\o_s^0$ sont injectifs tandis que les \mors $\omeh_s$ 
et $\omeh_s^0$ sont surjectifs. 
On ne peut cependant rien dire en g\'en\'eral concernant les morphismes 
$\o_s^1$ et $\omeh_s^1$.

\bigskip

\soussection{Sous-groupes de Levi} 
Soit $\Lb$ un \levi $F$-stable de $\Gb$ et soit $\Lb^*$ un 
\levi $F^*$-stable de $\Gb^*$ dual de $\Lb$. Supposons que 
$s \in \Lb^{*F^*}$. Le morphisme injectif $C_{\Lb^*}(s) \injto C_{\Gb^*}(s)$ 
induit un morphisme injectif $A_{\Lb^*}(s) \injto A_{\Gb^*}(s)$. En effet, 
le noyau du morphisme naturel $C_{\Lb^*}(s) \to A_{\Gb^*}(s)$ est 
$C_{\Gb^*}^\circ(s) \cap \Lb^*$. Mais, puisque $\Lb^*=C_{\Gb^*}(\Zb(\Lb^*)^\circ)$, 
on a $C_{\Gb^*}^\circ(s) \cap \Lb^*= C_{C_{\Gb^*}^\circ(s)}(\Zb(\Lb^*)^\circ)$. 
Le r\'esultat d\'ecoule alors de ce que le centralisateur d'un tore 
dans un groupe connexe est connexe \cite[corollaire 11.12]{borel}. 
Il est d'autre part facile de voir que le diagramme 
\equat\label{als ags zg}
\diagram
A_{\Lb^*}(s) \rto \dto & \ZC(\Lb)^\we \dto \\
A_{\Gb^*}(s) \rto & \ZC(\Gb)^\we,
\enddiagram
\endequat
est commutatif. De plus, toutes les applications sont injectives. 
En effet, l'injectivit\'e de l'application verticale de droite r\'esulte 
de la proposition \ref{ZG ZL} et l'injectivit\'e de l'autre application 
verticale a \'et\'e discut\'ee ci-dessus. 

En prenant les points fixes sous $F$ et $F^*$ 
et en dualisant, on obtient un autre diagramme commutatif 
\equat\label{commutativite als}
\diagram
H^1(F,\ZC(\Gb)) \dto_{h_\Lb^1} \rto^{\omeh_s^0} & (A_{\Gb^*}(s)^{F^*})^\we\dto^\Res\\
H^1(F,\ZC(\Lb)) \rto_{\omeh_{\Lb,s}^0}  & (A_{\Lb^*}(s)^{F^*})^\we.
\enddiagram
\endequat

\bigskip

\soussection{Centralisateurs de sous-tores} 
Fixons maintenant un groupe fini $H$ d'automorphismes de $X(\Tb_1^*) \otimes_\ZM \RM$ 
tel que, si $h \in H$ et $\a \in \Phi_1$, alors $h(\a)$ soit un multiple 
r\'eel d'une racine de $\Phi_1$. Cela d\'efinit une action de $H$ 
par permutation sur les racines que nous noterons $h * \a$~: elle 
est d\'efinie ainsi
$$h*\a \in \Phi_1 \cap \RM_+^\times h(\a).$$
Notons $Y^H$ le $\RM$-espace vectoriel des points fixes de 
$H$ dans son action sur $Y(\Tb_1^*) \otimes_\ZM \RM$. Posons
$$\Phi_H=\{\a \in \Phi_1~|~\forall~v \in Y^H,~<\a,v>_{\Tb_1^*} = 0 \}.$$
Alors $\Phi_H$ est le syst\`eme de racines relativement \`a $\Tb_1^*$ 
d'un \levi $\Lb^*$ de $\Gb^*$. Soit $\Lb$ le \levi de $\Gb$ 
dont le syst\`eme de coracines est $\Phi_H$. 
Nous noterons $W_\Lb$ le groupe de Weyl de 
$\Lb^*$ relativement \`a $\Tb^*$, $W_\Lb(s)$ le centralisateur de $s$ dans 
$W_\Lb$, $W_\Lb^\circ(s)$ le groupe de Weyl de 
$C_{\Lb^*}^\circ(s)$ relativement \`a $\Tb_1^*$. Alors

\bigskip

\proposition{lll}
{\it Supposons que $H$ stabilise $\Phi_s^+$ (pour l'action $*$). Alors~:
\begin{itemize}
\itemth{a} $C_{\Lb^*}^\circ(s)=\Tb_1^*$, \cad $W_\Lb^\circ(s) =1$.

\itemth{b} $W_\Lb(s)^H \subset A_{\Gb^*}(s)^H$. 

\itemth{c} $W_\Lb(s)^H$ commute avec $W^\circ(s)^H$. 
\end{itemize}}

\bigskip

\proof (a) Soit $\a \in \Phi_H \cap \Phi_s^+$ et soit $y \in Y(\Tb_1^*)$. Alors 
$<\a,\sum_{h \in H} h(y) >_{\Tb_1^*} =0$ par d\'efinition de $\Phi_H$. 
Mais, 
$$<\a,\sum_{h \in H} h(y) >_{\Tb_1^*} = <\sum_{h \in H} h(\a),y >_{\Tb_1^*}.$$
Donc, $\sum_{h \in H} h(\a)=0$, ce qui est impossible car $H$ stabilise 
$\Phi_s^+$ (dans son action $*$). Donc $\Phi_H \cap \Phi_s = \vide$, 
ce qui montre (a).

\medskip

Avant de continuer, introduisons quelques notations. Il existe $v \in Y(\Tb_1^*)$ 
tel que $<\a,v>_{\Tb_1^*} > 0$ pour tout $\a \in \Phi_s^+$. Posons 
$v_0 = \sum_{h \in H} h(v)$. 
Alors $v_0 \in Y^H$ et $<\a,v_0>_{\Tb_1^*} > 0$ pour tout $\a \in \Phi_s^+$ car 
$H$ stabilise $\Phi_s^+$ (via l'action $*$). 

\medskip

(b) Soit $w \in W_\Lb(s)^H$ et soit $\a \in \Phi_s^+$. 
Posons $f(\a)=\sum_{h \in H} h(\a)$. Alors $<f(\a),v_0>_{\Tb_1^*} > 0$. 
D'autre part, $f(w(\a))=w(f(\a))$ et, puisque $w(v_0)=v_0$, on a 
$<f(w(\a)),v_0>_{\Tb_1^*} > 0$. Cela montre que $<w(\a),v_0 >_{\Tb_1^*} > 0$ et donc 
que $w(\a) \in \Phi_s^+$. En d'autres termes, $w \in A_{\Gb^*}(s)$.

\medskip

(c) D'apr\`es (b), le groupe $W_\Lb(s)^H$ normalise $W^\circ(s)^H$. 
et $W_\Lb(s)^H \cap W^\circ(s)^H = 1$. 
D'autre part $W^H$ normalise $W_\Lb$, donc $W(s)^H$ normalise $W_\Lb(s)^H$. 
D'o\`u (c).\fin

\bigskip

\soussection{\'El\'ements semi-simples cuspidaux\label{sous semi cus}} 
Un \'el\'ement semi-simple $s \in \Gb^*$ (\resp $s \in \Gb^{*F^*}$) 
est dit {\it g\'eom\'etriquement cuspidal} (\resp {\it rationnellement 
cuspidal}) s'il existe un \'el\'ement $a \in A_{\Gb^*}(s)$ 
(\resp $a \in A_{\Gb^*}(s)^{F^*}$) tel que 
$\omega_s(a) \in \ZC_\cus^\we(\Gb)$. S'il est n\'ecessaire de 
pr\'eciser le groupe ambiant, nous parlerons d'\'el\'ements 
g\'eom\'etriquement ou rationnellement $\Gb^*$-cuspidaux. 
La d\'efinition pr\'ec\'edente est la m\^eme que 
\cite[d\'efinition 1.4.1]{bonnafe torsion}. 
Nous commen\c{c}ons cette sous-section par une caract\'erisation des 
\'el\'ements semi-simples g\'eom\'etriquement cuspidaux. 

\medskip

\proposition{caracterisation semi simple cuspidal}
{\it Soit $a \in A_{\Gb^*}(s)$. Alors $\o_s(a) \in \ZC_\cus^\we(\Gb)$ 
si et seulement si $a \not\in A_{\Lb^*}(s)$ pour tout \levi 
propre $\Lb^*$ de $\Gb^*$ contenant $s$ (o\`u on rappelle que 
$A_{\Lb^*}(s)$ peut \^etre vu naturellement comme un sous-groupe 
de $A_{\Gb^*}(s)$).}

\medskip

\proof S'il existe un sous-groupe de Levi propre $\Lb^*$ de $\Gb^*$ 
contenant $s$ tel que $a \in A_{\Lb^*}(s)$, il r\'esulte de 
la commutativit\'e du diagramme \ref{als ags zg} que $\o_s(a)$ 
n'appartient pas \`a $\ZC_\cus^\we(\Gb)$. 

R\'eciproquement, supposons que $a \not\in A_{\Lb^*}(s)$ pour 
tout \levi propre $\Lb^*$ de $\Gb^*$ contenant $s$. Nous devons 
montrer que $\o_s(a) \in \ZC_\cus^\we(\Gb)$. Pour cela, 
compte tenu de \cite[1.4.6 et 1.4.7]{bonnafe torsion}, 
nous pouvons supposer que $\Gb$ est semi-simple, 
simplement connexe et quasi-simple. D'apr\`es 
\cite[proposition 3.14 (b)]{bonnafe quasi}, 
$a$ peut-\^etre vu comme un automorphisme du diagramme de Dynkin affine 
de $\Gb$. L'hypoth\`ese implique que $a$ agit transitivement sur le 
diagramme de Dynkin affine (sinon $a$ appartient \`a un sous-groupe 
parabolique propre de $W$). Un examen de la classification des 
syst\`emes de racines montre que cela ne peut arriver que lorsque 
$\Gb$ est de type $A_n$. Dans ce cas, $a$ est d'ordre $n+1$ et 
$\ZC(\Gb)$ est aussi d'ordre $n+1$. Par cons\'equent, $\o_s(a)$
est injectif et le r\'esultat est imm\'ediat.\fin

\bigskip

Nous rappelons certaines des propri\'et\'es des \'el\'ements semi-simples 
cuspidaux d\'emontr\'ees dans 
\cite[\SEC 1.4]{bonnafe torsion}. Rappelons que $s$ est dit 
{\it isol\'e} (\resp {\it quasi-isol\'e}) si $C_{\Gb^*}^\circ(s)$ 
(\resp $C_{\Gb^*}(s)$) n'est contenu dans aucun sous-groupe 
de Levi propre de $\Gb^*$.

\bigskip

\proposition{semisimple cuspidal prop}
{\it Soit $s$ un \'el\'ement semi-simple g\'eom\'etriquement cuspidal 
de $\Gb^*$. Alors~:
\begin{itemize}
\itemth{a} Toutes les composantes quasi-simples de $\Gb^*$ sont de type $A$.

\itemth{b} Si $\Lb^*$ est un \levi propre de $\Gb^*$ contenant $s$, 
alors le morphisme injectif $A_{\Lb^*}(s) \injto A_{\Gb^*}(s)$ 
n'est pas surjectif.

\itemth{c} Si $s \in \Gb^{*F^*}$ est de plus rationnellement cuspidal, alors 
si $\Lb^*$ est un \levi $F^*$-stable propre de $\Gb^*$ contenant $s$, 
alors le morphisme injectif $A_{\Lb^*}(s)^{F^*} \injto A_{\Gb^*}(s)^{F^*}$ 
n'est pas surjectif.

\itemth{d} $s$ est quasi-isol\'e et r\'egulier.

\itemth{e} L'application $\omega_s : A_{\Gb^*}(s) \to \ZC(\Gb)^\we$ est 
un isomorphisme.
\end{itemize}}

\bigskip

\proof (a) d\'ecoule de la proposition \ref{cuspidal prop} (c). 
(b) d\'ecoule de la commutativit\'e du diagramme \ref{als ags zg}. 
(c) d\'ecoule de la proposition \ref{caracterisation semi simple cuspidal}. 
Le fait que $s$ est quasi-isol\'e d\'ecoule de (b). La preuve de la 
r\'egularit\'e de $s$ est faite dans \cite[lemme 3.2.9]{bonnafe mackey}. 
D'o\`u (d). Pour (e), voir \cite[proposition 1.4.9]{bonnafe torsion}.\fin

\bigskip

Fixons maintenant $a \in A_{\Gb^*}(s)$ et posons 
$$\Lb_{s,a}^* = C_{\Gb^*}(((\Tb_1^*)^a)^\circ).$$
C'est un \levi de $\Gb^*$ contenant $\Tb_1^*$. Soit $\Lb_{s,a}$ 
un \levi de $\Gb$ dual de $\Lb_{s,a}^*$. Notons que 
$a \in A_{\Lb_{s,a}^*}(s)$.

\bigskip

\proposition{cuspidal las}
{\it Si $a \in A_{\Gb^*}(s)$, alors 
$\omega_{\Lb_{s,a},s}(a) \in \ZC_\cus^\we(\Lb_{s,a})$. Donc 
$s$ est g\'eom\'etriquement cuspidal dans $\Lb_{s,a}^*$.}

\bigskip

\proof Pour montrer la proposition \ref{cuspidal las}, on peut travailler 
dans $\Lb_{s,a}$, \cad que l'on peut supposer que $\Lb_{s,a}=\Gb$. 
L'hypoth\`ese signifie donc que que $((\Tb_1^*)^a)^\circ=\Zb(\Gb^*)^\circ$, 
\cad que $a$ 
est un \'el\'ement cuspidal \cite[d\'efinition 3.1.1]{geck livre} 
de $W$. Par suite, si $\Mb$ est un \levi propre de $\Gb$ contenant $\Tb_1$ 
et si $\Mb^*$ est un \levi de $\Gb^*$ contenant $\Tb_1^*$ dual de $\Mb$, alors 
$a \not\in W_\Mb(s)$. La proposition \ref{caracterisation semi simple cuspidal} 
montre alors que $\omega_s(a) \in \ZC_\cus^\we(\Gb)$.\fin

\bigskip

\corollaire{ah ah ah}
{\it Soit $a \in A_{\Gb^*}(s)$. Alors~:
\begin{itemize}
\itemth{a} $C_{\Lb_{s,a}^*}^\circ(s)=\Tb_1^*$, \cad $W_{\Lb_{s,a}}^\circ(s) =1$.

\itemth{b} $W_{\Lb_{s,a}}(s)^a \subset A_{\Gb^*}(s)$. 

\itemth{c} $W_{\Lb_{s,a}}(s)^a$ commute avec $W^\circ(s)^a$. 

\itemth{d} $N_W(W_{\Lb_{s,a}})=W^a.W_{\Lb_{s,a}}$.

\itemth{e} $N_{W(s)}(W_{\Lb_{s,a}}) = W(s)^a$.
\end{itemize}}

\bigskip

\proof Le groupe $<a>$ stabilise $\Phi_s^+$. Par cons\'equent, 
(a), (b) et (c) d\'ecoulent de la proposition \ref{lll}. 
Montrons (d). Soit $w \in N_W(W_{\Lb_{s,a}})$. Par construction, 
$a$ est un \'el\'ement cuspidal \cite[d\'efinition 3.1.1]{geck livre} 
de $W_{\Lb_{s,a}}$. 
Par suite, $waw^{-1} \in W_{\Lb_{s,a}}$ et est conjugu\'e \`a $a$ sous $W$. 
Donc, d'apr\`es \cite[th\'eor\`eme 3.2.11]{geck livre}, il existe 
$x \in W_{\Lb_{s,a}}$ tel que 
$waw^{-1}=xax^{-1}$. Donc $w \in W^a.W_{\Lb_{s,a}}$. Cela montre 
que $N_W(W_{\Lb_{s,a}}) \subset W^a.W_{\Lb_{s,a}}$. L'inclusion r\'eciproque 
est imm\'ediate.

Montrons pour finir (e). Soit $w \in W(s)$ normalisant $W_{\Lb_{s,a}}$. 
Quitte \`a multiplier $w$ par un \'el\'ement de $A_{\Gb^*}(s)$ 
(et en utilisant (d)), on peut supposer que $w \in W^\circ(s)$. 
D'apr\`es (d), on a $awa^{-1} w^{-1} \in W_{\Lb_{s,a}}$. Mais, d'autre part, 
$awa^{-1} w^{-1} \in W^\circ(s)$. Donc, d'apr\`es (a), $awa^{-1} w^{-1}=1$.\fin

\bigskip

\noindent{\sc Remarque - } \`A ce jour, la preuve de 
\cite[th\'eor\`eme 3.2.11]{geck livre} n\'ecessite la classification des 
groupes de Coxeter finis et de leurs classes cuspidales, ce qui est 
d\'esagr\'eable. Il faut noter que l'analogue de 
\cite[th\'eor\`eme 3.2.11]{geck livre} est faux en g\'en\'eral pour les 
groupes de r\'eflexions complexes.\finl

\newpage

{\Large \part{Induction et restriction de Lusztig, 
s\'eries de Lusztig\label{chapitre lusztig}}}

\bigskip

Nous rappelons ici la construction de l'induction et de la restriction 
de Lusztig ainsi que la d\'efinition de s\'eries de Lusztig g\'eom\'etriques 
et rationnelles. Un des buts de ce chapitre est de fournir une preuve 
compl\`ete de la disjonction des s\'eries de Lusztig rationnelles 
(en partant de la disjonction des s\'eries g\'eom\'etriques des 
groupes \`a centre connexe). Cette preuve est largement esquiss\'ee 
dans \cite{luirr} et \cite[chapitre 14]{dmbook} mais dans ces deux cas, 
elle n'est pas tout-\`a-fait compl\`ete. Dans la section \ref{section cars T} 
nous rappelons quelques propri\'et\'es de la dualit\'e entre 
caract\`eres lin\'eaires d'un tore et \'el\'ements semi-simples du 
dual. Nous reprenons notamment un lemme d'Asai 
\cite[th\'eor\`eme 2.1.1]{asai} sur les caract\`eres centraux associ\'es 
\`a des \'el\'ements semi-simples g\'eom\'etriquement mais non 
rationnellement conjugu\'es. Dans la section \ref{section induction restriction}, 
nous \'etudions les actions du centre sur les applications de Lusztig. 
Nous rappelons aussi dans quels cas la formule de Mackey est connue. 
La section \ref{section series} est consacr\'ee \`a la disjonction des 
s\'eries de Lusztig ainsi qu'\`a quelques-unes de leurs propri\'et\'es 
(caract\`ere central, compatibilit\'e \`a l'induction de Lusztig...). 

\bigskip

\section{Caract\`eres lin\'eaires de tores maximaux\label{section cars T}}~

\medskip

\soussection{Dualit\'e}
Soit $\nabla(\Gb,F)$ l'ensemble des couples $(\Tb,\th)$ o\`u $\Tb$ est un \tor 
$F$-stable de $\Gb$ et $\th : \Tb^F \to \qlb^\times$ est un caract\`ere lin\'eaire. 
Soit $\nabla^*(\Gb,F)$ l'ensemble des couples $(\Tb^*,s)$ o\`u $\Tb^*$ est 
un \tor $F^*$-stable de $\Gb^*$ et $s \in \Tb^{*F^*}$. Le choix des morphismes 
$\imath$ et $\jmath$ d\'efinis dans \SEC\ref{sous groupes} 
induit une bijection \cite[proposition 13.13]{dmbook}
\equat\label{binabla}
\nabla(\Gb,F)/\Gb^F \longto \nabla^*(\Gb,F)/\Gb^{*F^*}.
\endequat
Si $(\Tb,\th) \in \nabla(\Gb,F)$ et $(\Tb^*,s) \in \nabla^*(\Gb,F)$, nous \'ecrivons 
$(\Tb,\th) \doublefleche{\Gb} (\Tb^*,s)$ pour dire que $(\Tb,\th)$ et $(\Tb^*,s)$ 
sont associ\'es par la bijection \ref{binabla}.

\bigskip

\lemme{tenseur T}
{\it Soient $(\Tbt,\thet) \in \nabla(\Gbt,F)$, $(\Tbt^*,\sti) \in \nabla^*(\Gbt,F)$ 
et $z \in (\Ker i^*)^{F^*}$. Voyons $\zha^\Gbt$ comme un \car 
lin\'eaire de $\Tbt^F$ par restriction depuis $\Gbt^F$. Alors 
$(\Tbt,\thet) \doublefleche{\Gbt} (\Tbt^*,\sti)$ \ssi
$(\Tbt,\thet \zha^\Gbt) \doublefleche{\Gbt} (\Tbt^*,\sti z)$.}

\bigskip

\proof Claire.\fin

\bigskip

\soussection{Restriction} Si $(\Tbt,\thet) \in \nabla(\Gbt,F)$, nous posons 
$$\RES_\Gb^\Gbt(\Tbt,\thet)=
(\Tbt \cap \Gb, \Res_{\Tbt^F \cap \Gb^F}^{\Tbt^F} \thet) \in \nabla(\Gb,F).$$
De m\^eme, si $(\Tbt^*,\sti) \in \nabla^*(\Gbt,F)$, nous posons 
$$\ORES_\Gb^\Gbt(\Tbt^*,\sti) = (i^*(\Tbt^*),i^*(\sti)) 
\in \nabla^*(\Gb,F).$$
Le lemme suivant se d\'emontre en revenant \`a la d\'efinition de la dualit\'e 
entre les tores.

\bigskip

\lemme{mise au point}
{\it \begin{itemize}
\itemth{a} Soient $(\Tbt,\thet)$ et $(\Tbt^*,\sti)$ deux \'el\'ements 
de $\nabla(\Gbt,F)$ et $\nabla^*(\Gbt,F)$ respectivement. On suppose 
que $(\Tbt,\thet) \doublefleche{\Gbt} (\Tbt^*,\sti)$. Alors 
$\RES_\Gb^\Gbt (\Tbt,\thet) \doublefleche{\Gb} \ORES_\Gb^\Gbt(\Tbt^*,\sti)$. 

\itemth{b} Soient $(\Tb,\th)$ et $(\Tb^*,s)$ deux \'el\'ements 
de $\nabla(\Gb,F)$ et $\nabla^*(\Gb,F)$ respectivement. 
On suppose que $(\Tb,\th) \doublefleche{\Gb} (\Tb^*,s)$. On pose 
$$\Tbt=\Tb.\Zb(\Gbt)\quad\quad\quad{\mathit{et}}
\quad\quad\quad\Tbt^*=i^{*-1}(\Tb^*)$$
et soit $\sti$ un \ele semi-simple de $\Gbt^{*F^*}$ tel que $i^*(\sti)=s$ 
(Un tel \ele $\sti$ existe d'apr\`es le corollaire \ref{surjecte}). 
Alors il existe une extension $\thet$ de $\th$ \`a $\Tbt^F$ telle que 
$(\Tbt,\thet) \doublefleche{\Gbt} (\Tbt^*,\sti)$.
\end{itemize}}

\bigskip

\soussection{Conjugaison g\'eom\'etrique et rationnelle} 
La d\'efinition suivante a \'et\'e pos\'ee par Deligne et Lusztig 
\cite[d\'efinition 5.5]{delu}.

\bigskip

\definition{conjugaison geo et ratio}
{\it Soient $(\Tb_1,\th_1)$ et $(\Tb_2,\th_2)$ deux \'el\'ement de $\nabla(\Gb,F)$ 
et soient $(\Tb_1^*,s_1)$ et $(\Tb_2^*,s_2)$ deux \'el\'ements de $\nabla^*(\Gb,F)$ 
tels que $(\Tb_k,\th_k) \doublefleche{\Gb} (\Tb_k^*,s_k)$ pour tout $k \in \{1,2\}$. 
On dit que $(\Tb_1,\th_1)$ et $(\Tb_2,\th_2)$ sont 
{\bfit g\'eom\'etriquement conjugu\'es} (\resp {\bfit sont dans la m\^eme 
s\'erie rationnelle}) si $s_1$ et $s_2$ sont g\'eom\'etriquement 
conjugu\'es (respectivement rationnellement conjugu\'es), \cad 
si ils sont conjugu\'es dans $\Gb^*$ (\resp $\Gb^{*F^*}$).}

\bigskip

Soit $s$ un \'el\'ement semi-simple de $\Gb^{*F^*}$. Nous noterons 
$(s)$, ou $(s)_{\Gb^*}$ s'il y a ambigu\"\i t\'e, la classe de conjugaison 
g\'eom\'etrique de $s$. De m\^eme, nous noterons $[s]$, 
ou $[s]_{\Gb^{*F^*}}$, la classe de conjugaison rationnelle de $s$. 
Soit $\nabla^*(\Gb,F,(s))$ 
(\resp $\nabla^*(\Gb,F,[s])$) l'ensemble des couples 
$(\Tb^*,s') \in \nabla^*(\Gb,F)$ 
tels que $s'$ est g\'eom\'etriquement (\resp rationnellement) conjugu\'e \`a $s$. 
Par dualit\'e, nous notons $\nabla(\Gb,F,(s))$ (\resp $\nabla(\Gb,F,[s])$) 
l'ensemble des couples $(\Tb,\th) \in  \nabla(\Gb,F)$ associ\'es  
aux couples appartenant \`a $\nabla^*(\Gb,F,(s))$ 
(\resp $\nabla^*(\Gb,F,[s])$) par la bijection 
\ref{binabla}. 

\bigskip

\corollaire{mise au point coro}
{\it Soit $\sti \in \Gbt^{F*^*}$ un \'el\'ement semi-simple et soit 
$s=i^*(\sti)$. Si $(\Tbt,\thet) \in \nabla(\Gbt,F,(\sti))$, alors 
$\RES_\Gb^\Gbt (\Tbt,\thet) \in \nabla(\Gb,F,[s])$.}

\bigskip

\soussection{Autres caract\'erisations des s\'eries g\'eom\'etriques et 
rationnelles} 
Concernant la conjugaison g\'eom\'etrique des couples 
appartenant \`a $\nabla(\Gb,F)$, le lemme suivant nous 
fournit une d\'efinition \'equivalente \cite[proposition 5.4]{delu}~:

\bigskip

\Lemme{Deligne-Lusztig}{conjugaison geometrique}
{\it Soient $(\Tb_1,\th_1)$ et $(\Tb_2,\th_2)$ deux \'el\'ements de $\nabla(\Gb,F)$. 
Alors $(\Tb_1,\th_1)$ et $(\Tb_2,\th_2)$ sont g\'eom\'etriquement conjugu\'es 
si et seulement si il existe un entier naturel non nul $n$ et un \'el\'ement 
$g \in \Gb^{F^n}$ tel que $\Tb_2=\lexp{g}{\Tb_1}$ et 
$$\th_2 \ci N_{F^n/F} = \lexp{g}{(\th_1 \ci N_{F^n/F})}$$
o\`u $N_{F^n/F} : \Tb_k^{F^n} \to \Tb_k^F$, $t \mapsto tF(t)\dots F^{n-1}(t)$ 
est la norme de $F^n$ \`a $F$ ($k=1$ o\`u $2$).}

\bigskip

Nous allons maintenant utiliser le groupe $\Gbt$ (et le fait que son centre est 
connexe) pour donner une autre d\'efinition des s\'eries rationnelles. 
Tout d'abord, notons que les s\'eries g\'eom\'etriques et rationnelles 
co\"{\i}ncident dans $\Gbt$ \`a cause du th\'eor\`eme de 
Steinberg (voir th\'eor\`eme \ref{cgs connexe})~:

\medskip

\proposition{geo=ratio} 
{\it Deux \'el\'ements semi-simples de $\Gbt^\OF$ 
sont g\'eom\'etriquement conjugu\'es si et 
seulement si ils sont rationnellement conjugu\'es.}

\bigskip

\corollaire{geo et ratio}
{\it Soient $(\Tbt_1,\thet_1)$ et $(\Tbt_2,\thet_2)$ deux \'el\'ements de 
$\nabla(\Gbt,F)$. Alors $(\Tbt_1,\thet_1)$ et $(\Tbt_2,\thet_2)$ 
sont g\'eom\'etriquement conjugu\'es si et seulement si 
ils appartiennent \`a la m\^eme s\'erie rationnelle.}

\bigskip

La prochaine proposition fournit, en termes du groupe $\Gbt$, 
un outil pratique pour d\'eterminer si deux \'el\'ements semi-simples 
de $\Gb^\OF$ sont rationnellement conjugu\'es.

\bigskip

\proposition{ratio g = geo gti}
{\it Soient $s_1$ et $s_2$ deux \'el\'ements semi-simples de $\Gb^\OF$. Alors 
les assertions suivantes sont \'equivalentes~:
\begin{itemize}
\itemth{1} $s_1$ et $s_2$ sont rationnellement conjugu\'es~;

\itemth{2} Il existe des \'el\'ements semi-simples rationnellement 
conjugu\'es $\sti_1$ et $\sti_2$ dans $\Gbt^{*F^*}$ tels que 
$i^*(\sti_i)=s_i$ pour tout $i \in \{1,2\}$~;

\itemth{3} Il existe des \'el\'ements semi-simples g\'eom\'etriquement 
conjugu\'es $\sti_1$ et $\sti_2$ dans $\Gbt^{*F^*}$ tels que 
$i^*(\sti_i)=s_i$ pour tout $i \in \{1,2\}$.
\end{itemize}}

\bigskip

\proof D'apr\`es la proposition \ref{geo=ratio}, (2) et (3) sont \'equivalents. 
Il est par ailleurs clair que (2) implique (1). Il nous reste donc \`a d\'emontrer 
que (1) implique (2). 

Supposons donc qu'il existe $g \in \Gb^\OF$ tel que $s_2=gs_1g^{-1}$. 
Puisque $i^* : \Gbt^\OF \to \Gb^\OF$ est surjective (voir corollaire 
\ref{surjecte}), il existe $\sti_1 \in \Gbt^\OF$ et $\gti \in \Gbt^\OF$ 
tels que $i^*(\sti_1)=s_1$ et $i^*(\gti)=g$. Posons maintenant 
$\sti_2=\gti \sti_1 \gti^{-1}$. Alors $\sti_1$ et $\sti_2$ 
sont des \'el\'ements semi-simples rationnellement conjugu\'es 
de $\Gbt^\OF$ et $i^*(\sti_k)=s_k$ pour tout $k \in \{1,2\}$.\fin

\bigskip

\corollaire{ratio T = geo T}
{\it Soient $(\Tb_1,\th_1)$ et $(\Tb_2,\th_2)$ deux \'el\'ements de $\nabla(\Gb,F)$. 
Soit $\Tbt_k=\Tb_k.\Zb(\Gbt)$ (pour $k\in \{1,2\}$). Alors les assertions 
suivantes sont \'equivalentes~:
\begin{itemize}
\itemth{1} $(\Tb_1,\th_1)$ et $(\Tb_2,\th_2)$ appartiennent \`a la m\^eme 
s\'erie rationnelle~;

\itemth{2} Il existe des extensions $\thet_1$ et $\thet_2$ de $\th_1$ et $\th_2$ 
respectivement (\`a $\Tbt_1^F$ et $\Tbt_2^F$ respectivement) telles que 
$(\Tbt_1,\thet_1)$ et $(\Tbt_2,\thet_2)$ sont g\'eom\'etriquement conjugu\'es.
\end{itemize}}

\bigskip

\proof Cela r\'esulte imm\'ediatement de la proposition \ref{ratio g = geo gti} 
et du lemme \ref{mise au point} (b).\fin

\bigskip

\soussection{Caract\`eres centraux}
La proposition suivante explique comment varient les caract\`eres 
lin\'eaires $\Res_{\Zb(\Gb)^F}^{\Tb^F} \th$ 
lorsque $(\Tb,\th)$ parcourt une s\'erie g\'eom\'etrique ou 
rationnelle.

\bigskip

\proposition{centre}
{\it Soient $(\Tb_1,\th_1)$ et $(\Tb_2,\th_2)$ deux \'el\'ements de $\nabla(\Gb,F)$. 
\begin{itemize}
\itemth{a} Si $(\Tb_1,\th_1)$ et $(\Tb_2,\th_2)$ sont g\'eom\'etriquement 
conjugu\'es, alors  
$$\Res_{\Zb(\Gb)^{\circ F}}^{\Tb^F_1} \th_1 = 
\Res_{\Zb(\Gb)^{\circ F}}^{\Tb^F_2} \th_2.$$ 

\itemth{b} Si $(\Tb_1,\th_1)$ et $(\Tb_2,\th_2)$ appartiennent \`a 
la m\^eme s\'erie rationnelle, alors  
$$\Res_{\Zb(\Gb)^F}^{\Tb^F_1} \th_1 = \Res_{\Zb(\Gb)^F}^{\Tb^F_2} \th_2.$$
\end{itemize}}

\bigskip

\proof (a) Si $(\Tb_1,\th_1)$ et $(\Tb_2,\th_2)$ 
sont g\'eom\'etriquement conjugu\'es, alors, 
d'apr\`es le lemme \ref{conjugaison geometrique}, 
il existe un entier naturel non nul $n$ et un \'el\'ement $g \in \Gb^{F^n}$ 
tel que $\Tb_2=\lexp{g}{\Tb_1}$ et 
$$\th_2 \ci N_{F^n/F} = \lexp{g}{(\th_1 \ci N_{F^n/F})}.$$
Soit $z \in \Zb(\Gb)^{\circ F}$. Puisque le groupe 
$\Zb(\Gb)^\circ$ est connexe, il existe $z' \in \Zb(\Gb)^{\circ F^n}$ tel que  
$N_{F^n/F}(z')=z$. Donc 
$$\th_2(z)=\th_2(N_{F^n/F}(z'))=\th_1(N_{F^n/F}(g^{-1}z'g))=
\th_1(N_{F^n/F}(z'))=\th_1(z)$$
car $z'$ est central. Cela montre (a).

\medskip

(b) Si $(\Tb_1,\th_1)$ et $(\Tb_2,\th_2)$ appartiennent \`a la 
m\^eme s\'erie rationnelle, alors, d'apr\`es le corollaire \ref{ratio T = geo T}, 
il existe 
$(\Tbt_1,\thet_1)$ et $(\Tbt_2,\thet_2)$ dans $\nabla(\Gbt,F)$ tels que  
$$\RES_\Gb^\Gbt(\Tbt_k,\thet_k)=(\Tb_k,\th_k)$$
($k \in \{1,2\}$) et tels que $(\Tbt_1,\thet_1)$ et $(\Tbt_2,\thet_2)$ 
sont g\'eom\'etriquement conjugu\'es. Par cons\'equent, d'apr\`es (a) 
(appliqu\'e au groupe $\Gbt$) et puisque $\Zb(\Gbt)$ est connexe, on a 
$$\Res_{\Zb(\Gbt)^F}^{\Tbt_1^F} \thet_1 = \Res_{\Zb(\Gbt)^F}^{\Tbt_2^F} \thet_2.$$
Donc (b) d\'ecoule de ce que $\Zb(\Gb)=\Zb(\Gbt) \cap \Gb$.\fin

\bigskip

Si $s$ est un \'el\'ement semi-simple de $\Gb^{*F^*}$, nous notons 
$\sha : \Zb(\Gb)^F \to {\overline{\QM}}_\ell^\times$ 
le caract\`ere lin\'eaire d\'efini par 
\equat\label{def sha}
\sha=\Res_{\Zb(\Gb)^F}^{\Tb^F} \th
\endequat 
pour tout $(\Tb,\th) \in \nabla(\Gb,F,[s])$~; $\sha$ est bien d\'efini 
d'apr\`es la proposition \ref{centre} (b). Nous notons d'autre part 
$\sha^\ci$ la restriction de $\sha$ \`a $\Zb(\Gb)^{\circ F}$~; le 
\car lin\'eaire $\sha^\circ$ peut aussi \^etre d\'efini par l'\'egalit\'e 
\equat\label{def sha0}
\sha^\ci=\Res_{\Zb(\Gb)^{\circ F}}^{\Tb^F} \th
\endequat
pour tout $(\Tb,\th) \in \nabla(\Gb,F,(s))$ (voir proposition \ref{centre} (a)).

\medskip

Soit maintenant $\a \in H^1(F^*,A_{\Gb^*}(s))$. Si $g_\a \in \Gb^*$ est tel 
que $g_\a^{-1}F(g_\a) \in C_{\Gb^*}(s)$ et repr\'esente $\a$, alors 
$s_\a = g_\a sg_\a^{-1} \in \Gb^{*F^*}$ est g\'eom\'etriquement conjugu\'e 
\`a $s$ et $[s_\a]$ ne d\'epend que de $\a$ et non pas 
du choix de $g_\a$. De plus, $(s_\a)_{\a \in H^1(F^*,A_{\Gb^*}(s))}$ est 
une famille de repr\'esentants des classes de $\Gb^{*F^*}$-conjugaison 
contenues dans $(s)_{\Gb^*}^{F^*}$. D'autre part, d'apr\`es la 
proposition \ref{centre} (a), $\sha_\a \sha^{-1}$ est un caract\`ere lin\'eaire 
du groupe $\ZC(\Gb)^F$. Il est donn\'e par la formule suivante~:

\bigskip

\Lemme{Asai}{asai} 
{\it Si $\a \in H^1(F^*,A_{\Gb^*}(s))$, alors $\sha_\a \sha^{-1} = \o_s^1(\a)$.}

\bigskip

\noindent{\sc Rappel - } Le morphisme 
$\o_s^1 : H^1(F^*,A_{\Gb^*}(s)) \to (\ZC(\Gb)^F)^\we$ a \'et\'e d\'efini 
en \SEC\ref{soussection ags}.\finl

\bigskip

\proof Le lemme \ref{asai} est d\'emontr\'e dans \cite[th\'eor\`eme 2.1.1]{asai} 
lorsque $F$ agit trivialement sur le groupe $\ZC(\Gb)$ (et lorsque $F$ est un 
endomorphisme de Frobenius). L'essentiel de sa preuve s'applique 
ici encore mais nous pr\'ef\'erons la rappeler dans ce cadre l\'eg\`erement 
plus g\'en\'eral.

\medskip

Soit $\Tbt^*$ un tore maximal $F^*$-stable de $\Gbt^*$ contenant $\sti$ 
et soit $\Tbt$ un \tor $F$-stable de $\Gbt$ dual de $\Tbt^*$. Posons 
$\Tb = \Tbt \cap \Gb$ et $\Tb^*=i^*(\Tbt^*)$. Soit $\a \in H^1(F^*,A_{\Gb^*}(s))$ 
et soit $a$ un \'el\'ement de $C_{\Gb^*}(s)$ representant $\a$. 
On peut supposer que $a$ normalise $\Tb^*$ et que  
$g_\a^{-1}F(g_\a)=a$. Alors $s_\a$ peut \^etre vu comme un \'el\'ement de 
$\Tb^*$~: plus pr\'ecis\'ement, $s_\a \in \Tb^{*aF^*}$. 

Choisissons un entier naturel non nul $n$ tel que $F^n$ soit un endomorphisme 
d\'eploy\'e de $\Tbt$ sur un corps fini \`a $q$ \'el\'ements et tel que 
$(a F^*)^n=F^{*n}$ sur $\Tbt^*$. 
On note encore par $a$ l'\'el\'ement de $N_\Gb(\Tb)$ correspondant \`a 
$a \in N_{\Gb^*}(\Tb^*)$. 

Il existe $\xti \in Y(\Tbt^*)=X(\Tbt)$ tel que 
$$\sti=N_{F^n/F}(\xti)(\tilde{\imath}(\frac{1}{q-1})).$$
On pose $x=\Res_\Tb^\Tbt \xti$. Alors 
$$s=N_{F^n/F}(x)(\tilde{\imath}(\frac{1}{q-1})).$$
D'autre part, il existe $x_\a \in Y(\Tb^*)=X(\Tb)$ tel que 
$$s=N_{F^n/aF}(x_\a)(\tilde{\imath}(\frac{1}{q-1})).$$
Alors $\sha=\Res_{\Zb^F}^\Tb x$ et $\sha_\a=\Res_{\Zb^F}^\Tb x_\a$. 

Soit $z \in \Zb(\Gb)^F$. Alors il existe $t \in \Tb$ tel que 
$t^{-1}F^n(t)=z$. On pose $u=F(t)$. 
Puisque $F(z)=z$ on a $u^{-1}F^n(u)=z$. Soit 
$t_1=t^{-1}F(t)$ et $t_\a=t^{-1}\lexp{a}{F(t)}$. 
Alors $N_{F^n/F}(t_1)=N_{F^n/aF}(t_\a)=z$ et $t_1 \in \Tb^{F^n}$ et 
$t_\a \in \Tb^{(aF)^n}=\Tb^{F^n}$. En particulier, 
$$\sha(z)=N_{F^n/F}(x)(t_1)\quad\quad\quad{\mathrm{et}}\quad\quad\quad \sha_\a(z)=
N_{F^n/F}(x_\a)(t_\a).$$
Mais, 
$$N_{F^n/F}(x)(\tilde{\imath}(\frac{1}{q-1}))=
N_{F^n/aF}(x_\a)(\tilde{\imath}(\frac{1}{q-1}))$$
donc 
$$\Res_{\Tb^{F^n}}^\Tb N_{F^n/F}(x)=\Res_{\Tb^{F^n}}^\Tb N_{F^n/aF}(x_\a).$$
Par cons\'equent, on a
\eqna
\sha_\a(z)\sha(z)^{-1} &=& N_{F^n/F}(x)(t_\a t_1^{-1})\\
&=&N_{F^n/F}(x)(aua^{-1}u^{-1}) \\
&=&(\lexp{a}{N_{F^n/F}}(\xti) - N_{F^n/F}(\xti))(u).
\endeqna
D'autre part,
$$\ph_s(a)=\ati\sti\ati^{-1}\sti^{-1}=(\lexp{a}{N_{F^n/F}}(\xti) - N_{F^n/F}(\xti))
(\tilde{\imath}(\frac{1}{q-1}))$$
donc, d'apr\`es \ref{definition omega}, on a 
\eqna
\o_s^1(\a)(z)&=&\o_s^1(\a)(u^{-1}F^n(u)) \\
&=&(\lexp{a}{N_{F^n/F}}(\xti) - N_{F^n/F}(\xti))
(\imath^{-1}(\frac{1}{q^n-1})),
\endeqna
ce qui est le r\'esultat attendu.\fin

\def\resr{\lexp{*}{R}}

\bigskip

\section{Induction et restriction de Lusztig\label{section induction restriction}}~

\medskip

\soussection{D\'efinitions\label{soussection definition induction}} 
Soit $\Pb$ un \para de $\Gb$ et supposons que 
$\Pb$ poss\`ede un \levic $F$-stable $\Lb$. Soit $\Ub$ le radical unipotent 
de $\Pb$. Posons, suivant Lusztig \cite{lufini},
$$\Yb_\Ub^\Gb=\{g \in \Gb~|~g^{-1}F(g) \in \Ub\}.$$
Alors $\Gb^F$ agit sur $\Yb_\Ub^\Gb$ par translations \`a gauche 
tandis que $\Lb^F$ agit par translations \`a droite. Par cons\'equent, 
$$H^*_c(\Yb_\Ub^\Gb)=\sum_{k \ge 0} (-1)^k H^k_c(\Yb_\Ub^\Gb,\qlb)$$
est un $\Gb^F$-module-$\Lb^F$ virtuel et 
$$H^*_c(\Yb_\Ub^\Gb)^\vee=\sum_{k \ge 0} (-1)^k H^k_c(\Yb_\Ub^\Gb,\qlb)^\vee$$
est un $\Lb^F$-module-$\Gb^F$ virtuel. Nous noterons 
$$\fonction{R_{\Lb \incl \Pb}^\Gb}{\ZM \Irr \Lb^F}{\ZM \Irr \Gb^F}{\L}{
H^*_c(\Yb_\Ub^\Gb) \otimes_{\qlb \Lb^F} \L}$$
$$\fonction{\resr_{\Lb \incl \Pb}^\Gb}{\ZM \Irr \Gb^F}{\ZM \Irr \Lb^F}{\G}{
H^*_c(\Yb_\Ub^\Gb) \otimes_{\qlb \Gb^F} \G}\leqno{\mathrm{et}}$$
les applications d'induction et de restriction de Lusztig respectivement. 
Elles s'\'etendent naturellement par lin\'earit\'e en applications entre espaces 
de fonctions centrales
$$R_{\Lb \incl \Pb}^\Gb : \Cent(\Lb^F) \longto \Cent(\Gb^F)$$
$$\resr_{\Lb \incl \Pb}^\Gb : \Cent(\Gb^F) \longto \Cent(\Lb^F).\leqno{\text{et}}$$
Ces deux applications sont adjointes l'une de l'autre par rapport 
aux produits scalaires $\langle,\rangle_{\Lb^F}$ et $\langle , \rangle_{\Gb^F}$. 

D'autre part, nous noterons 
$$\fonction{Q_{\Lb \incl \Pb}^\Gb}{\Gb^F_\uni \times \Lb^F_\uni}{\qlb}{(u,v)}{
\Tr((u,v), H^*_c(\Yb_\Ub^\Gb))}$$
la fonction de Green associ\'ee \`a la donn\'ee $(\Lb,\Pb,\Gb)$. 

\bigskip

\exemple{exemple harish} {\sc Induction de Harish-Chandra - } 
Si $\Pb$ est $F$-stable, alors $\Ub$ agit par translation 
\`a droite sur $\Yb_\Ub^\Gb$ et $\Yb_\Ub^\Gb/\Ub \simeq \Gb^F/\Ub^F$. 
Par cons\'equent, 
$$H^k_c(\Yb_\Ub^\Gb) \simeq \begin{cases} 
                            0 & {\text{si}}~k \not= 2 \dim \Ub \\
			    \qlb[\Gb^F/\Ub^F] & {\text{si}}~k=2 \dim \Ub.
			    \end{cases}$$
Par suite, $R_{\Lb \incl \Pb}^\Gb$ est le reflet, sur les groupes 
de Grothendieck, d'un vrai foncteur entre les cat\'egories de modules, 
foncteur qui sera toujours not\'e $R_{\Lb \incl \Pb}^\Gb$. 
De m\^eme pour $\resr_{\Lb \incl \Pb}^\Gb$. Ces deux foncteurs sont 
appel\'es respectivement {\it induction} et {\it restriction de Harish-Chandra}. 
Pratiquement toutes les formules d\'emontr\'ees dans ce chapitre 
au sujet des applications de Lusztig ont un sens en termes de modules 
lorsque l'on est en pr\'esence d'induction ou de restriction de 
Harish-Chandra (voir par exemple \ref{tauzg rlg}, \ref{tauzg resrlg}, 
propositions \ref{rlg G Gtilde} et \ref{tenseur}...).\finl

\bigskip

\soussection{Actions de ${\boldsymbol{\Zb(\Gb)^F}}$ et ${\boldsymbol{H^1(F,\ZC(\Gb))}}$} 
L'action de $\Zb(\Gb)^F$ par translation \`a gauche (ou \`a droite) 
sur $\Yb_\Ub^\Gb$ commute aux actions de $\Gb^F$ et $\Lb^F$. Par cons\'equent, 
si $z \in \Zb(\Gb)^F$, alors on a
\equat\label{tzg rlg}
t_z^\Gb \mathop{\circ} R_{\Lb \incl \Pb}^\Gb = 
R_{\Lb \incl \Pb}^\Gb \mathop{\circ} t_z^\Lb
\endequat
et
\equat\label{tzg resrlg}
t_z^\Lb \mathop{\circ} \resr_{\Lb \incl \Pb}^\Gb = 
\resr_{\Lb \incl \Pb}^\Gb \mathop{\circ} t_z^\Gb.
\endequat

Posons $\Lbt=\Lb.\Zb(\Gbt)$. 
Soit $a \in H^1(F,\ZC(\Gb))$ et soit $\lti_a \in \Lbt^F$ tel que 
$\s_\Lb^\Gb(a)=\lti_a \Lb^F \Zb(\Gbt)^F$. Alors la conjugaison 
par $\lti_a$ induit un automorphisme de $\Yb_\Ub^\Gb$, ce qui implique 
que
\equat\label{tauzg rlg}
\t_a^\Gb \mathop{\circ} R_{\Lb \incl \Pb}^\Gb = 
R_{\Lb \incl \Pb}^\Gb \mathop{\circ} \t_{h_\Lb^1(a)}^\Lb
\endequat
et
\equat\label{tauzg resrlg}
\t_{h_\Lb^1(a)}^\Lb \mathop{\circ} \resr_{\Lb \incl \Pb}^\Gb = 
\resr_{\Lb \incl \Pb}^\Gb \mathop{\circ} \t_a^\Gb.
\endequat

\bigskip

\exemple{exemple res}
Soient $\z_\Lb \in H^1(F,\ZC(\Gb))^\we$ et 
$\z \in H^1(F,\ZC(\Gb))^\we$ et soit $\l \in \Cent(\Lb^F)_{\z_\Lb}$ 
et $\g \in \Cent(\Gb^F)_\z$. Alors, d'apr\`es \ref{tauzg rlg}, on a 
$$R_{\Lb \incl \Pb}^\Gb \l \in \Cent(\Gb^F)_{\z_\Lb \circ h_\Lb^1}.$$
D'autre part, si $\z = \z_\Lb \circ h_\Lb^1$, alors, d'apr\`es \ref{tauzg resrlg}, 
on a 
$$\lexp{*}{R}_{\Lb \incl \Pb}^\Gb  \g \in \Cent(\Lb^F)_{\z_\Lb}.$$
De m\^eme, si $\Ker \z \notincl \Ker h_\Lb$, alors 
$$\lexp{*}{R}_{\Lb \incl \Pb}^\Gb \g = 0.~\SS{\square}$$

\bigskip

\soussection{Restriction de $\Gbt$ \`a $\Gb$}
Notons $\Pbt$ l'unique \para de $\Gbt$ tel que $\Pb=\Pbt \cap \Gb$ et soit 
$\Lbt$ l'unique \levic de $\Pbt$ tel que $\Lb=\Lbt \cap \Gb$. 
Alors $\Lbt$ est $F$-stable et $\Ub$ est le radical unipotent de $\Pbt$. 
De plus,
\equat\label{yug}
\Yb_\Ub^\Gbt = \Gbt^F \times_{\Gb^F} \Yb_\Ub^\Gb = 
\Yb_\Ub^\Gb \times_{\Lb^F} \Lbt^F.
\endequat
Par suite, on a, pour tout $k \in \NM$, un isomorphisme de 
$\Gbt^F$-modules-$\Lb^F$
\equat\label{premier}
H_c^k(\Yb_\Ub^\Gbt) \simeq \qlb \Gbt^F \otimes_{\qlb \Gb^F} 
H_c^k(\Yb_\Ub^\Gb)
\endequat
ainsi qu'un isomorphisme de $\Gb^F$-modules-$\Lbt^F$ 
\equat\label{second}
H_c^k(\Yb_\Ub^\Gbt) \simeq 
H_c^k(\Yb_\Ub^\Gb) \otimes_{\qlb \Lb^F} \qlb \Lbt^F.
\endequat
On en d\'eduit la proposition suivante~:

\bigskip

\proposition{rlg G Gtilde}
{\it On a~:
$$\Ind_{\Gb^F}^{\Gbt^F} \ci R_{\Lb \incl \Pb}^\Gb = 
R_{\Lbt \incl \Pbt}^\Gbt \ci 
\Ind_{\Lb^F}^{\Lbt^F},\leqno{\quad({\mathrm{a}})}$$
$$\Res_{\Lb^F}^{\Lbt^F} \ci \lexp{*}{R}_{\Lbt \incl \Pbt}^\Gbt 
= \lexp{*}{R}_{\Lb \incl \Pb}^\Gb \ci \Res_{\Gb^F}^{\Gbt^F},
\leqno{\quad({\mathrm{a}}^*)}$$
$$\Ind_{\Lb^F}^{\Lbt^F} \ci \lexp{*}{R}_{\Lb \incl \Pb}^\Gb = 
\lexp{*}{R}_{\Lbt \incl \Pbt}^\Gbt \ci \Ind_{\Gb^F}^{\Gbt^F},
\leqno{\quad({\mathrm{b}})}$$
$$\Res_{\Gb^F}^{\Gbt^F} \ci R_{\Lbt \incl \Pbt}^\Gbt = R_{\Lb \incl \Pb}^\Gb 
\ci \Res_{\Lb^F}^{\Lbt^F},\leqno{\quad({\mathrm{b}}^*)}$$
$$Q_{\Lbt \incl \Pbt}^\Gbt(u,v) = \sum_{g \in [\Gbt^F/\Gb^F]} 
Q_{\Lb \incl \Pb}^\Gb(\lexp{g}{u},v) = 
\sum_{l \in [\Lbt^F/\Lb^F]} Q_{\Lb \incl \Pb}^\Gb(u,\lexp{l}{v}).
\leqno{\quad({\mathrm{c}})}$$}

\bigskip

\proof Nous d\'emontrerons ici seulement la formule (a), 
les autres d\'ecoulant d'arguments similaires. 
Notons aussi que $({\mathrm{a}}^*)$ et $({\mathrm{b}}^*)$ 
sont des formules adjointes de (a) et (b).
Soit $\L$ un $\Lb^F$-module. Alors, d'apr\`es \ref{premier}, 
on a des isomorphismes de $\Gbt^F$-modules
\eqna
\Ind_{\Gb^F}^{\Gbt^F} \bigl(H^k_c(\Yb_\Ub^\Gb) \otimes_{\qlb\Lb^F} \L\bigr) 
&\simeq & \qlb\Gbt^F \otimes_{\qlb\Gb^F} \bigl(H^k_c(\Yb_\Ub^\Gb) 
\otimes_{\qlb\Lb^F} \L\bigr)\\
&\simeq & \bigl(\qlb\Gbt^F \otimes_{\qlb\Gb^F} H^k_c(\Yb_\Ub^\Gb)\bigr) 
\otimes_{\qlb\Lb^F} \L\\
&\simeq& H^k_c(\Yb_\Ub^\Gbt) \otimes_{\qlb\Lb^F} \L \\
&\simeq & \bigl(H^k_c(\Yb_\Ub^\Gbt) \otimes_{\qlb\Lbt^F} \qlb\Lbt^F\bigr)
\otimes_{\qlb\Lb^F} \L\\
&\simeq & H^k_c(\Yb_\Ub^\Gbt) \otimes_{\qlb\Lbt^F}
\bigl(\qlb\Lbt^F \otimes_{\qlb\Lb^F} \L\bigr)\\
&\simeq & H^k_c(\Yb_\Ub^\Gbt) \otimes_{\qlb\Lbt^F} 
\Ind_{\Lb^F}^{\Lbt^F} \L
\endeqna
et (a) en r\'esulte.\fin

\bigskip

Nous rappelons aussi la

\bigskip

\proposition{tenseur}
{\it Soit $\tau : \Lbt^F/\Lb^F \to \qlb^\times$ un caract\`ere lin\'eaire. 
Alors $\t$ peut \^etre vu comme un caract\`ere lin\'eaire de $\Gbt^F/\Gb^F
\simeq \Lbt^F/\Lb^F$ et, en utilisant cette identification, on a 
\begin{itemize}
\itemth{a} Si $\lamt \in \Class(\Lbt^F)$, alors  
$R_{\Lbt \incl \Pbt}^\Gbt(\lamt \otimes \t) = 
R_{\Lbt \incl \Pbt}^\Gbt(\lamt) \otimes \t$.

\itemth{b} Si $\gamt \in \Class(\Gbt^F)$, alors 
$\lexp{*}{R}_{\Lbt \incl \Pbt}^\Gbt(\gamt \otimes \t) = 
\lexp{*}{R}_{\Lbt \incl \Pbt}^\Gbt(\gamt) \otimes \t$.
\end{itemize}}

\bigskip

\soussection{Dualit\'e d'Alvis-Curtis} 
Notons $\PC(\D_0)$ l'ensemble des parties de $\D_0$. Alors 
$\phi_0$ agit sur $\PC(\D_0)$ et on note $\PC(\D_0)^{\phi_0}$ 
l'ensemble des parties $\phi_0$-stables de $\D_0$. 
Avec ces notations, on peut d\'efinir
$$D_\Gb = \sum_{I \in \PC(\D_0)^{\phi_0}} \eta_I ~R_{\Lb_I \incl \Pb_I}^\Gb 
~\circ~ \resr_{\Lb_I \incl \Pb_I}^\Gb.$$
Alors $D_\Gb : \ZM \Irr \Gb^F \to \ZM \Irr \Gb^F$ est une 
involution isom\'etrique \cite[proposition 8.10 et corollaire 8.14]{dmbook} 
appel\'ee {\it dualit\'e d'Alvis-Curtis}. 
Elle s'\'etend en une application lin\'eaire $\Cent(\Gb^F) \to \Cent(\Gb^F)$ 
toujours not\'e $D_\Gb$. 

\bigskip

\soussection{Formule de Mackey}
Soient $\Qb$ un \para de $\Gb$ et soit $\Mb$ un \levi $F$-stable de $\Qb$. 
Notons $\SC_\Gb(\Lb,\Mb)$ l'ensemble des $g \in \Gb$ tels que 
$\Lb \cap \lexp{g}{\Mb}$ contient un \tor de $\Gb$. Nous noterons 
$\D_{\Lb \incl \Pb, \Mb \incl \Qb}^\Gb$ l'application lin\'eaire 
$\Cent(\Mb^F) \to \Cent(\Lb^F)$ d\'efinie par 
$$\D_{\Lb \incl \Pb, \Mb \incl \Qb}^\Gb= 
\lexp{*}{R}_{\Lb \incl \Pb}^\Gb \circ R_{\Mb \incl \Qb}^\Gb - 
\sum_{g \in [\Lb^F\backslash \SC_\Gb(\Lb,\Mb)^F/\Mb^F]} 
R_{\Lb \cap \lexp{g}{\Mb} \incl \Lb \cap \lexp{g}{\Qb}}^\Lb 
\circ \lexp{*}{R}_{\Lb \cap \lexp{g}{\Mb} \incl \Pb \cap 
\lexp{g}{\Mb}}^{\lexp{g}{\Mb}} \circ (\ad g)_{\Mb^F}.$$
Ici, $(\ad g)_{\Mb^F} : \Cent(\Mb^F) \to Cent(\lexp{g}{\Mb}^F)$ est 
l'application induite par la conjugaison par $g$. Nous 
dirons que {\it ``la formule de Mackey a lieu dans $\Gb$''} si, 
pour tout sous-groupe r\'eductif connexe $\Gb'$ de $\Gb$ 
de m\^eme rang, pour tous \paras $\Pb'$ et $\Qb'$ de $\Gb'$ et 
pour tous compl\'ements de Levi $F$-stable $\Lb'$ et $\Mb'$ de $\Pb'$ et 
$\Qb'$ respectivement, on a $\D_{\Lb' \incl \Pb', \Mb' \incl \Qb'}^{\Gb'}=0$. 
Il est conjectur\'e que la formule de Mackey est valide sans 
hypoth\`ese. Pour l'instant, elle n'est connue que dans les cas suivants~:

\bigskip

\Theoreme{Formule de Mackey}{theo mackey}
{\it \begin{itemize}
\itemth{a} Supposons que l'une des conditions suivantes est satisfaite. 
\begin{itemize}
\itemth{a1} $\Pb$ et $\Qb$ sont $F$-stables. 

\itemth{a2} $\Lb$ ou $\Mb$ est un tore maximal de $\Gb$.
\end{itemize}
Alors $\D_{\Lb \incl \Pb, \Mb \incl \Qb}^\Gb = 0$. 

\itemth{b} Supposons que l'une des conditions suivantes est satisfaite. 
\begin{itemize}
\itemth{b1} $\d=1$ et $q \not= 2$.

\itemth{b2} $\d=1$ et $\Gb$ ne contient pas de composante 
quasi-simple de type $E_6$, $E_7$ ou $E_8$.

\itemth{b3} $\d=2$ et $\Gb$ est de type $B_2$, $G_2$ ou $F_4$.
\end{itemize}
Alors la formule de Mackey a lieu dans $\Gb$.
\end{itemize}}

\bigskip

\proof (a1) est d\^u \`a Deligne \cite[th\'eor\`eme 2.5]{luspa}, (a2) est d\^u \`a 
Deligne et Lusztig \cite[th\'eor\`eme 7]{delu2} et (b3) et (b4) et (b5) 
sont montr\'es dans \cite{ced jean}.\fin

\bigskip

De m\^eme que pour le th\'eor\`eme \ref{theo mackey}, il est conjectur\'e 
que le corollaire suivant reste vrai sans hypoth\`ese.

\bigskip

\corollaire{coro independance}
{\it Soit $\Pb_0$ un \para de $\Gb$ dont $\Lb$ est un compl\'ement de Levi. 
Supposons l'une des conditions suivantes v\'erifi\'ees~:
\begin{itemize}
\itemth{1} $\Lb$ est un tore maximal de $\Gb$.

\itemth{2} $\Pb$ et $\Pb_0$ sont $F$-stables.

\itemth{3} $\d=1$ et $q \not= 2$.

\itemth{4} $\d=1$ et $\Gb$ ne contient pas de composante 
quasi-simple de type $E_6$, $E_7$ ou $E_8$.

\itemth{5} $\d=2$ et $\Gb$ est de type $B_2$, $G_2$ ou $F_4$.
\end{itemize}
$$R_{\Lb \incl \Pb}^\Gb=R_{\Lb \incl \Pb_0}^\Gb,\leqno{\mathit{Alors}}$$
$$D_\Gb \circ R_{\Lb \incl \Pb}^\Gb = \e_\Gb\e_\Lb 
R_{\Lb \incl \Pb}^\Gb \circ D_\Lb$$
$$D_\Lb \circ \resr_{\Lb \incl \Pb}^\Gb = \e_\Gb\e_\Lb 
\resr_{\Lb \incl \Pb}^\Gb \circ D_\Gb.\leqno{\mathit{et}}$$}

\bigskip

\noindent{\sc Notation - } Si $\l$ est un \car virtuel de $\Lb^F$, nous noterons 
$W_{\Gb^F}(\Lb,\l)$ le groupe $N_{\Gb^F}(\Lb,\l)/\Lb^F$ et $\EC(\Gb^F,\Lb,\l)$ 
l'ensemble des \cars \irrs de $\Gb^F$ apparaissant dans le \car 
virtuel $R_{\Lb \incl \Pb}^\Gb \l$. A priori, cet ensemble 
pourrait d\'ependre du choix de $\Pb$. Bien s\^ur, il n'en d\'epend 
pas si l'une des hypoth\`eses du corollaire \ref{coro independance} 
est satisfaite. Dans la suite, nous n'emploierons cette notation que 
lorsque cet ensemble ne d\'epend pas du choix de $\Pb$ ou bien 
lorsque le choix du sous-groupe parabolique sera \'eclair\'e par 
le contexte.\finl

\bigskip

\noindent{\sc Remarque - } Si $\Tb$ est un \tor $F$-stable de $\Gb$ et 
si $\Bb$ est un sous-groupe de Borel contenant $\Tb$, nous noterons 
$R_\Tb^\Gb$ et $\lexp{*}{R}_\Tb^\Gb$ les applications 
$R_{\Tb \incl \Bb}^\Gb$ et $\lexp{*}{R}_{\Tb \incl \Bb}^\Gb$. De 
m\^eme, nous noterons $Q_\Tb^\Gb$ la fonction de Green 
$Q_{\Tb \incl \Bb}^\Gb$ respectivement.\finl

\bigskip

\soussection{Fonctions absolument cuspidales} 
Une fonction centrale $\g : \Gb^F \to \qlb$ est dite {\it absolument cuspidale} 
si, pour tout \levi $F$-stable propre de $\Gb$ et pour tout 
\para $\Pb$ de $\Gb$ dont $\Lb$ est un compl\'ement de Levi, 
on a $\resr_{\Lb \incl \Pb}^\Gb \g=0$ 
(voir \cite[d\'efinition du \SEC 3.1]{bonnafe mackey} ou 
\cite[\SEC 4.2]{bonnafe torsion}). Nous noterons $\Cus(\Gb^F)$ 
le $\qlb$-espace vectoriel des fonctions absolument cuspidales sur $\Gb^F$. 
D'apr\`es \ref{tzg rlg} et \ref{tauzg rlg}, $\Cus(\Gb^F)$ est stable 
sous les actions des groupes $\Zb(\Gb)^F$ et $H^1(F,\ZC(\Gb))$. 
Si la formule de Mackey a lieu dans $\Gb$, alors 
\equat
\Cent(\Gb^F)=\mathop{\oplus}_{\Lb \in [\LC(\ZC(\Gb))^F/\Gb^F]} 
R_\Lb^\Gb \Cus(\Lb^F).
\endequat
Rappelons que $\LC(\ZC(\Gb))$ est l'ensemble des sous-groupes de Levi 
de $\Gb$ (voir \SEC\ref{sous cuspidal levi})~; $\LC(\ZC(\Gb))^F$ est 
alors l'ensemble des \levi $F$-stables de $\Gb$, sur lequel 
le groupe $\Gb^F$ agit par conjugaison. Remarquons que, d'apr\`es 
\ref{tauzg rlg}, l'action de $H^1(F,\ZC(\Gb))$ sur $\Cent(\Gb^F)$ 
stabilise $\Cus(\Gb^F)$. 

\bigskip

\exemple{ZC absolument}
Soit $\z \in \ZC_\cus^\we(\Gb)$ et supposons que $\z$ est $F$-stale. 
Soit $\g \in \Cent(\Gb^F)_\z$. Alors, d'apr\`es l'exemple \ref{exemple res},
$\g$ est absolument cuspidale.\finl

\bigskip

\section{S\'eries de Lusztig g\'eom\'etriques et rationnelles\label{section series}}~

\medskip

\soussection{D\'efinitions} 
Si $(\Tb,\th) \doublefleche{\Gb} (\Tb^*,s)$, nous noterons 
$R_{\Tb^*}^\Gb(s)$ le caract\`ere (virtuel) de Deligne-Lusztig 
$R_\Tb^\Gb(\th)$. 
Fixons un \'el\'ement semi-simple $s \in \Gb^{*F*}$. Nous appellerons 
{\it s\'erie de Lusztig g\'eom\'etrique} (\resp {\it rationnelle}) 
associ\'ee \`a $s$ et nous noterons $\EC(\Gb^F,(s))$ 
(\resp $\EC(\Gb^F,[s])$) l'ensemble des caract\`eres 
irr\'eductibles de $\Gb^F$ apparaissant dans un $R_\Tb^\Gb(\th)$, 
o\`u $(\Tb,\th) \in \nabla(\Gb,F,(s))$ (\resp 
$(\Tb,\th) \in \nabla(\Gb,F,[s])$). En d'autres termes, 
$$\EC(\Gb^F,(s))=\bigcup_{(\Tb,\th) \in \nabla(\Gb,F,(s))} \EC(\Gb^F,\Tb,\th)$$
$$\EC(\Gb^F,[s])=\bigcup_{(\Tb,\th) \in \nabla(\Gb,F,[s])} \EC(\Gb^F,\Tb,\th).
\leqno{\text{et}}$$

\bigskip

\remarques{elementaire series} (a) Il est clair que 
$\nabla(\Gb,F,(s))=\DS{\bigcup_{[t] \incl (s)}} \nabla(\Gb,F,[t])$, donc 
$$\EC(\Gb^F,(s))=\DS{\bigcup_{[t] \incl (s)}} \EC(\Gb^F,[t]).$$

\medskip

(b) Si $\sti$ est un \'el\'ement semi-simple de $\Gbt^{*F^*}$, alors, 
d'apr\`es le corollaire \ref{geo et ratio}, on a 
$\EC(\Gbt^F,(\sti))=\EC(\Gbt^F,[\sti])$. 

\medskip

(c) Les s\'eries de Lusztig g\'eom\'etriques ou rationnelles sont stables 
sous l'action de $H^1(F,\ZC(\Gb))$. En effet, si $a \in H^1(F,\ZC(\Gb))$ 
et si $(\Tb,\th) \in \nabla(\Gb,F)$, alors $\t_a^\Gb(R_\Tb^\Gb(\th))=R_\Tb^\Gb(\th)$ 
d'apr\`es la formule \ref{tauzg rlg}. 

\medskip

(d) Si $\g \in \EC(\Gb^F,(s))$ (\resp $\g \in \EC(\Gb^F,[s])$) et si 
$z \in \Zb(\Gb)^{\circ F}$ (\resp $z \in \Zb(\Gb)^F$), alors 
$$t_z^\Gb \g = \sha^\circ(z) \g$$
$$t_z^\Gb \g = \sha(z) \g\quad).\leqno{\text{(respectivement}}$$
En effet, soit $\g \in \Irr \Gb^F$ et notons $\l$ le caract\`ere 
lin\'eaire de $\Zb(\Gb)^F$ tel que $t_z^\Gb \g=\l(z)\g$ pour tout 
$z \in \Zb(\Gb)^F$. D'apr\`es les formules \ref{tzg rlg} et \ref{def sha}, 
on a $t_z^\Gb R_\Tb^\Gb(\th) = \sha(z) R_\Tb^\Gb(\th)$ pour tout 
$(\Tb,\th) \in \nabla(\Gb,F,[s])$. Donc, si $\l \not= \sha$, 
$\g$ et $R_\Tb^\Gb(\th)$ sont orthogonaux.\finl

\bigskip

\soussection{Partition en s\'eries g\'eom\'etriques} 
La preuve du th\'eor\`eme suivant est d\^ue \`a Deligne et Lusztig 
\cite[th\'eor\`eme 6.2]{delu}.

\bigskip

\Theoreme{Deligne-Lusztig}{fondamental geo}
{\it \begin{itemize}
\itemth{a} Soit $s$ un \'el\'ement semi-simple de 
$\Gb^{*F^*}$, soit $(\Tb,\th) \in \nabla(\Gb,F,(s))$, soit $\Ub$ 
le radical unipotent d'un \borel de $\Gb$ contenant $\Tb$ et soit $k$ 
un entier naturel. Alors toute composante irr\'eductible de 
$H^k_c(\Yb_\Ub^\Gb) \otimes_{\qlb \Tb^F} \th$ appartient \`a 
la s\'erie de Lusztig g\'eom\'etrique 
$\EC(\Gb^F,(s))$. 

\itemth{b} On a 
$$\Irr \Gb^F=\coprod_{(s)} \EC(\Gb^F,(s)),$$
o\`u $(s)$ parcourt l'ensemble des classes de conjugaison 
g\'eom\'etriques d'\'el\'ements semi-simples de $\Gb^{*F^*}$.
\end{itemize}}

\bigskip

Le th\'eor\`eme \ref{fondamental geo} combin\'e au corollaire \ref{geo et ratio} 
fournit imm\'ediatement le corollaire suivant.

\bigskip

\corollaire{series gtilde}
{\it $\Irr \Gbt^F=\coprod_{[\sti]} \EC(\Gb^F,[\sti])$.}

\bigskip

Si $s$ est un \'el\'ement semi-simple de $\Gb^{*F^*}$, nous noterons 
$\Cent(\Gb^F,(s))$ le sous-$\qlb$-espace vectoriel de $\Cent(\Gb^F)$ engendr\'e 
par $\EC(\Gb^F,(s))$. Le th\'eor\`eme \ref{fondamental geo} (b) 
montre que 
\equat
\Cent(\Gb^F)= \mathop{\oplus}_{(s)}^\perp \Cent(\Gb^F,(s)).
\endequat

\bigskip

\soussection{Partition en s\'eries rationnelles} 
Dans cette sous-section, nous montrons qu'il est possible de remplacer 
``s\'erie g\'eom\'etrique'' par ``s\'erie rationnelle'' dans l'\'enonc\'e 
du th\'eor\`eme de Deligne-Lusztig pr\'ec\'edent. Ce r\'esultat 
est bien connu~: il \'etait annonc\'e dans \cite[7.3]{luirr} et Lusztig 
y donnait une indication pour la preuve. Une esquisse plus compl\`ete 
peut aussi \^etre trouv\'ee dans \cite[14.50]{dmbook}. 
La preuve (compl\`ete) que nous donnons ici ne pr\'etend \`a aucune 
originalit\'e~: il s'agit juste de suivre les indications des pr\'ec\'edents 
auteurs. Nous l'avons incluse car elle s'inscrit bien dans le cadre des 
m\'ethodes qui seront d\'evelopp\'ees tout au long de cet article.
Nous commen\c{c}ons par des rappels \'el\'ementaires.

\bigskip

\proposition{rtgs}
{\it Soient $(\Tbt^*,\sti) \in \nabla^*(\Gbt,F)$ et $z \in (\Ker i^*)^{F^*}$. 
Posons $(\Tb,s)=\lexp{*}{\RES}_\Gb^\Gbt (\Tbt^*,\sti)$. Alors
\begin{itemize}
\itemth{a} $\Res_{\Gb^F}^{\Gbt^F} 
R_{\Tbt^*}^\Gbt (\sti) = R_{\Tb^*}^\Gb(s)$.

\itemth{b} $R_{\Tbt^*}^\Gbt (\sti z)=
R_{\Tbt^*}^\Gbt (\sti) \zha^\Gbt$.
\end{itemize}}

\bigskip

\proof (a) r\'esulte du lemme \ref{mise au point} et de la proposition 
\ref{rlg G Gtilde} tandis que (b) d\'ecoule du lemme \ref{tenseur} 
et de la proposition \ref{tenseur T}.\fin

\bigskip

\corollaire{serie et lineaire}
{\it Soit $\sti$ un \'el\'ement semi-simple de $\Gbt^{*F^*}$ et soit 
$z \in (\Ker i^*)^{F^*}$. Alors l'application
$$\fonctio{\EC(\Gbt^F,[\sti])}{\EC(\Gbt^F,[\sti z])}{\g}{\g\zha}$$
est bijective.}

\bigskip

Le lien entre les s\'eries rationnelles de $\Gb^F$ 
et les s\'eries g\'eom\'etriques (ou rationnelles) de $\Gbt^F$ 
sont donn\'es par la proposition suivante.

\bigskip

\proposition{restriction series}
{\it Soit $\sti$ un \'el\'ement semi-simple de $\Gbt^{*F^*}$ et 
soit $s=i^*(\sti)$. Alors~:
\begin{itemize}
\itemth{a} Si $\gamt \in \EC(\Gbt^F,[\sti])$ et si $\g$ est une 
composante \irr de la restriction de $\gamt$ \`a $\Gb^F$, alors 
$\g \in \EC(\Gb^F,[s])$.

\itemth{b} Soit $\g \in \EC(\Gb^F,[s])$. Alors il existe $\gamt \in \EC(\Gbt^F,[\sti])$ 
tel que $\g$ est une composante \irr de la restriction de $\gamt$ \`a 
$\Gb^F$.
\end{itemize}}

\bigskip

\proof (a) Soit $W$ le groupe de Weyl de $\Gbt$ relatif \`a  
$\Tbt_0$. Pour tout $w \in W$, nous choisissons un \tor $F$-stable 
$\Tbt_w$ de $\Gbt$ de type $w$ par rapport \`a $\Tbt_0$. On pose aussi 
$$n_w=R_{\Tbt_w}^\Gbt(1_{\Tbt_w^F})(1)$$
o\`u $1_{\Tbt_w^F}$ est le \car trivial de $\Tbt_w^F$. Alors, d'apr\`es  
\cite[corollaire 2.11]{lucbms}, on a 
$$\sum_{\gamt \in \Irr \Gbt^F} \gamt(1)\gamt=\frac{1}{|W|} \sum_{w \in W} 
\Big(n_w \sum_{\thet \in (\Tbt_w^F)} R_{\Tbt_w}^\Gbt(\thet) \Big).$$
En projetant orthogonalement cette \'egalit\'e sur $\qlb\EC(\Gbt^F,[\sti])$ 
(voir corollaire \ref{series gtilde}), on obtient 
$$\sum_{\gamt \in \EC(\Gbt^F,[\sti])} \gamt(1)\gamt=\frac{1}{|W|} \sum_{w \in W} 
\Big(n_w \sum_{{\SS{\thet \in (\Tbt_w^F)}} \atop { 
\SS{(\Tbt_w,\thet) \in \nabla(\Gbt,F,[\sti])}}} R_{\Tbt_w}^\Gbt(\thet) \Big).$$
Donc, si $\gamt \in \EC(\Gbt^F,[\sti])$ et si $\g$ est une 
composante \irr de la restriction de $\gamt$ \`a $\Gb^F$, 
alors il existe $(\Tbt,\thet) \in \nabla(\Gbt,F,[\sti])$ tel que $\g$ 
soit une composante du \car virtuel $\Res_{\Gb^F}^{\Gbt^F} R_\Tbt^\Gbt(\thet)$. 
Il r\'esulte du corollaire \ref{mise au point coro} 
et de la proposition \ref{rtgs} que 
$\g \in \EC(\Gb^F,[s])$. Cela montre (a).

\medskip

(b) Soit $\g \in \EC(\Gb^F,[s])$. Alors il existe un \tor $F^*$-stable 
$\Tb^*$ de $\Gb^*$ contenant $s$ tel que $\g$ est une composante \irr de 
$R_{\Tb^*}^\Gb(s)$. Soit $\Tbt^*=i^{*-1}(\Tb^*)$. Alors, d'apr\`es la 
proposition \ref{rtgs}, 
$\g$ est une composante \irr de la restriction de $R_{\Tbt^*}^\Gbt(\sti)$ 
\`a $\Gb^F$. En particulier, il existe une composante 
\irr $\gamt$ de $R_{\Tbt^*}^\Gbt(\sti)$ telle que $\g$ 
soit une composante \irr de la restriction de $\gamt$ \`a $\Gb^F$. Mais 
$\gamt \in \EC(\Gbt^F,[\sti])$ par d\'efinition, donc (b) est d\'emontr\'e.\fin

\bigskip

\Theoreme{Lusztig}{fondamental ratio}
{\it \begin{itemize}
\itemth{a} Soit $s$ un \'el\'ement semi-simple de $\Gb^{*F^*}$, soit 
$(\Tb,\th) \in \nabla(\Gb,F,[s])$, soit $\Ub$ le radical unipotent 
d'un \borel de $\Gb$ contenant $\Tb$ et soit $k$ un entier naturel. 
Alors toute composante \irr du $\Gb^F$-module 
$H^k_c(\Yb_\Ub^\Gb)\otimes_{\qlb\Tb^F} \th$ appartient \`a $\EC(\Gb^F,[s])$.

\itemth{b} On a
$$\Irr \Gb^F = \coprod_{[s]} \EC(\Gb^F,[s]),$$
o\`u $[s]$ parcourt les classes de conjugaison rationnelles 
d'\'el\'ements semi-simples de $\Gb^{*F^*}$.
\end{itemize}}

\bigskip

\proof (a) Soit $\Tbt=\Tb.\Zbt$. Alors $\Ind_{\Gb^F}^{\Gbt^F} \g$ est un 
sous-$\Gbt^F$-module of
$$\Ind_{\Gb^F}^{\Gbt^F} H^k_c(\Yb_\Ub^\Gb)
\otimes_{\qlb\Tb^F} \th \simeq H^k_c(\Yb_\Ub^\Gbt) 
\otimes_{\qlb\Tbt^F} \Ind_{\Tb^F}^{\Tbt^F} \th$$
(voir l'isomorphisme \ref{premier}) donc il existe une extension 
$\thet$ de $\th$ \`a $\Tbt^F$et une composante \irr $\gamt$ du $\Gbt^F$-module 
$H^k_c(\Yb_\Ub^\Gbt) \otimes_{\qlb\Tbt^F} \thet$ tels que $\g$ 
soit une composante irr\'eductible de la restriction de $\gamt$ \`a $\Gb^F$. 
D'apr\`es le corollaire \ref{mise au point coro}, 
on a $(\Tbt,\thet) \in \nabla(\Gbt^F,[\tti])$ pour un \'el\'ement semi-simple 
$\tti$ de $\Gbt^{*F^*}$ tel que $i^*(\tti)=s$. Donc,  
par le th\'eor\`eme \ref{fondamental geo} (a), $\gamt \in \EC(\Gbt^F,[\tti])$. 
Il r\'esulte alors de la proposition \ref{restriction series} 
(a) que $\g \in \EC(\Gb^F,[s])$.

\medskip

(b) Soient $s_1$ et $s_2$ deux \'el\'ements semi-simples de $\Gb^{*F^*}$ tels que 
$\EC(\Gb^F,[s_1])$ et $\EC(\Gb^F,[s_2])$ ont un \'el\'ement en commun, disons $\g$. 
Pour $k=1$ ou $2$, soit $(\Tb_k,\th_k) \in \nabla(\Gb,F,[s_k])$ tel que $\g$ 
soit une composante \irr de $R_{\Tb_k}^\Gb(\th_k)$. Soit 
$\Ub_k$ le radical unipotent d'un \borel de $\Gb$ contenant 
$\Tb_k$ et soit $\Tbt_k$ le \tor $F$-stable de $\Gbt$ contenant $\Tb_k$. 
Alors il existe un entier naturel $n_k$ tel que  
$\g$ soit une composante \irr de 
$H^{n_k}_c(\Yb_{\Ub_k}^\Gb) \otimes_{\qlb \Tb_k^F} \th_k$. 
Si $\gamt$ est une composante \irr de $\Ind_{\Gb^F}^{\Gbt^F} \g$, 
alors, gr\^ace aux isomorphismes \ref{premier} et \ref{second}, 
il existe un \car lin\'eaire $\thet_k : \Tbt_k^F \to \qlb^\times$ 
tel que $\gamt$ est une composante \irr du $\Gbt^F$-module 
$H^{n_k}_c(\Yb_{\Ub_k}^\Gbt) \otimes_{\qlb \Tbt_k^F} \thet_k$. Mais 
$(\Tbt_k,\thet_k) \in \nabla(\Gb,F,[\sti_k])$ 
pour un \'el\'ement semi-simple $\sti_k \in \Gbt^{*F^*}$ tel que $i^*(\sti_k)=s_k$. 
Donc $\gamt \in \EC(\Gbt^F,[\sti_1]) \cap \EC(\Gbt^F,[\sti_2])$. 
Par cons\'equent, d'apr\`es le corollaire \ref{series gtilde} 
et le th\'eor\`eme \ref{fondamental geo} (a) appliqu\'e \`a $\Gbt^F$, 
on obtient que $\sti_1$ et $\sti_2$ sont conjugu\'es sous $\Gbt^{*F^*}$. 
Cela implique que $[s_1]=[s_2]$.\fin

\bigskip

Si $s$ est un \'el\'ement semi-simple de $\Gb^{*F^*}$, nous noterons 
$\Cent(\Gb^F,[s])$ le sous-$\qlb$-espace vectoriel de $\Cent(\Gb^F)$ engendr\'e 
par $\EC(\Gb^F,[s])$. Le th\'eor\`eme \ref{fondamental geo} (b) 
montre que 
\equat
\Cent(\Gb^F)= \mathop{\oplus}_{[s]}^\perp \Cent(\Gb^F,[s]).
\endequat

\bigskip

\soussection{Induction de Lusztig et s\'eries de Lusztig rationnelles} 
Soit $\Pb$ un \para de 
$\Gb$ et supposons que $\Pb$ admet un \levi $F$-stable $\Lb$. Soit 
$\Lb^*$ un \levi $F^*$-stable de $\Gb^*$ dual de $\Lb$. Alors  
le foncteur de Lusztig $R_{\Lb \incl \Pb}^\Gb$ pr\'eserve les s\'eries 
de Lusztig~:

\bigskip

\Theoreme{Lusztig}{rlg series}
{\it Soient $s$ un \'el\'ement semi-simple de $\Lb^{*F^*}$, $\l \in 
\EC(\Gb^F,[s]_{\Lb^{*F^*}})$ et $\g \in \EC(\Gb^F,\Lb,\l)$. 
Alors $\g \in \EC(\Gb^F,[s]_{\Gb^{*F^*}})$.}

\bigskip

\proof Soit $(\Tb,\th) \in \nabla(\Lb,F,[s]_{\Lb^{*F^*}})$ 
tel que $\g$ est une composante 
\irr de $R_\Tb^\Gb(\th)$. Notons qu'alors 
$(\Tb,\th) \in \nabla(\Gb,F,[s]_{\Gb^{*F^*}})$. 
Soit $\Bb$ un \borel de $\Lb$ contenant  
$\Tb$. Nous notons $\Vb$ et $\Ub$ les radicaux unipotents de $\Bb$ et $\Pb$
respectivement. Alors il existe un entier naturel $k$ tel que $\g$ 
soit une composante \irr du $\Gb^F$-module 
$H^k_c(\Yb_\Ub^\Gb) \otimes_{\qlb\Lb^F} \l$ 
et il existe un entier $k'$ such that 
$\l$ soit une composante \irr du $\Lb^F$-module 
$H^{k'}_c(\Yb_\Vb^\Lb) \otimes_{\qlb\Tb^F} \th$. Par suite, $\g$ 
est une composante \irr du $\Gb^F$-module
$$H^k_c(\Yb_\Ub^\Gb) \otimes_{\qlb\Lb^F} 
\Big(H^{k'}_c(\Yb_\Vb^\Lb) \otimes_{\qlb\Tb^F} \th \Big)
\simeq \Big(H^k_c(\Yb_\Ub^\Gb) \otimes_{\qlb\Lb^F} 
H^{k'}_c(\Yb_\Vb^\Lb) \Big) \otimes_{\qlb\Tb^F} \th.$$
Mais, par la formule de K\"unneth et d'apr\`es 
\cite[preuve de 11.5]{dmbook}, $H^k_c(\Yb_\Ub^\Gb) 
\otimes_{\qlb\Lb^F} H^{k'}_c(\Yb_\Vb^\Lb)$ est un sous-$\Gb^F$-module-$\Tb^F$ de  
$H^{k+k'}_c(\Yb_{\Vb\Ub}^\Gb)$, donc $\g$ est une composante 
\irr de $H^{k+k'}_c(\Yb_{\Vb\Ub}^\Gb)\otimes_{\qlb\Tb^F} \th$, 
ce qui montre que $\g \in \EC(\Gb^F,[s]_{\Gb^{*F^*}})$ 
d'apr\`es le th\'eor\`eme \ref{fondamental ratio} (a).\fin

\bigskip

\corollaire{rlg res series}
{\it Soit $s$ un \'el\'ement semi-simple de $\Gb^{*F^*}$, soit 
$\g \in \EC(\Gb^F,[s]_{\Gb^{*F^*}})$ 
et soit $\l$ une composante \irr de $\lexp{*}{R}_{\Lb \incl \Pb}^\Gb \g$. 
Alors $\l \in \EC(\Lb^F,[t]_{\Lb^{*F^*}})$ pour un \'el\'ement semi-simple 
$t \in \Lb^{*F^*}$ qui est $\Gb^{*F^*}$-conjugu\'e \`a $s$.}

\bigskip

\soussection{Stabilisateurs de \cars \irrs de ${\boldsymbol{\Gb^F}}$} 
Le r\'esultat suivant a \'et\'e montr\'e par Lusztig 
\cite[proposition 10]{luznc}. Pour cela, il a tout d'abord 
r\'eduit le probl\`eme au cas o\`u $\Gb$ est quasi-simple. 
Puisque ce r\'esultat est \'evident lorsque $\Gbt^F/\Gb^F.\Zb(\Gbt)^F$ 
est cyclique, il ne lui restait \`a traiter que le cas des groupes 
de type $D_{2n}$. Un d\'elicat argument de comptage lui 
a alors permis de conclure. Il serait plus satisfaisant d'avoir une 
preuve plus directe, mais nous en sommes incapables. Notons que cet 
argument de comptage est pr\'esent\'e en d\'etails dans 
\cite[chapitre 16]{cabanes}.

\bigskip

\Theoreme{Lusztig}{restriction}
{\it Soit $\gamt$ un \car \irr de $\Gbt^F$. Alors la 
restriction de $\gamt$ \`a $\Gb^F$ est sans multiplicit\'e.}

\bigskip

\noindent{\sc Question - } Le groupe $\Gbt^F/\Gb^F$ est un $p'$-groupe et, 
si $s$ est un $p'$-'\'el\'ement de $\Gbt^F$, alors $\Gbt^F=C_{\Gbt^F}(s).\Gb^F$ 
car, $s$ \'etant semi-simple, il est contenu dans un tore 
maximal. Il est alors naturel de se poser la question suivante~:

\smallskip

\begin{quotation}
\noindent{\it Soit $G$ un groupe fini et soit $N$ un sous-groupe distingu\'e 
de $G$. On suppose que $G/N$ est un $p'$-groupe et que $G=C_G(g) N$ pour 
tout $p'$-\'el\'ement $g \in G$. Est-ce que la restriction \`a $N$ 
d'un caract\`ere irr\'eductible de $G$ est toujours sans multiplicit\'e~?}
\end{quotation}

\smallskip

\noindent D'apr\`es ce qui pr\'ec\`ede, une r\'eponse positive \`a cette question 
fournirait une preuve du th\'eor\`eme \ref{restriction} qui n'utilise 
pas la classification des groupes r\'eductifs finis.\finl

\bigskip

Le morphisme 
$\omeh_s^0 \circ \s_\Gb^{-1} : \Gbt^F/\Gb^F.\Zb(\Gbt)^F \to (A_{\Gb^*}(s)^{F^*})^\we$ 
est surjectif. Nous noterons $\Gbt^F(s)$ le sous-groupe de $\Gbt^F$ tel que 
$\Gbt^F/\Gbt^F(s)$ soit isomorphe \`a $(A_{\Gb^*}(s)^{F^*})^\we$ via 
ce morphisme. Bien s\^ur, $\Gbt^F(s)$ contient $\Gb^F.\Zb(\Gbt)^F$. 

\bigskip

\corollaire{stabilisateur}
{\it Soient $s$ un \'el\'ement semi-simple de $\Gb^{*F^*}$ et 
$\g \in \EC(\Gb^F,[s])$. Soit $\Gbt^F(\g)$ le stabilisateur de $\g$ dans $\Gbt^F$. 
Alors $\Gbt^F(\g)$ contient $\Gbt^F(s)$. En d'autres termes, $\Gbt^F(s)$ 
(ou $\Ker \omeh_s^0$) agit trivialement sur $\EC(\Gb^F,[s])$.}

\bigskip

\proof Soit $\sti$ un \'el\'ement semi-simple de $\Gbt^{*F^*}$ tel que 
$s=i^*(\sti)$. D'apr\`es le th\'eor\`eme \ref{restriction} 
et la proposition \ref{restriction series} (b), 
il existe $\gamt \in \EC(\Gbt^F,[\sti])$ tel que 
$$\langle \Res_{\Gb^F}^{\Gbt^F} \gamt,\g \rangle_{\Gb^F}=1.$$
Donc, par la th\'eorie de Clifford et en utilisant l'isomorphisme 
$(\Ker i^*)^{F^*} \simeq (\Gbt^F/\Gb^F)^\we$, 
il suffit de montrer l'assertion suivante~:

\medskip

\begin{quotation}
{\it Si $z \in (\Ker i^*)^{F^*}$ v\'erifie $\gamt \otimes \zha=\gamt$,  
alors $z \in \Im \ph_s$.}
\end{quotation}

\medskip

Soit donc $z \in (\Ker i^*)^{F^*}$ tel que $\gamt \otimes \zha=\gamt$. Alors  
$\gamt \in \EC(\Gbt^F,[\sti]) \cap \EC(\Gbt^F,[\sti z])$ d'apr\`es 
le corollaire \ref{serie et lineaire}.
Donc $\sti$ et $\sti z$ sont $\Gbt^{*F^*}$-conjugu\'es 
ce qui montre que $z$ appartient \`a l'image de $\ph_s$ (voir lemme \ref{ags}).\fin

\bigskip

Si $a \in A_{\Gb^*}(s)^{F^*}$, rappelons que $\o_s^0(a)$ est un caract\`ere 
lin\'eaire de $H^1(F,\ZC(\Gb))$~; posons
$$\Cent(\Gb^F,(s),a) = \{\g \in \Cent(\Gb^F,(s))~|~\forall~z \in H^1(F,\ZC(\Gb)),~
\t_z^\Gb \g = \o_s^0(a)(z) \g\}$$
$$\Cent(\Gb^F,[s],a) = \{\g \in \Cent(\Gb^F,[s])~|~\forall~z \in H^1(F,\ZC(\Gb)),~
\t_z^\Gb \g = \o_s^0(a)(z) \g\}.\leqno{\text{et}}$$
Le corollaire \ref{stabilisateur} montre que 
\equat\label{decomposition action h1 geo}
\Cent(\Gb^F,(s))=\mathop{\oplus}^\perp_{a \in A_{\Gb^*}(s)^{F^*}} \Cent(\Gb^F,(s),a)
\endequat
et
\equat\label{decomposition action h1 ratio}
\Cent(\Gb^F,[s])=\mathop{\oplus}^\perp_{a \in A_{\Gb^*}(s)^{F^*}} \Cent(\Gb^F,[s],a).
\endequat
Nous verrons plus tard que, si $a \in A_{\Gb^*}(s)^{F^*}$, alors 
$\Cent(\Gb^F,[s],a) \not= \vide$ (voir \ref{action h1 rho}).

\bigskip

\soussection{Fonctions absolument cuspidales} 
Posons $\Cus(\Gb^F,(s))=\Cus(\Gb^F) \cap \Cent(\Gb^F,(s))$ et 
$\Cus(\Gb^F,[s])=\Cus(\Gb^F) \cap \Cent(\Gb^F,[s])$. De plus, 
si $a \in A_{\Gb^*}(s)^{F^*}$, nous poserons 
$\Cus(\Gb^F,(s),a)=\Cus(\Gb^F) \cap \Cent(\Gb^F,(s),a)$ et 
$\Cus(\Gb^F,[s],a)=\Cus(\Gb^F) \cap \Cent(\Gb^F,[s],a)$.
On a alors, d'apr\`es le th\'eor\`eme \ref{rlg series},
\equat
\Cus(\Gb^F)=\mathop{\oplus}_{(s)}^\perp \Cus(\Gb^F,(s))
=\mathop{\oplus}_{[s]}^\perp \Cus(\Gb^F,[s]).
\endequat
De plus, puisque $\Cus(\Gb^F)$ est stable par l'action de 
$H^1(F,\ZC(\Gb))$, il r\'esulte de \ref{decomposition action h1 geo} 
et \ref{decomposition action h1 ratio}
\equat
\Cus(\Gb^F,(s))=\mathop{\oplus}_{a \in A_{\Gb^*}(s)^{F^*}}^\perp 
\Cus(\Gb^F,(s),a)
\endequat
et
\equat
\Cus(\Gb^F,[s])=\mathop{\oplus}_{a \in A_{\Gb^*}(s)^{F^*}}^\perp 
\Cus(\Gb^F,[s],a).
\endequat

\newpage

\newpage

{\Large \part{Th\'eorie de Harish-Chandra\label{chapitre harish-chandra}}}

\bigskip

Soit $\Lbt$ un \levi $F$-stable d'un \para $F$-stable $\Pbt$ de $\Gbt$ et 
soit $\lamt$ un caract\`ere cuspidal de $\Lbt^F$. Notons $\Lb=\Lbt \cap \Gb$ 
et $\Pb=\Pbt \cap \Gb$ et posons $\l=\Res_{\Lb^F}^{\Lbt^F} \lamt$. 
Alors toutes les composantes irr\'eductibles de $\l$ sont cuspidales. 
Dans la section \ref{hc section}, nous montrons que ces composantes 
irr\'eductibles ne sont pas conjugu\'ees sous $N_{\Gb^F}(\Lb)$, \cad 
qu'elles d\'efinissent des s\'eries de Harish-Chandra diff\'erentes. 
L'ingr\'edient principal est un th\'eor\`eme de M. Geck \cite[page 400]{G}. 
Nous \'etudions dans la section \ref{section endo} 
les alg\`ebres d'endomorphismes des induits de ces caract\`eres 
cuspidaux pour obtenir un param\'etrage 
global des composantes irr\'eductibles de $R_\Lb^\Gb \l$.

\bigskip

\section{Autour d'un th\'eor\`eme de M. Geck\label{hc section}}~

\medskip

\soussection{Rappels} 
Un \car \irr $\g$ de $\Gb^F$ est dit {\it cuspidal} si, pour tout 
\para $F$-stable $\Pb$ de $\Gb$ et pour tout \levic $F$-stable de $\Pb$, 
on a $\resr_{\Lb \incl \Pb}^\Gb \g = 0$. Les faits suivants se d\'eduisent 
imm\'ediatement des propositions \ref{rlg G Gtilde} et \ref{tenseur}.

\bigskip

\proposition{restriction cuspidal}
{\it Soit $\gamt \in \Irr \Gbt^F$, soit $\g$ une composante 
\irr de la restriction de $\gamt$ \`a $\Gb^F$ et soit $\t \in (\Gbt^F/\Gb^F)^\we$. 
Les assertions suivantes sont \'equivalentes~:
\begin{itemize}
\itemth{1} $\gamt$ est cuspidal.

\itemth{2} $\gamt \otimes \t$ est cuspidal.

\itemth{3} $\g$ est cuspidal.
\end{itemize}}

\bigskip

Si $\Lb$ est un \levic $F$-stable d'un \para $F$-stable $\Gb$ de $\Pb$ et si $\l$ 
est un \car \irr cuspidal de $\Lb^F$, alors 
$\EC(\Gb^F,\Lb,\l)$ est appel\'e une {\it s\'erie de Harish-Chandra}. 
La proposition \ref{restriction cuspidal} montre qu'il existe un lien entre 
les s\'eries de Harish-Chandra de $\Gb^F$ et celles de $\Gbt^F$. 
Ce lien est pour l'essentiel l'objet de ce chapitre. 

\bigskip

\soussection{Notations} 
Dans ce chapitre, nous fixons un \para $F$-stable $\Pbt$ de $\Gbt$ et un \levic   
$F$-stable $\Lbt$ de $\Pbt$. Posons $\Pb=\Pbt \cap \Gb$ 
et $\Lb=\Lbt \cap \Gb$. Soit $\Ub$ le radical unipotent de $\Pbt$ ou de $\Pb$. 
Nous fixons un \car \irr {\it cuspidal} $\lamt$ de $\Lbt^F$ et nous 
posons 
$$\l=\Res_{\Lb^F}^{\Lbt^F} \lamt.$$
Nous fixons aussi une composante \irr $\l_1$ de $\l$. Pour tout 
$x \in \Lbt^F/\Lbt^F(\l_1)$, nous notons $\l_x=\lexp{x}{\l_1}$. Alors, d'apr\`es le 
th\'eor\`eme de Lusztig \ref{restriction}, on a 
$$\l=\sum_{x \in \Lbt^F/\Lbt^F(\l_1)} \lexp{x}{\l_1} = 
\sum_{x \in \Lbt^F/\Lbt^F(\l_1)} \l_x.$$
Puisque le groupe $\Lbt^F/\Lb^F$ est ab\'elien, on a  
$\Lbt^F(\l_x)=\Lbt^F(\l_1)$ pour tout $x \in \Lbt^F/\Lbt^F(\l_1)$.

Nous notons $\lamt^+$ (\resp $\l_x^+$, $x \in \Lbt^F/\Lbt^F(\l_1)$) 
le \car \irr de $\Pbt^F$ (\resp $\Pb^F$) obtenu en composant le \car $\lamt$ 
(\resp $\l_x$) avec le morphisme surjectif $\Pbt^F \to \Lbt^F$ 
(\resp $\Pb^F \to \Lb^F$). 
Nous fixons un $\Pbt^F$-module $\Mti$ ayant $\lamt^+$ comme caract\`ere. 
Nous notons $M$ la restriction de $\Mti$ \`a $\Pb^F$ et $M_x$ le sous-$\Pb^F$-module 
\irr de $M$ ayant $\l_x^+$ comme caract\`ere. Alors 
$$\Ind_{\Pbt^F}^{\Gbt^F}\Mti=\{\fti : \Gbt^F \to \Mti~|~\forall y \in \Pbt^F, 
~\forall g \in \Gbt^F, ~\fti(yg)=y.\fti(g)\},$$
$$\Ind_{\Pb^F}^{\Gb^F}M=\{f : \Gbt^F \to M~|~\forall y \in \Pb^F, 
~\forall g \in \Gb^F, ~f(yg)=y.f(g)\}$$
et des descriptions similaires sont valides pour $\Ind_{\Pb^F}^{\Gb^F} M_x$. Alors 
\equat\label{1 decompo}
\Ind_{\Pb^F}^{\Gb^F} M = \mathop{\oplus}_{x \in \Lbt^F/\Lbt^F(\l_1)} 
\Ind_{\Pb^F}^{\Gb^F} M_x. 
\endequat
Le $\Gbt^F$-module $\Ind_{\Pbt^F}^{\Gbt^F} \Mti$ a pour \car 
$R_{\Lbt \incl \Pbt}^\Gbt \lamt$. De m\^eme, les $\Gb^F$-modules 
$\Ind_{\Pb^F}^{\Gb^F} M$ et $\Ind_{\Pb^F}^{\Gb^F} M_x$ 
ont pour \cars $R_{\Lb \incl \Pb}^\Gb \l$ et $R_{\Lb \incl \Pb}^\Gb \l_x$ 
respectivement.

Soit $f \in \Ind_{\Pb^F}^{\Gb^F} M$. Notons $\fti$ la fonction $\Gbt^F \to \Mti=M$ 
definie comme suit. Soit $\gti \in \Gbt^F$. Alors il existe $\xti \in \Pbt^F$ et  
$g \in \Gb^F$ tels que $\gti=\xti g$. On pose $\fti(\gti)=\xti.f(g)$~: 
remarquons que 
$\fti(\gti)$ ne d\'epend pas du choix de $\xti \in \Pbt^F$ et $g \in \Gb^F$ tels que 
$\gti=\xti g$. Alors la restriction de $\fti$ \`a $\Gb^F$ est \'egale \`a $f$ et 
il est facile de v\'erifier que l'application 
\equat\label{1 iso}
\fonctio{\Ind_{\Pb^F}^{\Gb^F} M}{\Res_{\Gb^F}^{\Gbt^F} 
\Ind_{\Pbt^F}^{\Gbt^F} \Mti}{f}{\fti}
\endequat
est un \iso de $\Gb^F$-modules.

Le groupe $W_{\Gbt^F}(\Lbt,\lamt)$ 
est naturellement isomorphe \`a $W_{\Gb^F}(\Lbt,\lamt)=N_{\Gb^F}(\Lbt,\lamt)/\Lb^F$. 
En g\'en\'eral, nous nous r\'ef\`ererons \`a ce dernier car il 
est plus adapt\'e \`a notre situation. Par exemple, il est 
naturellement un sous-groupe de $W_{\Gb^F}(\Lb,\l)$.

\bigskip

\soussection{Un th\'eor\`eme de M. Geck} Dans \cite[corollaire 2]{G}, Geck \'enonce 
le r\'esultat suivant~:
{\it le \car $R_{\Lb \incl \Pb}^\Gb(\l_1)$ a une composante \irr 
de multiplicit\'e $1$}. Cependant, pour l'obtenir, 
il a en fait d\'emontr\'e le r\'esultat plus fort suivant \cite[page 400]{G}~:

\bigskip

\Theoreme{Geck}{geck}
{\it Le \car $R_{\Lb \incl \Pb}^\Gb \l$ a une composante \irr de multiplicit\'e $1$.}

\bigskip

\noindent{\sc Esquisse de la preuve de Geck - } Soient 
$\gamt$ et $\gamt'$ deux composantes \irrs de $R_{\Lbt \incl \Pbt}^\Gbt \lamt$ 
dont la restriction \`a $\Gb^F$ ont une composante \irr commune. 
Alors, par la th\'eorie de Clifford et par le th\'eor\`eme \ref{restriction}, 
ces deux restrictions co\"{\i}ncident. En particulier, $\gamt(1)=\gamt'(1)$. 
Pour prouver le th\'eor\`eme \ref{geck}, il est donc suffisant de montrer 
qu'il existe une composante \irr $\gamt$ de $R_{\Lbt \incl \Pbt}^\Gbt \lamt$, 
apparaissant avec la multiplicit\'e $1$ et telle que 
$\gamt(1)\not=\gamt'(1)$ pour toute autre composante \irr $\gamt'$ de 
$R_{\Lbt \incl \Pbt}^\Gbt \lamt$. En fait, Geck montre que le degr\'e 
g\'en\'erique de la signature de l'alg\`ebre de Hecke associ\'ee 
\`a la donn\'ee $(\Gbt,\Lbt,\lamt)$ a une $p$-partie sup\'erieure 
\`a celle des degr\'es g\'en\'eriques des autres \cars \irrs de 
l'alg\`ebre de Hecke \cite[th\'eor\`eme 1]{G}, ce qui permet de conclure.\fin

\bigskip 

\noindent{\sc Remarque - } Lusztig \cite{lubook} 
a montr\'e que $R_{\Lbt \incl \Pbt}^\Gbt\lamt$ 
poss\`ede une composante \irr de multiplicit\'e $1$. Donc le th\'eor\`eme de 
Geck g\'en\'eralise celui de Lusztig. Mais il faut bien noter que 
la preuve de Geck utilise le r\'esultat de Lusztig.\finl

\medskip

Cette version plus forte \ref{geck} du th\'eor\`eme de Geck 
est n\'ecessaire pour montrer le r\'esultat (tr\`es utile) suivant~:

\bigskip

\corollaire{stabilo}
{\it \begin{itemize}
\itemth{a} Pour tout $x \in \Lbt^F/\Lbt^F(\l_1)$, le \car 
$R_{\Lb \incl \Pb}^\Gb \l_x$ 
poss\`ede une composante \irr de multiplicit\'e $1$.

\itemth{b} Si $x$ et $y$ sont deux \'el\'ements distincts de $\Lbt^F/\Lbt^F(\l_1)$, 
alors les \cars cuspidaux $\l_x$ et  
$\l_y$ ne sont pas conjugu\'es sous $N_{\Gb^F}(\Lb)$. De fa\c{c}on \'equivalente, 
$$\langle R_{\Lb \incl \Pb}^\Gb \l_x, R_{\Lb \incl \Pb}^\Gb \l_y \rangle_{\Gb^F}=0.$$

\itemth{c} On a $W_{\Gb^F}(\Lb,\l)=W_{\Gb^F}(\Lb,\l_x)$ pour tout $x \in 
\Lbt^F/\Lbt^F(\l_1)$.

\itemth{d} Si $\g$ est une composante \irr de $R_{\Lb \incl \Pb}^\Gb \l$, 
alors $\Gbt^F(\g)$ est contenu dans $\Gb^F.\Lbt^F(\l_1)$.
\end{itemize}}

\bigskip

\proof (a) r\'esulte du th\'eor\`eme \ref{geck} 
et de l'\'egalit\'e~:
$$R_{\Lb \incl \Pb}^\Gb \l = \sum_{x \in \Lbt^F/\Lbt^F(\l_1)} 
\lexp{x}{R}_{\Lb \incl \Pb}^\Gb \l_1.$$
En fait, (a) est l'\'enonc\'e original de Geck \cite[corollaire 2]{G}.

\medskip

Prouvons maintenant (b) et (c). Les groupes
$N_{\Gb^F}(\Lb)$ et $\Lbt^F$ agissent sur $\Irr \Lb^F$ et ces deux actions 
commutent. Donc, si $x$ et $y$ appartiennent \`a $\Lbt^F/\Lbt^F(\l_1)$ et si 
$n \in N_{\Gb^F}(\Lb)$ est tel que $\lexp{n}{\l_x}=\l_y$, alors 
$n \in N_{\Gb^F}(\Lb,\l)$. Par cons\'equent, prouver (b) et (c) est \'equivalent 
\`a prouver que $N_{\Gb^F}(\Lb,\l)$ agit trivialement sur 
$E=\{\l_x~|~x \in \Lbt^F/\Lbt^F(\l_1)\}$. Puisque les actions de 
$N_{\Gb^F}(\Lb,\l)$ et
$\Lbt^F$ sur $E$ commutent, toutes les orbites de $N_{\Gb^F}(\Lb,\l)$ sur $E$ 
ont le m\^eme cardinal, disons $n$. Mais alors $n$ divise la multiplicit\'e 
de toute composante \irr de $R_{\Lb \incl \Pb}^\Gb \l$. 
Donc, d'apr\`es le th\'eor\`eme \ref{geck}, $n=1$ et (b) et (c) en r\'esultent.

\medskip

(d) Soit $g \in \Gbt^F$ tel que $\lexp{g}{\g}=\g$. On peut supposer que 
$g \in \Lbt^F$. Par hypoth\`ese, il existe $x \in \Lbt^F/\Lbt^F(\l_1)$ tel que 
$\g$ est une composante \irr de $R_{\Lb \incl \Pb}^\Gb \l_x$. Mais $\lexp{g}{\g}=\g$ 
donc $\g$ est aussi une composante \irr de $R_{\Lb \incl \Pb}^\Gb \lexp{g}{\l_x}$. 
Donc, d'apr\`es (b), $g \in \Lbt^F(\l_x)=\Lbt^F(\l_1)$.\fin

\bigskip

\soussection{Une extension centrale de 
${\boldsymbol{W_{\Gb^F}(\Lbt,\lamt)}}$\label{sub W'}} 
Posons 
$$W_{\Gb^F}^\pr(\Lbt,\lamt)=\{(w,\th) \in W_{\Gb^F}(\Lbt)\times (\Lbt^F/\Lb^F)^\we~|~
\lexp{w}{\lamt}=\lamt \otimes \th\}.$$
Notons qu'en fait 
$W_{\Gb^F}^\pr(\Lbt,\lamt) \incl W_{\Gb^F}(\Lb,\l) \times (\Lbt^F/\Lb^F)^\we$.
Nous identifierons $W_{\Gb^F}(\Lbt,\lamt)$ avec le sous-groupe 
$W_{\Gb^F}(\Lbt,\lamt) \times \{1\}$ de $W_{\Gb^F}^\pr(\Lbt,\lamt)$. 
De m\^eme, nous identifierons $(\Lbt^F/\Lbt^F(\l_1))^\we$ avec 
le sous-groupe $\{1\} \times (\Lbt^F/\Lbt^F(\l_1))^\we$ de 
$W_{\Gb^F}^\pr(\Lbt,\lamt)$. La premi\`ere projection
$$\fonctio{W_{\Gb^F}^\pr(\Lbt,\lamt)}{W_{\Gb^F}(\Lb,\l)}{(w,\th)}{w}$$
est surjective (d'apr\`es le th\'eor\`eme \ref{restriction} et la th\'eorie de 
Clifford) et son noyau est \'egal \`a $(\Lbt^F/\Lbt^F(\l_1))^\we$. 
De plus, $(\Lbt^F/\Lbt^F(\l_1))^\we$ est central dans 
$W_{\Gb^F}^\pr(\Lbt,\lamt)$.
La deuxi\`eme projection 
$$\fonctio{W_{\Gb^F}^\pr(\Lbt,\lamt)}{(\Lbt^F/\Lb^F)^\we}{(w,\th)}{\th}$$
n'est pas surjective en g\'en\'eral et son noyau 
est $W_{\Gb^F}(\Lbt,\lamt)$.
Nous notons $\Lbt^F(\Gb,\l_1)$ le sous-groupe de $\Lbt^F$ contenant $\Lb^F$ 
tel que l'image de cette deuxi\`eme projection soit $(\Lbt^F/\Lbt^F(\Gb,\l_1))^\we$. 
En fait, le groupe $\Lbt^F(\Gb,\l_1)$ est contenu dans $\Lbt^F(\l_1)$. 
On a des isomorphismes canoniques
\equat\label{W'/W}
W_{\Gb^F}^\pr(\Lbt,\lamt)/W_{\Gb^F}(\Lbt,\lamt) \simeq (\Lbt^F/\Lbt^F(\Gb,\l_1))^\we
\endequat
et
\equat\label{W/W}
W_{\Gb^F}(\Lb,\l)/W_{\Gb^F}(\Lbt,\lamt) \simeq (\Lbt^F(\l_1)/\Lbt^F(\Gb,\l_1))^\we.
\endequat
R\'esumons tout ceci dans le diagramme suivant, dans lequel toutes les suites 
horizontales ou verticales sont exactes et tous les carr\'es sont 
commutatifs~:
$$\diagram
&&1\dto & 1 \dto & \\
&& W_\Gb(\Lbt,\lamt) \rdouble \ddto& W_\Gb(\Lbt,\lamt) \ddto & \\
&&&&\\
1 \rto & (\Lbt^F/\Lbt^F(\l_1))^\we \rto \dddouble & 
W_\Gb'(\Lbt,\lamt) \rto \ddto
& W_\Gb(\Lb,\l) \ddto \rto & 1 \\
&&&&\\
1 \rto & (\Lbt^F/\Lbt^F(\l_1))^\we \rto & 
(\Lbt^F/\Lbt^F(\Gb,\l_1))^\we \rto \dto
& (\Lbt^F(\l_1)/\Lbt^F(\Gb,\l_1))^\we \rto \dto& 1 \\
&& 1 & 1 & \\
\enddiagram$$

\bigskip

\remarque{abelianite} Les \isos \ref{W'/W} et \ref{W/W} entra\^{\i}nent 
que les groupes $W_{\Gb^F}^\pr(\Lbt,\lamt)/W_{\Gb^F}(\Lbt,\lamt)$ 
et $W_{\Gb^F}(\Lb,\l)/W_{\Gb^F}(\Lbt,\lamt)$ sont ab\'eliens.\finl 

\bigskip

\proposition{aaa}
{\it On a~:
\begin{itemize}
\itemth{a} $(\Lbt^F/\Lbt^F(\Gb,\l_1))^\we = \{ \th \in (\Lbt^F/\Lb^F)^\we~|~
(R_{\Lbt \incl \Pbt}^\Gbt \lamt) \otimes \th = 
R_{\Lbt \incl \Pbt}^\Gbt \lamt\}$. 

\itemth{b} Si $\g \in \EC(\Gb^F,\Lb,\l)$, alors 
$\Gbt^F(\g)$ contient $\Gb^F.\Lbt^F(\Gb,\l_1)$.
\end{itemize}}

\bigskip

\proof (a) Soit $\th \in (\Lbt^F/\Lbt^F(\Gb,\l_1))^\we$. Alors il existe 
$w \in W_{\Gb^F}(\Lbt)$ tel que $(w,\th) \in W_{\Gb^F}^\pr(\Lbt,\lamt)$. 
Par cons\'equent, 
$\lexp{w}{\lamt}=\lamt \otimes \th$. Donc, d'apr\`es la proposition \ref{tenseur}, 
on a 
$(R_{\Lbt \incl \Pbt}^\Gbt \lamt)\otimes \th = R_{\Lbt \incl \Pbt}^\Gbt \lamt$. 
R\'eciproquement, soit $\th \in (\Lbt^F/\Lb^F)^\we$ tel que 
$(R_{\Lbt \incl \Pbt}^\Gbt \lamt)\otimes \th = R_{\Lbt \incl \Pbt}^\Gbt \lamt$. 
Alors, d'apr\`es la proposition \ref{tenseur}, on a $R_{\Lbt \incl \Pbt}^\Gbt 
(\lamt \otimes \th) = R_{\Lbt \incl \Pbt}^\Gbt \lamt$. 
De plus, $\lamt$ et $\lamt \otimes \th$ 
sont cuspidaux d'apr\`es la proposition \ref{restriction cuspidal}.  
Donc il existe $w \in W_{\Gb^F}(\Lb)$ tel que 
$\lexp{w}{\lamt}=\lamt \otimes \th$. En d'autres termes, 
$(w,\th) \in W_{\Gb^F}^\pr(\Lbt,\lamt)$ donc $\th \in (\Lbt^F/\Lbt^F(\Gb,\l_1))^\we$.

\medskip

(b) Soit $\gamt$ une composante \irr de $R_{\Lbt \incl \Pbt}^\Gbt \lamt$ 
telle que $\g$ soit une composante \irr de la restriction de $\gamt$ \`a $\Gb^F$. 
D'apr\`es le th\'eor\`eme \ref{restriction} et la th\'eorie de Clifford, 
(b) est \'equivalent \`a l'\'enonc\'e suivant~: si $\th \in 
(\Lbt^F/\Lb^F)^\we$ est tel que $\gamt \otimes \th=\gamt$, alors 
$\th \in (\Lbt^F/\Lbt^F(\Gb,\l_1))^\we$. Mais, si 
$\gamt \otimes \th=\gamt$, alors $\gamt$ est une composante \irr de 
$R_{\Lbt \incl \Pbt}^\Gbt \lamt$ ainsi que de 
$R_{\Lbt \incl \Pbt}^\Gbt (\lamt \otimes \th)$. Par suite, 
$R_{\Lbt \incl \Pbt}^\Gbt \lamt = R_{\Lbt \incl \Pbt}^\Gbt (\lamt \otimes \th)$
et donc $\th \in (\Lbt^F/\Lbt^F(\Gb,\l_1))^\we$ d'apr\`es (a).\fin

\bigskip

\corollaire{s}
{\it Si $\g \in \EC(\Gb^F,\Lb,\l)$, alors 
$\Gb^F.\Lbt^F(\Gb,\l_1) \incl \Gbt^F(\g) \incl \Gb^F.\Lbt^F(\l_1)$.}

\bigskip

\section{Alg\`ebres d'endomorphismes\label{section endo}}~

\medskip 

\soussection{Description} 
Nous noterons dans ce chapitre $\HCt$ l'alg\`ebre d'endomorphismes 
du $\Gbt^F$-module $\Ind_{\Pbt^F}^{\Gbt^F} \Mti$ et par $\HC$ l'alg\`ebre 
d'endomorphismes du $\Gb^F$-module $\Ind_{\Pb^F}^{\Gb^F} M$. Alors $\HCt$ est, 
via l'\iso \ref{1 iso}, une sous-alg\`ebre de $\HC$. Pour tout 
$x \in \Lbt^F/\Lbt^F(\l_1)$, 
nous notons $\HC_x$ l'alg\`ebre d'endomorphismes du $\Gb^F$-module 
$\Ind_{\Pb^F}^{\Gb^F} M_x$. Alors, d'apr\`es le corollaire \ref{stabilo} (b), on a 
$$\HC=\prod_{x \in \Lbt^F/\Lbt^F(\l_1)} \HC_x$$
$$\Irr \HC = \coprod_{x \in \Lbt^F/\Lbt^F(\l_1)} \Irr \HC_x.\leqno{\mathrm{et~donc}}$$
Cela nous donne un morphisme d'alg\`ebres 
$\HCt \to \HC_x$ pour tout $x \in \Lbt^F/\Lbt^F(\l_1)$. 

\medskip

D'apr\`es \cite[lemme 6.5]{HL} et d'apr\`es le corollaire \ref{stabilo}, 
le \car \irr $\l_1$ de $\Lb^F$ s'\'etend en en un \car \irr $\n_1$ de 
$N_{\Gb^F}(\Lb,\l)$.

\bigskip

\lemme{stabi}
{\it \begin{itemize}
\itemth{a} Le stabilisateur de $\n_1$ dans $\Lbt^F$ est $\Lbt^F(\Gb,\l_1)$.

\itemth{b} Le stabilisateur, dans $\Lbt^F$, 
de la restriction de $\n_1$ \`a $N_{\Gb^F}(\Lbt,\lamt)$ 
est $\Lbt^F(\l_1)$.
\end{itemize}}

\bigskip

\proof Soit $\Nti=N_{\Gb^F}(\Lb,\l).\Lbt^F$. 
Alors, d'apr\`es la formule de Mackey, on a 
$$\Res_{\Lbt^F}^{\Nti} 
\Ind_{N_{\Gb^F}(\Lbt,\lamt)}^{\Nti} \n_1 = \Ind_{\Lb^F}^{\Lbt^F} \l_1.$$
Le \car $\Ind_{\Lb^F}^{\Lbt^F} \l_1$ est sans multiplicit\'e, donc le \car 
$\Ind_{N_{\Gb^F}(\Lbt,\lamt)}^{\Nti} \n_1$ est sans multiplicit\'e. 
Soit $\nut$ une de ses composantes \irrs telle que $\Res_{\Lbt^F}^\Nti \nut$ 
a $\lamt$ comme composante irr\'eductible. Alors
$$\Res_{\Lbt^F}^\Nti \nut = \sum_{\t \in (\Lbt^F(\l_1)/\Lbt^F(\Gb,\l_1))^\we} 
\lamt \otimes \taut$$
o\`u, pour chaque $\t \in (\Lbt^F(\l_1)/\Lbt^F(\Gb,\l_1))^\we$, $\taut$ 
d\'esigne une extension de $\t$ \`a $\Lbt^F$~: alors 
le \car $\lamt \otimes \taut$ ne d\'epend que de $\t$ et non du choix 
de $\taut$. Cette d\'ecomposition a lieu car l'orbite de $\lamt$ sous cette action  
de $\Nti$ est $\{\lamt \otimes \taut~|~\t \in (\Lbt^F(\l_1)/\Lbt^F(\Gb,\l_1))^\we\}$.

Soit $\th$ un \car lin\'eaire de 
$\Lbt^F/\Lb^F \simeq \Nti/N_{\Gb^F}(\Lbt,\lamt)$. Alors, par la th\'eorie de Clifford, 
l'\'enonc\'e (a) du lemme \ref{stabi} est \'equivalent \`a l'assertion suivante~:

\medskip

\begin{quotation}{\it $\nut \otimes \th=\nut$ \ssi $\Lbt^F(\Gb,\l_1) \incl 
\Ker \th$.}\end{quotation}

\medskip

Supposons tout d'abord que $\nut \otimes \th=\nut$. Alors  
\eqna
(\Res_{\Lbt^F}^\Nti \nut) \otimes \th &=& \Res_{\Lbt^F}^\Nti (\nut \otimes \th) \\
&=& \Res_{\Lbt^F}^\Nti \nut
\endeqna
ce qui implique que $\Lbt^F(\Gb,\l_1) \incl \Ker \th$.

R\'eciproquement, supposons que $\Lbt^F(\Gb,\l_1) \incl \Ker \th$. Alors 
$\Res_{\Lbt^F}^\Nti (\nut \otimes \th) =\Res_{\Lbt^F}^\Nti \nut$ donc 
$\nut \otimes \th$ est une composante \irr de $\Ind_{N_{\Gb^F}(\Lb,\l)}^\Nti \n_1$ 
et $\Res_{\Lbt^F}^\Nti (\nut \otimes \th)$ a $\lamt$ composante irr\'eductible.
Ceci implique que $\nut \otimes \th=\nut$. D'o\`u (a).

\medskip

(b) d\'ecoule d'un argument similaire.\fin

\bigskip

Pour tout $x \in \Lbt^F/\Lb^F(\l_1)$, nous choisissons 
un repr\'esentant $\xti$ de $x$ dans $\Lbt^F$. 
Alors $\lexp{\xti}{\n_1}$ est une extension de $\l_x$ \`a $N_{\Gb^F}(\Lb,\l)$. 
Nous posons 
$$\n=\sum_{x \in \Lbt^F/\Lbt^F(\l_1)} \lexp{\xti}{\n_1}.$$
Alors $\n$ est une extension de $\l$ \`a $N_{\Gb^F}(\Lb,\l)$. Pour chaque  
$x \in \Lbt^F/\Lbt^F(\l_1)$, nous notons $\m_x$ la restriction de 
$\lexp{\xti}{\n_1}$ \`a $N_{\Gb^F}(\Lbt,\lamt)$.

\bigskip

\remarque{dependre de x} 
Soit $x \in \Lbt^F/\Lbt^F(\l_1)$. Alors, si $\Lbt^F(\Gb,\l_1) \not= \Lbt^F(\l_1)$, 
le caract\`ere $\lexp{\xti}{\n_1}$ depend du choix de $\xti$ d'apr\`es le 
lemme \ref{stabi} (a), tandis que le \car $\m_x$ n'en 
d\'epend pas d'apr\`es le lemme \ref{stabi} (b). Ceci justifie la notation.\finl

\bigskip

Posons 
$$\m=\sum_{x \in \Lbt^F/\Lbt^F(\l_1)} \m_x$$
de sorte que $\m$ est la restriction de $\n$ \`a $N_{\Gb^F}(\Lbt,\lamt)$. 
Par la formule de Mackey, on a 
$$\Res_{\Lbt^F}^{N_{\Gbt^F}(\Lbt,\lamt)} 
\Ind_{N_{\Gb^F}(\Lbt,\lamt)}^{N_{\Gbt^F}(\Lbt,\lamt)} \m_1 
= \Ind_{\Lb^F}^{\Lbt^F} \l_1$$
donc il existe une unique composante \irr $\mut$ de 
$\Ind_{N_{\Gb^F}(\Lbt,\lamt)}^{N_{\Gbt^F}(\Lbt,\lamt)} \m_1$ telle que 
$$\Res_{\Lbt^F}^{N_{\Gbt^F}(\Lbt,\lamt)}\mut=\lamt.$$ 
D'apr\`es le lemme \ref{stabi} (b), on a 
\equat\label{resres}
\m=\Res_{N_{\Gb^F}(\Lbt,\lamt)}^{N_{\Gbt^F}(\Lbt,\lamt)} \mut.
\endequat

\medskip

Fixons une repr\'esentation $\sigt : N_{\Gbt^F}(\Lbt,\lamt) \to \Gb\Lb(\Mti)$ 
ayant pour \car $\mut$ et une repr\'esentation $\s : N_{\Gb^F}(\Lb,\l) \to 
\Gb\Lb(M)$ ayant pour \car $\n$ \'etendant toutes deux les representations de 
$\Lbt^F$ et $\Lb^F$ sur $\Mti$ et $M$ respectivement. Alors 
\equat\label{resresres}
\Res_{N_{\Gb^F}(\Lbt,\lamt)}^{N_{\Gb^F}(\Lb,\l)} \s = 
\Res_{N_{\Gb^F}(\Lbt,\lamt)}^{N_{\Gbt^F}(\Lbt,\lamt)} \sigt
\endequat
d'apr\`es \ref{resres}. Pour tout $x \in \Lbt^F/\Lbt^F(\l_1)$, notons 
$\s_x$ la restriction de $\s$ aux sous-espace $M_x$ de $M$.

Pour tout $w \in W_{\Gb^F}(\Lb,\l)$, nous notons $\dot{w}$ un representant de 
$w$ dans $N_{\Gb^F}(\Lb,\l)$ et nous d\'efinissons, pour tous 
$f \in \Ind_{\Pb^F}^{\Gb^F} M_1$ et $g \in \Gb^F$,
$$T_w(f)(g)=\sum_{u \in \Ub^F}\s_1(\dot{w})f(\dot{w}^{-1}ug).$$
Alors, d'apr\`es \cite[proposition 3.9]{HL}, 
$(T_w)_{w \in W_{\Gb^F}(\Lb,\l)}$ est une base de $\HC_1$.
De m\^eme, nous d\'efinissons, pour tous $w \in W_{\Gb^F}(\Lbt,\lamt)$,  
$\fti \in \Ind_{\Pbt^F}^{\Gbt^F} \Mti$ et $\gti \in \Gbt^F$,
$$\Tti_w(\fti)(\gti)=\sum_{u \in \Ub^F}\sigt(\dot{w})\fti(\dot{w}^{-1}u\gti).$$
Alors, encore d'apr\`es \cite[proposition 3.9]{HL}, 
$(\Tti_w)_{w \in W_{\Gb^F}(\Lbt,\lamt)}$ est une base de $\HCt$. 

D'apr\`es \cite[corollaire 5.4]{HL}, il existe des \isos de $\qlb$-algebres 
\equat\label{iso 1}
\HC_1 \simeq \qlb[W_{\Gb^F}(\Lb,\l)]
\endequat
et
\equat\label{iso tilde}
\HCt\simeq \qlb[W_{\Gb^F}(\Lbt,\lamt)].
\endequat
De plus, l'image de $\Tti_w$ dans $\HC_1$ est $T_w$ d'apr\`es \ref{resresres} donc, 
par un argument de d\'eformation, les \isos ci-dessus peuvent \^etre choisis 
de sorte que le diagramme 
\equat
\diagram\label{diagram 1}
\HCt \rrto^{\sim\quad\quad} \dto && \qlb[W_{\Gb^F}(\Lbt,\lamt)] \dto \\
\HC_1 \rrto^{\sim\quad\quad} && \qlb[W_{\Gb^F}(\Lb,\l)]
\enddiagram
\endequat
soit commutatif. Nous ferons bien s\^ur ce choix-l\`a par la suite.

\bigskip

\remarque{racine carree} Une fois choisie une racine carr\'ee de $q$ dans $\qlb$
et une fois choisis $\n_1$ et $\nut$, 
alors les bijections entre ensembles de \cars \irrs induites par les \isos 
\ref{iso 1} et \ref{iso tilde}
sont uniquement d\'etermin\'ees (voir \cite[th\'eor\`eme 4.8]{HL2} 
et \cite[th\'eor\`eme 2.9]{benson}). 

\bigskip

\def\Ret{{\tilde{\Re}}}

\soussection{Param\'etrage des \cars dans une s\'erie de 
Harish-Chandra\label{soussection parametrage}} 
D'apr\`es le corollaire \ref{stabilo} (b), on a 
$$\EC(\Gb^F,\Lb,\l)=\coprod_{x \in \Lbt^F/\Lbt^F(\l_1)} \EC(\Gb^F,\Lb,\l_x).$$
Les isomorphismes \ref{iso 1} et \ref{iso tilde} induisent des bijections 
$$\fonctio{\Irr W_{\Gb^F}(\Lb,\l)}{\EC(\Gb^F,\Lb,\l)}{\eta}{R_\eta}$$
$$\fonctio{\Irr
W_{\Gbt^F}(\Lbt,\lamt)}{\EC(\Gbt^F,\Lbt,\lamt)}{\ch}{\Rti_\ch.}\leqno{\mathrm{et}}$$

\bigskip

\theoreme{harish chandra 1}
{\it Avec les notations ci-dessus, on a~:
\begin{itemize}
\itemth{a} $R_{\Lb \incl \Pb}^\Gb \l_1 = 
\DS{\sum_{\eta \in \Irr W_{\Gb^F}(\Lb,\l)}} \eta(1) R_\eta$. 

\itemth{b} Soit $\eta \in \Irr W_{\Gb^F}(\Lb,\l)$ et 
$\ch \in \Irr W_{\Gb^F}(\Lbt,\lamt)$. Alors 
$$\langle R_\eta , \Res_{\Gb^F}^{\Gbt^F} \Rti_\ch \rangle_{\Gb^F} =
\langle \Res_{W_{\Gb^F}(\Lbt,\lamt)}^{W_{\Gb^F}(\Lb,\l)} \eta, \ch 
\rangle_{W_{\Gb^F}(\Lbt,\lamt)}.$$

\itemth{c} Soit $l \in \Lbt^F(\l_1)$ et soit $\x$ le \car lin\'eaire 
de $W_{\Gb^F}(\Lb,\l)$ associ\'e \`a $l$ via l'\iso 
\ref{W/W}. 
Soit $\eta$ un \car \irr de $W_{\Gb^F}(\Lb,\l)$. Alors 
$$\lexp{l}{R_\eta}=R_{\eta \otimes \x}.$$
\end{itemize}}

\bigskip

\proof (a) est imm\'ediat.

(b) Faisons ici agir $\Gb^F$ et $\Gbt^F$ \`a droite sur $M$ et $\Mti$ 
respectivement. 
Soient $V_\eta$ (\resp $\Vti_\chi$) un $\HC_1$-module (\resp $\HCt$-module) 
\irr ayant pour \car $\eta$ (\resp $\chi$) \`a travers 
les isomorphismes \ref{iso 1} et \ref{iso tilde}. On voit 
$V_\eta$ comme un $\HC$-module. Posons 
$M_\eta=M^* \otimes_\HC V_\eta$ et $\Mti_\ch=\Mti^* \otimes_\HCt \Vti_\ch$. 
Alors $M_\eta$ (\resp $\Mti_\chi$) est un $\Gb^F$-module (\resp $\Gbt^F$-module) 
\`a droite \irr ayant pour \car $R_\eta$ (\resp $\Rti_\ch$). 
On a 
\eqna
\langle R_\eta,\Res_{\Gb^F}^{\Gbt^F} \Rti_\ch \rangle_{\Gb^F} 
&=&\dim_\qlb \Hom_{\Gb^F}(M_\eta, \Res_{\Gb^F}^{\Gbt^F} \Mti_\chi)\\
&=&\dim_\qlb M_\eta^* \otimes_{\qlb\Gb^F} \Res_{\Gb^F}^{\Gbt^F} \Mti_\chi \\
&=&\dim_\qlb V_\eta^* \otimes_\HC (M \otimes_{\qlb\Gb^F} \Res_{\Gb^F}^{\Gbt^F} \Mti^*) 
\otimes_{\HCt} \Vti_\chi.
\endeqna
Or, $\Res_{\Gb^F}^{\Gbt^F} \Mti^*=M^*$ et $M \otimes_{\qlb\Gb^F} M^*=\HC$. 
Par cons\'equent, 
\eqna
\langle R_\eta,\Res_{\Gb^F}^{\Gbt^F} \Rti_\ch \rangle_{\Gb^F} 
&=&\dim_\qlb V_\eta^* \otimes_\HC \HC \otimes_{\HCt} \Vti_\chi \\
&=&\dim_\qlb V_\eta^* \otimes_\HCt \Vti_\chi.
\endeqna
La commutativit\'e du diagramme \ref{diagram 1} montre que ce dernier 
terme est \'egal \`a 
$$\langle \Res_{W_{\Gb^F}(\Lbt,\lamt)}^{W_{\Gb^F}(\Lb,\l)} \eta, \ch 
\rangle_{W_{\Gb^F}(\Lbt,\lamt)},$$ 
ce qui montre (b).

\medskip

(c) Puisque $l \in \Lbt^F(\l_1)$, l'\auto $\sigt(l)$ de $M$ stabilise $M_1$. 
Soit $f \in \Ind_{\Pb^F}^{\Gb^F} M_1$. Posons 
$$\fonction{\o_l(f)}{\Gb^F}{M_1}{g}{\sigt(l)(f(l^{-1}gl)).}$$
Alors $\o_l(f) \in \Ind_{\Pb^F}^{\Gb^F} M_1$ et 
$$\o_l : \Ind_{\Pb^F}^{\Gb^F} M_1 \longto \Ind_{\Pb^F}^{\Gb^F} M_1$$
est un \iso de $\qlb$-espaces vectoriels. De plus, si 
$g \in \Gb^F$, alors $l^{-1}gl$ agit sur 
$\Ind_{\Pb^F}^{\Gb^F} M_1$ comme $\o_l^{-1} g \o_l$. 

Si $e_\eta$ d\'esigne l'idempotent primitif central de $\HC_1$ associ\'e \`a $\eta$, 
alors (c) est \'equivalent \`a l'\'egalit\'e 
$$\o_l e_\eta \o_l^{-1}=e_{\eta \otimes \x}.$$
Par un argument de d\'eformation, montrer cette \'egalit\'e revient 
\`a montrer que 
$$\o_l T_w \o_l^{-1} = \t_w(l)^{-1} T_w$$
pour tout $w \in W_{\Gb^F}(\Lb,\l)$, o\`u $\t_w$ est le \car lin\'eaire de 
$\Lbt^F(\l_1)$ associ\'e \`a $l$ par l'\iso 
$$W_{\Gb^F}(\Lb,\l)/W_{\Gb^F}(\Lbt,\lamt) \simeq (\Lbt^F(\l_1)/\Lbt^F(\Gb,\l_1))^\we.$$
Donc soient $f \in \Ind_{\Pb^F}^{\Gb^F} M_1$ et $g \in \Gb^F$. 
Alors
\eqna
(\o_l T_w \o_l^{-1} f)(g) &=& \sigt(l)(T_w \o_l^{-1} f)(l^{-1}gl) \\
&=&\DS{\sum_{u \in \Ub^F}} \sigt(l)\s_1(\dot{w}).(\o_l^{-1} f)(\dot{w}^{-1}ul^{-1}gl) \\
&=&\DS{\sum_{u \in \Ub^F}} \sigt(l)\s_1(\dot{w})\sigt(l)^{-1}.f(l\dot{w}^{-1}ul^{-1}g) \\
&=&\DS{\sum_{u \in \Ub^F}} \sigt(l)\s_1(\dot{w})\sigt(l)^{-1}.f(l\dot{w}^{-1}l^{-1}ug) \\
&=&\DS{\sum_{u \in \Ub^F}} \sigt(l)\s_1(\dot{w})\sigt(l)^{-1}.f(l\dot{w}^{-1}l^{-1}
\dot{w}\dot{w}^{-1}ug) \\
&=&\DS{\sum_{u \in \Ub^F}} \sigt(l)\s_1(\dot{w})\sigt(l)^{-1}\s_1(l\dot{w}^{-1}l^{-1}
\dot{w}).f(\dot{w}^{-1}ug) \\
&=&\sigt(l)\s_1(\dot{w})\sigt(l)^{-1}\s_1(l\dot{w}^{-1}l^{-1}
\dot{w})\s_1(\dot{w})^{-1}.(T_w f)(g).
\endeqna
D'apr\`es \ref{resresres}, on a $\s_1(l\dot{w}^{-1}l^{-1}\dot{w})=\sigt(l)
\sigt(\dot{w}^{-1}l^{-1}\dot{w})$. Donc 
$$(\o_l T_w \o_l^{-1} f)(g)=\sigt(l)\s_1(\dot{w})\s_1(\dot{w}^{-1}l^{-1}
\dot{w})\s_1(\dot{w})^{-1}.(T_w f)(g).$$
De plus, par d\'efinition de $\t_w$, la repr\'esentation
$$\fonctio{\Lbt^F(\l_1)}{\Gb\Lb(M_1)}{l}{\s_1(\dot{w})\s_1(\dot{w}^{-1}l
\dot{w})\s_1(\dot{w})^{-1}}$$
admet le m\^eme \car que la repr\'esentation 
$$\fonctio{\Lbt^F(\l_1)}{\Gb\Lb(M_1)}{l}{\t_w(l)\sigt(l)}$$
et leur restriction \`a $\Lb^F$ est \'egale \`a $\s_1$. 
Donc on a 
$$\s_1(\dot{w})\s_1(\dot{w}^{-1}l\dot{w})\s_1(\dot{w})^{-1}=\t_w(l)\sigt(l)$$
pour tout $l \in \Lbt^F(\l_1)$. D'o\`u le r\'esultat.\fin

\bigskip

\remarque{variation} 
Une fois choisi l'\iso d'alg\`ebres $\HCt \simeq \qlb[W_{\Gb^F}(\Lbt,\lamt)]$, 
l'\iso $\HC_1 \simeq \qlb[W_{\Gb^F}(\Lb,\l)]$ est d\'etermin\'e 
\`a un \car lin\'eaire pr\`es de $W_{\Gb^F}(\Lb,\l)/W_{\Gb^F}(\Lbt,\lamt)$ 
c'est-\`a-dire, d'apr\`es le th\'eor\`eme 
\ref{harish chandra 1} (c), \`a conjugaison pr\`es par un \'el\'ement de 
$\Lbt^F(\l_1)$.\finl

\bigskip

\soussection{Param\'etrage de 
${\boldsymbol{\EC(}}\Gb^{\boldsymbol{F}},\Lb,{\boldsymbol{\l)}}$} 
Soit maintenant $\eta$ un \car \irr de $W_{\Gb^F}^\pr(\Lbt,\lamt)$. Alors la 
restriction de $\eta$ au sous-groupe central $(\Lbt^F/\Lbt^F(\l_1))^\we$ de 
$W_{\Gb^F}^\pr(\Lbt,\lamt)$ est un multiple d'un \car lin\'eaire 
de $(\Lbt^F/\Lbt^F(\l_1))^\we$. Donc 
$\eta$ definit un \'el\'ement $x_\eta$ de $\Lbt^F/\Lbt^F(\l_1)$. 
Choisissons un relev\'e $\xti_\eta$ de $x_\eta$ dans $\Lbt^F/\Lbt^F(\Gb,\l_1)$ 
et soit 
$\x_\eta$ le \car lin\'eaire de $W^\pr_{\Gb^F}(\Lbt,\lamt)$ 
associ\'e \`a $\xti_\eta$ par l'isomorphisme 
$\bigl(W_{\Gb^F}^\pr(\Lbt,\lamt)/W_{\Gb^F}(\Lbt,\lamt)\bigr)^\we \simeq 
\Lbt^F/\Lbt^F(\Gb,\l_1)$. 
Alors $\eta \otimes \x_\eta^{-1}$ est un \car \irr de $W_{\Gb^F}^\pr(\Lbt,\lamt)$
et $(\Lbt^F/\Lbt^F(\l_1))^\we$ est contenu dans son noyau. 
Donc il peut \^etre vu comme un \car \irr de $W_{\Gb^F}(\Lb,\l)$. On d\'efinit alors 
\equat\label{definition Reta}
R_\eta=\lexp{\xti_\eta}{R_{\eta \otimes \x_\eta^{-1}}}.
\endequat
Remarquons que $R_\eta$ ne d\'epend pas du choix du repr\'esentant $\xti_\eta$ 
de $x_\eta$ (voir th\'eor\`eme \ref{harish chandra 1} (c)). 
De plus, si $\eta$ contient $(\Lbt^F/\Lbt^F(\l_1))^\we$ dans son noyau, alors 
$\eta$ peut \^etre vu comme un caract\`ere irr\'eductible de $W_{\Gb^F}(\Lb,\l)$ 
et le caract\`ere $R_\eta$ d\'efini par \ref{definition Reta} 
co\"\i ncide avec le caract\`ere $R_\eta$ d\'efini au d\'ebut de la sous-section 
\ref{soussection parametrage}. S'il y a ambigu\"\i t\'e, ce 
\car \irr sera not\'e $R_\eta^\Gb$. 
Le th\'eor\`eme suivant est une cons\'equence de cette discussion~:

\bigskip

\theoreme{harish}
{\it \begin{itemize} 
\itemth{a}L'application
$$\fonctio{\Irr W^\pr_{\Gb^F}(\Lbt,\lamt)}{\EC(\Gb^F,\Lb,\l)}{\eta}{R_\eta}$$
est bijective. 

\itemth{b} On a 
$$R_{\Lb \incl \Pb}^\Gb \l = 
\sum_{\eta \in \Irr W_{\Gb^F}^\pr(\Lbt,\lamt)} \eta(1) R_\eta.$$

\itemth{c} Soit $\eta$ et $\ch$ deux \cars \irrs de 
$W_{\Gb^F}^\pr(\Lbt,\lamt)$ et $W_{\Gb^F}(\Lbt,\lamt)$ respectivement. Alors 
$$\langle R_\eta , \Res_{\Gb^F}^{\Gbt^F} \Rti_\ch \rangle_{\Gb^F} =
\langle \Res_{W_{\Gb^F}(\Lbt,\lamt)}^{W_{\Gb^F}^\pr(\Lbt,\lamt)} \eta, \ch 
\rangle_{W_{\Gb^F}(\Lbt,\lamt)}.$$
 
\itemth{d} Soit $\eta \in \Irr W_{\Gb^F}^\pr(\Lbt,\lamt)$. Alors $R_\eta \in 
\EC(\Gb^F,\Lb,\l_{x_\eta})$.

\itemth{e} Soient $l \in \Lbt^F$ et $\eta \in \Irr W_{\Gb^F}^\pr(\Lbt,\lamt)$. Notons 
$\x_l$ le \car lin\'eaire de $W_{\Gb^F}^\pr(\Lbt,\lamt)$ associ\'e \`a $l$ via 
l'\iso $\bigl(W_{\Gb^F}^\pr(\Lbt,\lamt)/W_{\Gb^F}(\Lbt,\lamt)\bigr)^\we 
\simeq \Lbt^F/\Lbt^F(\Gb,\l_1)$ 
induit par \ref{W'/W}. Alors
$$\lexp{l}{R_\eta}=R_{\eta \otimes \x_l}.$$
\end{itemize}}

\bigskip

\corollaire{coro harish}
{\it Soit $\eta$ un \car \irr de $W_{\Gb^F}^\pr(\Lb,\l)$. Alors 
\begin{itemize}
\itemth{a} La restriction de $\eta$ \`a $W_{\Gb^F}(\Lbt,\lamt)$ 
est sans multiplicit\'e.

\itemth{b} Soit $\Lbt^F_\eta$ le sous-groupe de $\Lbt^F$ contenant $\Lbt^F(\Gb,\l_1)$ 
tel que  
$$\Lbt^F_\eta/\Lbt^F(\Gb,\l_1) \simeq \{\x \in 
(W_{\Gb^F}^\pr(\Lbt,\lamt)/W_{\Gb^F}(\Lbt,\lamt))^\we~|~ \eta \otimes \x = \eta\}.$$
(Rappelons que l'on a un \iso 
$(W_{\Gb^F}^\pr(\Lbt,\lamt)/W_{\Gb^F}(\Lbt,\lamt))^\we\simeq \Lbt^F/\Lbt^F(\Gb,\l_1)$.) 
Alors 
$$\Gbt^F(R_\eta)=\Gb^F.\Lbt^F_\eta.$$
\end{itemize}}

\bigskip

\proof (a) r\'esulte du th\'eor\`eme \ref{harish} (c) et du th\'eor\`eme 
\ref{restriction}. (b) d\'ecoule du th\'eor\`eme \ref{harish} (d).\fin

\bigskip

\remarques{remarques harish} 
(a) Gardons les notations du corollaire \ref{coro harish} (b). Le sous-groupe 
$(\Lbt^F/\Lbt^F(\l_1))^\we$ est central donc tout \car lin\'eaire $\x$ de 
$W_{\Gb^F}^\pr(\Lbt,\lamt)/W_{\Gb^F}(\Lbt,\lamt)$ v\'erifiant $\eta\otimes \x=\eta$ 
doit contenir ce sous-groupe dans son noyau. Donc 
$\Lbt^F(\Gb,\l_1) \incl 
\Lbt^F_\eta \incl \Lbt^F(\l_1)$ (comparer avec le corollaire \ref{s}).

\medskip

\tete{b} Comme dans la remarque suivant le th\'eor\`eme 
\ref{harish chandra 1}, le param\'etrage de 
$\EC(\Gb^F,\Lb,\l)$ donn\'e dans le th\'eor\`eme \ref{harish} (a) 
est bien d\'efini \`a tensorisation pr\`es par un \car lin\'eaire de 
$W_{\Gb^F}^\pr(\Lbt,\lamt)/W_{\Gb^F}(\Lbt,\lamt)$ 
(une fois l'\iso $\HCt \simeq \qlb[W_{\Gb^F}(\Lbt,\lamt)]$ choisi) 
donc il est bien d\'efini \`a conjugaison pr\`es par un \'el\'ement de 
$\Lbt^F$. Pour fixer pr\'ecis\'ement ce param\'etrage, 
il faut choisir une composante \irr de la restriction 
de $\Rti_1$ \`a $\Gb^F$ et associer \`a cette composante \irr 
le \car trivial de $W_{\Gb^F}^\pr(\Lbt,\lamt)$~: 
c'est \'equivalent \`a choisir une extension $\n_1$ 
de la composante \irr $\m_1$ de la restriction de $\nut$ \`a 
$N_\Gb^F(\Lbt,\lamt)$ telle que $\Res_{\Lb^F}^{N_{\Gb^F}(\Lbt,\lamt)}
\m_1=\l_1$.\finl

\bigskip

\soussection{Action de ${\boldsymbol{H^1(F,\ZC(\Gb))}}$} 
Soit $z \in H^1(F,\ZC(\Gb))$. Notons $\x_z$ le \car lin\'eaire du groupe  
$W_{\Gb^F}^\pr(\Lbt,\lamt)/W_{\Gb^F}(\Lbt,\lamt)$ image de $z$ 
par la suite de morphismes 
$$H^1(F,\ZC(\Gb)) \longmapright{\sim} \Lbt^F/\Lb^F.\Zb(\Gbt)^F \longto 
\Lbt^F/\Lbt^F(\Gb,\l_1) \longmapright{\sim} 
(W_{\Gb^F}^\pr(\Lbt,\lamt)/W_{\Gb^F}(\Lbt,\lamt))^\we.$$
Ici, le dernier isomorphisme est le dual de \ref{W'/W}. 
On peut voir $\x_z$ comme un \car lin\'eaire de $W_{\Gb^F}^\pr(\Lbt,\lamt)$. 
La proposition suivante est imm\'ediate.

\bigskip

\proposition{tauzg harish}
{\it Soient $z \in H^1(F,\ZC(\Gb))$ et $\eta \in \Irr W_{\Gb^F}^\pr(\Lbt,\lamt)$. 
Alors $\t_z^\Gb R_\eta=R_{\eta \x_z}$.}

\bigskip

\soussection{Induction de Harish-Chandra}
Pour une preuve de l'analogue du th\'eor\`eme suivant 
pour la s\'erie de Harish-Chandra $\EC(\Gb^F,\Lb,\l_1)$, 
voir \cite{HL2}. Le th\'eor\`eme ci-dessous en d\'ecoule facilement.

\bigskip

\theoreme{induit harish}
{\it Soit $\Mb$ un \levi $F$-stable de $\Gb$ tel que 
$\Lb \incl \Mb$. On choisit une bijection 
$$\fonctio{\Irr W_{\Mb^F}^\pr(\Lbt,\lamt)}{\EC(\Lb^F,\Lb,\l)}{\eta}{R_\eta^\Mb}$$
telle que $\langle R_{\Mb \incl \Qb}^\Gb R_1^\Mb, R_1^\Gb \rangle_{\Gb^F} \not = 0$. 
Alors 
$$\langle R_{\Mb \incl \Qb}^\Gb R_\eta^\Mb , R_\z^\Gb \rangle_{\Gb^F} = 
\langle \Ind_{W_{\Mb^F}^\pr(\Lbt,\lamt)}^{W_{\Gb^F}^\pr(\Lbt,\lamt)} \eta, \z 
\rangle_{W_{\Gb^F}^\pr(\Lbt,\lamt)}$$
pour tous \cars \irrs $\eta$ et $\z$ de $W_{\Mb^F}^\pr(\Lbt,\lamt)$ 
et $W_{\Gb^F}^\pr(\Lbt,\lamt)$ respectivement.}

\bigskip

\remarque{choix M} Choisir une bijection 
$$\fonctio{\Irr W_{\Mb^F}^\pr(\Lbt,\lamt)}{\EC(\Mb^F,\Lb,\l)}{\eta}{R_\eta^\Mb}$$
est \'equivalent \`a choisir une extension $\n_1^\pr$ de $\l_1$ \`a 
$N_{\Mb^F}(\Lb,\l)$ et une racine carr\'ee de $q$
(voir \cite[th\'eor\`eme 4.8]{HL2} et \cite[th\'eor\`eme 2.9]{benson}). 
Si nous demandons \`a cette bijection de v\'erifier 
$$\langle R_{\Mb \incl \Qb}^\Gb R_1^\Mb, R_1^\Gb \rangle_{\Gb^F} \not = 0,$$
alors nous devons choisir pour $\n_1^\pr$ la restriction de $\n_1$ et, 
pour la racine carr\'ee de $q$, la m\^eme que celle choisie pour le groupe $\Gb^F$.\finl

\bigskip

\newpage

\def\gelfand{{\mathrm{GG}}}
\def\Uni{{\mathrm{Uni}}}
\def\dlm{{\mathrm{dlm}}}

{\Large \part{Autour des caract\`eres de Gelfand-Graev\label{chapitre gelfand}}}

\bigskip

Le but de ce chapitre est d'\'etudier les caract\`eres de Gelfand-Graev 
ainsi que leurs composantes irr\'eductibles, appel\'es {\it caract\`eres 
r\'eguliers}. Dans les groupes \`a centre non connexe, il peut y avoir plusieurs 
caract\`eres de Gelfand-Graev (ils sont param\'etr\'es par $H^1(F,\ZC(\Gb))$) 
et ce ne sont pas forc\'ement des fonctions uniformes, contrairement 
\`a ce qui se passe dans les groupes \`a centre connexe \cite[\SEC 10]{delu}. 
Dans la section \ref{section gelfand}, nous rappelons la d\'efinition 
et les premi\`eres propri\'et\'es de ces caract\`eres, comme 
les r\'esultats de Digne, Lehrer et Michel \cite{DLM2} ou l'auteur 
\cite[partie II]{bonnafe action} sur la restriction de Lusztig. 
Dans la section \ref{section regulier}, nous rappelons comment 
sont param\'etr\'es les caract\`eres r\'eguliers ou semi-simples 
(un {\it caract\`ere semi-simple} est, au signe pr\`es, le dual d'Alvis-Curtis 
d'un caract\`ere r\'egulier). La section \ref{section semi cuspidal} est 
consacr\'ee \`a l'\'etude des s\'eries de Harish-Chandra au-dessus 
d'un caract\`ere semi-simple cuspidal en adaptant l'\'etude 
faite au chapitre pr\'ec\'edent \`a ce cas particulier. 
Dans la derni\`ere section de ce chapitre, nous \'etudions les combinaisons 
lin\'eaires d'induits de Lusztig de caract\`eres semi-simples.

\bigskip 

\section{Caract\`eres de Gelfand-Graev\label{section gelfand}}~

\medskip

\soussection{\'El\'ements unipotents r\'eguliers\label{unireg soussection}} 
Si $g \in \Gb$, alors $\dim C_\Gb(g) \ge \dim \Tb_0$. 
Un \'el\'ement $g$ de $\Gb$ est dit {\it r\'egulier} si $\dim C_\Gb(g)=\dim \Tb_0$. 
L'ensemble des \'el\'ements r\'eguliers de $\Gb$ forme un ouvert dense 
de $\Gb$ (voir \cite[th\'eor\`eme 1.3 (a)]{steinberg}). 

Concentrons-nous maintenant sur les \'el\'ements unipotents r\'eguliers. 
Tout d'abord, il existe des \'el\'ements unipotents r\'eguliers 
\cite[th\'eor\`eme 3.1]{steinberg}. Ils sont tous conjugu\'es dans $\Gb$ 
\cite[th\'eor\`eme 3.3]{steinberg}~: notons $\UCB_\reg^\Gb$ 
la classe de conjugaison des \'el\'ements unipotents r\'eguliers de $\Gb$. 
Un \'el\'ement unipotent $u$ de $\Gb$ 
est r\'egulier si et seulement si il est contenu dans un seul sous-groupe de 
Borel \cite[corollaire 3.12 (b)]{steinberg}. 
Un \'el\'ement unipotent 
$u \in \Ub_0$ est r\'egulier \ssi $u \not\in \Ub_I$ 
pour toute partie non vide $I$ de $\D_0$ \cite[lemme 3.2]{steinberg}. 
Notons $\Ub_{0,\reg}$ l'ensemble 
des \'el\'ements unipotents r\'eguliers de $\Ub_0$~: c'est 
un ouvert dense de $\Ub_0$. 

Soit $\Ub_1=\prod_{\a \in \Phi_0^+ - \D_0} \Ub_\a$, o\`u $\Phi_0^+$ est le 
syst\`eme de racines positives de $\Phi_0$ associ\'e \`a $\D_0$. Alors 
$\Ub_1$ est le groupe d\'eriv\'e de $\Ub_0$. De plus, l'action de $\Ub_1$ par 
translation (\`a droite ou \`a gauche) sur $\Ub_0$ stabilise $\Ub_{0,\reg}$. 
Il est de plus clair que $\Tbt_0$ agit transitivement sur $\Ub_{0,\reg}/\Ub_1$. 
En particulier, $(\Ub_{0,\reg}/\Ub_1)^F=\Ub_{0,\reg}^F/\Ub_1^F$ est non vide. 
D'autre part, le stabilisateur de tout \'el\'ement de $\Ub_{0,\reg}/\Ub_1$ dans 
$\Tbt_0$ est $\Zb(\Gbt)$. 
Ce dernier \'etant connexe, on en d\'eduit que $\Tbt_0^F$ agit 
transitivement sur $\Ub_{0,\reg}^F/\Ub_1^F$. 

\bigskip

\proposition{prop uni regulier}
{\it Soit $u \in \Ub_{0,\reg}$. Alors~:
\begin{itemize}
\itemth{a} $C_\Gb(u)=\Zb(\Gb). C_{\Ub_0}(u)$.

\itemth{b} Si $p$ est bon pour $\Gb$, alors $C_{\Ub_0}(u)$ est connexe.

\itemth{c} L'application $\Ub_0 \to u\Ub_1$, $x \mapsto xux^{-1}$ est surjective.

\itemth{d} Si $u'$ est un autre \'el\'ement de $\Ub_{0,\reg}$, alors 
il existe $b \in \Bb_0$ tel que $u'=bub^{-1}$.

\itemth{e} $\Bb_0$ est l'unique \borel de $\Gb$ contenant $u$.
\end{itemize}}

\bigskip

\proof (a) d\'ecoule du fait que le stabilisateur, dans $\Tb_0$, d'un \'el\'ement 
de $\Ub_{0,\reg}/\Ub_1$ est \'egal \`a $\Zb(\Gb)$. Pour (b), voir 
\cite[3.7]{springer steinberg}. Montrons maintenant (c). Soit $f : \Ub_0 \to u\Ub_1$, 
$x \mapsto xux^{-1}$. C'est un morphisme de vari\'et\'e. L'image de 
$f$ est une orbite sous l'action d'un groupe unipotent, donc 
c'est une sous-vari\'et\'e ferm\'ee de $u\Ub_1$. De plus, la dimension 
des fibres de $f$ est toujours \'egale \`a 
$$\dim C_{\Ub_0}(u)=\dim C_\Gb(u)-\dim \Zb(\Gb)=|\D_0|.$$ 
Par cons\'equent, la dimension de l'image de $f$ est 
$$\dim \Ub_0 - \dim C_{\Ub_0}(u)= |\Phi_0^+|-|\D_0|=\dim \Ub_1.$$
Puisque $u\Ub_1$ est irr\'eductible, l'image de $f$ est bien $u\Ub_1$. 
(d) d\'ecoule du fait qu'il n'y a qu'une classe de conjugaison 
d'\'el\'ements unipotents r\'eguliers, 
que chaque \'el\'ement unipotent r\'egulier est contenu dans un seul sous-groupe 
de Borel et que le normalisateur de $\Bb_0$ dans $\Gb$ est $\Bb_0$ 
\cite[th\'eor\`eme 11.16]{borel}. 
(e) est clair.\fin

\bigskip

\soussection{Caract\`eres r\'eguliers de ${\boldsymbol{\Ub_0^F}}$\label{sous chi}}
Un caract\`ere lin\'eaire $\psi : \Ub_0^F \to {\overline{\QM}}_\ell^\times$ 
est dit {\it r\'egulier} s'il contient $\Ub_1^F$ dans son noyau et si 
$\Res_{\Ub_I^F}^{\Ub_0^F} \psi \not= 1$ 
pour toute partie stricte $\phi_0$-stable $I$ de $\D_0$ 
(le lecteur peut ais\'ement v\'erifier que cette d\'efinition co\"\i ncide 
avec \cite[d\'efinition 2.3]{DLM1}). Notons $(\Ub_0^F)^\we_\reg$ 
l'ensemble des \cars lin\'eaires r\'eguliers de $\Ub_0^F$. 

Fixons maintenant et ce jusqu'\`a la fin de cet article un \car lin\'eaire 
non trivial $\chi_1 : \fp \to \qlb^\times$. Si $n$ est un entier naturel 
non nul, nous noterons 
$\ch_n : \FM_{p^n} \to \qlb^\times$, $x \mapsto \ch_1(\Tr_{\FM_{p^n}/\fp}(x))$. 
Avec ces choix, on obtient \cite[\SEC 2]{DLM2} une application 
bijective $\Tbt_0^F$-\'equivariante 
$\Ub_{0,\reg}^F/\Ub_1^F \longmapright{\sim} (\Ub_0^F)^\we_\reg$. 
Cela montre que le groupe $\Tbt_0^F$ agit transitivement 
sur l'ensemble des \cars r\'eguliers de $\Ub_0^F$ 
et que le stabilisateur d'un \car r\'egulier dans $\Tbt_0^F$ 
est $\Zb(\Gbt)^F$. Par cons\'equent, 
$\Tbt_0^F/\Tb_0^F.\Zb(\Gbt)^F \simeq H^1(F,\ZC(\Gb))$ agit librement 
sur l'ensemble des $\Tb_0^F$-orbites dans $(\Ub_0)_\reg^\wedge$. 

On appelle {\it \car de Gelfand-Graev} de $\Gb^F$ tout \car de la forme 
$\Ind_{\Ub_0^F}^{\Gb^F} \psi$, o\`u $\psi$ est un \car r\'egulier de $\Ub_0^F$. 
Notons $\Uni_\reg(\Gb^F)$ l'ensemble des classes de 
$\Gb^F$-conjugaison d'\'el\'ements unipotents r\'eguliers de $\Gb^F$ 
et $\gelfand(\Gb^F)$ l'ensemble des \cars de Gelfand-Graev de $\Gb^F$. 
La derni\`ere remarque du pr\'ec\'edent paragraphe montre 
que $H^1(F,\ZC(\Gb))$ agit transitivement sur $\gelfand(\Gb^F)$ 
(en particulier, il n'y a qu'un \car de Gelfand-Graev 
dans $\Gbt^F$).

Si $[u] \in \Uni_\reg(\Gb^F)$, alors $[u] \cap \Ub_0^F$ est 
contenu dans $\Ub_{0,\reg}^F$ et son image dans $\Ub_{0,\reg}^F/\Ub_1^F$ 
est une $\Tb_0^F$-orbite. On peut donc lui associer un 
\car de Gelfand-Graev $\G_u^\Gb$~: l'application
$$\fonctio{\Uni_\reg(\Gb^F)}{\gelfand(\Gb^F)}{[u]}{\G_u^\Gb}$$ 
est surjective et $H^1(F,\ZC(\Gb))$-\'equivariante. 
Nous verrons plus tard qu'elle est bijective. 

\bigskip 

\begin{quotation}
\noindent{\it Dor\'enavant, et ce jusqu'\`a la fin de cet article, 
nous fixons un \'el\'ement unipotent r\'egulier $u \in \Ub_{0,\reg}^F$, 
nous notons $\psi$ le \car r\'egulier de $\Ub_0^F$ associ\'e et nous posons
$$\G^\Gb=\Ind_{\Ub_0^F}^{\Gb^F} \psi.$$}
\end{quotation}

\bigskip

Une fois fix\'e $\psi$, l'ensemble des $\Tb_0^F$-orbites 
de \cars r\'eguliers de $\Ub_0^F$ est en bijection naturelle 
avec $\Tbt_0^F/\Tb_0^F.\Zb(\Gbt)^F 
\simeq H^1(F,\ZC(\Gb))$. 
Soient $z \in H^1(F,\ZC(\Gb))$ et $\tti_z \in \Tbt_0^F$ tels que 
$\s_{\Tb_0}^\Gb(\tti_z \Tb_0^F.\Zb(\Gbt)^F)=z$. Posons  
$\psi_z=\psi \ci (\ad \tti_z)^{-1}$. Alors $\{\psi_z~|~z \in H^1(F,\ZC(\Gb))\}$ 
est un ensemble de repr\'esentants des $\Tb_0^F$-orbites de 
\cars lin\'eaires r\'eguliers de $\Ub_0^F$. 
On d\'efinit alors
\equat
\G^\Gb_z=\Ind_{\Ub_0^F}^{\Gb^F} \psi_z.
\endequat
Remarquons que $\G^\Gb_z$ ne d\'epend que de $z$ et que 
\equat
\G^\Gb_z=\G^\Gb \ci (\ad \tti_z)^{-1} = \t_z^\Gb \G^\Gb.
\endequat
Donc les \cars de Gelfand-Graev de $\Gb^F$ sont les $\G^\Gb_z$ o\`u $z$ parcourt 
$H^1(F,\ZC(\Gb))$.

\bigskip

\soussection{Restriction de Harish-Chandra}  
Le th\'eor\`eme suivant a \'et\'e montr\'e par Digne, Lehrer et Michel 
\cite[th\'eor\`eme 2.9]{DLM1} (voir aussi \cite[preuve du lemme 3.6.1]{asai} 
pour un \'enonc\'e moins fort).

\bigskip

\Theoreme{Digne-Lehrer-Michel}{restriction gelfand}
{\it Soit $\Pb$ un \para $F$-stable de $\Gb$ et soit $\Lb$ un \levic 
$F$-stable de $\Pb$. 
Alors $\lexp{*}{R}_{\Lb \incl \Pb}^\Gb \G^\Gb$ est un \car de Gelfand-Graev de 
$\Lb^F$. Plus pr\'ecis\'ement, si $\Bb_0 \incl \Pb$ et $\Tb_0 \incl \Lb$, 
alors la restriction de $\psi$ 
\`a $\Ub_0^F \cap \Lb^F$ est un \car r\'egulier de $\Ub_0^F \cap \Lb^F$ et 
$$\lexp{*}{R}_{\Lb \incl \Pb}^\Gb \G^\Gb=
\Ind_{\Ub_0^F \cap \Lb^F}^{\Lb^F} (\Res_{\Ub_0^F \cap \Lb^F}^{\Ub_0^F} \psi).$$}

\bigskip

\noindent{\sc Notation - } Sous les hypoth\`eses et notations 
du th\'eor\`eme \ref{restriction gelfand}, on pose 
\equat
\G^\Lb = \lexp{*}{R}_{\Lb \incl \Pb}^\Gb \G^\Gb.
\endequat
Remarquons que $\G^\Lb$ ne d\'epend pas de $\Pb$. 
Si on pose $\G^\Lb_z=\t^\Lb_z\G^\Lb$ pour tout 
$z \in H^1(F,\Zb(\Lb))$, alors, d'apr\`es \ref{tauzg rlg}, on a 
\equat
\G^\Lb_{h_\Lb^1(z)} = \lexp{*}{R}_{\Lb \incl \Pb}^\Gb \G^\Gb_z
\endequat
pour tout $z \in H^1(F,\Zb(\Gb))$.\finl

\bigskip

\soussection{Dualit\'e d'Alvis-Curtis}
Les r\'esultats de cette section sont d\^us eux aussi \`a Digne, Lehrer 
et Michel. 

\bigskip

\Theoreme{Digne-Lehrer-Michel}{theo dual gelfand}
{\it \begin{itemize} 
\itemth{a} Si $g \in \Gb^F$ est tel que $D_\Gb \G^\Gb(g) \not= 0$, 
alors $g$ est un \'el\'ement unipotent r\'egulier.

\itemth{b} Si $z$ et $z'$ sont deux \'el\'ements de $H^1(F,\ZC(\Gb))$, alors
$$\langle D_\Gb \G^\Gb_z, \G_{z'}^\Gb \rangle_{\Gb^F}=\eta_\Gb |\Zb(\Gb)^F| \d_{z,z'}.$$
\end{itemize}}

\bigskip

\proof voir \cite[th\`eor\`eme 2.12 (i) et (ii)]{DLM1}.\fin

\bigskip

\Corollaire{Digne-Lehrer-Michel}{coro different gelfand}
{\it La famille $(\G_z^\Gb)_{z \in H^1(F,\ZC(\Gb))}$ est libre. 
En particulier, si $z$ et $z'$ sont deux \'el\'ements 
distincts de $H^1(F,\ZC(\Gb))$, alors $\G_z^\Gb\not= \G_{z'}^\Gb$.}

\bigskip

\soussection{Restriction de Lusztig\label{sous reslg}}~
Nous fixons dans cette sous-section un sous-groupe de Levi $F$-stable $\Lb$ 
de $\Gb$. 

\bigskip

\begin{quotation}
\noindent{\bf Hypoth\`ese :} {\it Nous supposerons dans cette sous-section, 
et seulement dans cette sous-section, que $p$ est un bon nombre premier 
pour $\Gb$.}
\end{quotation}

\bigskip

Dans \cite[\SEC 2]{bonnafe cras}, l'auteur a d\'efini une application 
$$\res_\Lb^\Gb : \Uni_\reg(\Gb^F) \longto \Uni_\reg(\Lb^F).$$
Rappelons sa d\'efinition (pour la preuve des faits utilis\'es dans la 
discussion suivante, nous nous r\'ef\'erons \`a \cite{bonnafe cras}).

Tout d'abord, si $\Lb$ est un \levic $F$-stable d'un \para $F$-stable $\Pb$ de $\Gb$, 
notons $\pi_{\Lb \incl \Pb} : \Pb \to \Lb$ la projection canonique et posons
$$\fonction{\r_\Lb^\Gb}{\Uni_\reg(\Gb^F)}{\Uni_\reg(\Lb^F)}{[g]_{\Gb^F}} 
{\pi_{\Lb \incl \Pb}([g]_\Gb^F \cap \Pb).}$$
Il est bien connu que $\r_\Lb^\Gb$ est bien d\'efinie et ne d\'epend pas 
du choix du \para $F$-stable ayant $\Lb$ comme sous-groupe de Levi 
\cite[proposition 5.3]{DLM1}. 

Revenons au cas g\'en\'eral. Soit $\Lb_1$ un \levi $F$-stable 
$\Lb$-d\'eploy\'e de $\Lb$ minimal tel que l'application 
$h_{\Lb_1}^\Lb$ soit un isomorphisme. Le groupe $\Lb_1$ \'etant 
cuspidal, il existe alors un \levi $F$-stable $\Gb$-d\'eploy\'e 
$\Lb_2$ qui lui est g\'eom\'etriquement conjugu\'e. Alors 
l'application 
$$\fonction{c_\Lb}{\Uni_\reg(\Lb_2^F)}{\Uni_\reg(\Lb_1^F)}{[l]_{\Lb_2^F}}
{[l]_{\Gb^F} \cap \Lb_1^F}$$
est bien d\'efinie et bijective. De plus, $\r_{\Lb_1}^\Lb$ est elle 
aussi bijective. On pose alors, comme dans \cite[page 279]{bonnafe cras}, 
$$\res_\Lb^\Gb = (\r_{\Lb_1}^\Lb)^{-1} ~\circ~ c_\Lb ~\circ~\r_{\Lb_2}^\Gb.$$

\medskip

Nous allons maintenant d\'efinir une autre application  
$\Uni_\reg(\Gb^F) \to \Uni_\reg(\Lb^F)$. Notons $x$ un \'el\'ement 
de $\Gb$ tel que $\Lb_1=\lexp{g}{\Lb_2}$. Puisque $\Lb_1$ et $\Lb_2$ sont 
$F$-stables, $F(g)g^{-1}$ appartient \`a $N_\Gb(\Lb_1)$. Notons 
$w_\Lb$ sa classe dans $W_\Gb(\Lb_1)$. Nous noterons ici 
$\ph_{\Lb_1}^\Gb : W_\Gb(\Lb_1) \to \ZC(\Lb_1)$ le morphisme de groupes 
not\'e $\ph_{\Lb_1,v_1}^\Gb$ dans \cite[partie I, corollaire 3.8]{bonnafe action}, 
o\`u $v_1$ est un \'el\'ement unipotent r\'egulier de $\Lb_1$. 
Il est \`a noter que ce morphisme a \'et\'e calcul\'e explicitement 
dans \cite[partie II, table 1]{bonnafe action}.
Identifions $\ZC(\Lb_1)$ et $\ZC(\Lb)$ via $h_{\Lb_1}^\Lb$ 
et notons $z_\Lb$ l'image dans $H^1(F,\ZC(\Lb))$ de 
l'\'el\'ement $\ph_{\Lb_1}^\Gb(w_\Lb)$ de $\ZC(\Lb)$. On pose alors
$$\dlm_\Lb^\Gb = \t_{z_\Lb}^\Lb ~\circ~ \res_\Lb^\Gb.$$
Nous rappelons les propri\'et\'es des applications $\res_\Lb^\Gb$ 
et $\dlm_\Lb^\Gb$. 

\bigskip

\proposition{prop prop res}
{\it On a~:
\begin{itemize}
\itemth{a} Les applications $\res_\Lb^\Gb$ et $\dlm_\Lb^\Gb$ ne d\'ependent 
pas des choix de $\Lb_1$, $\Lb_2$ et $g$ effectu\'es ci-dessus.

\itemth{b} Si $z \in H^1(F,\ZC(\Gb))$, alors 
$\res_\Lb^\Gb \circ \t_z^\Gb = \t_{h_\Lb(z)}^\Lb \circ \res_\Lb^\Gb$ et 
$\dlm_\Lb^\Gb \circ \t_z^\Gb = \t_{h_\Lb(z)}^\Lb \circ \dlm_\Lb^\Gb$.

\itemth{c} Les applications $\res_\Lb^\Gb$ et $\dlm_\Lb^\Gb$ sont surjectives. 

\itemth{d} Si $\Mb$ est un \levi $F$-stable de $\Lb$, alors
$\res_\Mb^\Gb = \res_\Mb^\Lb \circ \res_\Lb^\Gb$ et 
$\dlm_\Mb^\Gb = \dlm_\Mb^\Lb \circ \dlm_\Lb^\Gb$.
\end{itemize}}

\bigskip

\proof Les assertions (a), (b) et (c) sont \'evidentes. 
Le fait que $\res_\Mb^\Gb = \res_\Mb^\Lb \circ \res_\Lb^\Gb$ 
est d\'emontr\'e dans \cite[proposition 7.5 (c)]{bonnafe regulier}. 
Il suffit alors d'appliquer \cite[partie II, corollaire 12.4]{bonnafe action} 
pour obtenir que $\dlm_\Mb^\Gb = \dlm_\Mb^\Lb \circ \dlm_\Lb^\Gb$.\fin

\bigskip

\noindent{\sc Notation - } Dor\'enavant, et ce jusqu'\`a la fin de cet article, 
nous noterons, lorsque $p$ est bon pour $\Gb$, 
$u_\Lb$ un repr\'esentant de la classe de $\Lb^F$-conjugaison 
$\dlm_\Lb^\Gb [u]_{\Gb^F}$. D'autre part, nous noterons $\G^\Lb$ 
le \car de Gelfand-Graev de $\Lb^F$ associ\'e \`a $[u_\Lb]_{\Lb^F}$.\finl

\bigskip

La conjecture suivante propose une g\'en\'eralisation du th\'eor\`eme 
de Digne-Lehrer-Michel sur la restriction de Harish-Chandra d'un 
\car de Gelfand-Graev. 

\bigskip

\noindent{\bf Conjecture $\GGB$ :} {\it Si $\Lb$ est un \levic $F$-stable 
d'un \para $\Pb$ de $\Gb$, alors 
$\resr_{\Lb \incl \Pb}^\Gb \G^\Gb=\e_\Gb\e_\Lb \G^\Lb$.}

\bigskip

\noindent{\bf Conjecture $\GGB'$ :} {\it Si $\Lb$ est un \levic $F$-stable 
d'un \para $\Pb$ de $\Gb$, alors 
$\resr_{\Lb \incl \Pb}^\Gb \g_u^\Gb=\g_{u_\Lb}^\Lb$.}

\bigskip

Nous dirons que {\it ``la conjecture $\GG$ (resp. $\GG'$) a lieu dans 
$\Gb$''} si, pour tout sous-groupe de Levi $\Mb$ de $\Gb$, pour tout 
\para $\Pb$ de $\Mb$ et pour tout \levic $F$-stable $\Lb$ de $\Pb$, on a 
$\resr_{\Lb \incl \Pb}^\Mb \G^\Mb=\e_\Mb\e_\Lb \G^\Lb$ 
(resp. $\resr_{\Lb \incl \Pb}^\Gb \g_{u_\Mb}^\Mb=\g_{u_\Lb}^\Lb$). 

\bigskip

\proposition{equivalence}
{\it Si la formule de Mackey a lieu dans $\Gb$, alors la conjecture 
$\GG$ a lieu dans $\Gb$ si et seulement si la conjecture $\GG'$ 
a lieu dans $\Gb$.}

\bigskip

\proof Cela r\'esulte facilement de \cite[propositions 2.1 et 2.5]{DLM2}.\fin

\bigskip

Lorsque le centre de $\Gb$ est connexe, ces conjectures ont 
\'et\'e montr\'ees par Digne-Lehrer-Michel \cite[proposition 5.4]{DLM1}~: 
en effet, dans ce cas, les fonctions centrales $\G^\Gb$ et $\g^\Gb_u$ sont 
des combinaisons lin\'eaires explicites de \cars de Deligne-Lusztig 
$R_\Tb^\Gb(\th)$, o\`u $(\Tb,\th) \in \nabla(\Gb,F)$. 
Il suffit alors de calculer $\resr_{\Lb \incl \Pb}^\Gb R_\Tb^\Gb(\th)$~: 
cela se fait en utilisant la formule de Mackey qui est valable 
dans ce cas (voir th\'eor\`eme \ref{theo mackey} (a2)).

Lorsque le centre de $\Gb$ n'est pas connexe, ces conjectures ne sont 
d\'emontr\'ees qu'en utilisant la th\'eorie des faisceaux-caract\`eres, 
ce qui restreint leur domaine de validit\'e 
(notamment \`a cause de l'emploi de \cite[th\'eor\`eme 1.14]{lugf}).

\bigskip

\theoreme{theo restriction gamma}
{\it Si $p$ est bon pour $\Gb$, si $F$ est un $\fq$-endomorphisme de Frobenius 
de $\Gb$ et si $q > q_0(\Gb)$, o\`u $q_0(\Gb)$ est un constante ne 
d\'ependant que de la donn\'ee radicielle associ\'ee \`a $\Gb$, alors 
les conjectures $\GG$ et $\GG'$ ont lieu dans $\Gb$.}

\bigskip

\noindent{\sc Remarque - } Dans \cite[th\'eor\`eme 3.7]{DLM2}, Digne-Lehrer-Michel 
ont montr\'e que $\resr_{\Lb \incl \Pb}^\Gb \G^\Gb$ est \'egal, 
au signe $\e_\Gb\e_\Lb$ pr\`es, \`a un \car de Gelfand-Graev de 
$\Lb^F$. En revanche, ils n'ont pas d\'etermin\'e lequel. 
En \'etudiant plus pr\'ecis\'ement l'alg\`ebre d'endomorphismes 
de l'induit d'un faisceau-caract\`ere cuspidal support\'e 
par la classe unipotente r\'eguli\`ere \cite{bonnafe action} 
et en int\'egrant cette information suppl\'ementaire 
dans la preuve de Digne-Lehrer-Michel, nous avons obtenu 
le th\'eor\`eme ci-dessus 
\cite[partie II, th\'eor\`eme 15.2]{bonnafe action}.\finl

\bigskip

\section{Caract\`eres r\'eguliers et \cars semi-simples\label{section regulier}}

\bigskip

\begin{quotation}
\noindent{\it Dor\'enavant, et ce jusqu'\`a la fin de cet article, nous fixons 
un \'el\'ement semi-simple $s \in \Gb^{*F^*}$. Nous fixons aussi un \'el\'ement 
semi-simple $\sti \in \Gbt^{*F^*}$ tel que $i^*(\sti)=s$ et nous posons 
$s'=i^{\prime *}(\sti) \in \Gb^{\prime *F^*}$.}
\end{quotation}

\bigskip

Nous reprenons les notations introduites dans la section 
\ref{section semisimple} ($\Tb_1^*$, $\Bb_1^*$, $\Phi_s$, $\phi_1$...). 
Pour tout $w \in W^\ci(s)$, nous choisissons un \tor $F^*$-stable $\Tb_w^*$ de 
$C_{\Gb^*}^\ci(s)$ de 
type $w$ par rapport \`a $\Tb_1^*$. 
Notons $\e : W \to \{1,-1\}$ le \car signature de $W$. Nous noterons 
$\e_s$ (\resp $\e_s^\circ$) sa restriction \`a $W(s)$ (\resp 
$W^\circ(s)$). Alors, 
si $w \in W^\circ(s)$, on a $\e_{\Tb_w^*} = \e(w) \e_{C_{\Gb^*}^\circ(s)}$ 
et donc, d'apr\`es le corollaire \ref{coro independance},  
\equat\label{dg rtgs}
D_\Gb R_{\Tb_w^*}^\Gb(s) = \e_\Gb \e_{C_{\Gb^*}^\circ(s)} \e(w)R_{\Tb_w^*}^\Gb(s).
\endequat

\bigskip

\remarque{rem epsilon ags}
Il se peut que la restriction du \car signature \`a $A_{\Gb^*}(s)$ 
soit non triviale, comme le montre le cas o\`u $\Gb=\Sb\Lb_2(\FM)$, 
$F : \Gb \to \Gb$ est un endomorphisme de Frobenius d\'eploy\'e 
et $s$ est l'unique \'el\'ement d'ordre $2$ de $\Tb_0^*$.\finl

\bigskip

Soient $\Tbt_1^*=i^{*-1}(\Tb_1^*)$. 
Alors $W$ est canoniquement isomorphe au groupe de Weyl de 
$\Gbt^*$ relativement \`a $\Tbt_1^*$ 
et, puisque $C_{\Gbt^*}(\sti)$ est connexe 
(voir th\'eor\`eme \ref{cgs connexe}), on a $W(\sti)=W^\circ(\sti)$ et 
$A_{\Gbt^*}(\sti)=\{1\}$. De plus, $W(\sti)=W^\ci(s)$ car 
$W^\circ(s)$ est le groupe de Weyl du syst\`eme de racines $\Phi_s$. 
Pour tout $w \in W(\sti)$, on pose $\Tbt_w^*=i^{*-1}(\Tb^*_w)$. 
D'apr\`es la proposition \ref{tenseur} (a), on a, pour tout 
$a \in A_{\Gb^*}(s)$, 
\equat\label{tenseur rtgs}
R_{\Tbt_w^*}^\Gbt(\sti) \otimes \widehat{\ph_s(a)} = 
R_{\Tbt_{a^{-1}wa}^*}^\Gbt(\sti)
\endequat
car $\sti\ph_s(a)=a\sti a^{-1}$. En particulier, 
\equat
R_{\Tb_{awa^{-1}}^*}^\Gb(s)=R_{\Tb_w^*}^\Gb(s).
\endequat
D'apr\`es la formule de Mackey (th\'eor\`eme \ref{theo mackey} (2)), on a, 
pour tous $w$, $w' \in W^\circ(s)$,
\equat\label{scalaire rtg s}
\langle R_{\Tb_w^*}^\Gb(s),R_{\Tb_{w'}^*}^\Gb(s) \rangle_{\Gb^F} = 
\begin{cases} 
|C_{W(s)}(w\phi_1)| & \text{si $w\phi_1$ et $w'\phi_1$ sont conjugu\'es sous 
$W(s)$,} \\
0 & \text{sinon.}
\end{cases}
\endequat

\bigskip

\bigskip

\soussection{D\'efinition} 
Un \car \irr de $\Gb^F$ est dit 
{\it r\'egulier} (\resp {\it semi-simple}) s'il est 
une composante \irr d'un \car de Gelfand-Graev de $\Gb^F$ 
(\resp du dual d'Alvis-Curtis d'un \car de
Gelfand-Graev de $\Gb^F$). Nous allons dans cette section param\'etrer 
les \cars r\'eguliers (et semi-simples) de $\Gb^F$ 
appartenant \`a $\EC(\Gb^F,[s])$ ou $\EC(\Gb^F,(s))$. Nous allons aussi 
\'etablir les premi\`eres propri\'et\'es (action de $\Gbt^F$, 
restriction de Harish-Chandra...). 

\medskip

Posons
\equat
\begin{array}{rcl}
\r_s=\r_s^\Gb&=& \DS{\frac{\e_\Gb\e_{C_{\Gb^*}^\ci(s)}}{|W^\ci(s)|} 
\sum_{w \in W^\ci(s)} R_{\Tb_w^*}^\Gb(s)}, \\
&&\\
\ch_s=\ch_s^\Gb&=& \DS{\frac{\e_\Gb\e_{C_{\Gb^*}^\ci(s)}}{|W^\ci(s)|} 
\sum_{w \in W^\ci(s)} 
\e(w) R_{\Tb_w^*}^\Gb(s)}. \\
\end{array}
\endequat
Remarquons que, d'apr\`es \ref{dg rtgs},  
\equat
\ch_s=\e_\Gb \e_{C_{\Gb^*}^\circ(s)} D_\Gb(\r_s).
\endequat

Alors, d'apr\`es le corollaire \ref{rtgs}, on a 
\equat\label{res chis}
\begin{array}{rcl}
\r_s& =& \Res_{\Gb^F}^{\Gbt^F} \r_\sti, \\
&&\\
\ch_s& =& \Res_{\Gb^F}^{\Gbt^F} \ch_\sti. \\
\end{array}
\endequat
D'apr\`es la proposition \ref{tenseur}, on a 
\equat\label{tenseur chis}
\begin{array}{rcl}
\r_\sti \otimes \zha &=& \r_{\sti z}, \\
&&\\
\ch_\sti \otimes \zha &=& \ch_{\sti z},
\end{array}
\endequat
pour tout $z \in (\Ker i^*)^{F^*}$.

\bigskip

\Theoreme{Deligne-Lusztig}{gel connexe}
{\it \begin{itemize}
\itemth{a} $\langle \ch_\sti,\G^\Gbt \rangle_{\Gbt^F} = 1$.

\itemth{b} $\r_\sti$ et $\ch_\sti$ sont des \cars \irrs de $\Gbt^F$ 
et ils appartiennent \`a $\EC(\Gbt^F,[\sti])$.

\itemth{c} Le \car de  Gelfand-Graev de $\Gbt^F$ a la d\'ecomposition suivante~:
$$\G^\Gbt = \sum_{[\sti]} \ch_\sti.$$
\end{itemize}}

\bigskip

\noindent{\sc Remarque - } Le th\'eor\`eme pr\'ec\'edent est d\'emontr\'e 
dans \cite[th\'eor\`eme 10.7]{delu} lorsque $F$ est un endomorphisme 
de Frobenius. Dans le contexte l\'eg\`erement plus g\'en\'eral dans 
lequel nous nous pla\c{c}ons, le r\'esultat reste valide. En effet, 
d'apr\`es \cite[page 161]{delu}, il faut seulement utiliser 
le fait que la formule de Mackey a lieu lorsque l'un des deux sous-groupes 
de Levi est un tore (voir th\'eor\`eme \ref{theo mackey} (a2)).\finl

\bigskip

\Corollaire{Asai}{gel non connexe}
{\it \begin{itemize}
\itemth{a} $\langle \ch_s , \G^\Gb \rangle_{\Gb^F}=1$.

\itemth{b} $\r_s$ et $\ch_s$ sont des \cars de $\Gb^F$ sans multiplicit\'e 
et toute composante \irr de $\r_s$ ou $\ch_s$ appartient \`a $\EC(\Gb^F,[s])$.

\itemth{c} Si $\ch_{s,1}$ est l'unique composante \irr commune \`a 
$\G^\Gb$ et $\ch_s$ $($voir $($a$))$, alors 
$$\G^\Gb=\sum_{[s]} \ch_{s,1}.$$
\end{itemize}}

\bigskip

\proof (a) r\'esulte du th\'eor\`eme \ref{gel connexe} (a), de \ref{res chis} 
et de la r\'eciprocit\'e de Frobenius. (b) r\'esulte du 
th\'eor\`eme \ref{gel connexe} (b), de (a), de \ref{res chis} 
et de la proposition \ref{restriction series} (a).
(c) d\'ecoule de (b) et du th\'eor\`eme \ref{gel connexe} (c).\fin

\bigskip

On pose $\r_{s,1}=\e_\Gb\e_{C_{\Gb^*}^\ci(s)}  D_\Gb(\ch_{s,1}) 
\in \EC(\Gb^F,[s])$. C'est une composante \irr de $\r_s$. S'il y a ambigu\"\i t\'e, 
on notera $\r_{s,1}^\Gb$ et $\ch_{s,1}^\Gb$ les \cars \irrs 
$\r_{s,1}$ et $\ch_{s,1}$ de $\Gb^F$ respectivement.

\bigskip

\Corollaire{Digne-Lehrer-Michel}{rlg res chis}
{\it Soit $\Lb$ un \levic $F$-stable d'un \para $F$-stable $\Pb$ de $\Gb$. 
Notons $\Lb^*$ un \levic $F^*$-stable d'un \para $F^*$-stable de 
$\Gb^*$ dual de $\Lb$. Alors 
$$\lexp{*}{R}_{\Lb \incl \Pb}^\Gb \r_{s,1}^\Gb = \sum_{[t]_{\Lb^{*F^*}} \incl 
[s]_{\Gb^{*F^*}}} \r_{t,1}^\Lb$$
$$\lexp{*}{R}_{\Lb \incl \Pb}^\Gb \ch_{s,1}^\Gb = \sum_{[t]_{\Lb^{*F^*}} \incl 
[s]_{\Gb^{*F^*}}} \ch_{t,1}^\Lb.\leqno{\mathit{et}}$$}

\bigskip

\proof La deuxi\`eme \'egalit\'e r\'esulte du th\'eor\`eme \ref{restriction gelfand}, 
du corollaire \ref{gel non connexe} (c), et du corollaire  
\ref{rlg res series}. La premi\`ere d\'ecoule de la seconde et de la relation 
de commutation entre l'induction de Harish-Chandra et la dualit\'e 
d'Alvis-Curtis \cite[th\'eor\`eme 8.11]{dmbook}.\fin

\bigskip

Si $\x \in (A_{\Gb^*}(s)^{F^*})^\we$, on pose 
$$\r_{s,\x}=\t_z^\Gb \r_{s,1}$$
$$\ch_{s,\x}=\t_z^\Gb\ch_{s,1},\leqno{\mathrm{et}}$$
o\`u $z \in H^1(F,\ZC(\Gb))$ est tel que $\omeh_s^0(z)=\x$ 
(d'apr\`es le corollaire \ref{stabilisateur}, les caract\`eres 
$\chi_{s,\x}$ et $\r_{s,\x}$ ne d\'ependent que de $\x$ 
et non du choix de $z$).
La proposition suivante d\'ecrit les caract\`eres r\'eguliers ou semi-simples 
appartenant \`a $\EC(\Gb^F,[s])$. 

\bigskip

\Proposition{Asai}{decomposition chis}
{\it Le stabilisateur de $\r_{s,1}$ (ou $\ch_{s,1}$) dans $\Gbt^F$ 
est \'egal \`a $\Gbt^F(s)$. Par cons\'equent, 
$$\r_s=\sum_{\x \in (A_{\Gb^*}(s)^{F^*})^\we} \r_{s,\x}$$
$$\ch_s=\sum_{\x \in (A_{\Gb^*}(s)^{F^*})^\we} \ch_{s,\x}.\leqno{\mathit{et}}$$}

\bigskip

\proof Seule la premi\`ere assertion n\'ecessite une preuve, 
la deuxi\`eme r\'esultant imm\'ediatement de la premi\`ere et du corollaire 
\ref{gel non connexe} (b). Mais, par la th\'eorie de Clifford, 
elle d\'ecoule de la formule \ref{tenseur rtgs}.\fin

\bigskip

\corollaire{inverse}
{\it \begin{itemize} 
\itemth{a} Si $z \in H^1(F,\Zb(\Gb))$ et si $\x \in (A_{\Gb^*}(s)^{F^*})^\we$, alors 
$$\t^\Gb_z \r_{s,\x}=\r_{s,\x\omeh^0_s(z)}$$
$$\t^\Gb_z \ch_{s,\x}=\ch_{s,\x\omeh^0_s(z)}.\leqno{\mathit{et}}$$

\itemth{b} $\ch_{s,\x}$ est une composante \irr de $\G_z^\Gb$ \ssi 
$\x=\omeh^0_s(z)$.
\end{itemize}}

\bigskip

S'il y a ambigu\"\i t\'e, nous noterons $\r_{s,\x}^\Gb$ et $\ch_{s,\x}^\Gb$ 
les \cars \irrs $\r_{s,\x}$ et $\ch_{s,\x}$ de $\Gb^F$ respectivement 
($\x \in (A_{\Gb^*}(s)^{F^*})^\we$).

\bigskip

\Corollaire{Digne-Lehrer-Michel}{rlg res chis 2}
{\it Soit $\Lb$ un \levic $F$-stable d'un \para $F$-stable $\Pb$ de $\Gb$. 
Soit $\Lb^*$ un \levic 
$F^*$-stable d'un \para $F^*$-stable de $\Gb^*$ dual de $\Lb$. 
Soit $\x \in (A_{\Gb^*}(s)^{F^*})^\we$. Pour tout 
$t \in \Lb^{*F^*}$ tel qu'il existe 
$g \in \Gb^{*F^*}$ v\'erifiant $\lexp{g}{s}=t$, on pose 
$\x_t=\Res_{A_{\Lb^*}(t)^{F^*}}^{A_{\Gb^*}(t)^{F^*}} \lexp{g}{\x}$~; 
le \car lin\'eaire $\x_t$ ne d\'epend pas du choix de $g$. Alors 
$$\lexp{*}{R}_{\Lb \incl \Pb}^\Gb \r_{s,\x}^\Gb = \sum_{[t]_{\Lb^{*F^*}} \incl 
[s]_{\Gb^{*F^*}}} \r_{t,\x_t}^\Lb$$
$$\lexp{*}{R}_{\Lb \incl \Pb}^\Gb \ch_{s,\x}^\Gb = \sum_{[t]_{\Lb^{*F^*}} \incl 
[s]_{\Gb^{*F^*}}} \ch_{t,\x_t}^\Lb.\leqno{\mathit{et}}$$}

\bigskip

\proof Cela r\'esulte du corollaire \ref{rlg res chis}, 
de la commutativit\'e du diagramme \ref{commutativite als} et de 
\ref{tauzg resrlg}.\fin

\bigskip

\proposition{etonnant non ?}
{\it On a $\r_s=\chi_s$ \ssi $C_{\Gb^*}^\circ(s)$ est un \tor de $\Gb^*$. 
Dans ce cas, il existe un \car lin\'eaire $\x_s$ d'ordre $2$ de 
$A_{\Gb^*}(s)^{F^*}$ tel que, pour tout $\x \in (A_{\Gb^*}(s)^{F^*})^\we$, on ait 
$$D_\Gb \r_{s,\x}=\e_\Gb\e_{C_{\Gb^*}^\ci(s)} \r_{s,\x\x_s}.$$}

\bigskip

\proof La premi\`ere assertion est imm\'ediate. Supposons donc 
que $C_{\Gb^*}^\circ(s)$ est un \tor de $\Gb^*$. Notons $\xi_s$ le 
\car lin\'eaire de $A_{\Gb^*}(s)^{F^*}$ tel que
$$D_\Gb \r_{s,1}=\e_\Gb\e_{C_{\Gb^*}^\ci(s)} \r_{s,\x_s}.$$
Puisque $\t_z^\Gb \circ D_\Gb = D_\Gb \circ \t_z^\Gb$ 
pour tout $z \in H^1(F,\ZC(\Gb))$, on en d\'eduit 
la formule donn\'ee dans la proposition \ref{etonnant non ?}. 
Pour finir, puisque $D_\Gb$ est une involution, on a $\xi_s^2=1$.\fin

\bigskip

\exemple{dg sln} 
Le caract\`ere lin\'eaire $\xi_s$ de la proposition \ref{etonnant non ?} 
peut \^etre non trivial. En effet, supposons ici que $\Gb=\Sb\Lb_2(\FM)$, 
que $F$ est l'endomorphisme de Frobenius d\'eploy\'e standard sur $\FM_q$ 
et que $p$ (ou $q$) est impair. 
Notons $s$ l'unique \'el\'ement de $\Tb_0^*$ d'ordre $2$ et notons 
$\th$ l'unique caract\`ere lin\'eaire de $\Tb_0^F$ d'ordre $2$. 
Alors $\r_s=R_{\Tb_0}^\Gb(\th)$ et $\lexp{*}{R}_{\Tb_0}^\Gb(\r_{s,1})=\th$ 
d'apr\`es le corollaire \ref{rlg res chis}. D'autre part, 
$D_\Gb= (R_{\Tb_0}^\Gb \circ \lexp{*}{R}_{\Tb_0}^\Gb)-\Id_{\ZM\Irr \Gb^F}$. 
Donc $D_\Gb(\r_{s,1})=\r_{s,\xi}$, o\`u $\x$ est l'unique 
caract\`ere non trivial de $A_{\Gb^*}(s)^{F^*} \simeq \ZM/2\ZM$.\finl

\bigskip

\section{Caract\`eres semi-simples ou r\'eguliers 
cuspidaux\label{section semi cuspidal}}~

\medskip

Nous allons ici \'etudier les s\'eries de Harish-Chandra associ\'ees \`a un 
\car cuspidal semi-simple (ou r\'egulier). 

\bigskip

\soussection{Caract\'erisation de la cuspidalit\'e}
Soit $\x \in (A_{\Gb^*}(s)^{F^*})^\we$. Si $\r_{s,\x}$ est cuspidal, alors 
$D_\Gb \r_{s,\x}=\eta_\Gb \r_{s,\x}$ et donc $D_\Gb \r_s = \eta_\Gb \r_s$. 
En particulier, $\r_s=\chi_s$. 
Donc un \car \irr cuspidal est semi-simple \ssi il est r\'egulier. 
Dor\'enavant, nous fixons un couple $(\Tbt_1,\thet_1) \in \nabla(\Gbt,F)$ 
tel que $(\Tbt_1,\thet_1) \doublefleche{\Gbt} (\Tbt_1^*,\sti)$ 
et nous posons $(\Tb_1,\th_1)=\RES_\Gb^\Gbt (\Tbt_1,\thet_1)$. 
Le lemme suivant pr\'ecise quand est-ce qu'un \car semi-simple 
est cuspidal.

\bigskip

\lemme{condition}
{\it Soit $\x \in (A_{\Gb^*}(s)^{F^*})^\we$. Le \car \irr $\r_{s,\x}$ de 
$\Gb^F$ est cuspidal \ssi $C_{\Gb^*}^\ci(s)$ est un \tor $F^*$-stable 
de $\Gb^*$ qui n'est contenu dans aucun \levi $F^*$-stable $\Gb^*$-d\'eploy\'e.}

\bigskip

\proof Remarquons que, d'apr\`es le corollaire \ref{restriction cuspidal}, 
$\r_{s,\x}$ est cuspidal \ssi le \car \irr $\r_\sti$ de $\Gbt^F$ 
l'est. Il est donc suffisant de montrer le th\'eor\`eme pour $\Gbt$.

\medskip

Si $\r_\sti$ est cuspidal, alors $D_\Gbt \r_\sti = \eta_\Gbt \r_\sti$ donc 
$\r_\sti=\ch_\sti$.  
Cela montre que $C_{\Gbt^*}(\sti)$ est un tore maximal de $\Gbt^*$ 
donc que $C_{\Gbt^*}(\sti)=\Tbt_1^*$. Dans ce cas, 
$$\r_\sti=\e_\Gbt\e_{\Tbt_1^*} R_{\Tbt_1^*}^\Gbt(\sti)$$
donc $C_{\Gbt^*}(\sti)$ n'est contenu dans aucun \levic $F^*$-stable d'un 
\para $F^*$-stable propre de $\Gbt^*$. 

\medskip

R\'eciproquement, supposons que $C_{\Gbt^*}(\sti)$ soit un \tor $F^*$-stable  
de $\Gbt^*$ qui n'est contenu dans aucun \levi $F^*$-stable $\Gb^*$-d\'eploy\'e. 
Alors $\r_\sti=\e_\Gbt\e_{\Tbt_1} R_{\Tbt_1}^\Gbt(\thet_1)$. 
Soit $\Lbt$ un \levic 
$F$-stable d'un \para $F$-stable propre $\Pbt$ de $\Gbt$ et soit 
$\Lbt^*$ un \levic $F^*$-stable 
d'un \para $F^*$-stable propre de $\Gbt^*$ dual de $\Lbt$. Alors, d'apr\`es 
la formule de Mackey (voir th\'eor\`eme \ref{theo mackey} (a2)), on a 
$$\lexp{*}{R}_{\Lbt \incl \Pbt}^\Gbt R_{\Tbt_1}^\Gbt(\thet_1)=0,$$
ce qui montre la cuspidalit\'e de $\rho_\sti$.\fin

\bigskip

\soussection{Groupe d'inertie} 
Notons $\Lbt_s$ (\resp $\Lbt_s^*$) le \levi $F$-stable (\resp $F^*$-stable) 
$\Gbt$-d\'eploy\'e (\resp $\Gbt^*$-d\'eploy\'e) contenant $\Tbt_1$ (\resp 
$\Tbt_1^*$) et minimal pour ces propri\'et\'es (voir remarque \ref{plus petit levi}). 
Alors, $\Lbt_s$ et $\Lbt_s^*$ sont duaux. On pose $\Lb_s = \Lbt_s \cap \Gb$ 
et $\Lb_s^*=i^*(\Lbt_s^*)$. 
Soit $\Pbt_s$ un \para $F$-stable de $\Gbt$ dont $\Lbt_s$ soit un compl\'ement de 
Levi. On pose $\Pb_s=\Pbt_s \cap \Gb$. 
En appliquant le corollaire \ref{coro borel tits} au groupe $C_{\Gbt^*}(\sti)$, 
on obtient que $C_{\Lbt_s^*}(\sti)=\Tbt_1^*$. Donc, d'apr\`es le lemme 
\ref{condition}, $\r_\sti^{\Lbt_s}$ est un \car \irr cuspidal de $\Lbt_s^F$. 
Donc, d'apr\`es le corollaire 
\ref{restriction cuspidal}, les \cars \irrs $\r_{s,\x}^{\Lb_s}$ de $\Lb_s^F$ 
sont cuspidaux pour tout $\x \in (A_{\Lb_s^*}(s)^{F^*})^\we$.
Le groupe $W^{\phi_1}$ stabilise 
$\Ker(F^*-q^{1/\d}, Y(\Tb_1^*) \otimes \QM(q^{1/\d}))$ 
donc il normalise $\Lb^*$. On a en fait le r\'esultat plus pr\'ecis suivant~:

\bigskip

\proposition{W}
{\it Le groupe $W_{\Gb^F}^\pr(\Lbt_s,\r_\sti^{\Lbt_s})$ est canoniquement 
isomorphe \`a $W(s)^{F^*}$.}

\bigskip

\proof 
On a 
$$\r_\sti^{\Lbt_s}=\pm R_{\Tbt_1}^\Gbt(\thet_1).$$
Le groupe $W(s)^{F^*}$ est isomorphe \`a $W_{\Gb^F}(\Tb_1,\th_1)$ et le groupe 
$W(\sti)^{F^*}$ est isomorphe \`a $W_{\Gb^F}(\Tbt_1,\thet_1)$. De plus, puisque 
$C_{\Lbt_s^*}(\sti)=\Tbt_1^*$, on a $W_{\Lb_s^F}(\Tbt_1,\thet_1)=\{1\}$. 

Le groupe $W_{\Gb^F}(\Tb_1,\th_1)$ normalise $\Lb_s$. Soit 
$w \in W_{\Gb^F}(\Tb_1,\th_1)$. Notons $\t_w$ le \car lin\'eaire 
$\lexp{w}{\thet_1}. \thet_1^{-1}$ de $\Tbt_1^F/\Tb_1^F \simeq \Lbt^F/\Lb^F$. Alors 
$$\lexp{w}{\r_\sti^{\Lbt_s}}=\r_\sti^{\Lbt_s} \otimes \t_w.$$
Donc, si on note $\bar{w}$ l'image de $w$ dans $W_{\Gb^F}(\Lbt)$, alors 
l'application 
$$\fonction{\a}{W_{\Gb^F}(\Tb_1,\th_1)}{
W_{\Gb^F}^\pr(\Lbt_s,\r_\sti^{\Lbt_s})}{w}{(\bar{w},\t_w)}$$
est un \mor de groupes bien d\'efini.

Montrons d'abord que $\a$ est injective. Soit $w \in W_{\Gb^F}(\Tb_1,\th_1)$ 
tel que $\a(w)=1$. Alors $\t_w=1$ donc $w \in W_{\Gb^F}(\Tbt_1,\thet_1)$. Mais 
$\bar{w}=1$ donc $w \in W_{\Lb_s^F}(\Tbt_1,\thet_1)$ \cad $w=1$.

Il reste \`a montrer que $\a$ est surjective. Soit  
$(w,\t) \in W_{\Gb^F}^\pr(\Lbt_s,\r_\sti^{\Lbt_s})$. Notons $\dot{w}$ 
un repr\'esentant 
de $w$ dans $N_{\Gb^F}(\Lbt_s)$. On a 
$$R_{\lexp{\dot{w}}{\Tbt_1}}^{\Lbt_s}(\lexp{\dot{w}}{\thet_1})=
R_{\Tbt_1}^{\Lbt_s}(\thet_1 \otimes \t)$$
donc il r\'esulte de la formule de Mackey (th\'eor\`eme \ref{theo mackey} (2)) 
qu'il existe $l \in \Lb^F$ tel que 
$(\lexp{\dot{w}}{\Tbt_1},\lexp{\dot{w}}{\thet_1})=
(\lexp{l}{\Tbt_1},\lexp{l}{(\thet_1} \otimes \t))$. 
Soit $\dot{w}_+=l^{-1}\dot{w}$. Alors $\dot{w}_+ \in N_{\Gb^F}(\Tb_1,\th_1)$ et, 
si on note $w_+$ son image dans $W_{\Gb^F}(\Tb_1,\th_1)$, alors 
$\a(w_+)=(w,\t)$.\fin

\bigskip

Soit $W_{\Lb_s}(s)$ le groupe de Weyl de $C_{\Lb_s^*}(s)$ relativement \`a 
$\Tb_1^*$. On a $W_{\Lb_s}(s)=A_{\Lb_s^*}(s)$ car 
$C_{\Lb_s^*}^\ci(s)=i^*(C_{\Lbt_s^*}(\sti))=\Tb_1^*$. Donc $A_{\Lb_s^*}(s)$ est un 
sous-groupe $F^*$-stable de $W(s)$.

\bigskip

\proposition{quelques remarques}
{\it \begin{itemize} 
\itemth{a} $A_{\Lb_s^*}(s)^{F^*}$ est contenu dans $A_{\Gb^*}(s)^{F^*}$.

\itemth{b} $A_{\Lb_s^*}(s)^{F^*}$ est central dans $W(s)^{F^*}$.

\itemth{c} $\Gbt^F(s)=\Gb^F.\Lbt_s^F(\Gb,\r_{s,1}^{\Lb_s})$.

\itemth{d} $\r_\sti^\Gbt \in \EC(\Gbt^F,\Lbt_s,\r_\sti^{\Lbt_s})$.
\end{itemize}}

\bigskip

\proof (a) d\'ecoule de la proposition \ref{lll} (b). (b) r\'esulte du fait 
que $A_{\Gb^*}(s)$ est ab\'elien, de (a) et de la proposition \ref{lll} (c). 
(c) r\'esulte de la proposition \ref{W}. (d) d\'ecoule de l'\'egalit\'e 
$$\langle R_{\Lbt_s \incl \Pbt_s}^\Gbt \r_\sti^{\Lbt_s} , 
\r_\sti^\Gbt \rangle_{\Gb^F} = 1,$$
qui a \'et\'e montr\'ee dans le corollaire \ref{rlg res chis}.\fin

\bigskip

\remarque{als k} 
Le sous-groupe $A_{\Lb_s^*}(s)^{F^*}$ de 
$W(s)^{F^*} \simeq W_{\Gb^F}^\pr(\Lbt,\r_\sti^\Lbt)$ 
est isomorphe au sous-groupe central $(\Lbt^F/\Lbt^F(\r_\sti))^\wedge$ 
d\'efini dans le \SEC\ref{sub W'}.\finl

\bigskip

\soussection{La s\'erie ${\boldsymbol{\EC(\Gb^F,\Lb_s,\r_s^{\Lb_s})}}$}
D'apr\`es le th\'eor\`eme \ref{harish} (a) et d'apr\`es la proposition \ref{W} 
on obtient des bijections 
\equat\label{bi tilde}
\fonctio{\Irr W(\sti)^{F^*}}{\EC(\Gbt^F,\Lbt_s,\r_\sti^{\Lbt_s})}{\ch}{\Rti_\ch[\sti]}
\endequat
et
\equat\label{bi G}
\fonctio{\Irr W(s)^{F^*}}{\EC(\Gb^F,\Lb_s,\r_s^{\Lb_s})}{\eta}{R_\eta[s]}.
\endequat
D'apr\`es \cite[chapitre 8]{lubook}, la bijection 
\ref{bi tilde} est bien d\'efinie une fois fix\'ee la convention suivante~: 
\equat\label{convention}
\Rti_1[\sti]=\r_\sti^\Gbt
\endequat
et, par  la remarque \ref{remarques harish} (b), la bijection 
\ref{bi G} est bien d\'efinie une fois fix\'ee la convention suivante~: 
\equat
R_1[s]=\r_{s,1}^\Gb.
\endequat
S'il y a ambigu\"\i t\'e, nous noterons $R_\eta^\Gb[s]$ le \car \irr 
$R_\eta[s]$ de $\Gb^F$ $(\eta \in \Irr W(s)^{F^*}$) et par 
$R^\Gbt_\ch[\sti]$ le \car \irr $\Rti_\ch[\sti]$ de $\Gbt^F$ 
($\ch \in \Irr W(\sti)^{F^*}$).

\bigskip

\remarque{inclusions series} D'apr\`es le th\'eor\`eme \ref{rlg series}, on a 
$$\EC(\Gbt^F,\Lbt_s,\r_\sti^{\Lbt_s}) 
\incl \EC(\Gbt^F,[\sti])$$
$$\EC(\Gb^F,\Lb_s,\r_s^{\Lb_s}) 
\incl \EC(\Gb^F,[s]).~\SS{\square}\leqno{\mathrm{et}}$$

\bigskip

Gr\^ace aux th\'eor\`emes \ref{harish} et \ref{induit harish}, 
on obtient~:

\bigskip

\theoreme{parametrage}
{\it \begin{itemize}
\itemth{a} 
Si $\eta$ et $\ch$ sont des \cars \irrs de $W(s)^{F^*}$ et 
$W^\ci(s)^{F^*}=W(\sti)^{F^*}$ 
respectivement, alors 
$$\langle R_\eta[s],\Res_{\Gb^F}^{\Gbt^F} \Rti_\ch[\sti] \rangle_{\Gb^F} =
\langle \Res_{W(\sti)^{F^*}}^{W(s)^{F^*}} \eta,\ch \rangle_{W(\sti)^{F^*}}.$$

\itemth{b} Si $\eta \in \Irr W(s)^{F^*}$ et $\x \in (A_{\Gb^*}(s)^{F^*})^\we$, alors 
$$\lexp{g_\x}{R_\eta[s]}=R_{\eta \otimes \x}[s].$$

\itemth{c} Soit $\Lb$ un \levic $F$-stable d'un \para $F$-stable $\Pb$ de $\Gb$ 
contenant $\Lb_s$. Fixons un \levi $F^*$-stable $\Lb^*$ d'un \para $F^*$-stable de 
$\Gb^*$ contenant $\Lb_s^*$ et tel que la $\Lb^F$-classe de conjugaison de $\Lb_s$ 
soit associ\'ee \`a la $\Lb^{*F^*}$-classe de conjugaison de 
$\Lb_s^*$. Notons $W_\Lb(s)$ le groupe de Weyl de $C_{\Lb^*}(s)$ relativement 
\`a $\Tb_1^*$. Alors 
$$\langle R_{\Lb \incl \Pb}^\Gb R_\eta^\Lb[s],R_\z^\Gb[s] \rangle_{\Gb^F} =
\langle \Ind_{W_\Lb(s)^{F^*}}^{W(s)^{F^*}} \eta,\z \rangle_{W(s)^{F^*}}$$
pour tous \cars \irrs $\eta$ et $\z$ de $W_\Lb(s)^{F^*}$ et $W(s)^{F^*}$ 
respectivement.
\end{itemize}}

\bigskip

Remarquons que l'assertion (c) du th\'eor\`eme pr\'ec\'edent 
\ref{parametrage} utilise le corollaire \ref{rlg res chis}.
 
\bigskip

\section{Caract\`eres semi-simples et fonctions absolument cuspidales}

\medskip

\soussection{Un exemple de fonction absolument cuspidale} 
Si $a \in A_{\Gb^*}(s)^{F^*}$, on pose
$$\rhodot_{s,a}=\rhodot_{s,a}^\Gb= \e_\Gb\e_{C_{\Gb^*}^\circ(s)} 
\sum_{\xi \in (A_{\Gb^*}(s)^{F^*})^\we} \xi(a)^{-1} \rho_{s,\xi}.$$
Il est facile de retrouver les caract\`eres irr\'eductibles $\r_{s,\xi}$ 
comme combinaisons lin\'eaires des $\rhodot_{s,a}$. En effet, si 
$\xi \in (A_{\Gb^*}(s)^{F^*})^\we$, on a
\equat\label{changement de base semi}
\r_{s,\xi}=\frac{\e_\Gb\e_{C_{\Gb^*}^\circ(s)}}{|A_{\Gb^*}(s)^{F^*}|} 
\sum_{a \in A_{\Gb^*}(s)^{F^*}} \xi(a) \rhodot_{s,a}.
\endequat
Par ailleurs, il r\'esulte du corollaire \ref{inverse} que
\equat\label{action h1 rho}
\rhodot_{s,a} \in \Cent(\Gb^F,[s],a).
\endequat
Si $a$ et $b$ sont deux \'el\'ements de $A_{\Gb^*}(s)^{F^*}$, un calcul 
\'el\'ementaire montre que 
\equat\label{scalaire rhodot}
\langle \rhodot_{s,a},\rhodot_{s,b} \rangle_{\Gb^F} = 
\begin{cases}
|A_{\Gb^*}(s)^{F^*}| & \text{si } a=b,\\
0 & \text{sinon.}
\end{cases}
\endequat
D'apr\`es l'exemple \ref{exemple res}, on a~:

\bigskip

\proposition{rhodot absolument cuspidal}
{\it Si $a \in A_{\Gb^*}(s)^{F^*}$ est tel que $\o_s(a) \in \ZC_\cus^\we(\Gb)$, 
alors $\rhodot_{s,a}$ est une fonction absolument cuspidale.}

\bigskip

\soussection{Restriction de Lusztig} Nous travaillerons sous 
l'hypoth\`ese suivante~:

\medskip

\begin{quotation}
\noindent{\bf Hypoth\`ese~:} 
{\it Nous supposerons jusqu'\`a la fin de ce chapitre que $p$ est bon pour $\Gb$.}
\end{quotation}

\medskip

Soit $\Lb$ un \levi $F$-stable de $\Gb$ et soit $\Lb^*$ un \levi 
$F^*$-stable de $\Gb^*$ dual de $\Lb$. 
L'hypoth\`ese entra\^\i ne que l'application $\res_\Lb^\Gb$ entre 
ensemble de classes unipotentes r\'eguli\`eres est bien d\'efinie (voir 
\SEC\ref{sous reslg}). En particulier, le caract\`ere de Gelfand-Graev 
$\G^\Lb$ est lui aussi bien d\'efini. Nous allons ici donner une formule 
pour la restriction de Lusztig des caract\`eres $\r_{s,\xi}^\Gb$. 

\bigskip

\proposition{theo restriction chis}
{\it Supposons que la formule de Mackey et la conjecture $(\GG)$ ont lieu 
dans $\Gb$. Soit $\x \in (A_{\Gb^*}(s)^{F^*})^\wedge$. 
Pour tout $t \in \Lb^{*F^*}$ tel qu'il existe 
$g \in \Gb^{*F^*}$ v\'erifiant $\lexp{g}{s}=t$, on pose 
$\x_t=\Res_{A_{\Lb^*}(t)^{F^*}}^{A_{\Gb^*}(t)^{F^*}} \lexp{g}{\x}$~; 
le \car lin\'eaire $\x_t$ ne d\'epend pas du choix de $g$. Alors 
$$\lexp{*}{R}_\Lb^\Gb \r_{s,\x}^\Gb = 
\e_\Gb\e_\Lb\e_{C_{\Gb^*}^\circ(s)} \sum_{[t]_{\Lb^{*F^*}} \incl 
[s]_{\Gb^{*F^*}}} \e_{C_{\Lb^*}^\circ(t)} \r_{t,\x_t}^\Lb$$
$$\lexp{*}{R}_\Lb^\Gb \ch_{s,\x}^\Gb = \e_\Gb\e_\Lb 
\sum_{[t]_{\Lb^{*F^*}} \incl [s]_{\Gb^{*F^*}}} \ch_{t,\x_t}^\Lb.
\leqno{\mathit{et}}$$}

\bigskip

\remarque{rem generalisation chis} 
La proposition \ref{theo restriction chis} g\'en\'eralise 
le corollaire \ref{rlg res chis 2} tout comme le 
th\'eor\`eme \ref{theo restriction gamma} g\'en\'eralisait 
le th\'eor\`eme \ref{restriction gelfand}.\finl

\bigskip

\corollaire{coro restriction rhosa}
{\it Supposons que la formule de Mackey et la conjecture $(\GG)$ ont lieu 
dans $\Gb$. Soit $a \in A_{\Gb^*}(s)^{F^*}$. 
Pour tout $t \in \Lb^{*F^*}$ tel qu'il existe 
$g \in \Gb^{*F^*}$ v\'erifiant $\lexp{g}{s}=t$, on note 
$a_t$ l'\'el\'ement $gag^{-1}$ de $A_{\Gb^*}(t)^{F^*}$~; 
l'\'el\'ement $a_t$ ne d\'epend pas du choix de $g$. Alors 
$$\lexp{*}{R_\Lb^\Gb} \rhodot_{s,a}^\Gb = \sum_{\SS{[t]_{\Lb^{*F^*}} \incl 
[s]_{\Gb^{*F^*}}} \atop \SS{a_t \in A_{\Lb^*}(t)^{F^*}}}
\frac{|A_{\Gb^*}(s)^{F^*}|}{|A_{\Lb^*}(t)^{F^*}|} ~\rhodot_{t,a_t}^\Lb.$$}

\bigskip

\soussection{Combinaisons lin\'eaires d'induits de caract\`eres semi-simples}
Soit $w \in W(s)$. On fixe un \'el\'ement $g_w \in \Gb^*$ tel 
que $g_w^{-1} F(g_w)$ normalise $\Tb_1^*$ et repr\'esente $w$. 
On pose alors $\Tb_w^* = \lexp{g_w}{\Tb_1^*}$ et $s_w=g_w s g_w^{-1}$. 
D'apr\`es le th\'eor\`eme de Lang, on peut choisir $g_w$ de sorte 
que $s_w = s_\a$, o\`u $\a$ d\'esigne la classe de $w$ dans 
$H^1(F^*,A_{\Gb^*}(s))$. C'est ce que nous ferons dans la suite. 
Il est \`a noter que le couple $(\Tb_w^*,s_w)$ est bien 
d\'efini \`a $\Gb^{*F^*}$-conjugaison pr\`es par $w$ 
(et m\^eme par la classe de $w$ dans $H^1(F^*,W(s))$~: en effet, le stabilisateur 
du couple $(\Tb_1^*,s)$ dans $\Gb^*$ est \'egal \`a l'image 
inverse de $W(s)$ dans $N_{\Gb^*}(\Tb_1^*)$). 

Fixons maintenant $a \in A_{\Gb^*}(s)^{F^*}$. Alors le 
\levi $F^*$-stable $\Lb_{s,a}^*$ a \'et\'e d\'efini 
dans \SEC\ref{sous semi cus}. Si $w \in W(s)^a$, 
on pose $\Lb_{s,a,w}^*=\lexp{g_w}{\Lb_{s,a}^*}$. 
Alors le couple $(\Lb_{s,a,w}^*,s_w)$ est bien d\'efini \`a 
$\Gb^{*F^*}$-conjugaison pr\`es par $w$ (et m\^eme par la classe 
de $w$ dans $H^1(F^*,W(s)^a)$~: en effet, le stabilisateur 
du couple $(\Lb_{s,a}^*,s)$ dans $\Gb^*$ est \'egal \`a l'image 
inverse de $W(s)^a$ dans $N_{\Gb^*}(\Tb_1^*)$ d'apr\`es le corollaire 
\ref{ah ah ah} (e)). 
De plus, puisque $a \in A_{\Lb_{s,a}^*}(s)^{F^*}$ par construction, 
on en d\'eduit que $a \in A_{\Lb_{s,a,w}^*}(s)^{F^*}$ pour tout 
$w \in W(s)^a$ (\`a travers le morphisme injectif naturel 
$A_{\Lb_{s,a,w}^*}(s) \injto A_{\Gb^*}(s)$). 
Notons que
\equat\label{centralisateur lsa}
C_{\Lb_{s,a,w}^*}^\circ(s_w)=\Tb_w^*
\endequat
(voir corollaire \ref{ah ah ah} (a)). 
Notons $\Lb_{s,a,w}$ un \levi $F$-stable de $\Gb$ dual de 
$\Lb_{s,a,w}^*$. Alors le couple $(\Lb_{s,a,w},\rhodot_{s_w,a}^{\Lb_{s,a,w}})$ 
est bien d\'efini \`a $\Gb^F$-conjugaison pr\`es par $w$ (et m\^eme 
par la classe de $w$ dans $H^1(F^*,W(s)^a)$. Donc la fonction 
$R_{\Lb_{s,a,w}}^\Gb \rhodot_{s_w,a}^{\Lb_{s,a,w}}$ 
est bien d\'efinie~: nous la noterons $\RC_{s,a,w}$. Elle appartient 
\`a $\qlb \EC(\Gb^F,(s))$. Si $\a$ d\'esigne la classe de $w$ dans 
$H^1(F^*,A_{\Gb^*}(s))$, alors, d'apr\`es le th\'eor\`eme \ref{rlg series} et 
\ref{action h1 rho}, on a 
\equat\label{serie RC}
\RC_{s,a,w} \in \Cent(\Gb^F,[s_\a],a).
\endequat
Si $f \in \Cent(W(s)^a \phi_1)$, on pose~:
$$\RC(s,a)_f=\frac{1}{|W(s)^a|} \sum_{w \in W(s)^a} f(w\phi_1) 
\RC_{s,a,w}.$$
D'apr\`es \ref{serie RC}, cela nous d\'efinit une application lin\'eaire 
$\RC(s,a) : \Cent(W(s)^a \phi_1) \to \Cent(\Gb^F,(s),a)$. 
S'il est n\'ecessaire de pr\'eciser le groupe ambiant, 
nous noterons $\RC_{s,a,w}^\Gb$ la fonction $\RC_{s,a,w}$ et 
$\RC(s,a)_f^\Gb$ la fonction $\RC(s,a)_f$. 

\bigskip

\remarque{translation RC} 
Si $\t \in H^1(F^*,A_{\Gb^*}(s))$, nous identifierons $\t$ \`a une fonction 
centrale sur $W(s)^a \phi_1$ de la fa\c{c}on suivante~: si $w \in W(s)^a$, 
l'image de $w\phi_1$ par cette fonction centrale est \'egale \`a $\t(\wba)$, 
o\`u $\wba$ est l'image de $w$ \`a travers la suite de morphismes 
$W(s)^a \to A_{\Gb^*}(s) \to H^1(F^*,A_{\Gb^*}(s))$. Avec cette notation, on a, 
pour tout $z \in \Zb(\Gb)^F$ et pour tout $f \in \Cent(W(s)^a \phi_1)$, 
$$t_z^\Gb \RC(s,a)_f = \sha(z)\RC(s,a)_{f \omeh_s^1(\zba)}.$$
Ici, $\zba$ d\'esigne l'image de $z$ dans $\ZC(\Gb)^F$. Pour montrer cela, 
il suffit de remarquer que, d'apr\`es le lemme \ref{asai} et 
d'apr\`es la remarque \ref{elementaire series} (d), on a 
$$t_z^\Gb \RC_{s,a,w} = \sha(z) \omeh_s^1(\zba)(\wba) \RC_{s,a,w}$$
pour tout $w \in W(s)^a$.\finl

\bigskip

\proposition{produit scalaire RC}
{\it Supposons que la formule de Mackey a lieu dans $\Gb$. 
Soit $a \in A_{\Gb^*}(s)^{F^*}$ et 
soient $w$ et $w'$ deux \'el\'ements de $W(s)^a$. Alors
$$\langle \RC_{s,a,w},\RC_{s,a,w'} \rangle_{\Gb^F} = 
\begin{cases}
|C_{W(s)^a}(w\phi_1)| & \text{si $w\phi_1$ et $w'\phi_1$ sont 
conjugu\'es sous $W(s)^a$,} \\
0 & \text{sinon.}
\end{cases}$$}

\bigskip

\proof D'apr\`es la formule de Mackey, et compte tenu de la 
proposition \ref{rhodot absolument cuspidal}, on a 
$$\langle \RC_{s,a,w},\RC_{s,a,w'} \rangle_{\Gb^F} = 
\sum_{n \in [\NC_{w,w'}^F/\Lb_{s,a,w'}^F]} \langle \rhodot_{s_w,a}^{\Lb_{s,a,w}}, 
\lexp{n}{\rhodot_{s_{w'},a}^{\Lb_{s,a,w'}}} \rangle_{\Lb_{s,a,w}^F},$$
o\`u $\NC_{w,w'}=\{n \in \Gb~|~\Lb_{s,a,w} = \lexp{n}{\Lb_{s,a,w'}}\}$. 
D'autre part, on a une bijection naturelle entre $[\NC_{w,w'}^F/\Lb_{s,a,w'}^F]$ 
et $[\NC_{w,w'}^{*F^*}/\Lb_{s,a,w'}^{*F^*}]$ (o\`u bien s\^ur 
$\NC_{w,w'}^*=\{n \in \Gb^*~|~\Lb_{s,a,w}^* = \lexp{n}{\Lb_{s,a,w'}^*}\}$) et, 
\`a travers cette bijection, on a
$$\langle \RC_{s,a,w},\RC_{s,a,w'} \rangle_{\Gb^F} = 
\sum_{n \in [\NC_{w,w'}^{*F^*}/\Lb_{s,a,w'}^{*F^*}]} 
\langle \rhodot_{s_w,a}^{\Lb_{s,a,w}}, 
\rhodot_{ns_{w'}n^{-1},a}^{\Lb_{s,a,w}} \rangle_{\Lb_{s,a,w}^F}.$$
En particulier, si les couples $(\Lb_{s,a,w}^*,s_w)$ et 
$(\Lb_{s,a,w'}^*,ns_{w'}n^{-1})$ ne sont pas conjugu\'es sous $\Gb^{*F^*}$ 
(\cad si $w\phi_1$ et $w'\phi_1$ ne sont pas conjugu\'es sous $W(s)^a$), alors 
$\langle \RC_{s,a,w},\RC_{s,a,w'} \rangle_{\Gb^F} = 0$. 
Nous pouvons donc supposer maintenant que $w=w'$. On a, dans ce cas, 
$$\langle \RC_{s,a,w},\RC_{s,a,w'} \rangle_{\Gb^F} = 
\sum_{n \in [N_{\Gb^{*F^*}}(\Lb_{s,a,w}^*)/\Lb_{s,a,w}^{*F^*}]} 
\langle \rhodot_{s_w,a}^{\Lb_{s,a,w}}, 
\rhodot_{ns_w n^{-1},a}^{\Lb_{s,a,w}} \rangle_{\Lb_{s,a,w}^F}.$$
Soit maintenant $n \in N_{\Gb^{*F^*}}(\Lb_{s,a,w}^*)$. Posons 
$\b_n= \langle \rhodot_{s_w,a}^{\Lb_{s,a,w}}, 
\rhodot_{ns_w n^{-1},a}^{\Lb_{s,a,w}} \rangle_{\Lb_{s,a,w}^F}$. Si 
$s_w$ et $ns_w n^{-1}$ ne sont pas $\Lb_{s,a,w}^{*F^*}$-conjugu\'es, 
alors $\b_n=0$. 
Si $s_w$ et $ns_w n^{-1}$ sont $\Lb_{s,a,w}^{*F^*}$-conjugu\'es, 
alors il existe un repr\'esentant de la classe de $n$ dans 
$N_{\Gb^{*F^*}}(\Lb_{s,a,w}^*)/\Lb_{s,a,w}^{*F^*}$ qui centralise $s_w$ et 
alors $\b_n=|A_{\Lb_{s,a,w}^*}(s_w)^{F^*}|$. Par suite, 
\eqna
\langle \RC_{s,a,w},\RC_{s,a,w'} \rangle_{\Gb^F} &=& 
|A_{\Lb_{s,a,w}^*}(s_w)^{F^*}| \times 
|\Bigl(N_{\Gb^{*F^*}}(\Lb_{s,a,w}^*) \cap C_{\Gb^*}(s_w)^{F^*}\Bigr)/
C_{\Lb_{s,a,w}^*}(s_w)^{F^*}|\\
&=& |\Bigl(N_{\Gb^{*F^*}}(\Lb_{s,a,w}^*) \cap C_{\Gb^*}(s_w)^{F^*}\Bigr)/
C_{\Lb_{s,a,w}^*}^\circ(s_w)^{F^*}|.
\endeqna
Or, $C_{\Lb_{s,a,w}^*}^\circ(s_w)=\Tb_w^*$ et, d'apr\`es le corollaire 
\ref{ah ah ah} (e), on a 
$\Bigl( N_{\Gb^*}(\Lb_{s,a}^*) \cap C_{\Gb^*}(s)\Bigr) / \Tb_1^* \simeq W(s)^a$. 
D'o\`u le r\'esultat.\fin

\bigskip

\corollaire{isometrie RC}
{\it Supposons que la formule de Mackey a lieu dans $\Gb$. 
Alors l'application $\RC(s,a) : \Cent(W(s)^a \phi_1) \to \Cent(\Gb^F,(s),a)$ 
est une isom\'etrie.}

\bigskip

\soussection{Induction de Lusztig\label{soussection wl}} 
Soit $\Lb$ un \levi $F$-stable de $\Gb$ et soit $\Lb^*$ un \levi $F^*$-stable 
de $\Gb^*$ dual de $\Lb$. On suppose que $\Lb^*$ contient un \'el\'ement 
$s' \in \Lb^{*F^*}$ g\'eom\'etriquement conjugu\'e \`a $s$. Le but de cette 
section est de d\'ecrire l'action de l'induction de Lusztig $R_\Lb^\Gb$ 
sur l'image de $\RC(s',a)^\Lb$, pour $a \in A_{\Lb^*}(s')^{F^*}$. 
Le r\'esultat d\'ecrit cette action en termes d'une induction 
tordue entre les groupes $W_\Lb(s')$ et $W(s)$. Avant 
d'exprimer ce r\'esultat, nous avons besoin de comparer ces 
deux groupes. 
On se fixe un \para $\Pb^*$ de $\Gb^*$ dont $\Lb^*$ est un 
sous-groupe de Levi et on note $\Vb^*$ le radical unipotent de $\Pb^*$. 

Fixons tout d'abord un \'el\'ement $g \in \Gb^*$ tel que $gsg^{-1}=s'$. 
Soit $\Bb_\Lb^*$ un \borel $F^*$-stable de $C_{\Lb^*}^\circ(s')$ et 
soit $\Tb_\Lb^*$ un tore maximal $F^*$-stable de $\Bb_\Lb^*$. 
Alors $\lexp{g^{-1}}{(\Bb_\Lb^* C_{\Vb^*}(s'))}$ est un 
\borel de $C_{\Gb^*}^\circ(s)$ et $\lexp{g^{-1}}{\Tb_\Lb^*}$ 
est un tore maximal de $\lexp{g^{-1}}{(\Bb_\Lb^* C_{\Vb^*}(s'))}$. Donc 
il existe $h \in C_{\Gb^*}^\circ(s)$ tel que 
$$(\Tb_\Lb^*,\Bb_\Lb^*C_{\Vb^*}(s'))=\lexp{gh}{(\Tb_1^*,\Bb_1^*)}.$$
Notons que $(gh) s (gh)^{-1} = s'$. 
Par suite, $(gh)^{-1}F^*(gh)$ normalise $\Tb_1^*$ et centralise $s$~: 
on note $w_\Lb$ sa classe dans $W(s)$. Puisque le couple $(\Tb_\Lb^*,\Bb_\Lb^*)$ 
est bien d\'efini \`a conjugaison pr\`es par un \'el\'ement de 
$C_{\Lb^*}^\circ(s')^{F^*}$, l'\'el\'ement $w_\Lb$ est bien d\'efini par la 
couple $(\Lb^*,s')$. En particulier, 
si on identifie $A_{\Lb^*}(s')$ avec le sous-groupe correspondant 
de $A_{\Gb^*}(s)$ (via la conjugaison par $gh$), alors $w_\Lb$ commute 
avec $A_{\Lb^*}(s')$. D'autre part, via la conjugaison par $gh$, nous 
verrons $W_\Lb(s')$ et $W_\Lb^\circ(s')$ comme des sous-groupes 
$w_\Lb F^*$-stables de $W(s)$ et $W^\circ(s)$ respectivement.

\bigskip

\proposition{induction RC}
{\it Soit $a \in A_{\Lb^*}(s')^{F^*}$ et identifions $a$ avec un \'el\'ement 
de $A_{\Gb^*}(s)^{F^*}$ comme ci-dessus. Alors le diagramme 
$$\diagram
\Cent(W_\Lb(s')^a w_\Lb \phi_1) \rrto^{{\RC(s',a)^\Lb}} 
\ddto_{{\Ind_{W_\Lb(s')^a w_\Lb \phi_1}^{W(s)^a\phi_1}}} && 
\Cent(\Lb^F,(s'),a) \ddto^{{R_\Lb^\Gb}} \\
&&\\
\Cent(W(s)^a\phi_1) \rrto^{{\RC(s,a)^\Gb}} && \Cent(\Gb^F,(s),a)
\enddiagram$$
est commutatif.}

\bigskip

\proof Si $w \in W_\Lb(s')^a$ (vu comme un sous-groupe de $W(s)^a$), 
nous fixons un \'el\'ement $l_w \in \Lb^{*F^*}$ 
tel que $l_w^{-1} F^*(l_w)$ appartienne au normalisateur de $\Tb_\Lb^*$ 
et repr\'esente $(gh)w(gh)^{-1}$. On a $\Tb_\Lb^*=(gh)\Tb_1^*(gh)^{-1}$ 
et remarquons que $(l_w gh)^{-1} F^*(l_w gh)$ 
normalise $\Tb_1^*$ et repr\'esente $w w_\Lb$. 
La proposition d\'ecoule alors facilement de cette observation, 
de la transitivit\'e de l'induction et de \cite[lemme 3.1.1]{bonnafe couro}.\fin 

\bigskip

\soussection{Transform\'es de Fourier de caract\`eres semi-simples} 
Si $A$ est un groupe ab\'elien fini et si $\ph$ est un automorphisme 
de $A$, on note $\MC(A,\ph)$ le groupe $(A^\ph)^\we \times H^1(\ph,A)$. 
Son dual $\MC(A,\ph)^\we$ est \'egal \`a $A^\ph \times H^1(\ph,A)^\we$. 
Si $(a,\t) \in \MC(A_{\Gb^*}(s),F^*)^\we$, on pose
$$\rhoh_{s,a,\t}=\rhoh_{s,a,\t}^\Gb = \frac{1}{|A_{\Gb^*}(s)^{F^*}|} 
\sum_{(\xi,\a) \in \MC(A_{\Gb^*}(s),F^*)} \t(\a) \xi(a)^{-1} 
\r_{s_\a,\xi}.$$
Ici, nous avons identifi\'e le groupe $A_{\Gb^*}(s_\a)$ avec 
le groupe $A_{\Gb^*}(s)$ (via la conjugaison par l'\'el\'ement $g_\a$ 
tel que $g_\a s g_\a^{-1} = s_\a$)~: cette identification ne change 
pas l'action du morphisme de Frobenius car $A_{\Gb^*}(s)$ est ab\'elien. 
Si $(\xi,\a) \in \MC(A_{\Gb^*}(s),F^*)$, alors
\equat\label{retrouvage}
\r_{s_\a,\xi} = \frac{1}{|A_{\Gb^*}(s)^{F^*}|} 
\sum_{(a,\t) \in \MC(A_{\Gb^*}(s),F^*)^\we} \t(\a)^{-1}\xi(a) \rhoh_{s,a,\t}.
\endequat

\bigskip

\proposition{proprietes rhoh}
{\it Soient $(a,\t)$ et $(a',\t')$ deux \'el\'ements de $\MC(A_{\Gb^*}(s),F^*)^\we$. 
Alors~:
\begin{itemize}
\itemth{a} $\rhoh_{s,a,\t} \in \Cent(\Gb^F,(s),a)$. 

\itemth{b} $\DS{\rhoh_{s,a,\t}=\frac{1}{|A_{\Gb^*}(s)^{F^*}|} 
\sum_{\a \in H^1(F^*,A_{\Gb^*}(s))} \t(\a) \rhodot_{s_\a,a}}$.

\itemth{c} $\langle \rhoh_{s,a,\t},\rhoh_{s,a',\t'} \rangle_{\Gb^F} = 
\begin{cases}
1 & \text{si }(a,\t)=(a',\t'), \\
0 & \text{sinon.}
\end{cases}$

\itemth{d} Si $z \in \Zb(\Gb)^F$, alors 
$t_z^\Gb \rhoh_{s,a,\t} = \sha(z) \rhoh_{s,a,\t\omeh_s^1(z)}$. 

\itemth{e} Si la formule de Mackey et la conjecture $(\GG)$ ont lieu 
dans $\Gb$, alors $$\rhoh_{s,a,\t} = \RC(s,a)_\t.$$ 
Ici, $\t$ est vu 
comme la fonction centrale sur $W(s)^a \phi_1$ qui envoie $w\phi_1$ 
sur $\t(\wba)$, o\`u $\wba$ d\'esigne l'image de $w$ dans $H^1(F^*,A_{\Gb^*}(s))$.
\end{itemize}}

\bigskip

\proof (a), (b) et (c) sont \'evidents. (d) se montre de la m\^eme mani\`ere 
que la premi\`ere \'egalit\'e de la remarque \ref{translation RC}. Montrons (e). 
Tout d'abord, d'apr\`es (a) et le corollaire \ref{isometrie RC}, 
on a $\langle \rhoh_{s,a,\t},\rhoh_{s,a,\t} \rangle_{\Gb^F} = 
\langle \RC(s,a)_\t,\RC(s,a)_\t \rangle_{\Gb^F} = 1$. Il nous reste 
\`a montrer que 
$$\langle \rhoh_{s,a,\t},\RC(s,a)_\t \rangle_{\Gb^F}=1.\leqno{(*)}$$ 
Soit $w \in W(s)^a$ et notons $\a$ la classe de $w$ dans $H^1(F^*,A_{\Gb^*}(s))$. 
Pour montrer $(*)$, il suffit de montrer que 
$$\langle \rhoh_{s,a,\t} , 
\t(\a) R_{\Lb_{s,a,w}}^\Gb \rhodot_{s_w,a}^{\Lb_{s,a,w}} \rangle_{\Gb^F} = 1.
\leqno{(**)}$$
Mais, d'apr\`es (b), d'apr\`es \ref{serie RC} et d'apr\`es 
le th\'eor\`eme \ref{rlg series}, on a, par adjonction, 
$$\langle \rhoh_{s,a,\t} , 
\t(\a) R_{\Lb_{s,a,w}}^\Gb \rhodot_{s_w,a}^{\Lb_{s,a,w}} \rangle_{\Gb^F} =
\frac{1}{|A_{\Gb^*}(s)^{F^*}|} \langle \t(\a) 
\lexp{*}{R}_{\Lb_{s,a,w}}^\Gb \rhodot_{s_\a,a} , 
\t(\a) \rhodot_{s_w,a}^{\Lb_{s,a,w}} \rangle_{\Gb^F}.$$
Par suite, d'apr\`es \ref{scalaire rhodot} et 
le corollaire \ref{coro restriction rhosa}, on a 
$$\langle \rhoh_{s,a,\t} , 
\t(\a) R_{\Lb_{s,a,w}}^\Gb \rhodot_{s_w,a}^{\Lb_{s,a,w}} \rangle_{\Gb^F} = 
\frac{1}{|A_{\Gb^*}(s)^{F^*}|} \times 
\frac{|A_{\Gb^*}(s)^{F^*}|}{|A_{\Lb_{s,a,w}^*}(s_w)^{F^*}|} 
\times |A_{\Lb_{s,a,w}^*}(s_w)^{F^*}| = 1,$$
ce qui montre $(**)$.\fin

\bigskip

\exemple{a=1}
Supposons dans cet exemple, et uniquement dans cet exemple, que $a=1$. 
Nous poserons alors $\RC(s)=\RC(s,1)$. D'autre part, $\Lb_{s,1,w}^*=\Tb_w^*$. 
Si $\eta$ est un caract\`ere irr\'eductible $F^*$-stable de $W(s)$ et 
si $\etat$ est une extension de $\eta$ \`a $W(s) \rtimes <\phi_1>$, alors 
$\RC(s)_\etat$ est un {\it caract\`ere fant\^ome} de $\Gb^F$. Tous les 
caract\`eres fant\^omes de $\Gb^F$ ne sont pas obtenus ainsi.\finl

\bigskip

\soussection{S\'eries rationnelles} 
Nous allons maintenant construire une isom\'etrie $\RC[s,a]$ de 
l'espace des fonctions centrales sur $W^\circ(s)^a\phi_1$ invariantes 
par l'action de $A_{\Gb^*}(s)^{F^*}$ vers $\Cent(\Gb^F,[s],a)$. 
Si $f \in \Cent(W^\circ(s)^a\phi_1)$, on pose
$$\RC[s,a]_f=\RC[s,a]_f^\Gb = \frac{1}{|A_{\Gb^*}(s)^{F^*}|.|W^\circ(s)^a|} 
\sum_{w \in W^\circ(s)^a} f(w\phi_1) \RC_{s,a,w}.$$
D'apr\`es \ref{serie RC}, on a $\RC_{s,a,w} \in \Cent(\Gb^F,[s],a)$ pour 
tout $w \in W^\circ(s)^a$. 
On a  donc d\'efini une application 
$\RC[s,a] : \Cent(W^\circ(s)^a\phi_1) \to \Cent(\Gb^F,[s],a)$ dont 
il est facile de v\'erifier que, si $f \in \Cent(W^\circ(s)^a\phi_1)$ 
et $b \in A_{\Gb^*}(s)^{F^*}$, alors 
\equat\label{invariance RC}
\RC[s,a]_f=\RC[s,a]_{\lexp{b}{f}}.
\endequat
En particulier, l'image de $\RC[s,a]$ est \'egale \`a l'image de 
sa restriction \`a $\bigl(\Cent(W^\circ(s)^a\phi_1)\bigr)^{A_{\Gb^*}(s)^{F^*}}$. 
Si $f$ et $g$ sont deux \'el\'ements de 
$\bigl(\Cent(W^\circ(s)^a\phi_1)\bigr)^{A_{\Gb^*}(s)^{F^*}}$, on pose 
$$\langle f,g \rangle_{s,a} = 
\frac{\langle f,g \rangle_{W^\circ(s)^a\phi_1}}{|A_{\Gb^*}(s)^{F^*}|}.$$
Alors $\langle , \rangle_{s,a}$ est un produit scalaire sur 
$\bigl(\Cent(W^\circ(s)^a\phi_1)\bigr)^{A_{\Gb^*}(s)^{F^*}}$. 

\bigskip

\proposition{isometrie R}
{\it Soit $a \in A_{\Gb^*}(s)^{F^*}$. Alors l'application 
$\RC[s,a] : \bigl(\Cent(W^\circ(s)^a\phi_1)\bigr)^{A_{\Gb^*}(s)^{F^*}} \to 
\Cent(\Gb^F,[s],a)$ est une isom\'etrie (pour les produits scalaires 
$\langle , \rangle_{s,a}$ et $\langle , \rangle_{\Gb^F}$).}

\bigskip

\proof Si $f \in \Cent(W(s)^a\phi_1)$, on note $\Res_{s,a}^\circ f$ sa restriction 
\`a $W^\circ(s)^a \phi_1$. Il est alors imm\'ediat que 
$\Res_{s,a}^\circ f$ est stable sous l'action de $A_{\Gb^*}(s)^{F^*}$. Cela 
nous d\'efinit donc une application 
$\Res_{s,a}^\circ : \Cent(W(s)^a \phi_1) \to 
\bigl(\Cent(W^\circ(s)^a\phi_1)\bigr)^{A_{\Gb^*}(s)^{F^*}}$. R\'eciproquement, si 
$f \in \bigl(\Cent(W^\circ(s)^a\phi_1)\bigr)^{A_{\Gb^*}(s)^{F^*}}$, 
on pose, pour $w \in W^\circ(s)^a\phi_1$ et $b \in A_{\Gb^*}(s)$, 
$$(\Ext_{s,a}^\circ f)(wb\phi_1) = \begin{cases}
f(c w c^{-1} \phi_1) & \text{si } b=c^{-1} F^*(c)\text{ pour un }c \in A_{\Gb^*}(s),\\
0 & \text{si }b\not=c^{-1} F^*(c)\text{ pour tout }c \in A_{\Gb^*}(s).
\end{cases}.$$
Il est \`a noter que la premi\`ere formule ne d\'epend pas du choix 
de $c$ car $f$ est invariante sous l'action de $A_{\Gb^*}(s)^{F^*}$. 
Il est alors clair que $\Ext_{s,a}^\circ f \in \Cent(W(s)^a\phi_1)$. On a 
donc d\'efini une application 
$\Ext_{s,a}^\circ : \bigl(\Cent(W^\circ(s)^a\phi_1)\bigr)^{A_{\Gb^*}(s)^{F^*}}
\to \Cent(W(s)^a\phi_1)$. De plus,
\equat
\Res_{s,a}^\circ \circ \Ext_{s,a}^\circ = 
\Id_{(\Cent(W^\circ(s)^a\phi_1))^{A_{\Gb^*}(s)^{F^*}}}.
\endequat
D'autre part, $\Ext_{s,a}^\circ$ est une isom\'etrie (pour les produits 
scalaires $\langle , \rangle_{s,a}$ et $\langle , \rangle_{W(s)^a\phi_1}$) 
et le diagramme 
\equat\label{RC RR}
\diagram 
\bigl(\Cent(W^\circ(s)^a\phi_1)\bigr)^{A_{\Gb^*}(s)^{F^*}} 
\ddto_{\Ext_{s,a}^\circ} \rrto^{\quad\RC[s,a]} && \Cent(\Gb^F,[s],a) \ddto\\
&& \\
\Cent(W(s)^a\phi_1) \rrto^{\RC(s,a)} && \Cent(\Gb^F,(s),a) \\
\enddiagram
\endequat
est commutatif. Les applications $\RC(s,a)$, $\Ext_{s,a}^\circ$ et 
l'injection $\Cent(\Gb^F,[s],a) \injto \Cent(\Gb^F,(s),a)$ \'etant 
des isom\'etries, on en d\'eduit que $\RC[s,a]$ est une isom\'etrie.\fin

\bigskip

\remarque{possible} Il \'etait possible de d\'emontrer directement en utilisant 
la proposition \ref{produit scalaire RC} que $\RC[s,a]$ est une isom\'etrie. 
Nous avons cependant voulu introduire les applications $\Res_{s,a}^\circ$ 
et $\Ext_{s,a}^\circ$ car elles nous seront utiles par la suite.\finl

\bigskip

\proposition{rhodot R}
{\it Si la formule de Mackey et la conjecture $(\GG)$ ont lieu dans $\Gb$, 
alors
$$\RC[s,a]_1 = \frac{1}{|A_{\Gb^*}(s)^{F^*}|} \rhodot_{s,a}.$$
Ici, $1$ est vu comme la fonction constante et \'egale \`a $1$.}

\bigskip

\proof Notons $\pi_{[s]} : \Cent(\Gb^F) \to \Cent(\Gb^F,[s])$ la projection 
orthogonale. Alors le diagramme 
\equat\label{RC R}
\diagram 
\Cent(W(s)^a\phi_1) 
\ddto_{\Res_{s,a}^\circ} \rrto^{\quad\RC(s,a)} && 
\Cent(\Gb^F,(s),a) \ddto_{\pi_{[s]}}\\
&& \\
\bigl(\Cent(W^\circ(s)^a\phi_1)\bigr)^{A_{\Gb^*}(s)^{F^*}} \rrto^{\RC[s,a]} && 
\Cent(\Gb^F,[s],a) \\
\enddiagram
\endequat
est commutatif. D'o\`u $\RC[s,a]_1 = \pi_{[s]} \RC(s,a)_1$. 
Le r\'esultat d\'ecoule alors de la proposition \ref{proprietes rhoh} (e).\fin

\bigskip

Nous concluons ce chapitre par un r\'esultat d\'ecrivant l'induction de 
Lusztig \`a travers les applications $\RC[s,a]$. Soit donc $\Lb$ un 
\levi $F$-stable de $\Gb$ et soit $\Lb^*$ un \levi $F^*$-stable 
de $\Gb^*$ dual de $\Lb$. On suppose que $s \in \Lb^{*F^*}$. 
Reprenons les notations de \SEC\ref{soussection wl} (en rempla\c{c}ant 
$s'$ par $s$), de sorte que $W_\Lb(s)$ est vu comme un sous-groupe 
$w_\Lb F^*$-stable de $W(s)$. Remarquons aussi que, puisque $s'=s$, 
on a $w_\Lb \in W^\circ(s)$. 

\bigskip

\proposition{induction RRR}
{\it Supposons que $a \in A_{\Lb^*}(s)^{F^*}$. Alors le diagramme suivant 
est commutatif~:
$$\diagram 
\Cent(W_\Lb^\circ(s)^a w_\Lb \phi_1) 
\xto[0,3]^{{|A_{\Lb^*}(s)^{F^*}| \RC[s,a]^\Lb}} 
\ddto_{\DS{\Ind_{W_\Lb^\circ(s)^a w_\Lb \phi_1}^{W^\circ(s)^a \phi_1}}} && &
\Cent(\Lb^F,[s]) \ddto_{\DS{R_\Lb^\Gb}} \\
&&& \\
\Cent(W^\circ(s)^a \phi_1) \xto[0,3]^{{|A_{\Gb^*}(s)^{F^*}|\RC[s,a]^\Gb}} 
&&& \Cent(\Gb^F,[s]). \\
\enddiagram$$}

\bigskip

\proof Le m\^eme argument que dans la preuve de la proposition 
\ref{induction RC} alli\'e encore \`a \cite[lemme 3.1.1]{bonnafe couro} 
prouve imm\'ediatement cette proposition.\fin

\newpage

{\Large \part{Faisceaux-caract\`eres\label{chapitre faisceaux}}}

\bigskip

L'objet de ce chapitre est l'\'etude de l'influence de la non connexit\'e 
du centre de $\Gb$ sur la th\'eorie des faisceaux-caract\`eres. 
Le centre de $\Gb$ agit sur les faisceaux-caract\`eres de deux fa\c{c}ons. 
La premi\`ere est la trace de l'action par conjugaison~: 
cette action induit une action sur chaque faisceau-caract\`ere 
par multiplication par un caract\`ere lin\'eaire de $\ZC(\Gb)$. 
Le calcul de ce caract\`ere lin\'eaire est classique mais 
nous le rappelons ici (voir proposition \ref{calcul zeta}). 
La deuxi\`eme action est l'action par translation~: le translat\'e 
d'un faisceau-caract\`ere par un \'el\'ement de $\Zb(\Gb)$ 
est encore un faisceau-caract\`ere. Cela induit une action 
de $\ZC(\Gb)$ par permutations de 
l'ensemble des (classes d'isomorphie de) faisceaux-caract\`eres.
Nous \'etudions cette action \`a travers le processus d'induction 
parabolique (voir th\'eor\`eme \ref{translation z}). 
Nous en profitons pour tirer quelques cons\'equences de ce th\'eor\`eme 
sur le param\'etrage ou sur les fonctions caract\'eristiques 
des faisceaux-caract\`eres. \`A partir de la section 
\ref{section unip reg}, nous nous consacrons aux faisceaux-caract\`eres 
apparaissant dans l'induit d'un faisceau-caract\`ere cuspidal 
dont le support rencontre la classe unipotente r\'eguli\`ere. 
Nous y \'etablissons un param\'etrage de tels faisceaux-caract\`eres 
s\'eries par s\'eries et obtenons une formule pour leurs fonctions 
caract\'eristiques comme combinaisons lin\'eaires d'induits 
de Lusztig de caract\`eres semi-simples (voir th\'eor\`eme \ref{PRINCIPAL}). 
Comme cons\'equence, nous obtenons que la fonction caract\'eristique 
d'un faisceau-caract\`ere, non n\'ecessairement cuspidal, dont le support 
rencontre la classe unipotente r\'eguli\`ere est une transform\'ee 
de Fourier de caract\`eres semi-simples (voir corollaire 
\ref{regulier caracteristique}). 

\bigskip

\section{Action de $\ZC(\Gb)$ sur les faisceaux-caract\`eres}

\medskip

Nous rappelons dans cette section comment sont construits les 
faisceaux-caract\`eres avant d'\'etudier l'action de $\ZC(\Gb)$. 

\bigskip

\soussection{Syst\`emes locaux kumm\'eriens} 
Fixons un \borel $\Bb$ de $\Gb$ ainsi qu'un tore maximal 
$\Tb$ de $\Bb$. Soit $\Ub$ le radical unipotent de $\Bb$. 
Nous fixons aussi un tore maximal $\Tb^*$ de $\Gb^*$ 
dual de $\Tb$. Nous identifierons le groupe de Weyl $W$ de 
$\Gb$ relativement \`a $\Tb$ avec celui de $\Gb^*$ relativement 
\`a $\Tb^*$. Nous ferons aussi l'identification $X(\Tb)=Y(\Tb^*)$. 

Notons $\SC(\Tb)$ l'ensemble des classes d'isomorphie 
de syst\`emes locaux kumm\'eriens sur $\Tb$. Le produit tensoriel 
munit $\SC(\Tb)$ d'une structure de groupe ab\'elien. 
D'autre part, le groupe $W$ agit naturellement sur $\SC(\Tb)$. 
Le choix des applications $\imath$, $\jmath$ et $\kappa$ construites 
dans \SEC\ref{sous groupes} permet de 
construire un isomorphisme $W$-\'equivariant de groupes ab\'eliens 
$\Tb^* \simeq \SC(\Tb)$, $s \mapsto \LC_s$. 
Nous allons rappeler sa d\'efinition~: si $s \in \Tb^*$, 
il existe $x \in X(\Tb)=Y(\Tb^*)$ et $n \in \NM^*$, premier 
\`a $p$, tels que $\tilde{\imath}_{\Tb^*}(x/n)=s$. 
On note $e_n : \FM^\times \to \FM^\times$, $z \mapsto z^n$. 
C'est un rev\^etement \'etale galoisien de groupe $\mub_n(\FM)$. 
Nous noterons $\XC_n$ le syst\`eme local sur $\FM^\times$ 
associ\'e \`a ce rev\^etement et au caract\`ere lin\'eaire 
$\kappa : \mub_n(\FM) \to \qlb^\times$. On a alors~:
\equat\label{ls}
\LC_s=x^* \XC_n.
\endequat
Ici, $x : \Tb \to \FM^\times$ est seulement vu comme un 
morphisme de vari\'et\'es. 

\bigskip

\soussection{Faisceaux-caract\`eres} 
Fixons maintenant un \'el\'ement $w$ de $W$ et un repr\'esentant 
$\wdo$ de $w$ dans $N_\Gb(\Tb)$. Nous noterons $\pi_\wdo : \Bb w\Bb \to \Tb$ 
l'unique application telle que, si $v$ et $v'$ appartiennent \`a $\Ub$ 
et $t \in \Tb$, alors $\pi_\wdo(v\wdo t v')=t$. 
C'est un morphisme de vari\'et\'es. Soient 
$$\Ybh_w=\{(g,h\Ub) \in \Gb \times \Gb/\Ub~|~h^{-1}gh \in \Bb w \Bb\}$$
$$\Ybt_w=\{(g,h\Bb) \in \Gb \times \Gb/\Bb~|~h^{-1}gh \in \Bb w \Bb\}.
\leqno{\text{et}}$$
Notons $\b_w : \Ybh_w \to \Ybt_w$ l'application canonique. 
Posons 
$$\fonction{\a_\wdo}{\Ybh_w}{\Tb}{(g,h\Ub)}{\pi_\wdo(h^{-1}gh)}$$
$$\fonction{\g_w}{\Ybt_w}{\Gb}{(g,h\Bb)}{g.}\leqno{\text{et}}$$ 
Alors $\a_\wdo$, $\b_w$ et $\g_w$ sont des morphismes de vari\'et\'es 
bien d\'efinis. Nous avons donc construit un diagramme 
$$\diagram
\Tb && \dot{\Yb}_w \llto_{\DS{\a_\wdo}} \rrto^{\DS{\b_w}} && 
\Ybt_w \rrto^{\DS{\g_w}} && \Gb.
\enddiagram$$
Le groupe $\Tb$ agit sur $\Ybh_w$ de la fa\c{c}on suivante~: si $t \in \Tb$ 
et si $(g,h\Ub) \in \Ybh_w$, on pose
$$t*(g,h\Ub)=(g,ht^{-1}\Ub).$$
Alors $\b_w$ est une fibration principale de groupe $\Tb$. 
D'autre part, le groupe $\Tb$ agit sur $\Tb$ de la fa\c{c}on suivante~: 
si $t$ et $t'$ appartiennent \`a $\Tb$, on pose
$$t*_w t'= \wdo^{-1} t \wdo t' t^{-1}.$$
Alors il est facile de v\'erifier que $\a_\wdo$ est $\Tb$-\'equivariante. 
De plus, le groupe $\Gb$ agit diagonalement sur $\Ybh_w$ et $\Ybt_w$ 
par conjugaison sur la premi\`ere coordonn\'ee et par translation \`a gauche 
sur la deuxi\`eme, et il agit 
sur $\Gb$ par conjugaison. Les morphismes $\b_w$ et $\g_w$ sont 
alors $\Gb$-\'equivariants.

Soit $s \in \Tb^*$ et 
supposons que $w$ v\'erifie $w(s)=s$. Alors, d'apr\`es 
\cite[2.2.2]{lucs}, $\LC_s$ est $\Tb$-\'equivariant pour l'action $*_w$. 
En particulier, $\a_\wdo^* \LC_s$ est un syst\`eme local 
$\Tb$-\'equivariant sur $\Ybh_w$. Par suite, il existe un unique 
(\`a isomorphisme pr\`es) syst\`eme local $\LCt_{w,s}$ sur $\Ybt_w$ 
tel que $\b_w^* \LCt_{w,s} \simeq \a_\wdo^* \LC_s$. De plus, puisque 
$\a_\wdo^* \LC_s$ est $\Gb$-\'equivariant, il en est de m\^eme 
de $\LCt_{w,s}$. Posons~:
$$K_{w,s}=K_{w,s}^\Gb = R(\g_w)_! \LCt_{w,s}.$$
On rappelle qu'un faisceau pervers irr\'eductible $A$ sur $\Gb$ est 
appel\'e un {\it faisceau-caract\`ere} s'il existe un triplet $(s,w,i)$ 
o\`u $s \in \Tb^*$ et $w \in W$ v\'erifient $w(s)=s$ 
et $i$ est un entier relatif tels que $A$ soit une composante 
du faisceau pervers $\lexp{p}{H}^i(K_{w,s})$. Nous noterons 
$\faisceau(\Gb)$ l'ensemble des classes d'isomorphie de 
faisceaux-caract\`eres sur $\Gb$. Il est \`a noter que 
tout faisceau-caract\`ere est $\Gb$-\'equivariant pour 
l'action par conjugaison. 

\medskip

Nous allons conclure cette sous-section par une construction 
explicite du syst\`eme local $\LCt_{w,s}$. Pour cela, \'ecrivons 
$s=\tilde{\imath}_{\Tb^*}(x/n)$ comme pr\'ec\'edemment. 
Dire que $w(s)=s$ \'equivaut \`a dire que $\l=w(x/n)-x/n \in X(\Tb)$. 
En d'autres termes, $w(x)-x=n\l$, avec $\l \in X(\Tb)$. 
Soit $\FM_\l$ le $\Bb$-module irr\'eductible 
$\FM$ sur lequel $\Bb$ agit via l'unique caract\`ere 
$\lamt : \Bb \to \FM^\times$ qui \'etend $\l$. 
Soit $\BCB_\l$ le fibr\'e en droite associ\'e \`a $\l$ (il est obtenu 
en quotientant par $\Bb$ la vari\'et\'e $\Gb \times \FM_\l$, $\Bb$ 
agissant diagonalement sur $\Gb \times \FM_\l$ 
par translations \`a droite sur la premi\`ere coordonn\'ee et 
par le caract\`ere $\lamt$ sur la deuxi\`eme). 
Si $(g,z) \in \Gb \times \FM$, nous noterons $g *_\l z$ sa classe dans $\BCB_\l$. 
Nous noterons $\BCB_\l^\times$ le compl\'ementaire de la section nulle 
dans $\BCB_\l$. Posons alors 
$$\Ybh_{w,x,n}=\{(g,h\Ub,z) \in \Gb \times \Gb/\Ub \times \FM^\times~|~
h^{-1}gh \in \Bb w \Bb \text{ et } z^n=x(\pi_\wdo(h^{-1}gh)) \}$$
$$\Ybt_{w,x,n}=\{(g,h *_\l z) \in \Gb \times \BCB_\l^\times~|~h^{-1}gh \in \Bb w \Bb 
\text{ et } z^n=x(\pi_\wdo(h^{-1}gh)) \}.\leqno{\text{et}}$$
Il est facile de voir que les vari\'et\'es $\Ybh_{w,x,n}$ et $\Ybt_{w,x,n}$ sont 
bien d\'efinies. Notons $\fha_{w,x,n} : \Ybh_{w,x,n} \to \Ybh_w$, 
$(g,h\Ub,z) \mapsto (g,h\Ub)$ et $\fti_{w,x,n} : \Ybt_{w,x,n} \to \Ybt_w$, 
$(g, h *_\l z) \mapsto (g,h\Bb)$. Posons aussi 
$\alpt_{w,x,n} : \Ybh_{w,x,n} \to \FM^\times$, $(g,h\Ub,z) \mapsto z$ 
et $\bett_{w,x,n} : \Ybh_{w,x,n} \to \Ybt_{w,x,n}$, 
$(g,h\Ub,z) \mapsto (g,h*_\l z)$. Alors le diagramme 
$$\diagram
\FM^\times \ddto_{\DS{e_n}} && \Ybh_{w,x,n} \llto_{\DS{\alpt_{w,x,n}}} 
\rrto^{\DS{\bett_{w,x,n}}} 
\ddto_{\DS{\fha_{w,x,n}}} && \Ybt_{w,x,n} \ddto_{\DS{\fti_{w,x,n}}} && \\
&&&&&& \\
\FM^\times && 
\Ybh_w \llto_{\DS{x \circ \a_\wdo}} \rrto^{\DS{\b_w}} && 
\Ybt_w \rrto^{\DS{\g_w}} && \Gb
\enddiagram$$
est commutatif. De plus, les carr\'es sont cart\'esiens et 
les morphismes $\fha_{w,x,n}$ et $\fti_{w,x,n}$ sont des rev\^etements 
galoisiens de groupe $\mub_n(\FM)$. Le lemme suivant est alors 
imm\'ediat~:

\bigskip

\lemme{description ltilde}
{\it $\LCt_{w,s}$ est le syst\`eme local sur $\Ybt_w$ associ\'e 
au rev\^etement \'etale $\fti_{w,x,n}$ et au caract\`ere lin\'eaire 
$\kappa : \mub_n(\FM) \to \qlb^\times$.}

\bigskip 

\remarque{action t}
Le groupe $\Tb$ agit sur $\Ybh_{w,x,n}$ comme suit~: si $t \in \Tb$ et 
si $(g,h\Ub,z) \in \Ybh_{w,x,n}$, alors on pose 
$$\lexp{t}{(g,h\Ub,z)}=(g,ht^{-1}\Ub,\l(t)z).$$
Il est alors facile de voir que $\fti_{w,x,n}$ est $\Tb$-\'equivariant 
et que $\Ybt_{w,x,n}$ est le quotient de $\Ybh_{w,x,n}$ par cette action de $\Tb$.\finl

\bigskip

\soussection{Action de ${\boldsymbol{\ZC(\Gb)}}$} 
Si $A$ est un faisceau pervers irr\'eductible $\Gb$-\'equivariant sur $\Gb$, 
l'action par conjugaison de $\Gb$  sur lui-m\^eme induit une action de $\Zb(\Gb)$ 
sur $A$. Cette action se factorise par le groupe connexe $\Zb(\Gb)^\circ$ et, 
puisque $A$ est irr\'eductible, cette action est donn\'ee par un 
caract\`ere lin\'eaire $\z_A : \ZC(\Gb) \to \qlb^\times$. La proposition 
suivante donne un moyen de calculer ce caract\`ere lin\'eaire lorsque 
$A$ est un faisceau-caract\`ere~:

\bigskip

\proposition{calcul zeta}
{\it Soit $A$ un faisceau-caract\`ere sur $\Gb$. 
Soient $s \in \Tb^*$, $w \in W$ et $i \in \ZM$ tels que 
$w(s)=s$ et $A$ soit une composante irr\'eductible de $\lexp{p}{H}^i(K_{w,s})$. 
Notons $\wba$ l'image de $w$ dans $A_{\Gb^*}(s)$. Alors
$$\z_A=\omega_s(\wba).$$}

\bigskip

\proof Pour cela, il suffit de calculer l'action de $\ZC(\Gb)$ sur 
$\LCt_{w,s}$. En effet, l'action de $\Gb$ sur $\Ybt_w$ induit une action 
trivale de $\Zb(\Gb)$, donc $\ZC(\Gb)$ agit sur le syst\`eme local 
$\LCt_{w,s}$ par multiplication par un caract\`ere lin\'eaire, 
qui ne peut \^etre que $\z_A$. 

On utilise pour cela la description de $\LCt_{w,s}$ donn\'ee par le 
lemme \ref{description ltilde} dont on reprend les notations ($x$, $n$, $\l$...). 
Le groupe $\Gb$ agit sur $\Ybt_{w,x,n}$ comme suit~: si $\g \in \Gb$ et 
si $(g,h *_\l z) \in \Ybt_{w,x,n}$, on pose
$$
\lexp{\g}{(g,h*_\l z)}=(\g g \g^{-1},\g h *_\l z).
$$
Alors $\fti_{w,x,n}$ est $\Gb$-\'equivariant et il suffit de regarder comment 
agit $\g \in \Zb(\Gb)$. On a, si $\g \in \Zb(\Gb)$, 
$$
\lexp{\g}{(g,h*_\l z)}=(g,h*_\l \l(\g) z).
$$
Or, on a $\l(\g)^n=x(w^{-1}(\g)) x(\g)^{-1} = 1$ car $\g$ est central, 
donc $\l(\g) \in \mub_n(\FM)$. Par suite, $\g$ agit sur 
$\LCt_{w,s}$ par multiplication par $\kappa(\l(\g))$. Mais, par d\'efinition 
de $\omega_s$, on a $\kappa(\l(\g))=\omega_s(\wba)(\bar{\g})$, o\`u 
$\bar{\g}$ d\'esigne la classe de $\g$ dans $\ZC(\Gb)$ 
(voir \ref{definition omega}).\fin

\bigskip

\soussection{S\'eries g\'eom\'etriques} 
Soit $s$ un \'el\'ement semi-simple de $\Gb^*$. Nous noterons 
$\faisceau(\Gb,(s))$ l'ensemble des (classes d'isomorphie de) 
faisceaux-caract\`eres $A$ sur $\Gb$ tels qu'il existe $i \in \ZM$, 
$t \in (s) \cap \Tb^*$ et $w \in W$ tels que 
$w(t)=t$ et $A$ soit une composante de $\lexp{p}{H}^i(K_{w,t})$. 

Il r\'esulte de \cite[proposition 11.2 (c)]{lucs} que
\equat\label{partition faisceaux}
\faisceau(\Gb)=\coprod_{(s)} \faisceau(\Gb,(s)).
\endequat
Bien s\^ur, les ensembles $\faisceau(\Gb)$ et $\faisceau(\Gb,(s))$ ne d\'ependent 
pas des choix de $\Tb$, $\Bb$ et $\Tb^*$. Nous utiliserons \`a loisir 
cette souplesse en fonction des questions que nous aborderons.

Si $a \in A_{\Gb^*}(s)$, notons $\faisceau(\Gb,(s),a)$ l'ensemble des 
(classes d'isomorphie de) faisceaux-caract\`eres $A$ sur $\Gb$ tels qu'il 
existe $i \in \ZM$, $t \in (s) \cap \Tb^*$ et $w \in W$ tels que 
$w(t)=t$, $\wba \sim a$ et $A$ soit une composante de $\lexp{p}{H}^i(K_{w,t})$. 
Ici, la notation $\wba \sim a$ signifie que la classe $\wba$ de $w$ dans 
$A_{\Gb^*}(t) \simeq A_{\Gb^*}(s)$ est \'egale \`a $a$, l'isomorphisme 
entre $A_{\Gb^*}(s)$ et $A_{\Gb^*}(t)$ \'etant induit par un 
\'el\'ement conjuguant $s$ en $t$. La proposition \ref{calcul zeta} 
montre que 
\equat
\faisceau(\Gb,(s)) = \coprod_{a \in A_{\Gb^*}(s)} \faisceau(\Gb,(s),a)
\endequat
et que
\equat
\faisceau(\Gb,(s),a)=\{A \in \faisceau(\Gb,(s))~|~\z_A=\o_s(a)\}.
\endequat

\bigskip

\section{Action de $\ZC(\Gb)$ sur $\faisceau(\Gb)$}

\bigskip

Si $s \in \Tb^*$ et $w \in W$ sont tels que $w(s)=s$ et si $z \in \Zb(\Gb)$, 
alors $(t_z^\Gb)^* K_{w,s} \simeq K_{w,s}$. Cela montre en particulier que, 
si $A \in \faisceau(\Gb,(s))$, alors 
\equat\label{translation geometrique}
(t_z^\Gb)^* A \in \faisceau(\Gb,(s)).
\endequat
De plus, si $z \in \Zb(\Gb)^\circ$, alors $(t_z^\Gb)^* A \simeq A$. 
Cela nous d\'efinit donc une action de $\ZC(\Gb)$ sur l'ensemble 
$\faisceau(\Gb)$ ainsi que sur toutes les s\'eries g\'eom\'etriques 
$\faisceau(\Gb,(s))$. 

Le but de cette section est de d\'ecrire 
cette action via le processus d'induction des faisceaux-caract\`eres. 
Nous nous restreindrons aux faisceaux-caract\`eres apparaissant dans 
l'induit de faisceaux-caract\`eres cuspidaux dont le support contient 
des \'el\'ements unipotents. Nous aurons pour cela besoin d'introduire 
des notions d\'evelopp\'ees dans \cite[partie I]{bonnafe action}. 

\bigskip

\soussection{Notations} 
Soit $\Lb$ un \levi de $\Gb$, soit $\Cb$ une classe unipotente et supposons 
que $\Cb$ supporte un syst\`eme local cuspidal $\EC$. Soit $\LC$ un syst\`eme 
local kumm\'erien sur $\Zb(\Lb)^\circ$. Posons $\FC=\LC \boxtimes \EC$. 
C'est un syst\`eme local cuspidal sur $\Sigb=\Zb(\Lb)^\circ \Cb$. 
Rappelons que l'existence d'un syst\`eme local cuspidal support\'e 
par $\Cb$ implique que $N_\Gb(\Lb)$ stabilise $\Cb$ et $\EC$ et donc 
stabilise $\Sigb$ et $\FC$ (voir \cite[th\'eor\`eme 9.2 (b)]{luicc}). 

Soient $\Zb(\Lb)^\circ_\reg = \{z \in \Zb(\Lb)^\circ~|~C_\Gb^\circ(z) = \Lb\}$ 
et $\Sigb_\reg=\Zb(\Lb)_\reg^\circ \Cb$. On note $\LC_\reg$ la restriction de 
$\LC$ \`a $\Zb(\Lb)^\circ_\reg$ et $\FC_\reg$ la restriction de $\FC$ 
\`a $\Sigb_\reg$. On a $\FC_\reg=\LC_\reg \boxtimes \EC$. 
Posons $W_\Gb(\Lb,\LC)=N_\Gb(\Lb,\LC)/\Lb$. 
Fixons maintenant $v \in \Cb$ et notons $\z$ le caract\`ere irr\'eductible 
du groupe $A_\Lb(v)$ associ\'e \`a sa repr\'esentation sur $\EC_v$ 
par monodromie. 
Fixons aussi $x \in X(\Zb(\Lb)^\circ)$ et $n \in \NM^*$, premier \`a $p$, 
tels que $\LC=x^* \XC_n$. 

Si $w \in W_\Gb(\Lb)$, alors $w \in W_\Gb(\Lb,\LC)$ si et seulement si 
$w(x/n)-x/n \in X(\Zb(\Lb)^\circ)$. Dans ce cas, nous notons $\wha$ la restriction 
de $w(x/n)-x/n$ \`a $\Zb(\Gb) \cap \Zb(\Lb)^\circ$. Alors il est facile de v\'erifier 
que $\wha$ est trivial sur $\Zb(\Gb)^\circ$, \cad que $\wha \in X(\Ker h_\Lb)$. 
En composant avec $\kappa$, nous verrons $\wha$ comme 
un \'el\'ement de $(\Ker h_\Lb)^\wedge$ (voir \ref{X irr}). L'application 
$$\fonction{\omega_\LC}{W_\Gb(\Lb,\LC)}{(\Ker h_\Lb)^\wedge}{w}{\wha}$$ 
est un morphisme de groupes ne d\'ependant que de $\LC$ et non 
pas du choix de $x$ et $n$. Nous noterons $W_\Gb^+(\Lb,\LC)$ 
le noyau de $\omega_\LC$. Par dualit\'e, on obtient une application 
surjective 
$$\omeh_\LC : \Ker h_\Lb \longto 
(W_\Gb(\Lb,\LC)/W_\Gb^+(\Lb,\LC))^\wedge.$$

\bigskip

\soussection{Induction} 
Soit $A=IC(\overline{\Sigb},\FC)[\dim \Sigb]$. C'est un faisceau-caract\`ere 
cuspidal sur $\Lb$. Nous allons rappeler ici la construction du faisceau 
pervers induit de $A$. Pour cela, posons
$$\Ybh=\Gb \times \Sigb_\reg,\qquad\Ybt=\Gb \times_\Lb \Sigb_\reg 
\quad\text{et}\quad \Yb=\bigcup_{g \in \Gb} g \Sigb_\reg g^{-1}.$$
Ici, $\Gb \times_\Lb \Sigb_\reg$ d\'esigne le quotient de $\Gb \times \Sigb_\reg$ 
par l'action diagonale de $\Lb$ par translation \`a droite sur le premier 
facteur et par conjugaison sur le deuxi\`eme. 
Notons $\a : \Ybh \to \Sigb_\reg$ la deuxi\`eme projection, 
$\b : \Ybh \to \Ybt$ la projection canonique et $\pi : \Ybt \to \Yb$ 
l'application telle que $\pi \circ \b(g,x)=gxg^{-1}$. 

Alors $\Yb$ est une sous-vari\'et\'e localement ferm\'ee lisse de $\Gb$, 
$\pi$ est un rev\^etement \'etale galoisien de groupe $W_\Gb(\Lb,\Sigb)$ 
et $\a$ et $\b$ sont des morphismes de vari\'et\'es (voir \cite[\SEC 3]{luicc} et 
\cite[\SEC 1]{bonnafe action}). Il existe alors un syst\`eme local 
$\FCt_\reg$ sur $\Ybt$ tel que $\a^* \FC_\reg \simeq \b^* \FCt_\reg$. 
Par cons\'equent, $\pi$ \'etant un rev\^etement \'etale, $\pi_* \FCt_\reg$ est un 
syst\`eme local sur $\Yb$. On a alors \cite[proposition 4.5]{luicc}
\equat\label{definition induit A}
\Ind_\Lb^\Gb A = IC(\overline{\Yb},\pi_* \FCt_\reg)[\dim \Yb].
\endequat

\bigskip

\soussection{Alg\`ebre d'endomorphismes} 
Comme dans \cite[\SEC 3]{bonnafe action}, posons 
$W_\Gb(\Lb,v)=N_\Gb(\Lb,v) /C_\Lb^\circ(v)$ et 
$W_\Gb^\circ(\Lb,v)=(N_\Gb(\Lb) \cap C_\Gb^\circ(v)) /C_\Lb^\circ(v)$. 
L'introduction du syst\`eme local $\LC$ nous conduit \`a consid\'erer 
les sous-groupes 
$W_\Gb(\Lb,v,\LC)=N_\Gb(\Lb,v,\LC)/C_\Lb^\circ(v)$ et 
$W_\Gb^\circ(\Lb,v,\LC)=(N_\Gb(\Lb,\LC) \cap C_\Gb^\circ(v))/C_\Lb^\circ(v)$. 
Rappelons que $W_\Gb(\Lb,v)=W_\Gb^\circ(\Lb,v) \times A_\Lb(v)$ 
(voir \cite[5.3]{bonnafe action}). D'autre part, $A_\Lb(v)$ stabilise $\LC$, 
donc est contenu dans $W_\Gb(\Lb,v,\LC)$. Par suite
\equat\label{produit w}
W_\Gb(\Lb,v,\LC) = W_\Gb^\circ(\Lb,v,\LC) \times A_\Lb(v).
\endequat
Puisque $W_\Gb^\circ(\Lb,v) \simeq W_\Gb(\Lb)$, on a 
$W_\Gb^\circ(\Lb,v,\LC) \simeq W_\Gb(\Lb,\LC)$. Par suite, si 
$w \in W_\Gb(\Lb,\LC)$, nous noterons $\wdo$ un repr\'esentant de $w$ choisi dans 
$N_\Gb(\Lb,\LC) \cap C_\Gb^\circ(v)$. 

Notons $\AC$ l'alg\`ebre d'endomorphismes du faisceau pervers semi-simple 
$\Ind_\Lb^\Gb A$. Nous allons construire, en suivant 
\cite[proposition 3.5 et th\'eor\`eme 9.2]{luicc} et 
\cite[\SEC 5 et 6]{bonnafe action}, un isomorphisme entre 
$\AC$ et l'alg\`ebre de groupe de $W_\Gb(\Lb,\LC) \simeq W_\Gb^\circ(\Lb,\LC)$. 

Soit $w \in W_\Gb(\Lb,\LC)$. Soit $\t_w'$ 
l'isomorphisme $\EC \mapright{\sim} (\INT \wdo)^* \EC$ qui induit l'identit\'e 
sur $\EC_v$. Soit $\s_w$ l'isomorphisme $\LC \mapright{\sim} (\INT \wdo)^* \LC$ 
qui induit l'identit\'e sur $\LC_1$. Alors $\th_w' = \s_w \boxtimes \t_w'$ 
est un isomorphisme $\FC \mapright{\sim} (\INT \wdo)^* \FC$. 
Il lui correspond \cite[\SEC 3.4]{luicc} un automorphisme 
$\Th_w'$ de $\Ind_\Lb^\Gb A$. 

Dans \cite[corollaires 6.2 et 6.7]{bonnafe action} 
a \'et\'e construit un caract\`ere lin\'eaire 
$\g_{\Lb,v,\z}^\Gb : W_\Gb^\circ(\Lb,v) \to \{1,-1\}$. Nous noterons 
$\t_w=\g_{\Lb,v,\z}(w)\t_w'$, 
$\th_w=\g_{\Lb,v,\z}^\Gb(w)\th_w'$ et $\Th_w=\g_{\Lb,v,\z}^\Gb(w) \Th_w'$. 
Compte tenu des choix qui ont \'et\'e faits, il est facile de 
v\'erifier que l'application
$$\fonction{\Th}{W_\Gb(\Lb,\LC)}{\AC}{w}{\Th_w}$$ 
est un isomorphisme d'alg\`ebres. 

Si $\eta$ est un caract\`ere irr\'eductible de $W_\Gb(\Lb,\LC)$, 
nous noterons $K_\eta$ le faisceau-caract\`ere, composant irr\'eductible de 
$\Ind_\Lb^\Gb A$, associ\'e au caract\`ere $\eta$ gr\^ace 
\`a l'isomorphisme de Lusztig 
$\Th$. En d'autres termes, $\Hom(K_\eta, \Ind_\Lb^\Gb A)$ 
est un $\AC$-module admettant $\eta$ comme caract\`ere. 
Notons que, si $z \in \Zb(\Gb) \cap \Zb(\Lb)^\circ$, 
alors $(t_z^\Lb)^* A \simeq A$, donc 
$(t_z^\Gb)^* \Ind_\Lb^\Gb A \simeq \Ind_\Lb^\Gb A$. Par suite, 
on obtient une action de $\Ker h_\Lb$ sur 
$\{K_\eta~|~\eta \in \Irr W_\Gb(\Lb,\LC)\}$, \cad 
une action de $\Ker h_\Lb$ sur $\Irr W_\Gb(\Lb,\LC)$.  
Nous pouvons maintenant \'enoncer et d\'emontrer le r\'esultat 
principal de cette section, \`a savoir la description de cette action. 

\bigskip

\theoreme{translation z}
{\it Soient $\eta \in \Irr W_\Gb^\circ(\Lb,v,\LC)$ et soit $z \in \Ker h_\Lb$. 
Notons $\zdo$ un repr\'esentant de $z$ dans $\Zb(\Gb) \cap \Zb(\Lb)^\circ$. Alors 
$$(t_\zdo^\Gb)^* K_\eta \simeq K_{\eta \omeh_\LC(z)}.$$}

\bigskip

\proof Nous reprenons ici les constructions de \cite[\SEC 5 et 6]{bonnafe action}. 
Mais nous devons tenir compte du syst\`eme local $\LC$ 
(qui, dans \cite{bonnafe action}, 
\'etait suppos\'e constant). Il nous faut donc les modifier l\'eg\`erement 
en introduisant le rev\^etement \'etale de $\Zb(\Lb)^\circ$ qui trivialise 
le syst\`eme local $\LC$. 

Soit 
$$\Zb_{x,n}=\{(z,\x) \in \Zb(\Lb)^\circ~|~x(z)=\x^n\}.$$
Alors la permi\`ere projection $p_1 : \Zb_{x,n} \to \Zb(\Lb)^\circ$ 
est un rev\^etement \'etale de groupe $\mub_n(\FM)$~: le syst\`eme 
local $\LC$ est celui associ\'e \`a ce rev\^etement et 
au caract\`ere $\kappa$ de $\mub_n(\FM)$. Notons $\Zb_{x,n,\reg}$ 
l'image inverse de $\Zb(\Lb)^\circ_\reg$ dans $\Zb_{x,n}$. Posons 
$$\Ybt_{x,n}' = \Gb/C_\Lb^\circ(v) \times \Zb_{x,n,\reg}.$$
Comme dans \cite[\SEC 3.A]{bonnafe action}, posons 
$$\Ybt'=\Gb/C_\Lb^\circ(v) \times \Zb(\Lb)^\circ_\reg$$ 
et notons $\fti : \Ybt' \to \Ybt$ l'application naturelle. 
Soit $\fti^+ : \Ybt_{x,n}' \to \Ybt$ l'application d\'efinie 
par composition de $\Id_{\Gb/C_\Lb^\circ(v)} \times p_1$.  
Notons $\pi^+ =\pi \circ \fti^+$. 

Le groupe $W_\Gb(\Lb,v,\LC) \times \mub_n(\FM)$ agit \`a droite sur 
$\Ybt_{x,n}'$ de la fa\c{c}on suivante~: si 
$(w,\x) \in W_\Gb(\Lb,v,\LC) \times \mub_n(\FM)$ et si 
$(gC_\Lb^\circ(v),z,\x') \in \Ybt_{x,n}'$, on note 
$\l_w = w(x/n)-x/n \in X(\Zb(\Lb)^\circ)$ et on pose
$$(gC_\Lb^\circ(v),z,\x')\cdot(w,\x)=(g\wdo C_\Lb^\circ(v),\wdo^{-1}z\wdo,
\l_w(z)\x\x').$$
Il est alors facile de v\'erifier, en utilisant \cite[\SEC 3.A]{bonnafe action}, 
que $W_\Gb(\Lb,v,\LC) \times \mub_n(\FM)$ agit librement sur $\Ybt_{x,n}'$. 

Fixons maintenant $z \in \Ker h_\Lb$ et notons $\zdo$ un repr\'esentant 
de $z$ dans $\Zb(\Lb)^\circ$. Soit $z_1 \in \Zb(\Lb)^\circ$ tel que 
$z_1^n = z$. Alors l'application 
$$\fonction{t_1}{\Ybt_{x,n}'}{\Ybt_{x,n}'}{(g C_\Lb^\circ(v),
z',\x)}{(gC_\Lb^\circ(v),zz',x(z_1)\x)}$$ 
est un automorphisme de vari\'et\'es v\'erifiant 
$\pi^+ \circ t_1=t_\zdo^\Gb \circ \pi^+$. 
Le th\'eor\`eme \ref{translation z} d\'ecoule imm\'ediatement 
de ces remarques et du fait que, si $w \in W_\Gb^\circ(\Lb,v,\LC)$, alors 
$$t_1^{-1}(t_1((gC_\Lb^\circ(v),z',\x) \cdot w)) = 
(g \wdo C_\Lb^\circ(v),\wdo^{-1} z' \wdo, \l_w(z_1^{-1} w^{-1} z_1 w) \x)$$
et $\kappa \circ \l_w(z_1^{-1} w^{-1} z_1 w)= \o_\LC(w)(z)=\omeh_\LC(z)(w)$.\fin

\bigskip

\soussection{Un analogue du th\'eor\`eme \ref{harish} pour 
les faisceaux-caract\`eres\label{soussection translation}} 
Posons
$$W_\Gb'(\Lb,\LC)=\{(w,\mu) \in W_\Gb(\Lb,\LC) \times \ZC(\Gb)^\we~|~
\omeh_\LC(w) = \Res_{\Ker h_\Lb}^{\ZC(\Gb)} \mu\}.$$
Alors l'application $W_\Gb^+(\Lb,\LC) \to W_\Gb'(\Lb,\LC)$, 
$w \mapsto (w,1)$ est un morphisme de groupe injectif qui nous permettra 
d'identifier $W_\Gb^+(\Lb,\LC)$ avec un sous-groupe de $W_\Gb'(\Lb,\LC)$. 
De plus, l'application $\o_\LC' : W_\Gb'(\Lb,\LC) \to \ZC(\Gb)^\we$, 
$(w,\mu) \mapsto \mu$ est un morphisme de groupes dont le noyau est 
$W_\Gb^+(\Lb,\LC)$. Nous noterons 
$\omeh_\LC' : \ZC(\Gb) \to (W_\Gb'(\Lb,\LC)/W_\Gb^+(\Lb,\LC))^\we$ 
le morphisme dual.

D'autre part, l'application $\ZC(\Lb)^\we \to W_\Gb'(\Lb,\LC)$, 
$\mu \mapsto (1,\mu \circ h_\Lb)$ est aussi un morphisme injectif 
de groupes~: nous identifierons $\ZC(\Lb)^\we$ avec le sous-groupe 
correspondant de $W_\Gb'(\Lb,\LC)$. Alors l'application 
$W_\Gb'(\Lb,\LC) \to W_\Gb(\Lb,\LC)$, $(w,\mu) \mapsto w$ est 
surjective et son noyau est $\ZC(\Lb)^\we$. En d'autres termes,
\equat\label{iso drole}
W_\Gb'(\Lb,\LC)/\ZC(\Lb)^\we \simeq W_\Gb(\Lb,\LC).
\endequat
Pour finir, notons que $\ZC(\Lb)^\we$ est central dans $W_\Gb'(\Lb,\LC)$. 
R\'esumons tout ceci dans le diagramme suivant, dans lequel toutes 
les suites verticales ou horizontales sont exactes et tous les carr\'es 
sont commutatifs~:
$$\diagram
&&&1\dto && 1 \dto && \\
&&& W_\Gb^+(\Lb,\LC) \rrdouble \ddto&& W_\Gb^+(\Lb,\LC) \ddto && \\
&&&&&&&\\
1 \rto & \ZC(\Lb)^\we \rrto \dddouble && W_\Gb'(\Lb,\LC) \rrto \ddto^{\DS{\o_\LC'}}
&& W_\Gb(\Lb,\LC) \ddto^{\DS{\o_\LC}} \rto & 1 \\
&&&&&&&\\
1 \rto & \ZC(\Lb)^\we \rrto && \ZC(\Gb)^\we \rrto \dto
&& (\Ker h_\Lb)^\we \rto \dto& 1 \\
&&& 1 && 1 && \\
\enddiagram$$

\medskip

Soit $\eta$ un caract\`ere irr\'eductible de $W_\Gb'(\Lb,\LC)$. 
Notons $z_\eta$ l'\'el\'ement de $\ZC(\Lb)$ (vu comme un caract\`ere 
lin\'eaire de $\ZC(\Lb)^\we$) par lequel $\ZC(\Lb)^\we$ agit 
sur la repr\'esentation de $W_\Gb'(\Lb,\LC)$ associ\'ee \`a $\eta$. 
Notons $\zti_\eta$ un \'el\'ement de $\ZC(\Gb)$ tel que 
$h_\Lb(\zti_\eta)=z_\eta$. Nous verrons $\zti_\eta$ comme un \car lin\'eaire de 
$\ZC(\Gb)^\we$, \cad comme un caract\`ere lin\'eaire de $W_\Gb'(\Lb,\LC)$. 
Posons 
\equat
K_\eta = (t_{\zti_\eta}^\Gb)^* K_{\eta \omeh_\LC'(\zti_\eta)^{-1}}.
\endequat
Remarquons tout d'abord que cette notation a un sens. Premi\`erement, 
$\eta \omeh_\LC'(\zti_\eta^{-1})$ est trivial sur $\ZC(\Lb)^\we$ donc peut 
\^etre vu comme un \car irr\'eductible de $W_\Gb(\Lb,\LC)$ d'apr\`es 
\ref{iso drole}. D'autre part, en vertu du th\'eor\`eme \ref{translation z}, 
le membre de droite ne d\'epend pas du choix de $\zti_\eta$. On 
a donc montr\'e le r\'esultat suivant~:

\bigskip

\proposition{decomposition faisceau}
{\it L'application $\eta \mapsto K_\eta$ est une bijection entre 
$\Irr W_\Gb'(\Lb,\LC)$ et l'ensemble des composantes irr\'eductibles 
de $\Ind_\Lb^\Gb (\DS{\mathop{\oplus}_{z \in \ZC(\Lb)}} (t_z^\Lb)^* A)$. 
De plus
$$\Ind_\Lb^\Gb \Bigl(\mathop{\oplus}_{z \in \ZC(\Lb)} (t_z^\Lb)^* A \Bigr)
= \mathop{\oplus}_{\eta \in \Irr W_\Gb'(\Lb,\LC)} K_\eta^{\oplus \eta(1)}.$$}

\bigskip

La proposition \ref{decomposition faisceau} sugg\`ere fortement 
qu'il doit exister un isomorphisme naturel entre l'alg\`ebre 
d'endomorphismes du faisceau pervers 
$\Ind_\Lb^\Gb (\DS{\mathop{\oplus}_{z \in \ZC(\Lb)}} (t_z^\Lb)^* A)$ 
et l'alg\`ere de groupes de $W_\Gb'(\Lb,\LC)$. Nous allons ici le construire. 
Pour cela, posons $A'=\DS{\mathop{\oplus}_{z \in \ZC(\Lb)}} (t_z^\Lb)^* A$ 
et notons $\LC'$ le syst\`eme local 
$\DS{\mathop{\oplus}_{z \in \ZC(\Lb)}} (t_z^\Lb)^* \LC$ sur $\Zb(\Lb)$. 

Soit $(w,\mu) \in W_\Gb'(\Lb,\LC)$. On a construit un isomorphisme 
$\s_w : \LC \mapright{\sim} (\INT \wdo)^* \LC$. Si 
$z \in \Zb(\Gb) \cap \Zb(\Lb)^\circ$, la preuve du th\'eor\`eme 
\ref{translation z} montre que l'action de $(\s_w)_z$ sur $\LC_z$ 
est $\omeh_\LC(w)(z) \Id_{\LC_z} = \mu(z) \Id_{\LC_z}$. 
Par suite, il existe un unique isomorphisme 
$\s_{w,\mu} : \LC' \mapright{\sim} (\INT \wdo)^* \LC'$ tel que, 
pour tout $z \in \Zb(\Gb)$, on ait $(\s_w)_z=\mu(z) \Id_{\LC_z'}$. 
Par tensorisation avec $\t_w$, on obtient un isomorphisme 
$\th_{w,\mu} : \FC' \mapright{\sim} (\INT \wdo)^* \FC'$, o\`u 
$\FC'=\LC' \boxtimes \EC$.

\`A travers le diagramme d'induction, $\th_{w,\mu}$ induit un 
automorphisme $\Th_{w,\mu}$ du faisceau pervers $\Ind_\Lb^\Gb A'$. Si 
on note $\AC'$ l'alg\`ebre d'endomorphisme de $\Ind_\Lb^\Gb A'$, alors~:

\bigskip

\theoreme{algebre}
{\it L'application
$$\fonction{\Th}{\qlb W_\Gb'(\Lb,\LC)}{\AC'}{(w,\mu)}{\Th_{w,\mu}}$$
est un isomorphisme d'alg\`ebres. Si $\eta$ est un caract\`ere irr\'eductible 
de $W_\Gb'(\Lb,\LC)$, alors la composante irr\'eductible de 
$\Ind_\Lb^\Gb A'$ associ\'ee \`a $\eta$ \`a travers l'isomorphisme $\Th$ 
est $K_\eta$.}

\bigskip

\section{Fonctions caract\'eristiques}

\medskip

Nous allons maintenant introduire dans ce chapitre l'isog\'enie $F$. 
Un faisceau pervers $A$ sur $\Gb$ sera dit {\it $F$-stable} 
s'il est isomorphe \`a $F^* A$. Nous noterons $\faisceau(\Gb)^F$ 
l'ensemble des (classes d'isomorphie de) faisceaux-caract\`eres 
$F$-stables. Si $A$ est un faisceau pervers $F$-stable sur $\Gb$ 
et si $\ph : A \mapright{\sim} F^* A$ est un isomorphisme, nous noterons 
$\XC_{A,\ph} : \Gb^F \to \qlb$ la {\it fonction caract\'eristique} 
de $A$, d\'efinie par 
$$\XC_{A,\ph}(g)=\sum_{i \in \ZM} (-1)^i \Tr(\ph_x,\HC_x^i A)$$
pour tout $g \in \Gb^F$. Bien s\^ur, $\XC_{A,\ph}$ d\'epend de $\ph$. 
Cependant, si $A$ est irr\'eductible, alors $\ph$ est bien d\'etermin\'ee 
\`a un scalaire pr\`es. En cons\'equence, la fonction $\XC_{A,\ph}$ est bien 
d\'etermin\'ee par $A$ \`a un scalaire pr\`es. 

\bigskip

\soussection{Cas classique} 
Reprenons les notations de la section pr\'ec\'edente ($\Lb$, $\LC$, $\EC$,\dots). 
Supposons donc maintenant que $\Lb$ est $F$-stable, que $\Tb$ est $F$-stable, 
que $F(v)=v$ et que $F^* \FC \simeq \FC$. Fixons un isomorphisme 
$\ph : F^* \FC \simeq \FC$. Cet isomorphisme s'\'etend en un isomorphisme 
$\ph^\# : F^* A \mapright{\sim} A$. 
Soit $g_w$ un \'el\'ement de $\Gb$ tel que $g_w^{-1} F(g_w)=\wdo^{-1}$. 
Posons 
$$\Lb_w=\lexp{g_w}{\Lb},\quad v_w =\lexp{g_w}{v},\quad \Cb_w =\lexp{g_w}{\Cb}, 
\quad \Sigb_w=\lexp{g_w}{\Sigb},$$
$$ \LC_w=(\ad g_w^{-1})^* \LC, \quad
\EC_w=(\ad g_w^{-1})^* \EC, \quad\FC_w=(\ad g_w^{-1})^* \FC\quad\text{et}\quad 
A_w=(\ad g_w^{-1})^* A.$$
Alors $\FC_w = \LC_w \boxtimes \EC_w$. Alors $\Lb_w$, $v_w$, $\Cb_w$, $\Sigb_w$, 
$\LC_w$, $\EC_w$, $\FC_w$ et $A_w$ sont $F$-stables et, suivant la construction 
de \cite[\SEC 9.3]{lugf}, 
on obtient un isomorphisme $\ph_w : F^*\FC_w \mapright{\sim} \FC_w$. 
Il s'\'etend en un isomorphisme $\ph_w^\# : F^* A_w \mapright{\sim} A_w$. 

Fixons maintenant un caract\`ere $F$-stable $\eta$ de $W_\Gb(\Lb,\LC)$. 
On note $\phi$ l'automorphisme de $W_\Gb(\Lb,\LC)$ induit par $F$. 
On choisit une extension $\tilde{\eta}$ de $\eta$ au produit semi-direct 
$W_\Gb(\Lb,\LC) \rtimes <\phi>$. Ce choix  
d'une extension (ainsi que celui de $\ph$) d\'etermine 
un isomorphisme $\ph_\etat : F^* K_\eta \mapright{\sim} K_\eta$. 
Il r\'esulte de \cite[partie II, 10.4.5 and 10.6.1]{lucs} et 
\cite[proposition 9.2]{lugf} que~:

\bigskip

\Theoreme{Lusztig}{calcul fc}
{\it Supposons $p$ presque bon pour $\Gb$ et $q$ assez grand. 
Avec les notations pr\'ec\'edentes, on a 
$$\XC_{K_\eta,\ph_\etat} = \frac{1}{|W_\Gb(\Lb,\LC)|} 
\sum_{w \in W_\Gb(\Lb,\LC)} \etat(w\phi) R_{\Lb_w}^\Gb \XC_{A_w,\ph_w}.$$}

\bigskip

\soussection{Translation par ${\boldsymbol{\ZC(\Gb)}}$} 
Nous allons \'etudier ici le comportement des fonctions 
caract\'eristiques vis-\`a-vis de la translation par un \'el\'ement du centre. 
Cela sera fait en termes du param\'etrage de la proposition 
\ref{decomposition faisceau}. 

Le syst\`eme local $\LC$ \'etant $F$-stable, il en est de m\^eme du 
syst\`eme local $\LC'$ sur $\Zb(\Lb)$. De m\^eme, le syst\`eme local 
$\FC'$ est $F$-stable. On fixe un isomorphisme $\ph' : F^* \FC' \simeq \FC'$ 
\'etendant $\ph$. 

Soit $(w,\mu) \in W_\Gb'(\Lb,\LC)$. Dans la sous-section 
\ref{soussection translation}, nous avons construit un isomorphisme 
$\th_{w,\mu} : \FC' \mapright{\sim} (\INT \wdo)^* \FC'$. 
Reprenons les notations de la pr\'ec\'edente sous-section 
et posons $\LC_w'=(\INT g_w^{-1})^* \LC'$,  
$\FC_w'=(\INT g_w^{-1})^* \FC'$ et 
$A_w'=(\INT g_w^{-1})^* A'$. En suivant encore 
\cite[\SEC 9.3]{lugf}, on obtient un isomorphisme 
$\ph_{w,\mu}' : F^* \FC_w' \mapright{\sim} \FC_w'$. 
Cet isomorphisme s'\'etend en un isomorphisme 
$\ph_{w,\mu}^{\prime\#} : F^* A_w' \mapright{\sim} A_w'$. 

Fixons maintenant un caract\`ere irr\'eductible $\eta$ de 
$W_\Gb'(\Lb,\LC)$. Alors $K_\eta$ est $F$-stable si et seulement 
si $\eta$ est $F$-stable. Notons $\etat$ une extension de $\eta$ 
au produit semi-direct $W_\Gb'(\Lb,\LC) \rtimes <\phi>$, o\`u 
$\phi$ est l'automorphisme de $W_\Gb'(\Lb,\LC)$ induit par $F$. 
Comme dans le cas classique, le choix de $\etat$ d\'etermine 
un isomorphisme $\ph_\etat : F^* K_\eta \mapright{\sim} K_\eta$. 
Le th\'eor\`eme suivant d\'ecoule presque imm\'ediatement du th\'eor\`eme 
de Lusztig pr\'ec\'edent.

\bigskip

\theoreme{ouf !}
{\it Supposons $p$ presque bon pour $\Gb$ et $q$ assez grand. On a
$$\XC_{K_\eta,\ph_\etat} = \frac{1}{|W_\Gb'(\Lb,\LC)|} 
\sum_{(w,\mu) \in W_\Gb'(\Lb,\LC)} \etat((w,\mu) \phi) 
R_{\Lb_w}^\Gb \XC_{A_w',\ph_{w,\mu}^{\prime \#}}.$$}

\bigskip

\proof Notons $z=z_\eta \in \ZC(\Lb)$. Puisque $\eta$ est $F$-stable, 
on a $z \in \ZC(\Lb)^F$. L'alg\`ebre d'endomorphisme du faisceau pervers  
$\Ind_\Lb^\Gb (t_z^\Lb)^* A$ s'identifie, via la construction pr\'ec\'edente, 
\`a la sous-alg\`ebre de $\qlb W_\Gb'(\Lb,\LC)$ \'egale \`a 
$\qlb W_\Gb'(\Lb,\LC) e_z$, o\`u $e_z$ est l'idempotent central 
$\frac{1}{|\ZC(\Lb)|} \sum_{\t \in \ZC(\Lb)^\we} \t(z)^{-1} \t$. 
Il suffit alors d'appliquer le th\'eor\`eme de Lusztig en remarquant 
que la restriction de $\XC_{A_w',\ph_{w,\mu}^{\prime \#}}$ \`a 
$z \overline{\Sigma}_w^F$ est \'egale \`a 
$$\frac{1}{|\ZC(\Lb)|} \sum_{\t \in \ZC(\Lb)^\we} \t(z)^{-1} 
\XC_{A_w',\ph_{w,\mu\t}^{\prime \#}}.~\SS{\blacksquare}$$

\bigskip

\section{\'El\'ements unipotents r\'eguliers\label{section unip reg}}

\medskip

\begin{quotation}
\noindent{\bf Hypoth\`ese~:} {\it Dor\'enavant, et ce 
jusqu'\`a la fin de cet article, nous supposerons que 
$p$ est bon pour $\Gb$.}
\end{quotation}

\bigskip

Nous nous int\'eressons ici aux faisceaux-caract\`eres apparaissant 
dans l'induit, \`a partir d'un \levi $\Lb$ de $\Gb$, 
de faisceaux-caract\`eres cuspidaux dont 
le support rencontre $\Zb(\Lb)\UCB_\reg^\Lb$. On rappelle que, 
puisque $p$ est suppos\'e bon pour $\Gb$, le groupe 
$A_\Lb(u_\Lb)$ est isomorphe \`a $\ZC(\Lb)$. Si $\z \in \ZC(\Lb)^\we$, 
nous noterons $\EC_\z$ le syst\`eme local $\Lb$-\'equivariant sur $\UCB_\reg^\Lb$ 
tel que l'action de $\ZC(\Lb)$ sur la fibre en $u_\Lb \in \UCB_\reg^\Lb$ 
se fasse par le caract\`ere $\z$. Si $z \in \ZC(\Lb)$, on 
notera $\zdo$ un repr\'esentant de $z$ dans $\Zb(\Lb)$. En d'autres 
termes, $z=\zdo \Zb(\Lb)^\circ$.

\bigskip

\soussection{Cuspidalit\'e} 
Fixons un syst\`eme local kumm\'erien 
$\LC$ sur $\Zb(\Gb)^\circ$, un \'el\'ement $z \in \Zb(\Gb)$ et 
un caract\`ere lin\'eaire $\z$ de $\ZC(\Gb)$. 
Posons $\FC=((t_z^\Gb)^* \LC) \boxtimes \EC_\z$. 

\bigskip

\proposition{cuspidal noyau}
{\it $\FC$ est un syst\`eme local cuspidal \ssi $\z \in \ZC_\cus^\we(\Gb)$.}

\bigskip

\proof voir \cite[proposition 1.2.2]{bonnafe torsion}.\fin

\bigskip

Soit $\faisceau_\reg^\cus(\Gb)$ l'ensemble des (classes d'isomorphie de)
faisceaux-caract\`eres cuspidaux dont le support rencontre 
$\Zb(\Gb)\UCB_\reg^\Gb$. Soit $\Cus_\reg(\Gb)$ un ensemble de repr\'esentants 
(modulo l'action naturelle de $\Gb$) des 
triplets $(s,a,\t)$, o\`u $s$ est un \'el\'ement semi-simple de $\Gb^*$, 
$a \in A_{\Gb^*}(s)$ est tel que $\omega_s(a) \in \ZC_\cus^\we(\Gb)$ 
et $\t \in A_{\Gb^*}(s)^\we$. Nous allons construire une bijection entre 
$\Cus_\reg(\Gb)$ et $\faisceau_\reg^\cus(\Gb)$. 

Soit $(s,a,\t) \in \Cus_\reg(\Gb)$. On peut supposer, et nous le ferons, 
que $s \in \Tb^*$. 
Par construction, l'\'el\'ement $s$ est g\'eom\'etriquement cuspidal et donc 
$\omega_s : A_{\Gb^*}(s) \to \ZC(\Gb)^\we$ est un isomorphisme 
(voir proposition \ref{cuspidal prop} (e)). En particulier, 
$\omeh_s : \ZC(\Gb) \to A_{\Gb^*}(s)^\we$ est aussi 
un isomorphisme. Posons $z=\omeh_s^{-1}(\t)$ et notons $\LC_{s,z}$ la 
restriction de $\LC_s$ \`a $z^{-1}=\zdo^{-1} \Zb(\Gb)^\circ$. Posons maintenant 
$$\FC_{s,a,\t}=\FC_{s,a,\t}^\Gb = \LC_{s,z} \boxtimes \EC_{\omega_s(a)}.$$
Notons que $\LC_{s,z} \simeq (t_\zdo^\Gb)^* \LC_{s,1}$. 
C'est un syst\`eme local $\Gb$-\'equivariant 
irr\'eductible cuspidal sur la vari\'et\'e lisse 
$z^{-1} \UCB_\reg^\Gb=\zdo^{-1}\Zb(\Gb)^\circ \UCB_\reg^\Gb$. Posons 
$$A_{s,a,\t}=A_{s,a,\t}^\Gb=IC(\overline{z^{-1} \UCB_\reg^\Gb},
\FC_{s,a,\t})[\dim \Zb(\Gb)^\circ \UCB_\reg^\Gb].$$
C'est un faisceau pervers $\Gb$-\'equivariant irr\'eductible sur $\Gb$. 

\bigskip

\lemme{debut Asat}
{\it Soit $(s,a,\t) \in \Cus_\reg(\Gb)$. 
\begin{itemize}
\itemth{a} $\z_{A_{s,a,\t}} = \omega_s(a)$.

\itemth{b} $A_{s,a,\t}$ est l'extension par z\'ero du syst\`eme 
local $\EC_{s,a,\t}$.

\itemth{c} $A_{s,a,\t} \in \faisceau(\Gb,(s))$.

\itemth{d} $A_{s,a,\t}$ est cuspidal. 
\end{itemize}}

\bigskip

\proof Soit $(s,a,\t) \in \Cus_\reg(\Gb)$. Rappelons que 
l'existence de $(s,a,\t)$ implique que toutes les composantes quasi-simples de 
$\Gb$ sont de type $A$. Posons $\z=\omega_s(a)$ et $z=\omeh_s^{-1}(\t)$. 

\medskip

(a) d\'ecoule du fait que 
$\ZC(\Gb)$ agit sur $\EC_\z$ via le caract\`ere lin\'eaire $\z$. 

\medskip

(b) Soit $x$ un \'el\'ement de l'adh\'erence de $z\UCB_\reg^\Gb$ 
n'appartenant pas \`a $z\UCB_\reg^\Gb$. Puisque toutes les somposantes 
quasi-simples de $\Gb$ sont de type $A$, ceci implique qu'il existe 
un \levi $\Lb$ propre de $\Gb$ contenant $x$. En particulier, 
$\Zb(\Lb)^\circ \incl C_\Gb^\circ(x)$. Donc $\Zb(\Gb) \cap \Zb(\Lb)^\circ$ 
agit trivialement sur la fibre en $x$ de $A_{s,a,\t}$. Par suite, 
si cette fibre est non nulle, 
la restriction de $\z$ \`a $\Ker h_\Lb$ est triviale, ce qui contredit 
la proposition \ref{cuspidal noyau}. Donc $(A_{s,a,\t})_x=0$. 

\medskip

(c) Rappelons que toutes les composantes quasi-simples de $\Gb$ 
sont de type $A$. La classification des \'el\'ements quasi-isol\'es 
r\'eguliers \cite{bonnafe quasi} montre que, quitte \`a conjuguer 
le triplet $(s,a,\t)$, on peut supposer que $s \in \Tb^*$ et que 
$a$ est un \'el\'ement de Coxeter standard de $W$. Nous allons 
calculer $K_{a,s}$ dans ces conditions. 

Tous les \'el\'ements de $\Bb a \Bb$ sont r\'eguliers et $\Bb a \Bb$ rencontre 
toutes les classes de conjugaison d'\'el\'ements r\'eguliers 
\cite[remarque 8.8]{steinberg}. La proposition \ref{calcul zeta} et 
l'argument du (b) montre que le support de $A_{s,a,\t}$ est 
contenu dans $\Zb(\Gb).\UCB_\reg^\Gb$. Notons 
$i : \Zb(\Gb)\UCB_\reg^\Gb \injto \Gb$, 
$\Ybt$ l'image inverse de $\Zb(\Gb)\UCB_\reg^\Gb$ dans $\Ybt_a$, 
$\g : \Ybt \to \Zb(\Gb)\UCB_\reg^\Gb$ la restriction de $\g_a$ et 
$\LCt$ la restriction de $\LCt_{a,s}$ \`a $\Ybt$. On a alors
$$K_{a,s} = i_! R\g_! \LCt.\leqno{(*)}$$
Fixons un \'el\'ement unipotent r\'egulier 
$x \in \Bb a \cap a \Bb^- \cap \UCB_\reg$.  
Soit 
$$\fonction{\ph}{\Gb/\Zb(\Gb) \times \Zb(\Gb)}{\Ybt}{(g\Zb(\Gb),t)}{(tgxg^{-1},
g\Bb).}$$
Nous allons montrer que $\ph$ est un 
morphisme de vari\'et\'e bijectif purement ins\'eparable. 
On a construit une action de $\Gb$ sur $\Ybt_a$. 
Le groupe $\Zb(\Gb)$ agit aussi sur $\Ybt_a$ par translation 
de la premi\`ere coordonn\'ee. Cette action conserve $\Ybt$. 
Cela munit $\Ybt$ d'une action de $\Gb \times \Zb(\Gb)$. On remarque 
alors que $\ph$ est $\Gb \times \Zb(\Gb)$ \'equivariant. Il suffit 
donc de montrer que $\ph$ est bijectif. En effet, 
cela montre que la vari\'et\'e $\Ybt$ est une orbite sous 
l'action d'un groupe alg\'ebrique, donc elle est lisse, donc elle est 
normale et un morphisme bijectif entre vari\'et\'es normales 
est purement ins\'eparable \cite[th\'eor\`eme 18.2]{borel}. 

Soient $(g\Zb(\Gb),t)$ et $(g'\Zb(\Gb),t')$ deux \'el\'ements de 
$\Gb/\Zb(\Gb) \times \Zb(\Gb)$ 
ayant m\^eme image par $\ph$. 
Alors la partie semi-simple de $tgxg^{-1}$ co\"\i ncide avec celle 
de $t'g'xg^{\prime -1}$, \cad $t=t'$. On a par cons\'equent 
$g^{-1}g' \in \Bb \cap C_\Gb(x)$. Mais, d'apr\`es 
\cite[corollaire 10.3]{bonnafe action}, 
$\Ub \cap C_\Gb(x)=1$ donc la partie unipotente de $g^{-1}g'$ est \'egale \`a $1$. 
Donc, puisque $x$ est un unipotent r\'egulier, on en d\'eduit que 
$g\Zb(\Gb)=g'\Zb(\Gb)$, ce qui montre l'injectivit\'e de $\ph$. 

Montrons maintenant la surjectivit\'e de $\ph$. Soit $(g,h\Bb) \in \Ybt$. 
Alors, par d\'efinition, il existe $t \in \Zb(\Gb)$ et $y \in \Gb$ tels que 
$g=t yxy^{-1}$. Posons $k=h^{-1} y$. Alors, 
par hypoth\`ese, $kxk^{-1} \in \Bb a \Bb \cap \UCB_\reg^\Gb$. 
Donc, d'apr\`es \cite[corollaire 10.3]{bonnafe action} et 
\cite[th\'eor\`eme 1.4]{steinberg}, 
il existe $b \in \Bb$ tel que $kxk^{-1}=bxb^{-1}$. En d'autres termes, 
$yxy^{-1}=hbxb^{-1}h^{-1}$. Donc $(g,h\Bb) = \ph(hb\Zb(\Gb),t)$.

Posons $\LCt' = \ph^* \LCt$. Puisque $\ph$ est bijectif et purement 
ins\'eparable, on a $\ph_* \LCt'=\LCt$ et donc 
$$K_{a,s}=i_! R\g_!'\LCt',$$ 
o\`u $\g' : \Gb/\Zb(\Gb) \times \Zb(\Gb) \to \Zb(\Gb)\UCB_\reg^\Gb$, 
$(g\Zb(\Gb),t) \mapsto tgxg^{-1}$. 
Alors $\LCt'=\ECt \boxtimes(\oplus_{z \in \ZC(\Gb)} \LC_{s,z})$, 
o\`u $\ECt$ est un syst\`eme local $\Gb$-\'equivariant irr\'eductible 
sur $\Gb/\Zb(\Gb)$. 
L'action de $\ZC(\Gb)$ sur $\ECt$ \'etant donn\'ee par $\z$, $\ECt'$ est
l'unique syst\`eme local sur $\Gb/\Zb(\Gb)$ sur lequel $\ZC(\Gb)$ 
agit par $\z$. Pour $z \in \ZC(\Gb)$, notons $i_z : z \UCB_\reg^\Gb \injto \Gb$. 
Notons $\d : \Gb/\Zb(\Gb) \to \UCB_\reg^\Gb$, $g\Zb(\Gb) \mapsto gvg^{-1}$ 
et $\EC=R\d_! \ECt$. Alors 
$$K_{a,s} = \mathop{\oplus}_{z \in \ZC(\Gb)} i_{z !} (\EC \boxtimes \LC_{s,z}).$$
Il nous reste \`a calculer $\EC$. En d\'ecomposant $\d$ en la suite 
de morphismes $\Gb/\Zb(\Gb) \mapright{\d'} \Gb/C_\Gb(u) \to \UCB_\reg$, 
on est ramen\'e au calcul de $R\d_!' \ECt$. Mais, puisque $\d'$ est 
un morphisme lisse dont les fibres sont isomorphes 
\`a $C_\Ub(u)$, qui est, comme vari\'et\'e alg\'ebrique, 
un espace affine de dimension $\rang_\sem(\Gb)$, on a 
$$\EC=\EC_\z[-2\rang_\sem(\Gb)].$$
On en d\'eduit que 
\equat
K_{a,s} = \mathop{\oplus}_{\t \in A_{\Gb^*}(s)^\we} A_{s,a,\t}[m],
\endequat
pour un $m \in \ZM$ que je n'ai pas envie de calculer. Cela montre (c).

\medskip

(d) est \'evident.\fin

\bigskip

\proposition{bijection faisceaux}
{\it L'application $\Cus_\reg(\Gb) \to \faisceau_\reg^\cus(\Gb)$, 
$(s,a,\t) \mapsto A_{s,a,\t}$ est bijective. De plus, 
$A_{s,a,\t} \in \faisceau(\Gb,(s))$.}

\bigskip

\proof Soit $(s,a,\t) \in \Cus_\reg(\Gb)$. Alors 
$\EC_{s,a,\t}$ est un syst\`eme local cuspidal sur $z \UCB_\reg^\we$ 
(voir \cite[d\'efinition 2.4]{luicc} et proposition \ref{cuspidal noyau}). 
Par suite, $A_{s,a,\t}$ est un faisceau-caract\`ere cuspidal 
\cite[\SEC 7]{lucs} dont le support rencontre $\Zb(\Gb)\UCB_\reg^\Gb$. 
Cela montre que l'application est bien d\'efinie.

Montrons qu'elle est surjective. Soit $A \in \faisceau_\reg^\cus(\Gb)$. 
D'apr\`es \cite[\SEC 7]{lucs}, il existe $z \in \ZC(\Gb)$, 
un syst\`eme local $\LC$ sur $z=\zdo\Zb(\Gb)^\circ$ et $\z \in \ZC(\Gb)_\cus^\we$ 
tel que $A=IC(\overline{z\UC_\reg^\Gb},\LC \boxtimes \FC_\z)[\rang \Gb]$. 
En effet, puisque $\Gb$ est de type $A$, les \'el\'ements semi-simples 
isol\'es sont centraux et les \'el\'ements unipotents distingu\'es 
sont r\'eguliers. Notons que $\z_A=\z$. 

Soit $s$ tel que $A \in \faisceau(\Gb,(s))$ (voir \ref{partition faisceaux}). 
Il existe $w \in W_\Gb(\Tb)$ tel que $w(s)=s$ et $A$ est une composante 
irr\'eductible de $\lexp{p}{H}^i(K_{w,s})$. Notons $a$ la classe de $w$ dans 
$A_{\Gb^*}(s)$. D'apr\`es la proposition \ref{calcul zeta}, 
on a $\z=\omega_s(a)$. Posons maintenant $\t=\omeh_s(z) \in A_{\Gb^*}(s)^\we$. 
Puisque $A \in \faisceau(\Gb,(s))$, la restriction de $\LC_s$ \`a $\Zb(\Gb)^\circ$ 
est \'egale \`a $t_{\zdo}^* \LC$, o\`u $t_\zdo : \Gb \to \Gb$, $g \mapsto \zdo g$ 
est la translation par $\zdo$. Donc $\LC_{s,z} = \LC$, ce qui montre 
que $A=A_{s,a,\t}$. 

Montrons maintenant qu'elle est injective. Soient $(s,a,\t)$ et 
$(s',a',\t')$ deux \'el\'ements de $\Cus_\reg(\Gb)$ tels que 
$A_{s,a,\t} \simeq A_{s',a',\t'}$. D'apr\`es \ref{partition faisceaux} 
et le lemme \ref{debut Asat} (c), $s$ et $s'$ sont conjugu\'es 
sous $\Gb^*$. On peut donc supposer qu'ils sont \'egaux. 
De plus, $\z_{A_{s,a,\t}}=\z_{A_{s,a',\t'}}$ donc, d'apr\`es la 
proposition \ref{calcul zeta}, on a $\omega_s(a)=\omega_s(a')$. 
Donc $a=a'$ car $\omega_s$ est injectif. Pour finir, les supports 
de $A_{s,a,\t}$ et $A_{s,a,\t'}$ sont \'egaux, ce qui implique 
que $\omeh_s^{-1}(\t)=\omeh_s^{-1}(\t')$, d'o\`u l'on d\'eduit 
que $\t=\t'$.\fin

\bigskip

\soussection{Induction} 
Fixons un \'el\'ement semi-simple $s \in \Tb^*$ et un \'el\'ement 
$a \in A_{\Gb^*}(s)$. Posons $\Lb_{s,a}=C_\Gb((\Tb^a)^\circ)$. 
Alors, d'apr\`es la proposition \ref{cuspidal las}, 
$\o_s(a) \in \ZC_\cus^\we(\Lb_{s,a})$. Notons 
$\faisceau_\reg(\Gb,(s),a)$ l'ensemble des faisceaux-caract\`eres 
apparaissant dans $\Ind_{\Lb_{s,a}}^\Gb 
(\DS{\mathop{\oplus}_{\t \in A_{\Lb_{s,a}^*}(s)^\we}} A_{s,a,\t}^\Lb)$. 
Il est \`a noter que $\faisceau_\reg(\Gb,(s),a)$ est 
contenu dans $\faisceau(\Gb,(s),a)$. Alors, 
d'apr\`es la proposition \ref{decomposition faisceau}, on a une bijection 
\equat\label{bij w s a}
\Irr W_\Gb'(\Lb_{s,a},\LC_{s,a}) \longmapright{\sim} \faisceau_\reg(\Gb,(s),a),
\endequat 
o\`u $\LC_{s,a}$ d\'esigne la restriction de $\LC_s$ \`a $\Zb(\Lb_{s,a})^\circ$. 
Il nous reste \`a d\'eterminer le groupe $W_\Gb'(\Lb_{s,a},\LC_{s,a})$. 
C'est fait dans la proposition suivante (comparer avec la proposition 
\ref{W}). 

\bigskip

\proposition{W faisceau}
{\it Si $s \in \Tb^*$ et $a \in A_{\Gb^*}(s)$, alors 
$W_\Gb'(\Lb_{s,a},\LC_{s,a})$ est canoniquement isomorphe 
\`a $W(s)^a \simeq A_{\Gb^*}(s) \ltimes W^\circ(s)^a$. 
A travers cet isomorphisme, on a $A_{\Lb_{s,a}^*}(s) \simeq \ZC(\Lb_{s,a})^\we$ 
et $W(s)^a/A_{\Lb_{s,a}^*}(s) \simeq W_\Gb(\Lb_{s,a},\LC_{s,a})$.}

\bigskip

\proof Soit $w \in W(s)^a$. Alors $w$ normalise $\Lb_{s,a}$ et $\LC_{s,a}$. 
Posons 
$$\fonction{\aleph}{W(s)^a}{W_\Gb'(\Lb_{s,a},\LC_{s,a})}{w}{(\wti,
\omega_s(\wba)),}$$
o\`u $\wti$ d\'esigne la classe de $w$ dans $W_\Gb(\Lb_{s,a})$ et $\wba$ 
la classe de $w$ dans $A_{\Gb^*}(s)$. Nous allons montrer que $\aleph$ est 
un isomorphisme. 

Le fait que l'application $\aleph$ est bien d\'efinie d\'ecoule imm\'ediatement 
de la construction de $\omega_s$ et de la d\'efinition de 
$W_\Gb'(\Lb_{s,a},\LC_{s,a})$. Soit $w \in W(s)^a$ tel que 
$\aleph(w)=(1,1)$. Alors $w \in W_{\Lb_{s,a}}(s)^a=A_{\Lb_{s,a}^*}(s)$ car 
$s$ est g\'eom\'etriquement cuspidal dans $\Lb_{s,a}^*$ donc r\'egulier 
(voir les propositions \ref{cuspidal las} et \ref{semisimple cuspidal prop}). 
Mais alors $\omega_s(w) = 1$ et donc $w=1$ d'apr\`es 
l'injectivit\'e de $\omega_s$. Cela montre l'injectivit\'e de $\aleph$. 

Il nous reste \`a montrer la surjectivit\'e. Soit 
$(w,\mu) \in W_\Gb'(\Lb_{s,a},\LC_{s,a})$. Soit $\wdo$ un repr\'esentant 
de $w$ dans $W$. Il r\'esulte de \ref{partition faisceaux} (appliqu\'e au 
groupe $\Lb_{s,a}$) que l'on peut supposer que $\wdo \in W(s)$. 
Il est alors facile de v\'erifier que 
\equat\label{omega traduction}
\o_{\LC_{s,a}}(w)=\Res_{\Ker h_\Lb}^{\ZC(\Gb)} \o_s(\wdo).
\endequat
Donc $\Res_{\Ker h_\Lb}^{\ZC(\Gb)} \mu=\Res_{\Ker h_\Lb}^{\ZC(\Gb)} \o_s(\wdo)$. 
Puisque $A_{\Lb_{s,a}^*}(s) \simeq \ZC(\Lb_{s,a})^\we$, il existe 
$a \in A_{\Lb_{s,a}^*}(s)$ tel que $\mu = \o_s(\wdo)\o_s(a)=\o_s(\wdo a)$. 
Quitte \`a changer de repr\'esentant de $w$ dans $W(s)^a$, 
on peut donc supposer que $\o_s(\wdo)=\mu$.  
Notons $b$ la classe de $\wdo$ dans $A_{\Gb^*}(s)$. 
Alors $w'=b^{-1} \wdo \in W^\circ(s)$ et normalise $\Lb_{s,a}$. 
Donc, d'apr\`es le corollaire \ref{ah ah ah} (e), $aw'=w'a$ et 
donc $a\wdo = \wdo a$. Par suite, $\wdo \in W(s)^a$ et $\aleph(\wdo)=(w,\mu)$, 
ce qui montre la surjectivit\'e de $\aleph$.\fin

\bigskip

\section{Fonctions caract\'eristiques}

\medskip

Le but de cette section est de calculer les fonctions caract\'eristiques des 
faisceaux-caract\`eres $F$-stables appartenant \`a $\faisceau_\reg(\Gb,(s),a)$, 
o\`u $s$ est un \'el\'ement semi-simple de $\Gb^*$ et $a \in A_{\Gb^*}(s)$. 
Tout d'abord, remarquons que, si $\faisceau(\Gb,(s))^F \not= \vide$, alors 
$(s)$ est une classe de conjugaison $F^*$-stable. Par cons\'equent, on 
peut supposer que $s \in \Gb^{*F^*}$, ce qui sera fait par la suite. 
D'autre part, si $\faisceau(\Gb,(s),a)^F \not= \vide$, alors on a, 
pour tout $A \in \faisceau(\Gb,(s),a)^F$, $\z_A = \o_s(a) \in (\ZC(\Gb)^\we)^F$. 
Puisque $\o_s$ est injectif, cela implique que $a \in A_{\Gb^*}(s)^{F^*}$. 

Par cons\'equent, nous ferons l'hypoth\`ese suivante~:

\bigskip

\begin{quotation}
\noindent{\bf Hypoth\`ese : } {\it Dans cette section, nous fixons 
un \'el\'ement semi-simple $s \in \Gb^{*F^*}$ et un \'el\'ement 
$a \in A_{\Gb^*}(s)^{F^*}$. Nous reprenons les notations des chapitres 
pr\'ec\'edents ($\Tb_1^*$, $\Tb_1$,...) et nous supposons que 
$\Tb=\Tb_1$ et $\Tb^*=\Tb_1^*$.}
\end{quotation}

\bigskip

Nous aurons d'autre part besoin de la notation suivante. Si 
$\z \in H^1(F,\ZC(\Gb))^\we$, nous posons 
$$\GC(\Gb,\z) = \eta_\Gb q^{-\frac{1}{2} \rang_\sem(\Gb)} 
\sum_{z \in H^1(F,\ZC(\Gb))} 
\z(z)^{-1} \sum_{t \in \Tb_0^F/\Zb(\Gb)^F} \psi_z(tut^{-1}).$$
Remarquons que $\GC(\Gb,\z)$ est \'egal \`a 
$\eta_\Gb q^{-\frac{1}{2} \rang_\sem(\Gb)}$ 
fois le scalaire not\'e $\s_{\z^{-1}}$ dans \cite[\SEC 2]{DLM2}. En particulier 
\cite[proposition 2.5]{DLM2}~:
\equat\label{racine 4}
\text{\it Si $\z \in \ZC_\cus^\we(\Gb)$ est $F$-stable, alors $\GC(\Gb,\z)$ est 
une racine quatri\`eme de l'unit\'e.}
\endequat

\bigskip

\remarque{calcul racine de l'unite}
Le calcul de $\GC(\Gb,\z)$ lorsque $\z \in \ZC_\cus^\we(\Gb)$ sera 
effectu\'e dans l'appendice B.\fin

\bigskip

\soussection{Cas cuspidal} 
Nous allons rappeler dans cette sous-section comment 
les transform\'ees de Fourier de caract\`eres semi-simples 
sont reli\'ees aux fonctions caract\'eristiques 
de faisceaux-caract\`eres cuspidaux dont le support 
rencontre $\Zb(\Gb) \UCB_\reg^\Gb$. 

\smallskip

Supposons dans cette sous-section, et uniquement dans cette sous-section, 
que $\omega_s(a) \in \ZC_\cus^\we(\Gb)$. Dans ce cas, on a
$$\faisceau(\Gb,(s),a)^F=\{A_{s,a,\t}~|~\t \in H^1(F^*,A_{\Gb^*}(s))^\we\}.$$
Posons $A_{s,a}=\oplus_{\t \in (A_{\Gb^*}(s))^\we} A_{s,a,\t}$. 
Alors $A_{s,a}$ est $F$-stable et il existe un unique isomorphisme 
$\ph_{s,a} : F^* A_{s,a} \mapright{\sim} A_{s,a}$ tel que, 
pour tout $z \in \Zb(\Gb)^F$, on ait 
$$(\ph_{s,a})_{zu}=\sha(z) q^{\frac{1}{2} \rang_\sem(\Gb)} \Id_{(A_{s,a})_{zu}}.$$
Si $\t \in H^1(F^*,A_{\Gb^*}(s))^\we$, notons $\ph_{s,a,\t}$ 
la restriction de $\ph_{s,a}$ en un isomorphisme 
$F^* A_{s,a,\t} \mapright{\sim} A_{s,a,\t}$. 
Il r\'esulte de \cite[th\'eor\`eme 6.2.2]{bonnafe torsion} que 
\equat\label{fonction caracteristique cuspidale}
\XC_{A_{s,a,\t},\ph_{s,a,\t}} = \GC(\Gb,\o_s(a)) \rhoh_{s,a,\t}.
\endequat

\bigskip

\soussection{Le r\'esultat} 
Revenons au cas g\'en\'eral, \cad ne supposons plus que 
$\o_s(a) \in \ZC_\cus^\we(\Gb)$. On note alors $A_{s,a}$ le 
faisceau pervers $\oplus_{\t \in (A_{\Lb_{s,a}^*}(s))^\we} A_{s,a,\t}^{\Lb_{s,a}}$ 
sur $\Lb_{s,a}$ et on note $\ph_{s,a} : F^* A_{s,a} \mapright{\sim} A_{s,a}$ 
l'isomorphisme tel que, pour tout $z \in \Zb(\Lb_{s,a})^F$, on ait 
$$(\ph_{s,a})_{zu_{\Lb_{s,a}}}=\sha(z) q^{\frac{1}{2} \rang_\sem(\Lb_{s,a})} 
\Id_{(A_{s,a})_{zu_{\Lb_{s,a}}}}.$$
Si $\t \in H^1(F^*,A_{\Lb_{s,a}^*}(s))^\we$, on note $\ph_{s,a,\t}$ la restriction 
de $\ph_{s,a}$ en un isomorphisme 
$F^* A_{s,a,\t}^{\Lb_{s,a}} \mapright{\sim} A_{s,a,\t}^{\Lb_{s,a}}$. 
D'apr\`es \ref{fonction caracteristique cuspidale} appliqu\'ee 
\`a $\Lb_{s,a}$, on a 
\equat\label{fonction caracteristique cuspidale L}
\XC_{A_{s,a,\t},\ph_{s,a,\t}} = \GC(\Lb_{s,a},\o_{\Lb_{s,a},s}(a)) 
\rhoh_{s,a,\t}^{\Lb_{s,a}}.
\endequat

Soit $\eta$ un caract\`ere irr\'eductible $F^*$-stable de $W(s)^a$. 
Soit $\etat$ une extension de $\eta$ \`a $W(s)^a \rtimes <\phi_1>$. 
Le choix de cette extension fixe un isomorphisme 
$\ph_{s,a,\etat} : F^* K(s,a)_\eta \mapright{\sim} K(s,a)_\eta$. 

\bigskip

\theoreme{PRINCIPAL}
{\it Supposons que $p$ est bon pour $\Gb$ et que $q$ est assez grand. 
Soit $\eta$ un caract\`ere irr\'eductible $F^*$-stable de $W(s)^a$. 
Soit $\etat$ une extension de $\eta$ \`a $W(s)^a \rtimes <\phi_1>$. 
Alors
$$\XC_{K(s,a)_\eta,\ph_{s,a,\etat}} = 
\GC(\Lb_{s,a},\o_{\Lb_{s,a},s}(a)) \RC(s,a)_\etat.$$}

\bigskip

\proof Soit $w \in W(s)$. Notons $\wdo$ un repr\'esentant de $w$ dans 
$N_\Gb(\Tb_1)$. Fixons $g_w \in \Gb$ tel que $g_w^{-1} F(g_w) = \wdo$. 
On note $\wba$ la classe de $w$ dans $W(s)/A_{\Lb_{s,a}^*}(s)$. 
On peut choisir la famille $(g_w)_{w \in W(s)}$ de sorte que, 
si $w \in W(s)$ et $b \in A_{\Lb_{s,a}^*}(s)$, on ait 
$\lexp{g_w}{\Lb_{s,a}}=\lexp{g_{wb}}{\Lb_{s,a}}$. 
Par suite, on peut poser
$$\Lb_\wba=\lexp{g_w}{\Lb_{s,a}}=\Lb_{s,a,w}.$$
Posons $\Sigb'=\Zb(\Lb_{s,a}). \Cb$. Reprenons maintenant les 
hypoth\`eses et notations du th\'eor\`eme 
\ref{ouf !} et supposons de plus que $\Cb$ soit la classe unipotente 
r\'eguli\`ere de $\Lb_{s,a}$, que $\LC'$ soit \'egal \`a la restriction 
de $\LC_s$ \`a $\Zb(\Lb_{s,a})$, que $\EC=\EC_{\o_{\Lb_{s,a},s}(a)}$, 
que $\FC' = \LC' \boxtimes \EC$ et que $\ph'=\ph_{s,a}$. Posons alors 
$$\Cb_\wba =\lexp{g_w}{\Cb}, \quad \Sigb_\wba'=\lexp{g_w}{\Sigb'},$$
$$ \LC_\wba'=(\ad g_w^{-1})^* \LC', \quad
\EC_\wba=(\ad g_w^{-1})^* \EC, \quad\FC_\wba'=(\ad g_w^{-1})^* \FC'
\quad\text{et}\quad A_\wba=(\ad g_w^{-1})^* A_{s,a}.$$
L'\'el\'ement $w$ d\'efinit quant \`a lui un isomorphisme 
$\ph_w : F^* A_{\wba} \mapright{\sim} A_\wba$ qui, lui, d\'epend de 
$w$ et pas seulement de $\wba$. 

Compte tenu du th\'eor\`eme \ref{ouf !}, il nous reste 
\`a montrer que 
$$\XC_{A_\wba,\ph_w} = \GC(\Lb_{s,a},\o_{\Lb_{s,a},s}(a)) \rho_{s_w,a}^{\Lb_\wba}.
\leqno{(*)}$$
Rappelons que $s_w = g_w s g_w^{-1}$. On pose 
$\Tb_w=\lexp{g_w}{\Tb_1}$ et $\LC_{s,w} = (\ad g_w^{-1})^* \LC_s$. 
En fait, $\LC_\wba'$ est la restriction de $\LC_{s,w}$ \`a $\Zb(\Lb_\wba)$. 
Par cons\'equent, on a, pour tous $z \in \Zb(\Lb_\wba)^F$ et 
$x \in \Sigb_\wba^{\prime F}$, 
$$\XC_{A_\wba,\ph_w} (zx)= \sha_w(z) \XC_{A_\wba,\ph_w}(x).$$
D'autre part, l'action de $\ZC(\Lb_\wba)$ par conjugaison 
sur $\FC_\wba'$ montre que, pour prouver $(*)$, il suffit de montrer que 
$$\XC_{A_\wba,\ph_w} (u_{\Lb_\wba})= 
\GC(\Lb_{s,a},\o_{\Lb_{s,a},s}(a)) \rho_{s_w,a}^{\Lb_\wba}(u_{\Lb_\wba}).
\leqno{(**)}$$
Mais, d'apr\`es \cite[proposition 2.5]{DLM2}, on a 
$\GC(\Lb_{s,a},\o_{\Lb_{s,a},s}(a))=\GC(\Lb_\wba,\o_{\Lb_\wba,s}(a))$. 
Donc il suffit de montrer, d'apr\`es 
\ref{fonction caracteristique cuspidale L}, que 
$\XC_{A_\wba,\ph_w}(u_{\Lb_\wba})= q^{\frac{1}{2} \rang_\sem \Lb_\wba}$, 
ce qui d\'ecoule de \cite[th\'eor\`eme 15.10]{bonnafe action}.\fin

\bigskip

Nous allons nous int\'eresser maintenant aux faisceaux-caract\`eres 
dont le support rencontre $\Zb(\Gb)\UCB_\reg^\Gb$. 
Soit $K \in \faisceau(\Gb,(s))$ dont le support rencontre 
$\Zb(\Gb)\UCB_\reg^\Gb$. Notons $a$ l'unique \'el\'ement de 
$A_{\Gb^*}(s)$ tel que $\z_K=\o_s(a)$ (voir la proposition 
\ref{calcul zeta}). Alors 
$K$ est une composante irr\'eductible de 
$\Ind_{\Lb_{s,a}}^\Gb A_{s,a}$. Notons $\eta_K$ le caract\`ere 
irr\'eductible de $W(s)^a$ correspondant. 

\bigskip

\lemme{caractere trivial}
{\it Soit $z \in \ZC(\Gb)$ et supposons que le support de $K$ contienne 
$z \UCB_\reg$. Alors $\eta_K=\omeh_s(z)$.}

\bigskip

\proof Quitte \`a translater $K$ par un \'el\'ement de $\Zb(\Gb)$ (c'est-\`a-dire, 
d'apr\`es le th\'eor\`eme \ref{translation z}, 
\`a multiplier $\eta_K$ par un caract\`ere lin\'eaire de $A_{\Gb^*}(s)$), 
on peut supposer que le support de $K$ contient $\UCB_\reg^\Gb$. 
En utilisant les constructions de \cite[partie I]{bonnafe action}, 
on s'aper\c{c}oit, en utilisant \cite[corollaire 6.7]{bonnafe action}, 
que l'on peut supposer que $\LC'=\qlb$. Dans ce cas, il 
d\'ecoule de la d\'efinition de $\Th$ et 
\cite[corollaire 6.2]{bonnafe action} que $\eta_K=1$.\fin

\bigskip

\corollaire{regulier caracteristique}
{\it Soit $K \in \faisceau(\Gb,(s))^F$ dont le support rencontre 
$z \UCB_\reg^\we$, pour un $z \in \ZC(\Gb)^F$. Posons $\t=\omeh_s^1(z)$ et 
notons $\taut$ l'extension de $\t$ \`a $A_{\Gb^*}(s) \rtimes <\phi_1>$ 
qui est triviale sur $<\phi_1>$. Alors 
$$\XC_{K,\ph_{s,a,\taut}} = \GC(\Lb_{s,a},\o_{\Lb_{s,a},s}(a)) 
\rhoh_{s_\a,a,\t}.$$}

\bigskip

\proof Cela r\'esulte imm\'ediatement du th\'eor\`eme \ref{PRINCIPAL}, 
du lemme \ref{caractere trivial} et de la proposition 
\ref{proprietes rhoh} (e).\fin

\bigskip

\newpage

{\Large \part{Groupes de type ${\boldsymbol{A}}$\label{chapitre a}}}

\bigskip

\begin{quotation}
\noindent{\bf Hypoth\`ese : }{ \it Dor\'enavant, et ce jusqu'\`a la fin de 
cet article, nous supposerons que toutes les composantes quasi-simples de 
$\Gb$ sont de type $A$. Nous supposerons aussi que $\d=1$, 
\cad que $F : \Gbt \to \Gbt$ est un endomorphisme de Frobenius.}
\end{quotation}

\bigskip

Rappelons que l'hypoth\`ese ci-dessus implique que la formule de Mackey a 
lieu dans $\Gb$ (voir th\'eor\`eme \ref{theo mackey}). 
Le but de ce chapitre est d'obtenir un param\'etrage des caract\`eres 
irr\'eductibles de $\Gb^F$ et de montrer, lorsque $q$ est assez grand, 
que la conjecture de Lusztig a lieu. 

\bigskip

\section{Description de $\Cent(\Gb^F,[s])$}

\medskip

\soussection{Structure de ${\boldsymbol{W(s) \rtimes <\phi_1>}}$} 
Notons ici $\Phi_{(1)}$,\dots, $\Phi_{(r)}$ les composantes 
irr\'eductibles de $\Phi_s$ et posons $\Phi_{(i)}^+=\Phi_s^+ \cap \Phi_{(i)}$ et 
$\D_{(i)}=\D \cap \Phi_{(i)}$. Notons $W_{(i)}$ le groupe de Weyl 
du syst\`eme de racines $\Phi_{(i)}$. Alors 
$W^\circ(s)=W_{(1)} \times \dots \times W_{(r)}$. 

Chaque $W_{(i)}$ est isomorphe \`a un groupe sym\'etrique et $A_{\Gb^*}(s)$ 
permute les $W_{(i)}$. Il est possible de choisir une famille d'isomorphismes 
$W_{(i)} \simeq \SG_{n_i}$ (o\`u $n_i$ est un entier naturel $\ge 2$ et 
$\sum_{i=1}^r n_i = n$) telle que $A_{\Gb^*}(s)$ agisse seulement par permutation 
des composantes. Une fois un tel choix d'isomorphismes effectu\'e, il existe 
$w_s \in W^\circ(s)$ tel que $w_s F^*$ (ou $w_s \phi_1$) agisse 
sur $W^\circ(s)$ seulement par permutation des composantes. 
Il n'est pas d\'efini de mani\`ere unique car le centre de $W^\circ(s)$ 
n'est pas forc\'ement trivial. Cependant, cette non unicit\'e ne peut se produire 
que lorsqu'il existe des $i$ tels que $n_i=2$. Si $n_i=2$, nous 
supposerons que la composante de $w_s$ dans $W_{(i)}$ est \'egale \`a $1$. 
Cela d\'efinit $w_s$ de fa\c{c}on unique. S'il 
est n\'ecessaire de pr\'eciser, nous le noterons $w_{\Gb,s}$. 

\bigskip

\lemme{ags commute avec sigmas}
{\it $w_s$ commute avec les \'el\'ements de $A_{\Gb^*}(s)$.}

\bigskip

\proof Soit $a \in A_{\Gb^*}(s)$. Alors $F^*(a) \in A_{\Gb^*}(s)$. 
D'autre part, $[w_s \phi_1,a]$ agit sur $W^\circ(s)$ seulement 
par permutation des composantes. Or, 
$[w_s\phi_1,a]=w_s ~\lexp{F^*(a)}{w_s} \in W^\circ(s)$. Par suite, 
$w_s ~\lexp{F^*(a)}{w_s}$ est central dans $W^\circ(s)$ et donc 
\'egal \`a $1$ compte tenu du choix pr\'ecis fait pour $w_s$.\fin

\bigskip

\soussection{Fonctions absolument cuspidales} 
Nous rappelons la description de l'espace des fonctions absolument 
cuspidales dans notre cas \cite[th\'eor\`eme 4.3.3]{bonnafe torsion}~:

\bigskip

\theoreme{theo abs cus}
{\it Si $a \in A_{\Gb^*}(s)^{F^*}$, alors 
$$\Cus(\Gb^F,[s],a)=\begin{cases} 
\qlb \rhodot_{s,a} & \text{si } \o_s(a) \in \ZC_\cus^\we(\Gb),\\
0 & \text{sinon.}
\end{cases}$$}

\bigskip

\corollaire{surjectivite R}
{\it Soit $a\in A_{\Gb^*}(s)^{F^*}$. Alors l'application 
$\RC[s,a] : \Bigl(\Cent(W^\circ(s)^a\phi_1)\Bigr)^{A_{\Gb^*}(s)^{F^*}} \to 
\Cent(\Gb^F,[s],a)$ est une isom\'etrie bijective (pour les produits scalaires 
$\langle , \rangle_{s,a}$ et $\langle , \rangle_{\Gb^F}$).}

\bigskip

\proof Le fait que $\RC[s,a]$ est une isom\'etrie a \'et\'e 
montr\'e dans la proposition \ref{isometrie R}. 
Il nous reste \`a montrer la surjectivit\'e de $\RC[s,a]$. Puisque 
la formule de Mackey a lieu dans $\Gb$ (voir th\'eor\`eme \ref{theo mackey}), 
on a 
$$\Cent(\Gb^F,[s],a) = 
\mathop{\oplus}_{\Lb \in \LC_{s,a}} R_\Lb^\Gb(\Cus(\Lb^F,[s],a)),$$
o\`u $\LC_{s,a}$ est l'ensemble des sous-groupes de Levi $F$-stables 
de $\Gb$ dont un dual $\Lb^*$ contient $s$ et tels que $a \in A_{\Lb^*}(s)^{F^*}$. 
Mais, d'apr\`es le th\'eor\`eme \ref{theo abs cus}, on a 
$$\Cent(\Gb^F,[s],a) = \mathop{\oplus}_{\Lb \in \LC_{s,a,\text{cus}}} 
\qlb R_\Lb^\Gb\rhodot_{s,a}^\Lb,$$
o\`u $\LC_{s,a,\text{cus}}$ est l'ensemble des $\Lb \in \LC_{s,a}$ tels que 
$\o_{\Lb,s}(a) \in \ZC_\cus^\we(\Lb)$. Il suffit alors de montrer que, 
si $\Lb^*$ est un \levi $F^*$-stable de $\Gb^*$ contenant $s$ 
et v\'erifiant que $a \in A_{\Lb^*}(s)$ et $\o_{\Lb,s}(a) \in \ZC_\cus^\we(\Lb)$, 
alors $\Lb^*$ est conjugu\'e sous $C_{\Gb^*}^\circ(s)$ \`a un 
$\Lb_{s,a,w}$ pour un $w \in W^\circ(s)^a$. 
Reprenons les notations du \SEC\ref{soussection wl} ($W_\Lb^\circ(s)$, 
$W_\Lb(s)$, $A_{\Lb^*}(s)$ et $w_\Lb$). Alors $w_\Lb$ commute avec $a$ 
et est le type du tore $C_{\Lb^*}^\circ(s)$ (voir corollaire \ref{ah ah ah} (a)). 
Par suite, on peut supposer que $C_{\Lb^*}^\circ(s)=\Tb_{w_\Lb}^*$. 
Puisque $\o_{\Lb,s}(a) \in \ZC_\cus^\we(\Lb)$, on a 
$\Lb^* = C_{\Gb^*}(((\Tb_{w_\Lb}^*)^a)^\circ)$, ce qui montre 
que $\Lb^*$ est conjugu\'e sous $C_{\Gb^*}^\circ(s)$ \`a $\Lb_{s,a,w_\Lb}$.\fin

\bigskip

\soussection{Une autre isom\'etrie} 
Fixons encore $a \in A_{\Gb^*}(s)^{F^*}$. L'\'el\'ement $w_s\phi_1$ 
agit sur $W^\circ(s)^a$ par permutation des composantes donc, d'apr\`es 
\ref{iso pia}, on a une isom\'etrie bijective naturelle entre 
$\Cent(W^\circ(s)^a \phi_1)=\Cent(W^\circ(s)^a w_s \phi_1)$ et 
$\Cent(W^\circ(s)^{<a,w_s\phi_1>})$. Cette isom\'etrie commute \`a l'action de 
$A_{\Gb^*}(s)^{F^*}$. De m\^eme, $a$ agissant par permutation des composantes 
de $W^\circ(s)^{w_s F^*}$, on a une isom\'etrie bijective naturelle entre 
$\Cent(W^\circ(s)^{<a,w_s\phi_1>})$ et $\Cent(W^\circ(s)^{w_s F^*} a)$ 
commutant \`a l'action de $A_{\Gb^*}(s)^{F^*}$. Nous noterons 
$$\s_{s,a} : \Cent(W^\circ(s)^{w_sF^*}a) \longto \Cent(W^\circ(s)^a\phi_1)$$
la composition de ces isom\'etries. 
Posons alors $R[s,a] = \RC[s,a] \circ \s_{s,a}$. Notons que, 
si $b \in A_{\Gb^*}(s)^{F^*}$ et si $f \in \Cent(W^\circ(s)^{w_sF^*}a)$, alors, 
d'apr\`es \ref{invariance RC}, on a 
\equat\label{invariance RR}
R[s,a]_{\lexp{b}{f}}=R[s,a]_f.
\endequat
L'application $\s_{s,a}$ se restreint en une isom\'etrie bijective 
toujours not\'ee 
$$\s_{s,a} : \Bigl(\Cent(W^\circ(s)^{w_sF^*}a)\Bigr)^{A_{\Gb^*}(s)^{F^*}}
\longto \Bigl(\Cent(W^\circ(s)^a\phi_1)\Bigr)^{A_{\Gb^*}(s)^{F^*}}$$
\`a condition de munir 
$\bigl(\Cent(W^\circ(s)^{w_sF^*}a)\bigr)^{A_{\Gb^*}(s)^{F^*}}$ du produit 
scalaire $\langle , \rangle_{s,a}' = 
\DS{\frac{1}{|A_{\Gb^*}(s)^{F^*}|}} \langle , \rangle_{W^\circ(s)^{w_sF^*}a}$. 
L'application $R[s,a] : 
\bigl(\Cent(W^\circ(s)^{w_sF^*}a)\bigr)^{A_{\Gb^*}(s)^{F^*}} \to \Cent(\Gb^F,[s],a)$ 
est alors une isom\'etrie bijective (pour les produits scalaires 
$\langle,\rangle_{s,a}'$ et $\langle , \rangle_{\Gb^F}$). 
Faisons l'identification isom\'etrique canonique 
$$\Cent(W(s)^{w_s F^*}) = \mathop{\oplus}_{a \in A_{\Gb^*}(s)^{F^*}}^\perp 
\Bigl(\Cent(W^\circ(s)^{w_sF^*}a)\Bigr)^{A_{\Gb^*}(s)^{F^*}}$$
et, \`a travers cette identification, 
posons 
$$R[s]=\mathop{\oplus}_{a \in A_{\Gb^*}(s)^{F^*}} R[s,a].$$
On a alors~:

\bigskip

\proposition{isometrie RR}
{\it L'application $R[s] : \Cent(W(s)^{w_s F^*}) \longto \Cent(\Gb^F,[s])$ 
est une isom\'etrie bijective.}

\bigskip

Si cela s'av\`ere n\'ecessaire, nous noterons $\s_{s,a}^\Gb$, $R[s,a]^\Gb$ 
et $R[s]^\Gb$ les applications $\s_{s,a}$, $R[s,a]$ et $R[s]$.

\bigskip

\soussection{Quelques propri\'et\'es de l'isom\'etrie ${\boldsymbol{R[s]}}$} 
Nous allons commencer par \'etudier l'action de $H^1(F,\ZC(\Gb))$ 
\`a travers cette isom\'etrie. Si $z \in H^1(F,\ZC(\Gb))$ et si 
$f \in \Cent(W(s)^{w_s F^*})$, alors 
\equat\label{tauzg RC}
\t_z^\Gb R[s]_f = R[s]_{f \omeh_s^0(z)},
\endequat
o\`u $\omeh_s^0(z)$ est vu comme un caract\`ere lin\'eaire de 
$W(s)^{w_s F^*}=W^\circ(s)^{w_s F^*} \rtimes A_{\Gb^*}(s)^{F^*}$. 

Nous allons maintenant \'etudier un cas particulier d'induction de 
Lusztig. Un sous-groupe de Levi $\Lb^*$ de $\Gb^*$ est dit $(s,\Gb^*)$-d\'eploy\'e 
s'il est $F^*$-stable et s'il contient $\Tb_{w_s}^*$. 
Un sous-groupe de Levi $\Lb$ de $\Gb$ est dit $(s,\Gb)$-d\'eploy\'e 
s'il est $F^*$-stable et s'il contient $\Tb_{w_s}$. 
Nous allons ici calculer l'induction de Lusztig 
$R_\Lb^\Gb : \Cent(\Lb^F,[s]) \to \Cent(\Gb^F,[s])$ 
lorsque $\Lb$ est $(s,\Gb)$-d\'eploy\'e en utilisant les isom\'etries 
$R[s]^\Lb$ et $R[s]^\Gb$. Mais avant cela, nous allons \'etudier 
quelques-unes des propri\'et\'es de ces sous-groupes de Levi. 

Soit $\Lb$ un \levi $(s,\Gb)$-d\'eploy\'e. Notons $\Lb^*$ un \levi 
$F^*$-stable de $\Gb^*$ contenant $\Tb_{w_s}^*$ tel que le triplet 
$(\Lb^*,\Tb_{w_s}^*,F^*)$ soit dual de $(\Lb,\Tb_{w_s},F)$. 
Par d\'efinition, $\Lb^*$ est $(s,\Gb^*)$-d\'eploy\'e. 
Comme dans \SEC\ref{soussection wl}, d\'efinissons 
un sous-groupe parabolique standard $W_\Lb^\circ(s)$ de $W^\circ(s)$ 
ainsi qu'un \'el\'ement $w_\Lb$ de $W^\circ(s)$. 
Notons $w_{\Lb,s}$ l'\'el\'ement de $W^\circ_\Lb(s)$ d\'efini 
comme $w_s$. Alors il est possible de choisir $w_\Lb$ et $w_{\Lb,s}$ 
tels que $w_s=w_{\Lb,s}w_\Lb$. Par suite 
$W_\Lb^\circ(s) w_\Lb \phi_1 =W_\Lb^\circ(s) w_s \phi_1$ et l'application 
$R[s]^\Lb$ est une isom\'etrie 
$\Cent(W_\Lb(s)^{w_s F^*}) \mapright{\sim} \Cent(\Lb^F,[s])$. 

\bigskip

\proposition{induit commute}
{\it Soit $\Lb$ un sous-groupe de Levi $(s,\Gb)$-d\'eploy\'e de $\Gb$. 
Alors le diagramme 
$$\diagram 
\Cent(W_\Lb(s)^{w_s F^*}) \rrto^{\DS{R[s]^\Lb}} 
\ddto_{\DS{\Ind_{W_\Lb(s)^{w_s F^*}}^{W(s)^{w_s F^*}}}} && 
\Cent(\Lb^F,[s]) \ddto_{\DS{R_\Lb^\Gb}} \\
&& \\
\Cent(W(s)^{w_s F^*}) \rrto^{\DS{R[s]^\Gb}} 
&& \Cent(\Gb^F,[s]) \\
\enddiagram$$
est commutatif.}

\bigskip

\proof Soit $a \in A_{\Lb^*}(s)^{F^*}$. 
Soit $f$ une fonction centrale sur $W_\Lb^\circ(s)^{w_s F^*} a$ 
invariante par l'action de $A_{\Lb^*}(s)^{F^*}$. Notons $f^\#$ son extension 
par $0$ en une fonction centrale sur $W_\Lb(s)^{w_s F^*}$. Il nous suffit de montrer 
que 
$$R_\Lb^\Gb R[s]_{f^\#}^\Lb = \bigl(R[s]^\Gb \circ 
\Ind_{W_\Lb(s)^{w_s F^*}}^{W(s)^{w_s F^*}}\bigr)( f^\#).$$ 
Posons 
$I = \Ind_{W_\Lb(s)^{w_s F^*}}^{W(s)^{w_s F^*}} f^\#$, 
$g = \Ind_{W_\Lb^\circ(s)^{w_s F^*}a}^{W^\circ(s)^{w_s F^*}a} f$ et notons 
$g^\#$ l'extension de $g$ par $0$ en une fonction (pas forc\'ement centrale) 
sur $W(s)^{w_s F^*}$. Notons tout de m\^eme que $g^\#$ est invariante 
par $W^\circ(s)^{w_s F^*}$-conjugaison. On a alors
$$f^\#=\frac{1}{|A_{\Lb^*}(s)^{F^*}|} 
\Ind_{W_\Lb^\circ(s)^{w_s F^*}a}^{W_\Lb(s)^{w_s F^*}} f$$
et donc 
$$I = \frac{1}{|A_{\Lb^*}(s)^{F^*}|}
\sum_{a \in A_{\Gb^*}(s)^{F^*}} \lexp{a}{g^\#}.$$
Compte tenu de \ref{invariance RR}, on a donc 
$$R[s]_I^\Gb = \frac{|A_{\Gb^*}(s)^{F^*}|}{|A_{\Lb^*}(s)^{F^*}|} 
R[s,a]_g^\Gb.$$
Il nous suffit donc de montrer que 
$$|A_{\Gb^*}(s)^{F^*}| R[s,a]_g^\Gb = |A_{\Lb^*}(s)^{F^*}| R_\Lb^\Gb R[s,a]_f^\Lb.$$
En d'autres termes, nous devons montrer que le diagramme 
$$\diagram 
\Cent(W_\Lb^\circ(s)^{w_s F^*}a) \xto[0,3]^{|A_{\Lb^*}(s)^{F^*}| {R[s,a]^\Lb}} 
\ddto_{\DS{\Ind_{W_\Lb^\circ(s)^{w_sF^*} a}^{W^\circ(s)^{w_s F^*}a}}} &&& 
\Cent(\Lb^F,[s],a) \ddto_{\DS{R_\Lb^\Gb}} \\
&&& \\
\Cent(W^\circ(s)^{w_s F^*}a) \xto[0,3]^{|A_{\Gb^*}(s)^{F^*}|{R[s,a]^\Gb}} 
&&& \Cent(\Gb^F,[s],a), \\
\enddiagram$$
est commutatif. Cela d\'ecoule de la commutativit\'e du diagramme 
\ref{diagramme induction tordue} et de la proposition \ref{induction RRR}.\fin

\bigskip

Nous terminons par une formule exprimant $R[s]_f$ dans un cas particulier. 
Si $\xi$ est un caract\`ere lin\'eaire de $A_{\Gb^*}(s)^{F^*}$, nous verrons 
$\xi$ comme une fonction centrale sur $W(s)^{w_s F^*}$ comme dans 
la formule \ref{tauzg RC}. 

\bigskip

\proposition{rho xi}
{\it Supposons que la conjecture $(\GG)$ a lieu dans $\Gb$. Soit 
$\xi \in (A_{\Gb^*}(s)^{F^*})^\we$. Alors
$$R[s]_\xi = \e_\Gb \e_{C_{\Gb^*}^\circ(s)} \rho_{s,\xi}.$$}

\bigskip

\proof Notons $1_{s,a}$ la fonction sur $W^\circ(s)^{w_s F^*} a$ 
constante et \'egale \`a $1$. Alors 
$\s_{s,a}(1_{s,a})$ est constante et \'egale \`a $1$. Puisque 
$\xi = \sum_{a \in A_{\Gb^*}(s)^{F^*}} \xi(a) 1_{s,a}$, on a 
$$R[s]_\xi = \sum_{a \in A_{\Gb^*}(s)^{F^*}} \xi(a) \RC[s,a]_{\s_{s,a}(1_{s,a})}.$$
Il suffit alors d'utiliser \ref{changement de base semi} et 
la proposition \ref{rhodot R}.\fin

\bigskip

\soussection{Caract\`eres irr\'eductibles de $\Gb^F$\label{caracteres GF}} 
Nous allons montrer ici que, si $\eta$ est un caract\`ere irr\'eductible 
de $W(s)^{w_s F^*}$, alors $\pm R[s]_\eta$ est un caract\`ere irr\'eductible de 
$\Gb^F$. Nous aurons cependant besoin de la conjecture $(\GG)$ ce qui, 
\`a l'heure actuelle, restreint le domaine de validit\'e de ce r\'esultat 
au cas o\`u $q$ est grand. Sachant que c'est une fonction centrale de norme $1$, 
il suffit de montrer que c'est un caract\`ere virtuel de $\Gb^F$. 
Pour cela, nous utiliserons le corollaire \ref{frobenius} 
ainsi que les propositions \ref{induit commute} et \ref{rho xi}. 

\bigskip

\theoreme{irreductibles G}
{\it Supposons que la conjecture $(\GG)$ a lieu dans $\Gb$. 
Soit $\eta \in \Cent(W(s)^{w_s F^*})$. Alors 
$R[s]_\eta \in \pm \EC(\Gb^F,[s])$ si et seulement si 
$\eta \in \pm \Irr W(s)^{w_s F^*}$.}

\bigskip

\proof Nous allons commencer par montrer le lemme suivant~:

\bigskip

\begin{quotation}
\lemme{existence L}
{\it Soit $W_1$ un sous-groupe parabolique $w_s F^*$-stable de $W^\circ(s)$ et 
soit $A_1$ son normalisateur dans $A_{\Gb^*}(s)^{F^*}$. Alors 
il existe un sous-groupe de Levi $(s,\Gb)$-d\'eploy\'e de $\Gb$ tel que 
$W_\Lb(s) = W_1 \rtimes A_1$.}

\bigskip

\proof On peut identifier le groupe 
$W_1 \rtimes A_1$ \`a un sous-groupe $F$-stable 
de $N_\Gb(\Tb_{w_s})/\Tb_{w_s}$. Posons alors 
$$\Lb=C_\Gb((\Tb_{w_s}^{W_1 \rtimes A_1})^\circ).$$
Alors $\Lb$ est $F$-stable et $(s,\Gb)$-d\'eploy\'e. D'autre part, 
d'apr\`es le lemme \ref{centralisateur tore}, 
on a $W_\Lb^\circ(s)=W_1$. Il est de plus clair 
que $W_\Lb(s)$ contient $W_1 \rtimes A_1$. Comme $A_1$ est le normalisateur 
de $W_\Lb^\circ(s)=W_1$ dans $A_{\Gb^*}(s)$, on en d\'eduit que 
$W_\Lb(s)=W_1 \rtimes A_1$.\finl
\end{quotation}

\bigskip

Il r\'esulte du corollaire \ref{frobenius}, du lemme \ref{existence L} et 
de la proposition \ref{rho xi} que, si $\eta \in \ZM\Irr W(s)^{w_s F^*}$, alors 
$R[s]_\eta$ est une combinaison lin\'eaire, \`a coefficients dans $\ZM$, 
de caract\`eres virtuels de la forme $R_\Lb^\Gb \rho$, o\`u $\Lb$ est 
un \levi $(s,\Gb)$-d\'eploy\'e de $\Gb$ et $\rho \in \EC(\Lb^F,[s])$ 
est un caract\`ere semi-simple. En particulier, $R[s]_\eta$ est un caract\`ere 
virtuel. Par suite, si $\eta$ est de plus irr\'eductible, $R[s]_\eta$ est un 
caract\`ere virtuel de norme $1$, donc c'est, au signe pr\`es, un caract\`ere 
irr\'eductible.\fin

\bigskip

Comme cons\'equence directe de la preuve du th\'eor\`eme pr\'ec\'edent, on obtient~:

\bigskip

\corollaire{caracterisation sem}
{\it Supposons que la conjecture $(\GG)$ a lieu dans $\Gb$. 
Notons $\SC$ l'ensemble des couples $(\Lb,\rho)$ 
ou $\Lb$ est un \levi $(s,\Gb)$-d\'eploy\'e de $\Gb$ et 
$\rho \in \EC(\Lb^F,[s])$ est un caract\`ere semi-simple. 
Soit $\g \in \Cent(\Gb^F,[s])$. Alors les assertions 
suivantes sont \'equivalentes~:
\begin{itemize}
\itemth{1} $\g \in \ZM\EC(\Gb^F,[s])$. 

\itemth{2} $\g \in \DS{\mathop{\oplus}_{(\Lb,\rho) \in \SC}} \ZM R_\Lb^\Gb \rho$.

\itemth{3} Pour tout $(\Lb,\rho) \in \SC$, 
$\langle \lexp{*}{R}_\Lb^\Gb \g , \rho \rangle_{\Lb^F} \in \ZM$.
\end{itemize}}

\bigskip

\corollaire{caracterisation reg}
{\it Supposons que la conjecture $(\GG)$ a lieu dans $\Gb$. 
Notons $\RC$ l'ensemble des couples $(\Lb,\chi)$ 
ou $\Lb$ est un \levi $(s,\Gb)$-d\'eploy\'e de $\Gb$ et 
$\chi \in \EC(\Lb^F,[s])$ est un caract\`ere r\'egulier. 
Soit $\g \in \Cent(\Gb^F,[s])$. Alors les assertions 
suivantes sont \'equivalentes~:
\begin{itemize}
\itemth{1} $\g \in \ZM\EC(\Gb^F,[s])$. 

\itemth{2} $\g \in \DS{\mathop{\oplus}_{(\Lb,\chi) \in \RC}} \ZM R_\Lb^\Gb \chi$.

\itemth{3} Pour tout $(\Lb,\chi) \in \RC$, 
$\langle \lexp{*}{R}_\Lb^\Gb \g , \chi \rangle_{\Lb^F} \in \ZM$.
\end{itemize}}

\bigskip

\soussection{Signes} 
Nous allons terminer cette section en expliquant comment d\'eterminer le signe 
intervenant dans le th\'eor\`eme \ref{irreductibles G}. Plus pr\'ecis\'ement, 
soit $\eta \in \Irr W(s)^{w_s F^*}$. D'apr\`es le th\'eor\`eme \ref{irreductibles G}, 
il existe un unique $\e_\eta \in \{1,-1\}$ tel que $\e_\eta R[s]_\eta$ 
soit un caract\`ere irr\'eductible de $\Gb^F$. Pour calculer $\e_\eta$, 
nous rappelons la construction des caract\`eres irr\'eductibles 
de $\Gbt^F$ par Lusztig et Srinivasan \cite{lusr} dans le cadre de la 
th\'eorie de Deligne-Lusztig. Si $\chi$ est un caract\`ere irr\'eductible 
$F^*$-stable de $W(\sti)=W^\circ(\sti)=W^\circ(s)$, nous noterons $\chih$ 
son {\it extension pr\'ef\'er\'ee} (au sens de \cite[chapitre IV, \SEC 17]{lucs}) 
au groupe $W(\sti) \rtimes <\phi_1>$. Rappelons aussi que l'isom\'etrie 
$\RC(\sti) : \Cent(W(\sti)\phi_1) \to \Cent(\Gbt^F,(\sti))$ 
a \'et\'e d\'efinie dans l'exemple \ref{a=1}.

\bigskip 

\Theoreme{Lusztig-Srinivasan}{lusztig srinivasan}
{\it Si $\chi$ est un carct\`ere irr\'eductible $F^*$-stable de 
$W(\sti)$, alors la fonction centrale 
$\e_\Gbt \e_{C_{\Gbt^*}(\sti)}\RC(\sti)_\chih$ est un \car \irr de 
$\Gbt^F$ appartenant \`a $\EC(\Gbt^F,(\sti))$. De plus, l'application 
$$\fonctio{(\Irr W(\sti))^{F^*}}{\EC(\Gbt^F,[\sti])}{\chi}{\e_\Gbt 
\e_{C_{\Gbt^*}(\sti)} \RC(\sti)_\chih}$$
est bijective.}

\bigskip

Notons $\IC(W^\circ(s),A_{\Gb^*}(s),F^*)$ l'ensemble des couples 
$(\chi,\xi)$ o\`u $\chi \in (\Irr W^\circ(s))^{F^*}$ et 
$\xi \in (A_{\Gb^*}(s,\chi)^{F^*})^\we$. Si 
$(\chi,\xi) \in \IC(W^\circ(s),A_{\Gb^*}(s),F^*)$, on note 
$\chi_\phi$ le caract\`ere irr\'eductible de $W^\circ(s)^{w_sF^*}$ 
correspondant \`a $\chi$ par la bijection \ref{bi alpha}. Si 
$A_{\Gb^*}(s,\chi)$ d\'esigne le stabilisateur de $\chi$ dans $A_{\Gb^*}(s)$, 
alors $A_{\Gb^*}(s,\chi)^{F^*}$ est le stabilisateur de $\chi_\phi$ 
dans $A_{\Gb^*}(s)^{F^*}$. Notons $\chit_\phi$ l'extension canonique 
de $\chi_\phi$ \`a $W^\circ(s)^{w_s F^*} \rtimes A_{\Gb^*}(s,\chi)^{F^*}$ 
(voir proposition \ref{extension canonique}). Posons alors
$$\eta_{\chi,\xi}=
\Ind_{W^\circ(s)^{w_s F^*} \rtimes A_{\Gb^*}(s,\chi)^{F^*}}^{W(s)^{w_sF^*}} 
(\chit_\phi \otimes \xi).$$
Le groupe $A_{\Gb^*}(s)^{F^*}$ agit sur $\IC(W^\circ(s),A_{\Gb^*}(s),F^*)$ 
par conjugaison sur la premi\`ere composante. 
D'apr\`es \ref{bijection car}, l'application
\equat\label{surjective IC}
\fonctio{\IC(W^\circ(s),A_{\Gb^*}(s),F^*)}{\Irr W(s)^{w_s F^*}}{(\chi,\xi)}{
\eta_{\chi,\xi}}
\endequat
est surjective et ses fibres sont les orbites sous $A_{\Gb^*}(s)^{F^*}$. 
Notons maintenant $\chit$ l'extension canonique de $\chi$ \`a 
$W^\circ(s) \rtimes <w_s \phi_1> = W^\circ(s) \rtimes <\phi_1>$. 

\bigskip

\proposition{somme}
{\it Soit $\chi$ un caract\`ere irr\'eductible $F^*$-stable de $W^\circ(s)$. 
Alors
$$\Res_{\Gb^F}^{\Gbt^F} \RC(\sti)_\chit = 
\sum_{\xi \in (A_{\Gb^*}(s,\chi)^{F^*})^\we} R[s]_{\eta_{\chi,\xi}}.$$}

\bigskip

\proof On a 
$$\sum_{\xi \in (A_{\Gb^*}(s,\chi)^{F^*})^\we} \eta_{\chi,\xi} 
=\Ind_{W^\circ(s)^{w_s F^*}}^{W(s)^{w_sF^*}} \chi_\phi.$$
Cette fonction centrale est nulle en dehors 
de $W^\circ(s)^{w_s F^*}$ et co\"\i ncide avec $\sum_{a \in A_{\Gb^*}(s)^{F^*}} 
\lexp{a}{\chi_\phi}$ sur $W^\circ(s)^{w_sF^*}$. Compte tenu de 
\ref{invariance RR}, on a 
$$\sum_{\xi \in (A_{\Gb^*}(s,\chi)^{F^*})^\we} R[s]_{\eta_{\chi,\xi}} 
= |A_{\Gb^*}(s)^{F^*}| \times R[s,1]_{\chi_\phi}.$$
Mais, par construction, $R[s,1]_{\chi_\phi}=\RC(s,1)_{\chit}$. 
Le r\'esultat d\'ecoule alors de la proposition \ref{rtgs} (a).\fin

\bigskip

D'apr\`es \cite[chapitre IV, \SEC 17]{lucs}, il existe un unique 
$\e_\chi \in \{1,-1\}$ tel que $\chih = \e_\chi \chit$ sur $W^\circ(s)\phi_1$. 
Ce signe est d\'etermin\'e \`a partir des deux exemples extr\^emes suivants~:

\bigskip

\exemple{deploye signe} 
Si $w_s=1$, alors $\e_\chi=1$ (voir corollaire \ref{canonique preferee}).\finl

\bigskip

\exemple{unitaire signe}
Si $C_{\Gb^*}^\circ(s)$ est quasi-simple et si $w_s \not= 1$, alors, d'apr\`es 
\cite[chapitre IV, \SEC 17.2 (a)]{lucs}, on a 
$\e_\chi = (-1)^{\ab_\chi}$, o\`u $\ab_\chi$ est le $\ab$-invariant 
attach\'e \`a $\chi$ (voir \cite[\SEC 16.2]{lucs}).\finl

\bigskip

Il r\'esulte alors de la proposition \ref{somme} et du th\'eor\`eme 
de Lusztig-Srinivasan que~:

\bigskip

\corollaire{signe eta}
{\it Si $(\chi,\xi) \in \IC(W^\circ(s),A_{\Gb^*}(s),F^*)$, alors 
$\e_{\eta_{\chi,\xi}}=\e_\Gb \e_{C_{\Gb^*}^\circ(s)} \e_\chi$.}

\bigskip

\section{Conjecture de Lusztig}

\medskip

Dor\'enavant, nous noterons $\EC'(\Gb^F,[s])$ l'ensemble  
$\{\e_\eta R[s]_\eta~|~\eta \in \Irr W(s)^{w_s F^*}\}$. Il est vraisemblable 
qu'en g\'en\'eral $\EC'(\Gb^F,[s])=\EC(\Gb^F,[s])$. Cependant, nous 
ne sommes pour l'instant capable de le prouver que lorsque 
$q$ est assez grand (voir le th\'eor\`eme \ref{irreductibles G}). 
Nous travaillerons n\'eanmoins avec cet ensemble. Nous posons 
$\EC'(\Gb^F,(s))=\coprod_{\a \in H^1(F^*,A_{\Gb^*}(s))} \EC'(\Gb^F,[s_\a])$.

\bigskip

\soussection{Familles} 
Reprenons les notations pr\'ec\'edant le lemme \ref{asai}.
Si $\a \in H^1(F^*,A_{\Gb^*}(s))$, soit $\d(\a)$ l'\'el\'ement de 
$A_{\Gb^*}(s)$ repr\'esent\'e par $g_\a^{-1}F(g_\a)$. En fait, l'application 
$\d : H^1(F^*,A_{\Gb^*}(s)) \to A_{\Gb^*}(s)$ est une section du 
morphisme canonique $A_{\Gb^*}(s) \to H^1(F^*,A_{\Gb^*}(s))$. Ce n'est 
en g\'en\'eral pas un morphisme de groupes. 

Le groupe $W(s_\a)$ est naturellement isomorphe \`a $W(s)$ (gr\^ace \`a la 
conjugaison par $g_\a$), mais nous devons le munir d'un \auto de Frobenius 
diff\'erent $\d(\a)F^*$. Par exemple, on a 
une application surjective (voir \ref{surjective IC}) 
$$\IC(W^\ci(s),A_{\Gb^*}(s),\d(\a)F^*) \longto \Irr W(s)^{w_s \a F^*}$$ 
qui induit une application surjective 
$$\IC(W^\ci(s),A_{\Gb^*}(s),\d(\a)F^*) \longto \EC'(\Gb^F,[s_\a]).$$
On en d\'eduit une troisi\`eme application surjective 
\equat\label{surjective e bis}
\coprod_{\a \in H^1(F^*,A_{\Gb^*}(s))} \IC(W^\ci(s),A_{\Gb^*}(s),\d(\a)F^*) 
\longto \EC'(\Gb^F,(s)).
\endequat
Les fibres de cette application sont des $A_{\Gb^*}(s)^{F^*}$-orbites (en effet, 
puisque $A_{\Gb^*}(s)$ est ab\'elien, on a $A_{\Gb^*}(s)^{aF^*}=A_{\Gb^*}(s)^{F^*}$ 
pour tout $a \in A_{\Gb^*}(s)$). 
Nous allons ici donner une autre description de cette surjection en termes 
de familles.

\medskip

Soit $\ch$ un \car \irr de $W^\ci(s)$. S'il existe $\a \in H^1(F^*,A_{\Gb^*}(s))$ 
tel que $\ch$ est $\d(\a)F^*$-stable, 
alors l'orbite de $\ch$ sous $A_{\Gb^*}(s)$ est $F^*$-stable. R\'eciproquement, 
si l'orbite de $\ch$ sous $A_{\Gb^*}(s)$ est $F^*$-stable, alors 
il existe $\a \in H^1(F^*,A_{\Gb^*}(s))$ et $a \in A_{\Gb^*}(s)$ tels que 
$\ch=\lexp{a^{-1}\d(\a)F^*(a)}{(\lexp{F^*}{\ch})}$. 
En particulier, $\lexp{a}{\ch}$ est $\d(\a)F^*$-stable. Nous allons 
donc utiliser les $A_{\Gb^*}(s)$-orbites $F^*$-stables 
de \cars \irrs pour regrouper les \'el\'ements de $\EC(\Gb^F,(s))$ en familles.
On note 
$$\m_\ch : H^1(F^*,A_{\Gb^*}(s,\ch)) \longto H^1(F^*,A_{\Gb^*}(s))$$
le \mor naturel de groupes.

Si $\ch$ est un \car \irr de $W^\ci(s)$, nous noterons $[\ch]$ son orbite 
sous $A_{\Gb^*}(s)$. On note $A_{\Gb^*}(s)\backslash \Irr W^\ci(s)$ l'ensemble 
de ces orbites. Fixons un \'el\'ement $F^*$-stable 
$[\ch] \in A_{\Gb^*}(s)\backslash \Irr W^\ci(s)$. Alors, 
d'apr\`es le calcul pr\'ec\'edent, on peut choisir 
$\ch$ de sorte que $\ch$ soit $\d(\a_\ch)F^*$-stable pour un 
$\a_\ch \in H^1(F^*,A_{\Gb^*}(s))$. 
Bien s\^ur, le couple $(\ch,\a_\ch)$ n'est pas uniquement d\'etermin\'e 
par $[\ch]$. De plus, l'\'el\'ement 
$\a_\ch$ n'est pas forc\'ement d\'etermin\'e par le choix de $\ch$. 

Soit maintenant $(\x,\a) \in \MC(A_{\Gb^*}(s,\ch),F^*)$. 
Alors il existe $a \in A_{\Gb^*}(s)$ tel que 
$$a^{-1}\d(\m_\ch(\a)\a_\ch) \d(\a_\ch)^{-1} F^*(a) \in A_{\Gb^*}(s,\ch)$$
et repr\'esente $\a$. On pose alors 
$$\ch_\a=\lexp{a}{\ch}.$$
Le \car $\ch_\a$ est d\'etermin\'e \`a $A_{\Gb^*}(s)^{F^*}$-conjugaison pr\`es 
et il est facile de v\'erifier que $\ch_\a$ est stable sous l'action de 
$\d(\m_\ch(\a)\a_\ch)F^*$. On pose alors 
\equat\label{definition rs}
R(s,\chi,\a_\chi)_{\a,\xi}= \e_\Gb \e_{C_{\Gb^*}^\circ(s)} \e_\chi 
R[s_{\m_\chi(\a)\a_\chi}]_{\eta_{\chi_\a,\xi}} \in \EC'(\Gb^F,(s)).
\endequat
Alors la fonction centrale $R(s,\chi,\a_\chi)_{\a,\xi}$ 
d\'epend uniquement du choix de $\ch$, $\a_\ch$ et $\d$~: 
rappelons que c'est un caract\`ere irr\'eductible, du moins 
lorsque $q$ est assez grand (voir le th\'eor\`eme \ref{irreductibles G} et 
le corollaire \ref{signe eta}). 

\bigskip 

\proposition{geo}
{\it Soient $(\x,\a)$ et $(\x',\a')$ deux \'el\'ements de 
$\MC(A_{\Gb^*}(s,\ch),F^*)$. Alors 
$R(s,\chi,\a_\chi)_{\a,\xi}=R(s,\chi,\a_\chi)_{\a',\xi'}$ 
\ssi $(\x,\a)=(\x',\a')$.}

\bigskip

\proof Soient $(\x,\a)$ et $(\x',\a')$ deux \'el\'ements de 
$\MC(A_{\Gb^*}(s,\ch),F^*)$ 
tels que $R(s,\chi,\a_\chi)_{\a,\xi}=R(s,\chi,\a_\chi)_{\a',\xi'}$. 
Alors, par construction, 
on sait que $\m_\ch(\a)=\m_\ch(\a')$, les \cars $\ch_\a$ et $\ch_{\a'}$ 
sont conjugu\'es sous $A_{\Gb^*}(s)^{F^*}$ et $\x=\x'$ (d'apr\`es la proposition 
\ref{surjective IC}). 
Soient $a$ et $a'$ deux \'el\'ements de $A_{\Gb^*}(s)$ 
tels que $a^{-1}\d(\m_\ch(\a)\a_\ch) \d(\a_\ch)^{-1} F^*(a)$ 
et $a^{\pr -1}\d(\m_\ch(\a')\a_\ch) \d(\a_\ch)^{-1} F^*(a')$ appartiennent \`a 
$A_{\Gb^*}(s,\ch)$ et repr\'esentent $\a$ et $\a'$ respectivement. Alors, puisque 
$\ch_\a=\lexp{a}{\ch}$ et $\ch_{\a'}=\lexp{a'}{\ch}$ sont 
conjugu\'es sous $A_{\Gb^*}(s)^{F^*}$, il existe $b \in A_{\Gb^*}(s,\ch)$ et 
$c \in A_{\Gb^*}(s)^{F^*}$ tels que $a'=abc$. Cela montre que $\a=\a'$.\fin

\bigskip

Si on note $\EC(\Gb^F,(s),[\chi])$ l'ensemble 
$\{R(s,\chi,\a_\chi)_{\a,\xi}~|~(\x,\a) \in \MC(A_{\Gb^*}(s,\ch),F^*)\}$, 
alors $\EC(\Gb^F,(s),[\chi])$ ne d\'epend que de $[\ch]$ 
et non du choix de $\ch$ et $\a_\ch$. On a alors 
\equat\label{partition familles}
\EC'(\Gb^F,(s)) =
\coprod_{[\ch] \in \bigl(A_{\Gb^*}(s)\backslash \Irr W^\ci(s)\bigl)^{F^*}}
\EC(\Gb^F,(s),[\chi])
\endequat
et la proposition \ref{geo} montre qu'il y a une bijection 
\equat
\diagram
\MC(A_{\Gb^*}(s,\ch),F^*) \rrto^\sim && \EC(\Gb^F,(s),[\chi])
\enddiagram
\endequat
pour tout $[\ch] \in \bigl(A_{\Gb^*}(s)\backslash \Irr W^\ci(s)\bigl)^{F^*}$.  
Cette bijection d\'epend du choix de $\ch$, $\a_\ch$, $\d$. 
Nous noterons $\Cent(\Gb^F,(s),[\chi])$ le sous-$\qlb$-espace 
vectoriel de $\Cent(\Gb^F,(s))$ engendr\'e par $\EC(\Gb^F,(s),[\chi])$. On a, 
d'apr\`es \ref{partition familles},
\equat\label{decomposition familles}
\Cent(\Gb^F,(s)) =\mathop{\oplus}^\perp_{[\ch] \in \bigl(A_{\Gb^*}(s)
\backslash \Irr W^\ci(s)\bigl)^{F^*}}  \Cent(\Gb^F,(s),[\chi])
\endequat

\bigskip

\soussection{Transformation de Fourier} 
Fixons $\chi \in \Irr W^\ci(s)\bigl)^{F^*}$ dont la 
$A_{\Gb^*}(s)$-orbite est $F^*$-stable et soit $\a_\chi$ un \'el\'ement de 
$H^1(F^*,A_{\Gb^*}(s))$ tel que $\chi$ soit $\d(\a_\chi)F^*$-stable. 
Si $(\xi,\a) \in \MC(A_{\Gb^*}(s,\chi),F^*)$, 
notons $R(s,\chi,\a_\chi)_{\a,\xi}$ la fonction centrale d\'efinie 
pr\'ec\'edemment gr\^ace \`a ce choix-ci via \ref{definition rs}. 

Si $(a,\t) \in \MC^\vee(A_{\Gb^*}(s,\chi),F^*)$, posons 
$$\Rha(s,\chi,\a_\chi)_{a,\t} = \frac{1}{|A_{\Gb^*}(s,\chi)^{F^*}|} 
\sum_{(\xi,\a) \in \MC(A_{\Gb^*}(s,\chi),F^*)} 
\xi(a)^{-1} \t(\a) R(s,\chi,\a_\chi)_{\a,\xi}.$$
Alors $(\Rha(s,\chi,\a_\chi)_{a,\t})_{(a,\t) \in \MC^\vee(A_{\Gb^*}(s,\chi),F^*)}$ 
est une base orthonormale de $\Cent(\Gb^F,(s),[\chi])$. 
La proposition suivante d\'ecrit l'action de $H^1(F,\ZC(\Gb))$ et de 
$\Zb(\Gb)^F$ sur $\Cent(\Gb^F,(s),[\chi])$ dans cette base. 

\bigskip

\proposition{proprietes familles}
{\it Soit $(a,\t) \in \MC^\vee(A_{\Gb^*}(s,\chi),F^*)$. Alors~: 
\begin{itemize}
\itemth{a} $\Rha(s,\chi,\a_\chi)_{a,\t} \in \Cent(\Gb^F,(s),a)$.

\itemth{b} Soit $z \in \Zb(\Gb)^F$. Notons $\t_z$ la restriction \`a 
$H^1(F^*,A_{\Gb^*}(s,\chi))$ du caract\`ere lin\'eaire $\omeh_s^1(\zba)$ de 
$H^1(F^*,A_{\Gb^*}(s))$, o\`u $\zba$ d\'esigne la classe de $z$ dans 
$\ZC(\Gb)^F$. Alors
$$t_z^\Gb \Rha(s,\chi,\a_\chi)_{a,\t} = \sha_{\a_\chi}(z) 
\Rha(s,\chi,\a_\chi)_{a,\t\t_z}.$$
\end{itemize}}

\bigskip

\proof La premi\`ere assertion d\'ecoule imm\'ediatement de la d\'efinition. 
La deuxi\`eme d\'ecoule du lemme \ref{asai} et de la remarque 
\ref{elementaire series} (d).\fin

\bigskip

\soussection{Conjecture de Lusztig} 
Fixons un caract\`ere $\chi$ de $W^\circ(s)$ et un \'el\'ement 
$\a_\chi$ de $H^1(F^*,A_{\Gb^*}(s))$ tels que $\lexp{\d(\a_\chi)F^*}{\chi}=\chi$. 
Soit $(a,\t) \in \MC^\vee(A_{\Gb^*}(s,\chi),F^*)$. 
Notons $\chi_a$ le caract\`ere irr\'eductible de $W^\circ(s)^a$ 
associ\'e \`a $\chi$ par la bijection \ref{bi alpha}. Son stabilisateur dans 
$A_{\Gb^*}(s) \rtimes <\phi_1>$ est 
$A_{\Gb^*}(s,\chi) \rtimes <\d(\a_\chi) \phi_1>$. On note $\chit_a$ 
l'extension canonique de $\chi_a$ \`a 
$W^\circ(s)^a \rtimes (A_{\Gb^*}(s,\chi) \rtimes <\d(\a_\chi) \phi_1>)$. 
On note $\taut_{\a_\chi}$ l'extension  
de $\t$ \`a $A_{\Gb^*}(s,\chi) \rtimes <\d(\a_\chi) \phi_1>$ telle que 
$\taut_{\a_\chi}(\d(\a_\chi)\phi_1)=1$. Pour finir, on pose
$$\eta_{\chi,a,\t}=\Ind_{W^\circ(s)^a \rtimes A_{\Gb^*}(s,\chi)}^{W(s)^a} 
(\chit_a \otimes \t)$$
$$\etat_{\chi,\a_\chi,a,\t}=
\Ind_{W^\circ(s)^a \rtimes (A_{\Gb^*}(s,\chi) \rtimes <\d(\a_\chi)\phi_1>}^{W(s)^a 
\rtimes <\phi_1>} 
(\chit_a \otimes \taut_{\a_\chi}).\leqno{\text{et}}$$
Alors $\eta_{\chi,a,\t}$ est un caract\`ere $F^*$-stable de $W(s)^a$ et 
$\etat_{\chi,\a_\chi,a,\t}$ est une extension de $\eta_{\chi,a,\t}$ 
\`a $W(s)^a \rtimes <\phi_1>$. 

\bigskip

\theoreme{conjecture de lusztig}
{\it Avec les notation pr\'ec\'edentes, on a
$$\Rha(s,\chi,\a_\chi)_{a,\t}=\e_\Gb \e_{C_{\Gb^*}^\circ(s)} \e_\chi 
\RC(s,a)_{\etat_{\chi,\a_\chi,a,\t}}.$$}

\bigskip

\proof Quitte \`a changer d'\'el\'ement semi-simple $F^*$-stable dans $(s)$, 
on peut supposer que $\a_\chi=1$ (rappelons qu'alors $\d(\a_\chi)=1$). 
D'autre part, le r\'esultat ne d\'epend pas du choix de la section $\d$, 
donc nous pourrons supposer que $\d(\m_1(\a)) \in A_{\Gb^*}(s,\chi)$ 
pour tout $\a \in A_{\Gb^*}(s,\chi)$. Ici, 
$\m_1 : A_{\Gb^*}(s) \to H^1(F^*,A_{\Gb^*}(s))$ est l'application canonique. 

Soit $\a \in A_{\Gb^*}(s,\chi)$. Notons $\bar{\a}$ la classe de $\a$ dans 
$H^1(F^*,A_{\Gb^*}(s,\chi))$. Alors il existe $b \in A_{\Gb^*}(s)$ 
tel que $b^{-1}\d(\m_1(\a)) F^*(b) = \a$ et on pose $\chi_\a=\lexp{b}{\chi}$. 
Alors $R[s_{\mu_1(\a)}]_{\eta_{\chi_\a,\xi}}$ ne d\'epend que 
de $\bar{\a}$ et est \'egal \`a 
$R[s_{\mu_\chi(\bar{\a})}]_{\eta_{\chi_{\bar{\a}},\xi}}$. 
Par suite, si on note $\chi_{\a,a}$ le caract\`ere irr\'eductible de 
$W^\circ(s)^a$ associ\'e \`a $\chi_\a$ par la bijection \ref{bi alpha}, on a
$$\Rha(s,\chi,\a_\chi)_{a,\t}=\frac{1}{|A_{\Gb^*}(s,\chi)|} 
\sum_{\a \in A_{\Gb^*}(s,\chi)} \t(\a)^{-1} \RC[s_{\m_1(\a)},a]_{\chit_{\a,a}}.$$
D'autre part, si $w \in W^\circ(s)^a$, on a 
$\chit_{\a,a}(w\d(\m_1(\a))\phi_1)=\chit_a(b^{-1} w b \a)$ 
et $R_{\Lb_{s_{\m_1(\a)},w,a}}^\Gb 
\rhodot_{s_{w\d(\m_1(\a))},a}^{\Lb_{s_{\m_1(\a)},w,a}} = 
R_{\Lb_{s,w\a,a}}^\Gb \rhodot_{s_{w\a},a}^{\Lb_{s,w\a,a}}$. 
Par suite, 
$$\Rha(s,\chi,\a_\chi)_{a,\t}=\frac{1}{|W^\circ(s)^a \rtimes A_{\Gb^*}(s,\chi)|} 
\sum_{w \in W^\circ(s)^a \rtimes A_{\Gb^*}(s,\chi)} 
(\chit_a \otimes \taut_1)(w\phi_1) 
R_{\Lb_{s,w,a}}^\Gb \rhodot_{s_w,a}^{\Lb_{s,w,a}},$$
ce qui montre le r\'esultat (voir aussi \cite[lemme 3.1.1]{bonnafe couro}).\fin

\bigskip

Le th\'eor\`eme \ref{conjecture de lusztig} et le th\'eor\`eme 
\ref{PRINCIPAL} montre la conjecture de Lusztig pour tous les groupes 
de type $A$ lorsque $q$ est assez grand~:

\bigskip

\Corollaire{Conjecture de Lusztig en type A}{essentiel}
{\it Avec les notation pr\'ec\'edentes, on a
$$\Rha(s,\chi,\a_\chi)_{a,\t}=\e_\Gb \e_{C_{\Gb^*}^\circ(s)} \e_\chi 
\GC(\Lb_{s,a},\o_{\Lb_{s,a},s}(a))^{-1} 
\XC_{K(s,a)_{\eta_{\chi,a,\t}},\etat_{\chi,\a_\chi,a,\t}}.$$}

\bigskip

\section{Le groupe sp\'ecial lin\'eaire\label{chapitre sln}}

\medskip

Nous pr\'ecisons ici, dans le cas du groupe sp\'ecial lin\'eaire, 
quelques-uns des r\'esultats obtenus dans le chapitre pr\'ec\'edent. 
Nous \'etablissons aussi une {\it d\'ecomposition de Jordan} des 
caract\`eres et sa compatibilit\'e avec l'induction de Lusztig 
(voir le diagramme \ref{jordan sln}).

\bigskip

\begin{quotation}
\noindent{\bf Hypoth\`ese : }{ \it Dor\'enavant, et ce jusqu'\`a la fin de 
cette section, nous fixons un groupe \`a centre connexe $\Gbt_\bullet$ 
de type $A_{n-1}$ muni d'un endomorphisme de Frobenius d\'eploy\'e $F$ sur $\fq$ 
et nous supposerons que $\Gbt$ est un \levi $F$-stable de $\Gbt_\bullet$.}
\end{quotation}

\bigskip

Notons avant de commencer que l'hypoth\`ese entra\^\i ne que $w_s =1$ 
et donc que $\e_\chi=1$ pour tout caract\`ere irr\'eductible $F^*$-stable 
$\chi$ de $W^\circ(s)$ (voir exemple \ref{deploye signe}). 
En d'autres termes, $F^*$ agit sur $W^\circ(s)$ (qui est un produit 
direct de groupes sym\'etriques) seulement par permutation des composantes. 

\bigskip

\soussection{Th\'eorie de Harish-Chandra}
Commen\c{c}ons par d\'ecrire les s\'eries de Harish-Chandra 
contenues dans $\EC(\Gb^F,[s])$. 

\bigskip

\proposition{serie gln}
{\it On a $\EC(\Gb^F,[s])=\EC(\Gb^F,\Lb_s,\rho_s^{\Lb_s})$.}

\bigskip

\proof D'apr\`es la proposition \ref{isometrie RR}, on a 
$|\EC(\Gb^F,[s])|=|\Irr W(s)^{F^*}|$. D'apr\`es \ref{bi G}, on 
a $|\EC(\Gb^F,\Lb_s,\rho_s^{\Lb_s})|=|\Irr W(s)^{F^*}|$. Puisque 
$\EC(\Gb^F,\Lb_s,\rho_s^{\Lb_s})$ est contenu dans $\EC(\Gb^F,[s])$, 
on a en fait l'\'egalit\'e de ces deux ensembles.\fin

\bigskip

D'apr\`es la proposition \ref{serie gln} et d'apr\`es \ref{bi G}, on 
obtient une bijection 
$$\fonctio{\Irr W(s)^{F^*}}{\EC(\Gb^F,[s])}{\eta}{R_\eta[s].}$$
En fait, cette bijection co\"\i ncide avec celle obtenue via l'isom\'etrie 
$R[s]$, du moins lorsque $q$ est assez grand~:

\bigskip

\theoreme{egalite parametrage}
{\it Si la conjecture $(\GG)$ a lieu dans $\Gb$, alors 
$R_\eta[s]=\e_\Gb \e_{C_{\Gb^*}^\circ(s)} R[s]_\eta$ 
pour tout $\eta \in \Irr W(s)^{F^*}$.}

\bigskip

\proof Soit $W_1$ un sous-groupe parabolique standard $F^*$-stable 
de $W^\circ(s)$ et soit $A_1$ son normalisateur dans $A_{\Gb^*}(s)^{F^*}$. 
Alors, par le m\^eme raisonnement que dans la preuve du lemme 
\ref{existence L}, on a $W_1 \rtimes A_1 = W_\Lb(s)$ pour 
$$\Lb=C_\Gb\bigl((\Tb_1^{W_1 \rtimes (A_1 \rtimes <\phi_1>)})^\circ\bigr).$$
Le groupe $\Lb$ est en fait $\Gb$-d\'eploy\'e. Donc, compte tenu 
du corollaire \ref{frobenius}, de la proposition \ref{induit commute} 
et du th\'eor\`eme \ref{parametrage} (c), il suffit de montrer le 
r\'esultat lorsque $\eta$ se factorise en un caract\`ere 
lin\'eaire de $A_{\Gb^*}(s)^{F^*}$. Cela d\'ecoule alors 
de \ref{convention}, du th\'eor\`eme \ref{parametrage} (b) et  
de la proposition \ref{rho xi}.\fin

\bigskip

La preuve du th\'eor\`eme \ref{egalite parametrage} montre 
\'egalement la version suivante plus pr\'ecise du corollaire 
\ref{caracterisation sem}. Il n'y a pas besoin d'hypoth\`ese 
sur $q$ car on peut passer par la th\'eorie de Harish-Chandra 
et le th\'eor\`eme \ref{parametrage} (c).

\bigskip

\corollaire{caracterisation sem sln}
{\it Notons $\SC_d$ l'ensemble des couples $(\Lb,\rho)$ 
ou $\Lb$ est un \levi $\Gb$-d\'eploy\'e de $\Gb$ dont un dual 
$\Lb^*$ dans $\Gb^*$ contient $s$ et 
$\rho \in \EC(\Lb^F,[s])$ est un caract\`ere semi-simple. 
Soit $\g \in \Cent(\Gb^F,[s])$. Alors les assertions 
suivantes sont \'equivalentes~:
\begin{itemize}
\itemth{1} $\g \in \ZM\EC(\Gb^F,[s])$. 

\itemth{2} $\g \in \DS{\mathop{\oplus}_{(\Lb,\rho) \in \SC_d}} \ZM R_\Lb^\Gb \rho$.

\itemth{3} Pour tout $(\Lb,\rho) \in \SC_d$, 
$\langle \lexp{*}{R}_\Lb^\Gb \g , \rho \rangle_{\Lb^F} \in \ZM$.
\end{itemize}}

\bigskip

\corollaire{caracterisation reg sln}
{\it Notons $\RC_d$ l'ensemble des couples $(\Lb,\chi)$ 
ou $\Lb$ est un \levi $\Gb$-d\'eploy\'e de $\Gb$ dont un dual 
$\Lb^{\!*}$ dans $\Gb^*$ contient $s$ et 
$\chi \in \EC(\Lb^F,[s])$ est un caract\`ere r\'egulier. 
Soit $\g \in \Cent(\Gb^F,[s])$. Alors les assertions 
suivantes sont \'equivalentes~:
\begin{itemize}
\itemth{1} $\g \in \ZM\EC(\Gb^F,[s])$. 

\itemth{2} $\g \in \DS{\mathop{\oplus}_{(\Lb,\chi) \in \RC_d}} \ZM R_\Lb^\Gb \chi$.

\itemth{3} Pour tout $(\Lb,\chi) \in \RC_d$, 
$\langle \lexp{*}{R}_\Lb^\Gb \g , \chi \rangle_{\Lb^F} \in \ZM$.
\end{itemize}}

\bigskip

\soussection{D\'ecomposition de Jordan}
D'apr\`es \cite[proposition 6.4.3]{bonnafe couro}, l'application $i^*$ induit 
une bijection entre $\EC(i^{*-1}(C_{\Gb^*}(s))^{F^*},1)$ et 
$\EC(C_{\Gb^*}(s)^{F^*},1)$. Par suite, 
d'apr\`es \cite[7.4.3]{bonnafe couro}, on a une bijection 
$$\fonctio{\Irr W(s)^{F^*}}{\EC(C_{\Gb^*}(s)^{F^*},1)}{\eta}{R_{s,\eta}.}$$
On obtient donc une bijection 
$$\fonction{\aleph_{\Gb,s}}{\EC(\Gb^F,[s])}{
\EC(C_{\Gb^*}(s)^{F^*},1)}{R_\eta[s]}{R_{s,\eta}}$$
appel\'ee {\it d\'ecomposition de Jordan} des caract\`eres de $\Gb^F$. 

Soit $\Lb$ un sous-groupe de Levi $F$-stable de $\Gb$ et soit $\Lb^*$ 
un sous-groupe de Levi $F^*$-stable de $\Gb^*$ dual de $\Lb$ et 
contenant $s$. D'apr\`es le th\'eor\`eme \ref{egalite parametrage}, d'apr\`es 
la proposition \ref{induction RRR} et d'apr\`es 
\cite[th\'eor\`eme 7.6.1]{bonnafe couro}, le diagramme 
\equat\label{jordan sln}
\diagram
\ZM\EC(\Lb^F,[s]) \rrto^{\DS{\aleph_{\Lb,s}}} \ddto_{\DS{\e_\Gb\e_\Lb R_\Lb^\Gb}} 
&& \ZM\EC(C_{\Lb^*}(s)^{F^*},1) 
\ddto^{\DS{\e_{C_{\Gb^*}^\circ(s)}\e_{C_{\Lb^*}^\circ(s)} 
R_{C_{\Lb^*}(s)}^{C_{\Gb^*}(s)}}} \\
&& \\
\ZM(\EC(\Gb^F,[s]) \rrto^{\DS{\aleph_{\Gb,s}}} 
&& \ZM\EC(C_{\Gb^*}(s)^{F^*},1) 
\enddiagram
\endequat
est commutatif lorsque la conjecture $(\GG)$ a lieu dans $\Gb$.

\bigskip

\section{Questions en suspens}

\medskip

\noindent{\bf 1.} Il serait int\'eressant d'\'etudier les questions 
abord\'ees dans ce chapitre (s\'eries de Harish-Chandra, d\'ecomposition 
de Jordan) dans le cas des groupes de type $A$ non n\'ecessairement 
d\'eploy\'es. Concernant la question de la d\'ecomposition de Jordan, 
il faudrait \'etablir l'analogue de la bijection \cite[7.4.3]{bonnafe couro}. 
La commutativit\'e de l'analogue du diagramme \ref{jordan sln} 
est alors purement formelle.

\medskip

\noindent{\bf 2.} \`A travers le th\'eor\`eme \ref{irreductibles G}, on obtient 
un param\'etrage des caract\`eres irr\'eductibles de $\Gb^F$ par 
les paires $(\chi,\xi)$ en utilisant seulement un caract\`ere de Gelfand-Graev 
de $\Gb^F$. Un autre param\'etrage par les paires $(\chi,\xi)$ a 
\'et\'e obtenu par Shoji \cite{shoji} (ou encore \cite{bonnafe torsion}) 
en utilisant les caract\`eres de Gelfand-Graev g\'en\'eralis\'es. 
Il serait int\'eressant de relier ces deux param\'etrages. Cela 
permettrait de relier les transform\'ees de Fourier introduites 
ici et les {\it caract\`eres fant\^omes} d\'efinis par Shoji 
dans \cite{shoji banff}.

\medskip

\noindent{\bf 3.} L'\'ecriture d'un algorithme effectif, \`a partir 
des r\'esultats de cet article, pour calculer la table de caract\`eres des 
groupes r\'eductifs connexes de type $A$ est maintenant 
th\'eoriquement possible. Cela reste tout de m\^eme un travail 
consid\'erable.

\newpage

{\Large \part*{Appendice A. Produits en couronne}}

\bigskip

Nous allons rappeler dans cet appendice quelques faits sur les 
caract\`eres de produits en couronne de groupe finis. Nous 
reprendrons essentiellement ce qui est fait dans \cite[\SEC 2]{bonnafe couro}. 

Dans cet appendice, $r$, $d_1$,\dots, $d_r$ d\'esigneront des entiers 
naturels non nuls. On se fixe des groupes finis $G_1$,\dots, $G_r$ et on pose 
$$G^\ci=\prod_{i=1}^r 
(\underbrace{G_1 \times \dots \times G_r}_{d_i~{\mathrm{fois}}}).$$
On fixe aussi un morphisme de groupes 
$A \to \SG_{d_1} \times \dots \times \SG_{d_r}$. 
Alors $A$ agit sur $G^\ci$ via ce \mor par permutations des composantes 
et on pose
$$G = G^\ci \rtimes A.$$

\bigskip

\section{Extension canonique}

\medskip

\soussection{D\'efinition} Si $\ch$ est un \car \irr de $G^\ci$, 
on note $A(\ch)$ son stabilisateur dans $A$ et $G(\ch)$ son stabilisateur dans 
$G$. On a $G(\ch)=G^\ci \rtimes A(\ch)$. 
La proposition suivante est classique et sa preuve peut par exemple 
\^etre trouv\'ee dans \cite[Proposition 2.3.1]{bonnafe couro}~:

\bigskip

\proposition{extension canonique}
{\it Soit $\ch$ un \car \irr de $G^\ci$. Alors il existe une unique extension 
$\chit$ de $\ch$ \`a $G$ telle que $\chit(\a)$ soit un entier naturel 
non nul pour tout $\a \in A(\ch)$.} 

\bigskip

Soit $\IC(G^\circ,A)$ l'ensemble des couples $(\chi,\xi)$ o\`u 
$\chi \in \Irr G^\circ$ et $\xi \in \Irr A(\chi)$. Le groupe $A$ agit 
par conjugaison sur $\IC(G^\circ,A)$ et, par la th\'eorie de Clifford, 
l'application 
\equat\label{bijection car}
\fonctio{\IC(G^\circ,A)}{\Irr G}{(\chi,\xi)}{\Ind_{G(\chi)}^G (\chit \otimes \xi)}
\endequat
induit une bijection entre $\IC(G^\circ,A)/A$ et $\Irr G$. 

\medskip

Le \car \irr $\chit$ de la proposition \ref{extension canonique} 
sera appel\'e l'{\it extension canonique} de $\ch$ \`a $G(\ch)$. 
Soit $a \in A$. 
Dans \cite[2.2]{bonnafe couro}, l'auteur a construit une application 
\equat\label{application pi}
\pi_a : G^\ci a\longto G^{\ci a}
\endequat
induisant une bijection bien d\'efinie entre l'ensemble des classes de conjugaison 
de $G^\circ \rtimes <a>$ contenues dans $G^\circ a$ 
et l'ensemble des classes de conjugaison de $(G^\circ)^a$ 
ainsi qu'une bijection 
\equat\label{bi alpha}
\fonctio{(\Irr G^\ci)^a}{\Irr((G^\circ)^a)}{\ch}{\ch_a.}
\endequat
Si $\chi \in (\Irr G^\circ)^a$, nous noterons $\chit_a$ la restriction de 
$\chit$ \`a $G^\circ a$. Alors $(\chit_a)_{\chi \in (\Irr G^\circ)^a}$ 
est une base orthonormale de $\Cent(G^\circ a)$. 
Ces deux applications satisfont la propri\'et\'e suivante~: si $\ch$ 
est un \car \irr de $G^\ci$ et si $a \in A(\chi)$, alors 
\equat\label{formule chi a}
\chit(w)=\ch_a(\pi_a(w))
\endequat
pour tout $w \in G^\ci a$. Rappelons aussi que $\pi_a(a)=1$ donc 
$\chit(a)=\ch_a(1)$. D'autre part, l'application $\pi_a$ induit 
une application lin\'eaire 
\equat\label{iso pia}
\fonction{\pi_a^*}{\Cent((G^\circ)^a)}{\Cent(G^\circ a)}{f}{f \circ \pi_a}
\endequat
et l'\'egalit\'e \ref{formule chi a} montre que c'est une isom\'etrie. 

\bigskip

\exemple{pi a} 
Nous rappelons ici la d\'efinition de ces deux applications dans un cas particulier 
dont le cas g\'en\'eral peut ais\'ement se d\'eduire par produit direct. 
Supposons que $r=1$ et posons $d=d_1$. Supposons aussi que 
$\lexp{a}{(w_1,\dots,w_d)}=(w_d,w_1,\dots,w_{d-1})$. Alors 
$G_1 \simeq (G^\circ)^a$ et, via cet isomorphisme, 
$$\pi_a(w_1,\dots,w_d)=w_1 \dots w_d$$
et
$$\fonctio{\Irr G_1}{(\Irr G^\circ)^a}{\ch}{
\underbrace{\ch \otimes \dots \otimes \ch}_{d~{\text{fois}}}}$$
est la bijection r\'eciproque de la bijection \ref{bi alpha}.\finl

\bigskip

\soussection{Induction} 
Fixons maintenant $a \in A$ et une famille 
$(H_{ij})_{1 \le i \le r, 1 \le j \le d_i}$, o\`u $H_{ij}$ est un sous-groupe de 
$G_i$. On pose 
$$H^\circ=\prod_{i=1}^r (H_{i1} \times \dots \times H_{id_i})$$
et on suppose que $H^\circ$ est $a$-stable. Alors la restriction de $\pi_a$ 
\`a $H^\circ a$ est l'analogue de $\pi_a$ pour le groupe $H^\circ$ 
et on a encore une isom\'etrie toujours not\'ee 
$\pi^*_a : \Cent((H^\circ)^a) \to \Cent(H^\circ a)$. 
D'autre part, il r\'esulte de \cite[lemme 3.2.1]{bonnafe couro} que le diagramme 
\equat\label{diagramme induction tordue}
\diagram
\Cent((H^\circ)^a) \rrto^{\DS{\pi_a^*}} \ddto_{\DS{\Ind_{(H^\circ)^a}^{(G^\circ)^a}}} && 
\Cent(H^\circ a) \ddto^{\DS{\Ind_{H^\circ a}^{G^\circ a}}} \\
&& \\
\Cent((G^\circ)^a) \rrto^{\DS{\pi_a^*}} && 
\Cent(G^\circ a)
\enddiagram
\endequat
est commutatif. 

\bigskip

\proposition{debuggage}
{\it On a~:
\begin{itemize}
\itemth{a} Si $\chi \in (\Irr G^\circ)^a$, alors
$$\Ind_{G^\circ a}^G \chit_a = \sum_{\xi \in \Irr G(\chi)/G^\circ} 
\overline{\xi(a)} \Ind_{G(\chi)}^G (\chit \otimes \xi).$$

\itemth{b} L'application $\Ind_{G^\circ a}^G$ a pour image 
l'espace des fonctions centrales sur $G$ qui s'annulent en dehors 
de $G^\circ [a]$, o\`u $[a]$ est la classe de conjugaison de $a$ dans $A$.
\end{itemize}}

\bigskip

\proof (a) Par la transitivit\'e de l'induction, on peut supposer, 
et nous le ferons, que $G(\chi)=G$. Notons $A'=<a>$ et $G'=G^\circ \rtimes A'$. 
Nous noterons $\chit'$ la restriction de $\chit$ \`a $G'$. 
D'apr\`es \ref{une formule tordue}, on a 
$$\Ind_{G^\circ a}^G \chit_a= \Ind_{G'}^G\Bigl(\sum_{\xi \in A^{\prime \wedge}} 
\xi(a)^{-1} (\chit \otimes \xi)\Bigr).$$
Par cons\'equent, 
$$\Ind_{G^\circ a}^G \chit_a=\chit \otimes 
\Bigl(\Ind_{A'}^A (\sum_{\xi \in A^{\prime \wedge}} \xi(a)^{-1}\xi)\Bigr).$$
Il ne reste donc plus qu'\`a montrer que 
$$\Ind_{A'}^A (\sum_{\xi \in A^{\prime \wedge}} \overline{\xi(a)}\xi)
=\sum_{\xi \in \Irr A} \overline{\xi(a)}\xi,$$
ce qui est \'evident. 

\medskip

(b) D'apr\`es les formules donn\'ees dans la preuve 
de \ref{une formule tordue}, l'application $\Ind_{G^\circ a}^{G'}$ a pour 
image l'espace des fonctions centrales sur $G'$ nulles en dehors 
de $G^\circ a$. (b) en d\'ecoule car toute classe de conjugaison 
de $G$ contenue dans $G^\circ [a]$ rencontre $G^\circ a$.\fin

\bigskip

\section{Produits en couronne de groupes sym\'etriques}

\medskip

Nous ferons dans cette section l'hypoth\`ese suivante~:

\bigskip

\begin{quotation}
\noindent{\bf Hypoth\`ese :} 
{\it 
Jusqu'\`a la fin de cette section, nous fixons une famille  
d'entiers naturels non nuls $(n_i)_{1 \le i \le r}$ et nous 
supposons que $G_i=\SG_{n_i}$, le groupe sym\'etrique de degr\'e $n_i$.}
\end{quotation}

\bigskip

Fixons un entier $i \in \{1,2,\dots,r\}$. Notons $\e_i : G_i \to \{1,-1\}$ 
la signature. 
Si $\l=(\l_1,\dots,\l_x)$ est une partition de $n_i$, nous noterons 
$G_{i,\l}$ le sous-groupe (parabolique) de $G_i$ isomorphe 
\`a $\SG_{\l_1} \times \dots \times \SG_{\l_x}$ (sous-groupe de Young). 
Nous noterons $\l^*$ la partition duale de $\l$. 
Nous noterons $b_\l$ le nombre de transpositions contenues dans $G_{i,\l^*}$~:
on a
$$b_\l=\sum_{j=1}^x (j-1)\l_j.$$
Par exemple, $b_{(n_i)}=0$. 
Notons $\chi_\l$ l'unique \car \irr commun \`a 
$\Ind_{G_{i,\l}}^{G_i} 1_{G_{i,\l}}$ et 
$\e_i \otimes \Ind_{G_{i,\l^*}}^{G_i} 1_{G_{i,\l^*}}$ 
(voir \cite[th\'eor\`eme 5.4.7]{geck livre}). 
Toujours d'apr\`es \cite[th\'eor\`eme 5.4.7]{geck livre}, on a 
\equat\label{induit sn}
\Ind_{G_{i,\l}}^{G_i} 1_{G_{i,\l}} = \chi_\l + 
\sum_{\SS{\mu \vdash n_i} \atop \SS{b_\mu < b_\l}} \b_{\l\m} \chi_\mu.
\endequat

\medskip

Revenons \`a notre groupe $G^\circ$. 
Nous noterons $\PC$ l'ensemble des familles 
$\lamb=(\l_{ij})_{1 \le i \le r, 1 \le j \le d_i}$ 
o\`u $\l_{ij}$ est une partition de $n_i$. 
Le groupe $A$ agit naturellement sur $\PC$ par permutations. 
Nous posons alors $b_\lamb=\sum_{i,j} b_{\l_{ij}}$,  
$$G_\lamb^\circ=\prod_{i=1}^r (G_{i,\l_{i1}} \times \dots \times G_{i,\l_{id_i}})$$
et nous notons $A_\lamb$ le stabilisateur de $\lamb$ dans $A$, \cad le normalisateur 
de $G_\lamb^\circ$ dans $A$. On pose 
$$G_\lamb=G_\lamb^\circ \rtimes A_\lamb.$$
Nous noterons $\chi_\lamb$ le \car \irr de $G^\circ$ d\'efini par
$$\chi_\lamb=\bigotimes_{i=1}^r(\chi_{\l_{i1}}\otimes\dots\otimes\chi_{\l_{i d_i}}).$$
Il est facile de v\'erifier que
$$A(\chi_\lamb)=A_\lamb.$$
Notons $\PC^+$ l'ensemble des couples $(\lamb,\xi)$ tels que 
$\lamb \in [\PC/A]$ et $\xi \in \Irr A_\lamb$. Posons maintenant 
$$\chi_{\lamb,\xi}^+=\Ind_{G^\circ \rtimes A_\lamb}^G (\chi_\lamb \otimes \xi)$$
$$\Pi_{\lamb,\xi}=\Ind_{G_\lamb}^G \xi.\leqno{\text{et}}$$
Alors l'application $\PC^+ \to \IC(G^\circ,A)$, $(\lamb,\xi) \mapsto (\chi_\lamb,\xi)$ 
induit une bijection entre $\PC^+$ et $\IC(G^\circ,A)/A$. Par cons\'equent, l'application 
$$\fonctio{\PC^+}{\Irr G}{(\lamb,\xi)}{\chi_{\lamb,\xi}^+}$$
est bijective. 

\bigskip

\proposition{triangulaire induit}
{\it Si $(\lamb,\xi) \in \PC^+$, alors
$$\Pi_{\lamb,\xi}=\chi_{\lamb,\xi}^+
+ \sum_{\SS{(\mub,\xi') \in \PC^+} \atop \SS{b_\mub < b_\lamb}} 
\b_{\lamb,\xi,\mub,\xi'} \chi_{\mub,\xi}^+,$$
o\`u $\b_{\lamb,\xi,\mub,\xi'} \in \NM$.}

\bigskip

\proof Tout d'abord, remarquons que, si le r\'esultat de la proposition 
est vrai lorsque $A=A_\lamb$, alors il est vrai dans le cas g\'en\'eral. 
Nous pouvons donc supposer, et nous le ferons, que $A=A_\lamb$. 

Soit $(\mub,\xi') \in \PC^+$ tel que 
$\langle \Pi_{\lamb,\xi},\chi_{\mub,\xi'}^+ \rangle_G \not= 0$. Alors 
$$\langle \Res_{G^\circ}^G \Pi_{\lamb,\xi},\Res_{G^\circ}^G \chi_{\mub,\xi'}^+ 
\rangle_{G^\circ} \not= 0.$$ 
Mais, par la formule de Mackey 
$$\Res_{G^\circ}^G \Pi_{\lamb,\xi} = \xi(1)
\sum_{a \in [A/A_\lamb]}  \Ind_{G_{\lexp{a}{\lamb}}^\circ}^{G^\circ} 
1_{G_{\lexp{a}{\lamb}}^\circ}$$
$$\Res_{G^\circ}^G \chi_{\mub,\xi'}^+ = \xi(1)\sum_{a \in [A/A_\mub]} 
\chi_{\lexp{a}{\mub}}.\leqno{\text{et}}$$
Par suite, il existe $a \in A$ tel que 
$\langle \Ind_{G_{\lexp{a}{\lamb}}^\circ}^{G^\circ} 1_{G_{\lexp{a}{\lamb}}^\circ}, 
\chi_\mub \rangle_{G^\circ} \not= 0$. 
D'apr\`es \ref{induit sn}, ceci implique que $\mub=\lexp{a}{\lamb}$ 
ou que $b_\mub < b_{\lexp{a}{\lamb}}=b_\lamb$. Dans le premier cas, on a alors 
$\lamb=\mub$ car $\lamb$ et $\mub$ parcourent $[\PC/A]$. Il nous reste donc \`a 
montrer que, si $\xi$ et $\xi'$ sont deux caract\`eres irr\'eductibles de 
$A_\lamb$, alors 
$$\langle \Pi_{\lamb,\xi},\chi_{\lamb,\xi'}^+ \rangle_G =\begin{cases}
    1 & \text{si } \xi=\xi' \\
    0 & \text{sinon.}
    \end{cases}\leqno{(*)}$$
Puisque $A=A_\lamb$, on a 
\eqna
\langle \Pi_{\lamb,\xi},\chi_{\lamb,\xi'}^+ \rangle_G 
&=& \langle \xi, \xi' \otimes \Res_{G_\lamb}^G \chit_\lamb \rangle_{G_\lamb} \\
&=& \DS{\frac{1}{|A_\lamb|} \sum_{a \in A_\lamb} \xi(a)\overline{\xi'(a)} 
\Bigl(\frac{1}{|G_\lamb^\circ|} }
\DS{\sum_{w \in G_\lamb^\circ}} \chit_\lamb(wa) \Bigr). \\
\endeqna
Il nous suffit donc de montrer que, pour tout $a \in A_\lamb$, 
$$\frac{1}{|G_\lamb^\circ|} 
\sum_{w \in G_\lamb^\circ} \chit_\lamb(wa) = 1.\leqno{(**)}$$
mais, 
\eqna
\DS{\frac{1}{|G_\lamb^\circ|} \sum_{w \in G_\lamb^\circ} \chit_\lamb(wa)}& = & 
\langle 1_{G_\lamb^\circ a} , \Res_{G_\lamb^\circ a}^G \chit_\lamb 
\rangle_{G_\lamb^\circ a} \\
&=& \langle \Ind_{G_\lamb^\circ a}^{G^\circ a} 1_{G_\lamb^\circ a} , 
\Res_{G^\circ a}^G \chit_\lamb \rangle_{G^\circ a} \\
&=& \langle \Ind_{(G_\lamb^\circ)^a}^{(G^\circ)^a} 1_{(G_\lamb^\circ)^a} , 
(\chi_\lamb)_a \rangle_{(G^\circ)^a} \\
&=&1,
\endeqna
l'avant-derni\`ere \'egalit\'e d\'ecoulant de la commutativit\'e 
du diagramme \ref{diagramme induction tordue}, la derni\`ere 
d\'ecoulant de la formule \ref{induit sn} appliqu\'ee au groupe 
$(G^\circ)^a$.\fin

\bigskip

Les deux corollaires suivants d\'ecoulent facilement (par une r\'ecurrence 
descendante sur $b_\lamb$) de la proposition \ref{triangulaire induit}.

\bigskip

\corollaire{frobenius}
{\it $\ZM\Irr G = 
\DS{\mathop{\oplus}_{(\lamb,\xi) \in \PC^+} \ZM \Pi_{\lamb,\xi}}$.}

\bigskip

\corollaire{brauer}
{\it Soit $\eta \in \Cent(G)$. Alors les conditions suivantes sont \'equivalentes~:
\begin{itemize}
\itemth{1} $\eta \in \ZM \Irr G$. 

\itemth{2} Pour tout $(\lamb,\xi) \in \PC^+$, 
$\langle \Res_{G_\lamb}^G \eta, \xi\rangle_{G_\lamb} \in \ZM$.
\end{itemize}}

\bigskip

\section{Extension canonique, extension pr\'ef\'er\'ee}

\bigskip

\def\coxeter{{\mathrm{Coxeter}}}

Soit $(W^\circ,S)$ un groupe de Coxeter cristallographique 
fini. Notons $\Aut_\coxeter(W^\circ,S)$ 
le groupe des automorphismes $\s$ de $W^\circ$ tels que $\s(S)=S$. Fixons 
un groupe fini $A$ et un morphisme de groupes $A \to \Aut_\coxeter(W^\circ,S)$. 
A travers ce morphisme, $A$ agit sur $W^\circ$ et on peut former le 
produit semi-direct $W=W^\circ \rtimes A$. 

Soit $\chi$ un caract\`ere irr\'eductible de $W^\circ$. Si $a \in A(\chi)$, 
Lusztig \cite[chapitre IV, \SEC 17]{lucs} a d\'efini une extension de $\chi$ \`a 
$W^\circ \rtimes <a>$, appel\'ee {\it extension pr\'ef\'er\'ee} 
dont une des propri\'et\'es est que la repr\'esentation sous-jacente est 
d\'efinie sur $\QM$. Nous voulons montrer ici la proposition suivante~:

\bigskip

\proposition{extension preferee}
{\it Il existe une unique extension de $\chi$ \`a 
$W^\circ \rtimes A(\chi)$ dont la restriction \`a tout sous-groupes 
$W^\circ \rtimes <a>$, o\`u $a$ parcourt $A(\chi)$, soit l'extension pr\'ef\'er\'ee 
de Lusztig.}

\bigskip

\proof L'unicit\'e de l'extension v\'erifiant les conditions de 
l'\'enonc\'e est \'evidente. Montrons en l'existence. 
Tout d'abord, nous pouvons travailler \`a conjugaison pr\`es par un 
\'el\'ement de $A$ en raison du lemme suivant (dont la preuve est imm\'ediate)~:

\bigskip

\begin{quotation}
\lemme{conjugaison preferee}
{\it Soient $a$ et $b$ deux \'el\'ements de $A$. Notons $\chit$ l'extension 
pr\'ef\'er\'ee de Lusztig \`a $W^\circ \rtimes <a>$. Alors $\lexp{b}{\chit}$ est 
l'extension pr\'ef\'er\'ee de Lusztig de $\lexp{b}{\chi}$ 
\`a $W^\circ \rtimes <bab^{-1}>$.}
\end{quotation}

\bigskip

Nous pouvons supposer, et nous le ferons, que $A=\Aut_\coxeter(W^\circ,S)$ 
et que le morphisme $A \to \Aut_\coxeter(W^\circ, S)$ est l'identit\'e. 
Par produit direct, on peut aussi supposer, et nous le ferons, que 
$A(\chi)$ agit transitivement sur l'ensemble des composantes irr\'eductibles 
de $W^\circ$. En d'autres termes, $W^\circ$ est un produit direct de 
$d$ copies d'un groupe de Coxeter cristallographique fini irr\'eductible 
$(W_1,S_1)$, on peut \'ecrire 
$\chi=\chi_1 \boxtimes \dots \boxtimes \chi_d$ o\`u, pour tout $i$, $\chi_i$ 
est un caract\`ere irr\'eductible de $W_1$, $A=(A_1)^d \rtimes \SG_d$, 
o\`u $A_1=\Aut_\coxeter(W_1,S_1)$ et l'image de $A(\chi)$ dans $\SG_d$ 
est un sous-groupe transitif de $\SG_d$. 

Soit $i \in \{1,2,\dots,d\}$. Puisque l'image de $A(\chi)$ dans $\SG_d$ est 
un sous-groupe transitif, il existe $\s \in \SG_d$ et 
$(b_1,\dots, b_d) \in (A_1)^d$ 
tels que $\s(1)=i$ et $\s.(b_1,\dots,b_d) \in A(\chi)$. En particulier, 
$\lexp{b_i}{\chi_i}=\chi_1$. Par cons\'equent, il existe 
$(a_1,\dots,a_d) \in (A_1)^d$ 
tels que, pour tout $i$, on ait $\chi_i = \lexp{a_i}{\chi_1}$. 

Notons $\t$ le cycle de longueur $d$ \'egal \`a $(1,2,\dots,d)$. Alors, 
quitte \`a remplacer $\chi$ par $\lexp{(a_1,\dots,a_d)}{\chi}$ 
(en utilisant le lemme \ref{conjugaison preferee}), on peut supposer, et 
nous le ferons, que $\chi_1=\dots=\chi_d$. En particulier, 
$A(\chi) = A_1(\chi_1)^d \rtimes \SG_d$. 

Supposons d\'emontr\'e le r\'esultat lorsque $d=1$. Notons alors 
$\chit_1$ l'extension de $\chi_1$ \`a $W_1 \rtimes A_1(\chi_1)$ telle que 
$\Res_{W_1 \rtimes <a>}^{W_1 \rtimes A(\chi_1)} \chit_1$ soit 
l'extension pr\'ef\'er\'ee de $\chi_1$ \`a $W_1 \rtimes A_1(\chi_1)$. 
Alors l'extension canonique de $\chit_1 \otimes \dots \otimes \chit_1$ 
\`a $W=(W_1 \rtimes A_1)^d \rtimes \SG_d$ v\'erifie les conditions de 
l'\'enonc\'e. 

On est donc ramen\'e au cas o\`u $d=1$, \cad au cas o\`u $W^\circ=W_1$ est 
irr\'eductible. On peut m\^eme supposer que $A(\chi_1) \not= 1$. 
Dans ce cas, \`a part lorsque $W_1$ est de type $D_4$, 
$A_1$ est cyclique d'ordre $2$ et le r\'esultat est \'evident. 
Supposons donc que $W^\circ$ est de type $D_4$. Alors $A \simeq \SG_3$. 
Deux cas peuvent se produire~:

$\bullet$ Si $3$ ne divise pas $|A(\chi)|$, alors $|A(\chi)|=2$ et 
le r\'esultat est \'evident.

$\bullet$ Si $3$ divise $|A(\chi)|$, alors, d'apr\`es 
\cite[proposition 3.2]{lubook}, $A(\chi)=A$ et il existe une extension 
$\chit'$ de $\chi$ \`a $W^\circ \rtimes A(\chi)$ dont la restriction \`a 
$W^\circ \rtimes <b>$ est l'extension 
pr\'ef\'er\'ee de Lusztig lorsque $b \in A$ est d'ordre $3$. 
Notons $\e$ le caract\`ere signature de $A\simeq \SG_3$. Il est 
alors clair que $\chit'$ ou $\chit' \otimes \e$ satisfait 
aux conditions de la proposition.\fin

\bigskip

L'extension de $\chi$ v\'erifiant les conditions de la proposition 
pr\'ec\'edente sera appel\'ee l'{\it extension pr\'ef\'er\'ee} de $\chi$ 
\`a $W^\circ \rtimes A(\chi)$. 

\bigskip

\corollaire{canonique preferee}
{\it Si $A$ agit sur $W^\circ$ seulement par permutation des composantes 
irr\'eductibles, alors l'extension pr\'ef\'er\'ee de $\chi$ \`a 
$W^\circ \rtimes A(\chi)$ est l'extension canonique.}

\bigskip

\corollaire{restriction preferee}
{\it Soit $A'$ un sous-groupe de $A$ et posons $W'=W^\circ \rtimes A'$. Soit 
$\chi \in \Irr W^\circ$ et notons $\chit$ (respectivement $\chit'$) l'extension 
canonique de $\chi$ \`a $W(\chi)$ (respectivement $W'(\chi)$). 
Alors $\chit'$ est la restriction de $\chit$.}

\newpage

{\Large \part*{Appendice B. Sommes de Gauss}}

\bigskip

Le but de cet appendice est de donner des formules pour les 
racines de l'unit\'e $\GC(\Gb,\z)$ lorsque $\z \in \ZC_\cus^\we(\Gb)$ 
est $F$-stable. D'apr\`es \cite[proposition 2.4]{DLM2}, 
$\GC(\Gb,\z)$ peut \^etre d\'ecrit par un produit de 
{\it sommes de Gauss}. Dans la section \ref{section gauss}, nous rappelons 
quelques propri\'et\'es g\'en\'erales des sommes de Gauss. Nous appliquons 
ces r\'esultats dans la section \ref{calcul constantes} pour calculer 
explicitement les racines de l'unit\'e.

\medskip

\noindent{\sc Notation - } Nous noterons $r$ l'entier naturel 
non nul tel que $q=p^r$. 

\bigskip

\section{Sommes de Gauss\label{section gauss}}

\medskip

Un caract\`ere additif $\chi_1 : \FM_p \to \qlb^\times$ non trivial 
a \'et\'e fix\'e dans le \SEC\ref{sous chi}. Si $s \in \NM^*$, on rappelle que 
$\chi_s : \FM_{p^s} \to \qlb^\times$ est d\'efini par 
$\chi_s = \chi_1 \circ \Tr_s$, o\`u $\Tr_s : \FM_{p^s} \to \FM_p$ 
est la trace. 
Si $s \in \NM^*$ et si $\th : \FM_{p^s}^\times \to \qlb^\times$ 
est un caract\`ere lin\'eaire, nous noterons 
$$\GC_s(\th) = \sum_{x \in \FM_{p^s}^\times} \th(x) \chi_s(x)$$ 
la {\it somme de Gauss} associ\'ee. La premi\`ere identit\'e 
sur les sommes de Gauss est bien connue~: 
si $\th$ est non trivial, alors 
\equat\label{norme gauss}
\GC_s(\th) \GC_s(\th^{-1}) = p^s \th(-1).
\endequat

Nous allons maintenant traiter un cas particulier intervenant 
dans le groupe sp\'ecial unitaire (voir la preuve de la proposition 
\ref{sln gauss}). Il s'agit des sommes de Gauss de la forme 
$\GC_{2r}(\th)$, o\`u $\th$ est un caract\`ere lin\'eaire 
non trivial de $\FM_{q^2}^\times$ tel que $\th^{q+1}=1$ 
(rappelons que $q=p^r$). Tout d'abord, notons qu'un tel caract\`ere est 
trivial sur $\FM_q^\times$ (surjectivit\'e de la norme). 
Il d\'ecoule alors de la formule 
\ref{norme gauss} que $\GC(\th)=\pm q$. Nous allons d\'eterminer exactement 
le signe. 

Pour cela, notons $\Tr : \FM_{q^2} \to \FM_q$ la trace. Elle est surjective 
et $\FM_q$-lin\'eaire. Donc il existe $\xi \in \FM_{q^2}^\times$ tel que 
$\Tr \xi = 0$, de sorte que $\Ker \Tr = \FM_q \xi$. 
Si $q$ est pair, alors $\Ker \Tr = \FM_q$ (donc on peut prendre $\xi=1$). 
Si $q$ est impair, $\xi$ est un \'el\'ement de $\FM_{q^2}^\times$ tel que 
$\xi^2$ appartient \`a $\FM_q$ mais n'est pas le carr\'e d'un 
\'el\'ement de $\FM_q$. Alors, si $\th$ est 
un caract\`ere lin\'eaire non trivial de $\FM_{q^2}^\times$ 
tel que $\th^{q+1}=1$, on a
\equat\label{waldspurger}
\GC_{2r}(\th)=\th(\xi) q.
\endequat

\medskip

\noindent{\sc Remarque - } On a toujours $\xi^2 \in \FM_q^\times$, 
donc $\th(\xi) \in \{1,-1\}$. De plus $\th(\xi)$ ne d\'epend pas 
du choix de $\xi$. 
Notons aussi que, si $q$ est pair, on a $\xi \in \FM_q^\times$ 
donc $\th(\xi)=1$. \finl

\medskip

\noindent{\sc Preuve de \ref{waldspurger} - } 
Nous remercions J.L. Waldspurger pour nous avoir pr\'esent\'e 
la formule \ref{waldspurger} ainsi que l'argument suivant. 
Tout d'abord, puisque $\th$ est trivial sur $\FM_q^\times$, on a 
$$\GC_{2r}(\th) = \frac{1}{q-1} \sum_{x \in \FM_{q^2}^\times} 
\sum_{y \in \FM_q^\times} \th(xy) \chi_r(\Tr(x)).$$
Par un changement de variable \'evident, on obtient
\eqna
\GC_{2r}(\th) &=& \DS{\frac{1}{q-1} 
\sum_{x \in \FM_{q^2}^\times} \sum_{y \in \FM_q^\times} \th(x) \chi_r(y\Tr(x))}\\
&=& \DS{\frac{1}{q-1} 
\sum_{x \in \FM_{q^2}^\times} \th(x) 
\Bigl(\sum_{y \in \FM_q^\times} \chi_r(y\Tr(x))\Bigr)}\\
&=& \DS{\frac{1}{q-1} 
\Bigl((q-1) \sum_{x \in \FM_{q^2}^\times \cap \Ker \Tr} \th(x) 
-\sum_{\SS{x \in \FM_{q^2}^\times}\atop \SS{\Tr x \not= 0}} \th(x)\Bigr)}\\
&=& \DS{\frac{1}{q-1} 
\Bigl(q \sum_{x \in \FM_{q^2}^\times \cap \Ker \Tr} \th(x) 
-\sum_{x \in \FM_{q^2}^\times} \th(x)\Bigr)}\\
\endeqna
Puisque $\th$ est non trivial, la deuxi\`eme somme est nulle. 
Puisque $\Ker \Tr = \FM_q \xi$ et que $\th$ est trivial sur $\FM_q^\times$, 
la premi\`ere somme vaut $(q-1) \th(\xi)$. Au bilan, il nous reste 
$\GC_{2r}(\th)=\th(\xi) q$.\fin

\bigskip

Nous terminons cette section par un cas particulier classique, d\^u \`a Gauss. 
Supposons ici $p$ impair. Soit $\LC_s : \FM_{p^s}^\times \to \{1,-1\}$ le 
{\it caract\`ere de Legendre}, \cad l'unique caract\`ere d'ordre $2$. 
La formule \ref{norme gauss} montre que $\l_s=p^{-s/2} \GC_s(\LC_s)$ 
est une racine quatri\`eme de l'unit\'e. Sa valeur d\'epend du choix de 
$\chi_1$ ainsi que du choix d'une racine carr\'ee de $p$ dans $\qlb$. 
Posons $i=\tilde{\jmath}(1/4)$ (rappelons que 
$\jmath : \QM \to \qlb^\times$ est le morphisme de groupe de noyau $\ZM$ 
d\'efini dans la sous-section \ref{sous groupes}). Alors $i$ est 
une racine primitive quatri\`eme de l'unit\'e. 
Si $p \equiv 1 \mod 4$, on prend $p^{1/2} = \GC_1(\LC_1)$. 
Si $p \equiv 3 \mod 4$, on prend $p^{1/2} = i^{-1} \GC_1(\LC_1)$. 
Il r\'esulte de \ref{norme gauss} que $p^{1/2}$ est alors une 
racine carr\'ee de $p$. Alors, d'apr\`es \cite[VII.356]{gauss}, on a
\equat\label{legendre}
\l_s = \begin{cases}
1 & \text{si } p \equiv 1 \mod 4, \\
i^s & \text{si } p \equiv 3 \mod 4.
\end{cases}
\endequat

\bigskip

\section{Calcul de $\GC(\Gb,\z)$\label{calcul constantes}}

\medskip

\soussection{Propri\'et\'es g\'en\'erales} 
Rappelons qu'il a \'et\'e fix\'e un 
morphisme injectif $\kappa : \FM^\times \injto \qlb^\times$ qui fournit, 
par restriction \`a $\FM_q^\times$ un caract\`ere lin\'eaire 
que l'on notera encore $\kappa$. 
Avant d'exprimer $\GC(\Gb,\z)$ 
sous forme d'un produit de $\GC_s(\kappa^m)$, nous aurons besoin de quelques 
notations. Notons $(\varpi_\a^\vee)_{\a \in \D_0}$ 
la $\QM$-base de $Y(\Tb_0/\Zb(\Gb)^\circ)$ duale de $\D_0$. 
Notons $\D_0/\phi_0$ l'ensemble des orbites de $\phi_0$ dans 
$\D_0$. Pour finir, si $\o \in \D_0/\phi_0$, notons $\a_\o \in \o$ 
un repr\'esentant et $r_\o$ l'entier naturel non nul tel 
que $F^{|\o|}(\a_\o)=p^{r_\o} \a_\o$. 

\medskip

Soit $\z \in H^1(F,\ZC(\Gb))^\we$. Fixons $\dot{\z} \in X(\Tb_0/\Zb(\Gb)^\circ)$ 
tel que $\z = \kappa \circ \Res_{\ZC(\Gb)}^{\Tb/\Zb(\Gb)^\circ} \dot{\z}$. Alors, 
d'apr\`es \cite[proposition 2.4]{DLM2}, on a 
$(q^{|\o|}-1) < \dot{\z},\varpi_{\a_\o}^\vee >_{\Tb_0} \in \ZM$ pour tout 
$\o \in \D_0/\phi_0$ (car $(F-1)(\dot{\z})$ est trivial sur $\ZC(\Gb)$) et 
\equat\label{dlm gauss}
\GC(\Gb,\z) = \eta_\Gb q^{-\frac{1}{2} \rang_\sem(\Gb)} 
\prod_{\o \in \D_0/\phi_0} 
\GC_{r_\o}(\kappa^{(p^{r_\o}-1)< \dot{\z},\varpi_{\a_\o}^\vee >_{\Tb_0}}).
\endequat
Il est imm\'ediat que le membre de droite ne d\'epend pas du choix 
de $\kappa$, ce qui est souhaitable. Nous allons maintenant 
rappeler quelques propri\'et\'es des nombres $\GC(\Gb,\z)$. 

Soit $\Gbh$ un groupe r\'eductif connexe muni d'un endomorphisme 
de Frobenius $F : \Gb \to \Gb$ d\'efini sur $\FM_q$ et soit 
$\pi : \Gbh \to \Gb$ est un morphisme isotypique d\'efini sur $\FM_q$.  
Alors $\pi$ induit un morphisme surjectif $\ZC(\Gbh) \to \ZC(\Gb)$ et, 
si on note alors $\hat{\z} = \z \circ \pi \in \ZC(\Gbh)^\we$, alors 
il d\'ecoule facilement de \ref{dlm gauss} que 
\equat\label{isotypique gauss}
\GC(\Gb,\z)=\GC(\Gbh,\hat{\z}).
\endequat

D'autre part, si $\Gb = \Gb_0 \times \dots \times \Gb_0$ ($k$ fois) 
et si $F$ permute transitivement les composantes $\Gb_0$, 
notons $\z_0 \in H^1(F^k,\ZC(\Gb_0))$ le caract\`ere correspondant 
\`a $\z$. On a alors, d'apr\`es \cite[preuve de la proposition 2.5]{DLM2}, 
\equat\label{couronne gauss}
\GC(\Gb,\z) = \GC(\Gb_0,\z_0). 
\endequat
Ici, $\GC(\Gb,\z_0)$ est calcul\'e en utilisant l'isog\'enie $F^k$. 

\bigskip

\soussection{Cas cuspidal}
Nous supposons maintenant que $\z \in \ZC_\cus^\we(\Gb)$ et que $\z$ 
est $F$-stable. Alors toutes les composantes irr\'eductibles de 
$\Gb$ sont de type $A$. Par suite, compte tenu de \ref{isotypique gauss} 
et \ref{couronne gauss}, il suffit de calculer $\GC(\Gb,\z)$ lorsque 
$\Gb=\Sb\Lb_n(\FM)$ et $F = \s^k \circ F_{\text{nat}}$, o\`u $k \in \{0,1\}$, 
$F_{\text{nat}} : \Gb \to \Gb$ est l'endomorphisme de Frobenius 
d\'eploy\'e sur $\FM_q$ 
et $\s : \Gb \to \Gb$, $g \mapsto J \lexp{t}{g}^{-1} J$. Ici, 
$J$ d\'esigne la matrice monomiale dont les coefficients sur la 
deuxi\`eme diagonale sont tous \'egaux \`a $1$. 

\medskip

\proposition{sln gauss}
{\it Supposons que $\Gb=\Sb\Lb_n(\FM)$ et que $F=\s^k \circ F_{\mathrm{nat}}$, o\`u 
$k \in \{0,1\}$. On pose $\e=(-1)^k$. Soit $\z \in \ZC_\cus^\we(\Gb)$ 
et supposons que $\z$ est $F$-stable. Alors $\z$ est d'ordre $n$, 
$n$ divise $q-\e$, et~:
$$\GC(\Gb,\z) = \begin{cases}
1 & \text{si $n$ est impair,} \\
-\l_r (-1)^{\frac{(q-\e)(n-2)}{8}} & \text{si $n$ est pair.}
\end{cases}$$}

\bigskip

\noindent{\sc Remarque - } 
La formule de la proposition pr\'ec\'edente 
est close gr\^ace \`a \ref{legendre}.\finl

\bigskip

\proof Le fait que $\z$ est d'ordre $n$ d\'ecoule de la table \ref{tabletable}. 
Par cons\'equent, $\z$ est injectif et $F$-stable, donc $F$ agit trivialement 
sur $\ZC(\Gb)$. Par cons\'equent, $n$ divise $q-\e$ car $F$ agit sur 
$\ZC(\Gb)$ par \'el\'evation \`a la puissance $\e q$. En particulier, $n$ 
est premier \`a $p$. Nous aurons besoin de quelques notations. 
On num\'erote les $n-1$ racines simples comme suit~:
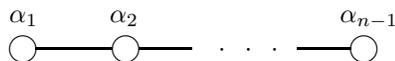
\begin{center}
\begin{picture}(150,40)
\put( 15, 10){\circle{10}}
\put( 20, 10){\line(1,0){29}}
\put(54, 10){\circle{10}}
\put(59, 10){\line(1,0){20}}
\put(89,  7){$\cdot$}
\put(99,  7){$\cdot$}
\put(109,  7){$\cdot$}
\put(119, 10){\line(1,0){20}}
\put(144, 10){\circle{10}}
\put(10, 20){$\a_1$}
\put(48, 20){$\a_2$}
\put(135, 20){$\a_{n-1}$}
\end{picture}
\end{center}
Notons pour simplifier $\varpi_j^\vee=\varpi_{\a_j}^\vee$ pour $1 \le j \le n-1$. 
Il existe alors un entier $k \in \ZM$ premier \`a $n$ tel que, 
pour tout $1 \le j \le n-1$, on ait 
$$< \dot{\z}, \varpi_j^\vee >_{\Tb_0} = \frac{kj}{n}.$$

\medskip

$\bullet$ Commen\c{c}ons par \'etudier le cas d\'eploy\'e (\cad $k=0$, ou encore 
$\e=1$). Alors $\phi_0$ est l'identit\'e et on a 
\eqna
\GC(\Gb,\z)&=&(-1)^{n-1} q^{-\frac{n-1}{2}} 
\DS{\prod_{j=1}^{n-1} \GC_r((\kappa^{(q-1)/n})^{kj})} \\
&=& \DS{(-1)^{n-1} q^{-\frac{n-1}{2}} \prod_{j=1}^{n-1}\GC_r((\kappa^{(q-1)/n})^j)},
\endeqna
la derni\`ere \'egalit\'e d\'ecoulant de ce que $k$ est premier \`a $n$. 
Notons $\g$ le caract\`ere lin\'eaire $\kappa^{(q-1)/n}$ de $\FM_q^\times$. 
En regroupant $j$ et $n-j$ et en utilisant \ref{norme gauss},  
on obtient~:
$$\GC(\Gb,\z)=\begin{cases}
\g(-1)^{1+2+ \dots + \frac{n-1}{2}} & \text{si $n$ est impair,}\\
-\l_r \g(-1)^{1+2+ \dots + \frac{n-2}{2}} & \text{si $n$ est pair.}
\end{cases}$$
Puisque $\kappa(-1)=-1$, on a $\g(-1)=(-1)^{(q-1)/n}$. Cela montre 
le r\'esultat annonc\'e lorsque $n$ est pair. Lorsque $n$ est impair, 
la simplification provient du fait que $\g(-1)^2=\g(-1)^n=1$ 
et donc que $\g(-1)=1$. 

\medskip

$\bullet$ \'Etudions maintenant le cas o\`u $k=1$, \cad $\e=-1$. 
Alors $\phi_0(\a_j)=\a_{n-j}$ pour tout $1 \le j \le n-1$. Deux cas se pr\'esentent~:
$$\GC(\Gb,\z) =  \begin{cases}
q^{-\frac{n-1}{2}}\DS{\prod_{j=1}^{(n-1)/2} \GC_{2r}((\kappa^{(q^2-1)/n})^{kj})} & 
\text{si $n$ est impair,} \\ 
-\l_r q^{-\frac{n-2}{2}}
\DS{\prod_{j=1}^{(n-2)/2} \GC_{2r}((\kappa^{(q^2-1)/n})^{kj})} & 
\text{sinon.} \\
\end{cases}$$
Notons $\g$ la restriction de $\kappa^{k(q^2-1)/n}$ \`a $\FM_{q^2}^\times$. 
Alors $\g$ est non trivial et $\g^{q+1}=1$ car $n$ divise $q+1=q-\e$. 
Notons, comme dans la preuve de \ref{waldspurger}, $\x$ un g\'en\'erateur 
du noyau de la trace $\Tr : \FM_{q^2} \to \FM_q$. D'apr\`es 
\ref{waldspurger}, on a 
$$\GC(\Gb,\z)=\begin{cases}
\g(\x)^{1+2+ \dots + \frac{n-1}{2}} & \text{si $n$ est impair,}\\
-\l_r \g(\x)^{1+2+ \dots + \frac{n-2}{2}} & \text{si $n$ est pair.}
\end{cases}$$
Lorsque $n$ est impair, alors $\g(\x)^2=\g(\x)^n=1$ (car $\x^2 \in \FM_q^\times$ 
et $\g$ est trivial sur $\FM_q^\times$) et donc $\g(\x)=1$, ce qui 
montre le r\'esultat attendu. Supposons maintenant $n$ pair. 
Alors $q$ et $k$ sont impairs et on peut prendre $\xi=\xi_0^{(q+1)/2}$, o\`u 
$\xi_0$ est un g\'en\'erateur de $\FM_{q^2}^\times$. On a alors 
$\g(\xi)=\kappa(\xi_0^{(q^2-1)/2})^{k(q+1)/n}$. Or, 
$\xi_0^{(q^2-1)/2}=-1$, donc $\g(\xi)=(-1)^{k(q+1)/n}=(-1)^{(q+1)/n}$. Cela 
termine la preuve de la proposition.\fin

\newpage


\end{document}